# NATURAL PRODUCT $\times_n$ ON MATRICES

W. B. Vasantha Kandasamy
Florentin Smarandache

**2012**

# NATURAL PRODUCT $\times_n$ ON MATRICES

# CONTENTS

# PREFACE

In this book the authors introduce a new product on matrices called the natural product. We see when two row matrices of $1 \times n$ order are multiplied, the product is taken component wise; for instance if $X = (x_1, x_2, x_3, \ldots, x_n)$ and $Y = (y_1, y_2, y_3, \ldots, y_n)$ then $X \times Y = (x_1y_1, x_2y_2, x_3y_3, \ldots, x_ny_n)$ which is also the natural product of X with Y. But we cannot find the product of a $n \times 1$ column matrix with another $n \times 1$ column matrix, infact the product is not defined. Thus if

$$X = \begin{bmatrix} x_1 \\ x_2 \\ \vdots \\ x_n \end{bmatrix} \text{ and } Y = \begin{bmatrix} y_1 \\ y_2 \\ \vdots \\ y_n \end{bmatrix}$$

under natural product

$$X \times_n Y = \begin{bmatrix} x_1 y_1 \\ x_2 y_2 \\ \vdots \\ x_n y_n \end{bmatrix}.$$

Thus by introducing natural product we can find the product of column matrices and product of two rectangular matrices of same order. Further this product is more natural which is just identical with addition replaced by multiplication on these matrices.

Another fact about natural product is this enables the product of any two super matrices of same order and with same type of partition. We see on supermatrices products cannot be defined easily which prevents from having any nice algebraic structure on the collection of super matrices of same type.

This book has eight chapters. The first chapter is introductory in nature. Polynomials with matrix coefficients are introduced in chapter two. Algebraic structures on these polynomials with matrix coefficients is defined and described in chapter three. Chapter four introduces natural product on matrices. Natural product on super matrices is introduced in chapter five. Super matrix linear algebra is introduced in chapter six. Chapter seven claims only after this notion becomes popular we can find interesting applications of them. The final chapter suggests over 100 problems some of which are at research level.

We thank Dr. K.Kandasamy for proof reading and being extremely supportive.

<div style="text-align: right;">
W.B.VASANTHA KANDASAMY  
FLORENTIN SMARANDACHE
</div>

**Chapter One**

# INTRODUCTION

In this chapter we only indicate as reference of those the concepts we are using in this book. However the interested reader should refer them for a complete understanding of this book.

In this book we define the notion of natural product in matrices so that we have a nice natural product defined on column matrices, m × n (m ≠ n) matrices. This extension is the same in case of row matrices.

We make use of the notion of semigroups and Smarandache semigroups refer [13].

Also the notion of semirings, Smarandache semirings, semi vector spaces and semifields are used, please refer [16].

Likewise S-rings, S-ideals, S-subrings are also used, refer [18].

The concept of polynomials with matrix coefficients are used. That is if

$$p(x) = \sum_{i=0}^{\infty} a_i x^i$$

where x is an indeterminate and if $a_i$ is a matrix (a square matrix or a row matrix of a column matrix or a m × n matrix m ≠ n), then p(x) is a polynomial in the variable x with matrix coefficients ('or' used in the mutually exclusive sense).

Suppose

$$p(x) = \begin{bmatrix} 3 \\ 2 \\ 0 \\ -1 \end{bmatrix} + \begin{bmatrix} -2 \\ 3 \\ 1 \\ 5 \end{bmatrix} x + \begin{bmatrix} 0 \\ 1 \\ 0 \\ 2 \end{bmatrix} x^3 + \begin{bmatrix} 7 \\ 0 \\ 1 \\ 0 \end{bmatrix} x^5$$

is a polynomial with column matrix coefficients.

We also introduce polynomial matrix coefficient semiring. We call usual matrices as simple matrices.

The super matrix concepts are used. If $X = (a_1 \ a_2 \ | \ a_3 \ a_4 \ | \ a_5)$, $a_i \in R$ (or Q or Z) then X is a super row matrix [8, 19].

If

$$Y = \begin{bmatrix} a_1 \\ \hline a_2 \\ \hline a_3 \\ a_4 \\ \hline a_5 \\ a_6 \end{bmatrix}, \ a_i \in R \text{ (or Q or Z)}$$

then Y is a super column matrix.

Let

$$M = \begin{bmatrix} x_1 & x_2 & x_3 & x_4 \\ \hline x_5 & x_6 & x_7 & x_8 \\ \hline x_9 & x_{10} & x_{11} & x_{12} \\ \hline x_{13} & x_{14} & x_{15} & x_{16} \end{bmatrix}$$

with $x_i \in R$ (or Q or Z); $1 \le i \le 16$ be a super square matrix.

$$P = \begin{bmatrix} a_1 & a_4 & a_7 & a_{10} & a_{13} & a_{16} \\ a_2 & a_5 & a_8 & a_{11} & a_{14} & a_{17} \\ a_3 & a_6 & a_9 & a_{12} & a_{15} & a_{18} \end{bmatrix}$$

is a super row vector.

$$S = \begin{bmatrix} a_1 & a_2 & a_3 & a_4 & a_5 & a_6 & a_7 & a_8 \\ a_9 & a_{10} & a_{11} & \ldots & \ldots & \ldots & \ldots & a_{16} \\ a_{17} & a_{18} & a_{19} & \ldots & \ldots & \ldots & \ldots & a_{24} \\ a_{25} & a_{26} & a_{27} & \ldots & \ldots & \ldots & \ldots & a_{32} \end{bmatrix}$$

is a super $4 \times 8$ matrix (vector).

Likewise

$$B = \begin{bmatrix} a_1 & a_2 & a_3 \\ \hline a_4 & a_5 & a_6 \\ a_7 & a_8 & a_9 \\ \hline a_{10} & a_{11} & a_{12} \\ \hline a_{13} & a_{14} & a_{15} \\ a_{16} & a_{17} & a_{18} \\ \hline a_{19} & a_{20} & a_{21} \end{bmatrix} \text{ with } a_i \in R \text{ (or Q or Z)}$$

is a super column vector [8, 19].

Also we use the notion of vector spaces, Smarandache vector spaces and Smarandache linear algebra [17].

**Chapter Two**

# POLYNOMIALS WITH MATRIX COEFFICIENTS

In this chapter we define polynomials in the variable x with coefficients from the collection of matrices of same order. We call such polynomials as matrix coefficient polynomials or polynomials with matrix coefficients. We first give some examples before we define operations on them.

***Example 2.1:*** Let $p(x) = (5, 3, 0, -3, 2) + (0, 1, 2, 3, 4)x + (7, 0, 1, 0, 1)x^2 + (-7, -9, 10, 0, 0)x^5 - (3, 2, 1, 2, 1)x^7$; we see $p(x)$ is a polynomial in the variable x with row matrix coefficients.

***Example 2.2:*** Let $m(x) = \begin{bmatrix} 1 \\ 2 \\ 3 \\ 0 \\ 0 \end{bmatrix} + \begin{bmatrix} -1 \\ -2 \\ 3 \\ 4 \\ 5 \end{bmatrix} x + \begin{bmatrix} 6 \\ 0 \\ 1 \\ 2 \\ 0 \end{bmatrix} x^2 + \begin{bmatrix} 1 \\ -1 \\ 2 \\ -2 \\ 3 \end{bmatrix} x^5$ be a column matrix coefficient polynomial or a polynomial with column matrix coefficients.

*Example 2.3:* Let

$$p(x) = \begin{bmatrix} 3 & 0 \\ -1 & 2 \end{bmatrix} + \begin{bmatrix} 1 & 0 \\ 0 & 2 \end{bmatrix} x^2 + \begin{bmatrix} 0 & 1 \\ 0 & 3 \end{bmatrix} x^3 + \begin{bmatrix} 1 & 0 \\ 4 & 0 \end{bmatrix} x^5 +$$

$$\begin{bmatrix} 1 & 4 \\ 0 & 0 \end{bmatrix} x^8 + \begin{bmatrix} 0 & 0 \\ 1 & 2 \end{bmatrix} x^9 + \begin{bmatrix} 0 & 1 \\ 5 & 0 \end{bmatrix} x^{10}$$

be a square matrix coefficient polynomial.

*Example 2.4:* Let

$$T(x) = \begin{bmatrix} 2 & 1 \\ 0 & 1 \\ 5 & 2 \end{bmatrix} + \begin{bmatrix} 1 & 0 \\ 1 & 1 \\ 1 & 0 \end{bmatrix} x + \begin{bmatrix} 1 & 2 \\ 0 & 3 \\ 4 & 0 \end{bmatrix} x^3 + \begin{bmatrix} 9 & 1 \\ 0 & 2 \\ 6 & 0 \end{bmatrix} x^7$$

be a polynomial with $3 \times 2$ matrix coefficient.

Now we define some operations on the collection.

**DEFINITION 2.1:** *Let*

$$V_R = \left\{ \sum_{i=0}^{\infty} a_i x^i \;\middle|\; a_i = (x_1, \ldots, x_n) \right.$$

*are $1 \times n$ row matrices, $x_i \in R$ (or Q or Z); $1 \leq i \leq n$} be the collection of all row matrix coefficient polynomials. $V_R$ is a group under addition.*

*For if $p(x) = \sum_{i=0}^{\infty} a_i x^i$ and $q(x) = \sum_{j=0}^{\infty} b_j x^j$ then*

*we define* $p(x) + q(x) = \sum_{i=0}^{\infty} a_i x^i + \sum_{j=0}^{\infty} b_j x^j$

$$= \sum_{i=0}^{\infty} (a_i + b_i) x^i$$

$0 = (0,...,0) + (0,...,0)x + ... + (0,... ,0)x^n$ $(n \in N)$ is defined as the row matrix coefficient zero polynomial.

Let $p(x) = \sum_{i=0}^{\infty} a_i x^i$ now $-p(x) = \sum_{i=0}^{\infty} -a_i x^i$ is defined as the inverse of the row matrix coefficient polynomial. Thus $(V_R, +)$ is an abelian group of infinite order.

*Example 2.5:* Let

$$V_R = \left\{ \sum_{i=0}^{\infty} a_i x^i \,\middle|\, a_i = (x_1, x_2, x_3, x_4) \text{ with } x_j \in Q; 1 \leq j \leq 4 \right\}$$

be the collection of row matrix coefficient polynomials $V_R$ is a group under addition.

For if $p(x) = (0, 2, 1, 0) + (7, 0, 1, 2)x + (1, 1, 1, 1)x^3 + (0, 1, 2, 0)x^5$ and

$q(x) = (7, 8, 9, 10) + (3, 1, 0, 7)x + (3,0,1, 4)x^3 - (4, 2, 3, 4)x^4 + (7, 1, 0, 0)x^5 + (1, 2, 3, 4)x^8$ are in $V_R$ then

$p(x) + q(x) = ((0, 2, 1, 0) + (7, 8, 9, 10)) + ((7, 0, 1, 2) + (3, 1, 0, 7))x + ((1, 1, 1, 1) + (3, 0, 1, 4))x^3 + ((0, 0, 0, 0) - (4, 2, 3, 4))x^4 + ((0, 1, 2, 0) + (7, 1, 0, 0))x^5 + (1, 2, 3, 4)x^8$

$= (7, 10, 10, 10) + (10, 1, 1, 9)x + (4, 1, 2, 5)x^3 - (4, 2, 3, 4)x^4 + (7, 2, 2, 0)x^5 + (1, 2, 3, 4)x^8$.

We see $-p(x) = (0, -2, -1, 0) + (-7, 0, -1, -2)x + (-1, -1, -1, -1)x^3 + (0, -1, -2, 0)x^5$ acts as the additive inverse of $p(x)$.

*Example 2.6:* Let

$$V_R = \left\{ \sum_{i=0}^{\infty} a_i x^i \;\middle|\; a_i = (x_1, x_2, x_3); x_j \in Z_{12}; 1 \leq j \leq 3 \right\}$$

be the collection of all row coefficient polynomials. $V_R$ is a group under modulo addition 12.

*Example 2.7:* Let

$$V_R = \left\{ \sum_{i=0}^{10} a_i x^i \;\middle|\; a_i = (d_1, d_2) \text{ with } d_j \in Q; 1 \leq j \leq 2 \right\}$$

be the row coefficient polynomial. $V_R$ is a group under addition.

*Example 2.8:* Let

$$V_R = \left\{ \sum_{i=0}^{5} a_i x^i \;\middle|\; a_i = (x_1, x_2); x_1, x_2 \in Z_{10} \right\}$$

be the row coefficient polynomial. $V_R$ is a finite group under addition.

We now can define other types of operations on $V_R$.

We see if $(1,1,1,1)x^3 - (0, 0, 8, 27) = p(x)$ then $(1,1,1,1)x^3 - (0, 0, 2, 3)^3 = p(x)$

$= [(1, 1, 1, 1)x - (0, 0, 2, 3)] [(1, 1, 1, 1)x^2 + (0, 0, 2, 3)x + (0, 0, 4, 9)]$.

For this we have to define another operation on $V_R$ called product.

Throughout this book $V_R$ will denote polynomial with row matrix coefficient in the variable x.

We know $V_R$ is a group under addition.

Now we can define product on $V_R$ as follows:

Let $p(x) = (0,1,2) + (3,4,0)x + (2,1,5)x^2 + (3,0,2)x^3$ and

$q(x) = (6,0,2) + (0,1,4)x + (3,1,0)x^2 + (1,2,3)x^4$ be in $V_R$.

We define product of $p(x)$ with $q(x)$ as follows.

$p(x) \times q(x)$

$= [(0,1,2) + (3,4,0)x + (2,1,5)x^2 + (3,0,2)x^3] \times [(6,0,2) + (0,1,4)x + (3,1,0)x^2 + (1,2,3)x^4]$

$= (0,1,2)(6,0,2) + (3,4,0)(6,0,2)x + (2,1,5)(6,0,2)x^2 + (3,0,2)(6,0,2)x^3 + (0,1,2)(0,1,4)x + (3,4,0)(0,1,4)x^2 + (2,1,5)(0,1,4)x^3 + (3,0,2)(0,1,4)x^4 + (0,1,2)(3,1,0)x^2 + (3,4,0)(3,1,0)x^3 + (2,1,5)(3,1,0)x^4 + (3,0,2)(3,1,0)x^5 + (0,1,2)(1,2,3)x^4 + (3,4,0)(1,2,3)x^5 + (2,1,5)(1,2,3)x^6 + (3,0,2)(1,2,3)x^7$

$= (0,0,4) + (18,0,0)x + (12,0,10)x^2 + (18,0,4)x^3 + (0,1,8)x + (0,4,0)x^2 + (0,1,20)x^3 + (0,0,8)x^4 + (0,1,0)x^2 + (9,4,0)x^3 + (6,1,0)x^4 + (9,0,0)x^5 + (0,2,6)x^4 + (3,8,0)x^5 + (2,2,15)x^6 + (3,0,6)x^7$

$= (0,0,4) + (18,1,8)x + (12,5,10)x^2 + (27,5,24)x^3 + (6,3,14)x^4 + (12,8,0)x^5 + (2,2,15)x^6 + (3,0,6)x^7$.

Now we see with componentwise product we see $V_R$ under product is a commutative semigroup.

We see $V_R$ has zero divisors.

Now we proceed onto give one or two examples.

*Example 2.9:* Let

$$V_R = \left\{ \sum_{i=0}^{5} a_i x^i \;\middle|\; a_i = (x_1, \ldots, x_8); x_i \in Q; 1 \le i \le 8 \right\}$$

be a semigroup of row matrix polynomials. $V_R$ is a monoid under product.

Now we see $(V_R, +, \times)$ is a commutative ring of polynomials with row matrix coefficients.

We give examples of them.

*Example 2.10:* Let

$$V_R = \left\{ p(x) = \sum_{i=0}^{\infty} a_i x^i \;;\; a_j = (x_1, x_2 \ldots, x_{18}); x_i \in R; 1 \le i \le 18 \right\}$$

be a ring of polynomials with row matrix coefficients.

Now we have shown examples of polynomial row matrix coefficients in the variable x.

*Example 2.11:* Let

$$V_C = \left\{ \sum_{i=0}^{\infty} a_i x^i \;\middle|\; a_j = \begin{bmatrix} x_1 \\ x_2 \\ x_3 \\ x_4 \\ x_5 \end{bmatrix} \text{ with } x_i \in Z; 1 \le i \le 5 \right\},$$

$V_C$ is a group under addition.

$$p(x) = \begin{bmatrix} 3 \\ 0 \\ 1 \\ 0 \\ 2 \end{bmatrix} + \begin{bmatrix} 1 \\ 0 \\ 0 \\ 2 \\ 0 \end{bmatrix} x + \begin{bmatrix} 0 \\ 1 \\ 2 \\ 3 \\ 0 \end{bmatrix} x^2 \text{ and}$$

$$q(x) = \begin{bmatrix} 4 \\ 2 \\ 1 \\ -4 \\ 0 \end{bmatrix} + \begin{bmatrix} 2 \\ 3 \\ 4 \\ 5 \\ 0 \end{bmatrix} x + \begin{bmatrix} 2 \\ 3 \\ 1 \\ 4 \\ 5 \end{bmatrix} x^2 + \begin{bmatrix} 2 \\ -1 \\ 1 \\ 0 \\ 4 \end{bmatrix} x^3 + \begin{bmatrix} 2 \\ -1 \\ 2 \\ 3 \\ 0 \end{bmatrix} x^4 \text{ be in } V_C.$$

$$p(x) + q(x) = \begin{bmatrix} 3 \\ 0 \\ 1 \\ 0 \\ 2 \end{bmatrix} + \begin{bmatrix} 4 \\ 2 \\ 1 \\ -4 \\ 0 \end{bmatrix} + \begin{bmatrix} 1 \\ 0 \\ 0 \\ 2 \\ 0 \end{bmatrix} + \begin{bmatrix} 2 \\ 3 \\ 4 \\ 5 \\ 0 \end{bmatrix} x + \left( \begin{bmatrix} 0 \\ 1 \\ 2 \\ 3 \\ 0 \end{bmatrix} + \begin{bmatrix} 2 \\ 3 \\ 1 \\ 4 \\ 5 \end{bmatrix} \right) x^2$$

$$+ \begin{bmatrix} 2 \\ -1 \\ 1 \\ 0 \\ 4 \end{bmatrix} x^3 + \begin{bmatrix} 2 \\ -1 \\ 2 \\ 3 \\ 0 \end{bmatrix} x^4$$

$$= \begin{bmatrix} 7 \\ 2 \\ 2 \\ -4 \\ 2 \end{bmatrix} + \begin{bmatrix} 3 \\ 3 \\ 4 \\ 7 \\ 0 \end{bmatrix} x + \begin{bmatrix} 2 \\ 4 \\ 3 \\ 7 \\ 5 \end{bmatrix} x^2 + \begin{bmatrix} 2 \\ -1 \\ 1 \\ 0 \\ 4 \end{bmatrix} x^3 + \begin{bmatrix} 2 \\ -1 \\ 2 \\ 3 \\ 0 \end{bmatrix} x^4 \text{ is in } V_C.$$

Thus $V_C$ is a commutative group under addition.

We see on $V_C$ we cannot define product for it is not defined.

Thus

$$V_C = \left\{ \sum_{i=0}^{\infty} a_i x^i \;\middle|\; a_i = \begin{bmatrix} x_1 \\ x_2 \\ \vdots \\ x_n \end{bmatrix} \right\}; x_i \in Q \text{ (or R or Z)} ; 1 \leq i \leq n\}$$

is an abelian group under addition with polynomials whose coefficients are column matrices.

Now $V_{n \times m}$ denotes the collection of all polynomials whose coefficients are n×m matrices. $V_{n \times m}$ is a group under addition.

Now if $m \neq n$ then on $V_{n \times m}$ we cannot define product. We will illustrate this situation by an example.

*Example 2.12:* Let

$$V_{5 \times 3} = \left\{ \sum_{i=0}^{\infty} a_i x^i \;\middle|\; a_j = \begin{bmatrix} x_1 & x_6 & x_{11} \\ x_2 & x_7 & x_{12} \\ x_3 & x_8 & x_{13} \\ x_4 & x_9 & x_{14} \\ x_5 & x_{10} & x_{15} \end{bmatrix} \text{ where } x_i \in R; 1 \leq i \leq 15 \right\}$$

be the group of polynomials under addition whose coefficients are 5×3 matrix.

*Example 2.13:* Let

$$V_{2 \times 4} = \left\{ \sum_{i=0}^{\infty} a_i x^i \;\middle|\; a_i = \begin{bmatrix} y_1 & y_2 & y_3 & y_4 \\ y_5 & y_6 & y_7 & y_8 \end{bmatrix} \right.$$

where $y_i \in R; 1 \leq i \leq 8\}$

be the group of polynomials under addition whose coefficients are $2 \times 4$ matrices ($a_{ij} \in R; 1 \leq i \leq n, 1 \leq j \leq m$).

Thus we can say

$$V_{n \times m} = \left\{ \sum_{i=0}^{\infty} a_i x^i \,\middle|\, a_i = \begin{bmatrix} a_{11} & a_{12} & \ldots & a_{1m} \\ a_{21} & a_{22} & \ldots & a_{2m} \\ \vdots & \vdots & & \vdots \\ a_{n1} & a_{n2} & \ldots & a_{nm} \end{bmatrix} \right\}$$

is the group of polynomials in the variable x with coefficients as n × m matrices. Clearly if n ≠ m we cannot define product on $V_{n \times m}$.

Now we can define product on $V_{n \times n}$, that is when n = m. We first illustrate this by an example.

*Example 2.14:* Let

$$V_{n \times n} = \left\{ \sum_{i=0}^{\infty} a_i x^i \,\middle|\, a_i = \begin{bmatrix} a_{11} & a_{12} & \ldots & a_{1n} \\ a_{21} & a_{22} & \ldots & a_{2n} \\ \vdots & \vdots & & \vdots \\ a_{n1} & a_{n2} & \ldots & a_{nn} \end{bmatrix} \right.$$

where $a_{ij} \in R$; $1 \leq i, j \leq n$}

be the group of polynomials under addition with n × n square matrix coefficients. We see on $V_{n \times n}$, one can define product. $V_{n \times n}$ is only a semigroup which is non commutative.

We will illustrate this situation by examples.

*Example 2.15:* Let

$$V_{3 \times 3} = \left\{ \sum_{i=0}^{\infty} a_i x^i \,\middle|\, a_i = \begin{pmatrix} x_1 & x_2 & x_3 \\ x_4 & x_5 & x_6 \\ x_7 & x_8 & x_9 \end{pmatrix} \text{ where } x_i \in Q; 1 \leq i \leq 9 \right\}$$

be the group of polynomials in the variable x with coefficients from 3 × 3 matrices.

We will show how addition in $V_{3\times 3}$ is carried out.

$$\text{Let } p(x) = \begin{pmatrix} 0 & 3 & -2 \\ 1 & 0 & 0 \\ 0 & 0 & 4 \end{pmatrix} + \begin{pmatrix} 2 & 1 & 0 \\ 3 & 0 & 2 \\ 1 & 2 & 3 \end{pmatrix} x^2 + \begin{pmatrix} 0 & 1 & 2 \\ 1 & 2 & 0 \\ 2 & 1 & 0 \end{pmatrix} x^3$$

and

$$q(x) = \begin{pmatrix} 1 & 2 & 1 \\ 0 & 1 & 3 \\ -6 & 1 & 2 \end{pmatrix} + \begin{pmatrix} 1 & 2 & 3 \\ 0 & 1 & 5 \\ -5 & 0 & 1 \end{pmatrix} x + \begin{pmatrix} -1 & 2 & 3 \\ -2 & 3 & 1 \\ -3 & 2 & 1 \end{pmatrix} x^2 +$$

$$\begin{pmatrix} 0 & 1 & 0 \\ 9 & 0 & 1 \\ 0 & 2 & 3 \end{pmatrix} x^3 \text{ be in } V_{3\times 3}.$$

$$p(x) + q(x) = \begin{pmatrix} 0 & 3 & -2 \\ 1 & 0 & 0 \\ 0 & 0 & 4 \end{pmatrix} + \begin{pmatrix} 1 & 2 & 1 \\ 0 & 1 & 3 \\ -6 & 1 & 2 \end{pmatrix} + \begin{pmatrix} 1 & 2 & 3 \\ 0 & 1 & 5 \\ -5 & 0 & 1 \end{pmatrix} x$$

$$+ \left[ \begin{pmatrix} 2 & 1 & 0 \\ 3 & 0 & 2 \\ 1 & 2 & 3 \end{pmatrix} + \begin{pmatrix} -1 & 2 & 3 \\ -2 & 3 & 1 \\ -3 & 2 & 1 \end{pmatrix} \right] x^2$$

$$+ \left[ \begin{pmatrix} 0 & 1 & 2 \\ 1 & 2 & 0 \\ 2 & 1 & 0 \end{pmatrix} + \begin{pmatrix} 0 & 1 & 0 \\ 9 & 0 & 1 \\ 0 & 2 & 3 \end{pmatrix} \right] x^3$$

$$= \begin{pmatrix} 1 & 5 & -1 \\ 1 & 1 & 3 \\ -6 & 1 & 6 \end{pmatrix} + \begin{pmatrix} 1 & 2 & 3 \\ 0 & 1 & 5 \\ -5 & 0 & 1 \end{pmatrix} x +$$

$$\begin{pmatrix} 1 & 3 & 3 \\ 1 & 3 & 3 \\ -2 & 4 & 4 \end{pmatrix} x^2 + \begin{pmatrix} 0 & 2 & 2 \\ 10 & 2 & 1 \\ 2 & 3 & 3 \end{pmatrix} x^3.$$

We see $V_{3\times 3}$ is an abelian group under addition.

*Example 2.16:* Let

$$V_{2\times 2} = \left\{ \sum_{i=0}^{\infty} a_i x^i \,\middle|\, a_i = \begin{pmatrix} x_1 & x_2 \\ x_3 & x_4 \end{pmatrix}; x_i \in R; 1 \leq i \leq 4 \right\}$$

be the semigroup of polynomials in the variable x with coefficients from the collection of all 2 × 2 matrices under product.

$$p(x) = \begin{pmatrix} 1 & 2 \\ 0 & 4 \end{pmatrix} + \begin{pmatrix} 0 & 1 \\ 2 & 3 \end{pmatrix} x + \begin{pmatrix} 1 & 2 \\ 3 & 0 \end{pmatrix} x^2 \text{ and}$$

$$q(x) = \begin{pmatrix} 0 & 1 \\ 2 & 0 \end{pmatrix} + \begin{pmatrix} 1 & 0 \\ 2 & 3 \end{pmatrix} x + \begin{pmatrix} 1 & 2 \\ 3 & 4 \end{pmatrix} x^3 \text{ be in } V_{3\times 3}.$$

Now

$$p(x) \cdot q(x) = \begin{pmatrix} 1 & 2 \\ 0 & 4 \end{pmatrix}\begin{pmatrix} 0 & 1 \\ 2 & 0 \end{pmatrix} + \begin{pmatrix} 0 & 1 \\ 2 & 3 \end{pmatrix}\begin{pmatrix} 1 & 0 \\ 2 & 0 \end{pmatrix} x +$$

$$\begin{pmatrix} 1 & 2 \\ 3 & 0 \end{pmatrix}\begin{pmatrix} 0 & 1 \\ 2 & 0 \end{pmatrix} x^2 + \begin{pmatrix} 1 & 2 \\ 0 & 4 \end{pmatrix}\begin{pmatrix} 1 & 2 \\ 2 & 3 \end{pmatrix} x$$

$$+ \begin{pmatrix} 0 & 1 \\ 2 & 3 \end{pmatrix}\begin{pmatrix} 1 & 0 \\ 2 & 3 \end{pmatrix} x^2 + \begin{pmatrix} 1 & 2 \\ 3 & 0 \end{pmatrix}\begin{pmatrix} 1 & 0 \\ 2 & 3 \end{pmatrix} x^3$$

$$+ \begin{pmatrix} 1 & 2 \\ 0 & 4 \end{pmatrix}\begin{pmatrix} 1 & 2 \\ 3 & 4 \end{pmatrix} x^3 + \begin{pmatrix} 0 & 1 \\ 2 & 3 \end{pmatrix}\begin{pmatrix} 1 & 2 \\ 3 & 4 \end{pmatrix} x^4$$

$$+ \begin{pmatrix} 1 & 2 \\ 3 & 0 \end{pmatrix} \begin{pmatrix} 1 & 2 \\ 3 & 4 \end{pmatrix} x^5 + \begin{pmatrix} 4 & 1 \\ 8 & 0 \end{pmatrix} + \begin{pmatrix} 2 & 0 \\ 6 & 2 \end{pmatrix} x$$

$$+ \begin{pmatrix} 4 & 1 \\ 0 & 3 \end{pmatrix} x^2 + \begin{pmatrix} 5 & 6 \\ 8 & 12 \end{pmatrix} x + \begin{pmatrix} 2 & 3 \\ 8 & 9 \end{pmatrix} x^2$$

$$+ \begin{pmatrix} 5 & 6 \\ 3 & 0 \end{pmatrix} x^3 + \begin{pmatrix} 7 & 10 \\ 12 & 16 \end{pmatrix} x^3 + \begin{pmatrix} 3 & 4 \\ 11 & 16 \end{pmatrix} x^4 + \begin{pmatrix} 7 & 10 \\ 3 & 6 \end{pmatrix} x^5$$

$$= \begin{pmatrix} 4 & 1 \\ 8 & 0 \end{pmatrix} + \begin{pmatrix} 7 & 6 \\ 14 & 14 \end{pmatrix} x + \begin{pmatrix} 6 & 4 \\ 8 & 12 \end{pmatrix} x^2 + \begin{pmatrix} 12 & 16 \\ 15 & 16 \end{pmatrix} x^3 +$$

$$\begin{pmatrix} 3 & 4 \\ 11 & 16 \end{pmatrix} x^4 + \begin{pmatrix} 7 & 10 \\ 3 & 6 \end{pmatrix} x^5.$$

This is the way product is defined. Thus $V_{2\times 2}$ is a semigroup under multiplication.

$V_{2\times 2}$ is a monoid and infact $V_{2\times 2}$ has zero divisors.

This is a polynomial ring.

*Example 2.17:* Let

$$V_{4\times 4} = \left\{ \sum_{i=0}^{\infty} a_i x^i \,\middle|\, a_i = \begin{pmatrix} a_{11} & a_{12} & a_{13} & a_{14} \\ a_{21} & a_{22} & a_{23} & a_{24} \\ a_{31} & a_{32} & a_{33} & a_{34} \\ a_{41} & a_{42} & a_{43} & a_{44} \end{pmatrix} \right.$$

where $a_{ij} \in R; 1 \leq i, j \leq 4\}$

be a group of polynomials in the variable x with $4 \times 4$ matrices as coefficients.

$V_{4\times 4}$ is a group under addition and $V_{4\times 4}$ is a semigroup under product ($V_{4\times 4}$, +, ×) is a ring which is non commutative. This ring has zero divisors and units and all p(x) of degree greater than or equal to one have no inverse.

*Example 2.18:* Let

$$V_{2\times 2} = \left\{ \sum_{i=0}^{\infty} a_i x^i \,\middle|\, a_i = \begin{pmatrix} a_{11} & a_{12} \\ a_{21} & a_{22} \end{pmatrix}; a_{ij} \in R; 1 \le i, j \le 2 \right\}$$

be the ring of polynomials with 2×2 matrix coefficients in the variable x. $V_{2\times 2}$ is non commutative and has zero divisors and no $p(x) \in V_{2\times 2}$, of degree greater than one has inverse. We cannot have idempotent in them.

We can differentiate and integrate these polynomials with matrix coefficients apart from finding roots in them.

Now we first illustrate this situation by some examples.

*Example 2.19:* Let

$$p(x) = \begin{bmatrix} 3 & 0 \\ 1 & 2 \end{bmatrix} + \begin{bmatrix} 2 & 6 \\ 1 & 5 \end{bmatrix} x + \begin{bmatrix} 7 & 0 \\ 0 & 8 \end{bmatrix} x^2 - \begin{bmatrix} 3 & 1 \\ 0 & 0 \end{bmatrix} x^3$$

$$+ \begin{bmatrix} 8 & 1 \\ 0 & 1 \end{bmatrix} x^4 - \begin{bmatrix} 0 & 4 \\ -2 & 0 \end{bmatrix} x^5$$

be a polynomial in matrix coefficients or matrix coefficient polynomial.

To find the derivative of p(x).

$$\frac{dp(x)}{dx} = 0 + \begin{bmatrix} 2 & 6 \\ 1 & 5 \end{bmatrix} + 2\begin{bmatrix} 7 & 0 \\ 0 & 8 \end{bmatrix}x - 3\begin{bmatrix} 3 & 1 \\ 0 & 0 \end{bmatrix}x^2$$

$$+ 4\begin{bmatrix} 8 & 1 \\ 0 & 1 \end{bmatrix}x^3 - 5\begin{bmatrix} 0 & 4 \\ -2 & 0 \end{bmatrix}x^4$$

$$= \begin{bmatrix} 2 & 6 \\ 1 & 5 \end{bmatrix} + \begin{bmatrix} 14 & 0 \\ 0 & 16 \end{bmatrix}x - \begin{bmatrix} 9 & 3 \\ 0 & 0 \end{bmatrix}x^2$$

$$+ \begin{bmatrix} 32 & 4 \\ 0 & 4 \end{bmatrix}x^3 - \begin{bmatrix} 0 & 20 \\ -10 & 0 \end{bmatrix}x^4.$$

We see $\frac{dp(x)}{dx}$ is again a matrix coefficient polynomial in the variable x.

We can find the second derivative of p(x).

Consider

$$\frac{d^2p(x)}{dx} = \begin{bmatrix} 14 & 0 \\ 0 & 16 \end{bmatrix} - 2\begin{bmatrix} 9 & 3 \\ 0 & 0 \end{bmatrix}x$$

$$+ 3\begin{bmatrix} 32 & 4 \\ 0 & 4 \end{bmatrix}x^2 - 4\begin{bmatrix} 0 & 20 \\ -10 & 0 \end{bmatrix}x^3$$

$$= \begin{bmatrix} 14 & 0 \\ 0 & 16 \end{bmatrix} - \begin{bmatrix} 18 & 6 \\ 0 & 0 \end{bmatrix}x + \begin{bmatrix} 96 & 12 \\ 0 & 12 \end{bmatrix}x^2 - \begin{bmatrix} 0 & 80 \\ -40 & 0 \end{bmatrix}x^3.$$

Clearly $\frac{d^2p(x)}{dx}$ also belongs to the collection of $2 \times 2$ matrix coefficient polynomials.

*Example 2.20:* Let

$$V_{2\times 4} = \left\{ \sum_{i=0}^{\infty} a_i x^i \;\middle|\; a_i = \begin{pmatrix} x_1 & x_2 & x_3 & x_4 \\ x_5 & x_6 & x_7 & x_8 \end{pmatrix} \right.$$

where $x_i \in R; 1 \le i \le 8\}$

be the $2 \times 4$ matrix coefficient polynomial.

Let $p(x) = \begin{pmatrix} 1 & 0 & 2 & 4 \\ 0 & 3 & 0 & 5 \end{pmatrix} + \begin{pmatrix} 3 & 1 & 5 & 2 \\ 0 & 4 & 0 & 5 \end{pmatrix} x$

$+ \begin{pmatrix} -3 & 4 & 2 & 4 \\ 0 & 0 & 0 & 3 \end{pmatrix} x^3 + \begin{pmatrix} 1 & -1 & 0 & 2 \\ 2 & 0 & -2 & 0 \end{pmatrix} x^4$

be a $2 \times 4$ matrix coefficient polynomial. To find the derivative of $p(x)$.

$$\frac{dp(x)}{dx} = \begin{pmatrix} 3 & 1 & 5 & 2 \\ 0 & 4 & 0 & 5 \end{pmatrix} + 3 \begin{pmatrix} -3 & 4 & 2 & 4 \\ 0 & 0 & 0 & 3 \end{pmatrix} x^2$$

$$+ 4 \begin{pmatrix} 1 & -1 & 0 & 2 \\ 2 & 0 & -2 & 0 \end{pmatrix} x^3$$

$$= \begin{pmatrix} 3 & 1 & 5 & 2 \\ 0 & 4 & 0 & 5 \end{pmatrix} + \begin{pmatrix} -9 & 12 & 6 & 12 \\ 0 & 0 & 0 & 9 \end{pmatrix}$$

$$+ \begin{pmatrix} 4 & -4 & 0 & 8 \\ 8 & 0 & -8 & 0 \end{pmatrix} x^3.$$

Clearly $\dfrac{dp(x)}{dx}$ is in $V_{2\times 4}$.

Consider

$$\frac{d^2p(x)}{dx^2} = 2\begin{pmatrix} -9 & 12 & 6 & 12 \\ 0 & 0 & 0 & 9 \end{pmatrix} x + 3\begin{pmatrix} 4 & -4 & 0 & 8 \\ 8 & 0 & -8 & 0 \end{pmatrix} x^2$$

$$= \begin{pmatrix} -18 & 24 & 12 & 24 \\ 0 & 0 & 0 & 18 \end{pmatrix} x + \begin{pmatrix} 12 & -12 & 0 & 24 \\ 24 & 0 & -24 & 0 \end{pmatrix} x^2.$$

We see $\dfrac{d^2p(x)}{dx^2} \in V_{2\times 4}$.

If we consider the third derivative of p(x);

$$\frac{d^3p(x)}{dx^3} = \begin{pmatrix} -18 & 24 & 12 & 24 \\ 0 & 0 & 0 & 18 \end{pmatrix} +$$

$$2\begin{pmatrix} 12 & -12 & 0 & 24 \\ 24 & 0 & -24 & 0 \end{pmatrix} x$$

$$= \begin{pmatrix} -18 & 24 & 0 & 24 \\ 0 & 0 & 0 & 18 \end{pmatrix} + \begin{pmatrix} 24 & -24 & 0 & 48 \\ 48 & 0 & -48 & 0 \end{pmatrix} x.$$

We see $\dfrac{d^3p(x)}{dx^3} \in V_{2\times 4}$.

Further the forth derivative.

$$\frac{d^4p(x)}{dx^4} = \begin{pmatrix} 24 & -24 & 0 & 48 \\ 48 & 0 & -48 & 0 \end{pmatrix} \in V_{2\times 4}.$$

However the fifth derivative $\dfrac{d^5p(x)}{dx^5}$ is zero.

*Example 2.21:* Let

$$V_R = \left\{ \sum_{i=0}^{\infty} a_i x^i \,\middle|\, a_i = (x_1, \ldots, x_6); x_i \in Z; 1 \le i \le 6 \right\}$$

be a row matrix coefficient polynomial.

$p(x) = (2,0,1,0,1,5) + (3,2,1,0,0,0)x + (0,1,0,2,0,4)x^2$
$\qquad + (0,-2,-3,0,0,0)x^3 + (8,0,7,0,1,0)x^5$ be in $V_R$.

To find the derivative of

$$p(x) = \frac{dp(x)}{dx}$$

$= 0 + (3,2,1,0,0,0) + 2(0,1,0,2,0,4)x + 3(0,-2,-3,0,0,0)x^2 + 5(8,0,7,0,1,0)x^4$

$= (3,2,1,0,0,0) + (0,2,0,4,0,8)x + (0,-6,-9,0,0,0)x^2 + (40,0,35,0,5,0)x^4$.

We see $\dfrac{dp(x)}{dx}$ is in $V_R$.

$\dfrac{d^2p(x)}{dx^2} = (0,2,0,4,0,8) + 2(0,-6,-9,0,0,0)x$
$\qquad\qquad + 4(40,0,35,0,5,0)x^3$

$= (0,2,0,4,0,8) + (0,-12,-18,0,0,0)x + (160,0,140,0,20,0)x^3$.

Clearly $\dfrac{d^2p(x)}{dx^2} \in V_R$.

*Example 2.22:* Let

$$V_C = \left\{ \sum_{i=0}^{\infty} a_i x^i \,\middle|\, a_i = \begin{bmatrix} x_1 \\ x_2 \\ x_3 \\ x_4 \end{bmatrix} \text{ where } x_i \in Q;\ 1 \le i \le 4 \right\}$$

be a 4 × 1 column matrix coefficient polynomial.

Let $p(x) = \begin{bmatrix} 2 \\ 0 \\ 4 \\ 0 \end{bmatrix} + \begin{bmatrix} 3 \\ 2 \\ 1 \\ -4 \end{bmatrix} x + \begin{bmatrix} 0 \\ 1 \\ 2 \\ 3 \end{bmatrix} x^3 + \begin{bmatrix} 4 \\ 5 \\ 2 \\ 1 \end{bmatrix} x^6$ belongs to $V_C$.

$$\frac{dp(x)}{dx} = \begin{bmatrix} 3 \\ 2 \\ 1 \\ -4 \end{bmatrix} + 3\begin{bmatrix} 0 \\ 1 \\ 2 \\ 3 \end{bmatrix} x^2 + 6\begin{bmatrix} 4 \\ 5 \\ 2 \\ 1 \end{bmatrix} x^5$$

$$= \begin{bmatrix} 3 \\ 2 \\ 1 \\ -4 \end{bmatrix} + \begin{bmatrix} 0 \\ 3 \\ 6 \\ 9 \end{bmatrix} x^2 + \begin{bmatrix} 24 \\ 30 \\ 12 \\ 6 \end{bmatrix} x^5 \in V_C.$$

$$\frac{d^2 p(x)}{dx^2} = 2\begin{bmatrix} 0 \\ 3 \\ 6 \\ 9 \end{bmatrix} x + 5\begin{bmatrix} 24 \\ 30 \\ 12 \\ 6 \end{bmatrix} x^4$$

$$= \begin{bmatrix} 0 \\ 6 \\ 12 \\ 18 \end{bmatrix} x + \begin{bmatrix} 120 \\ 150 \\ 60 \\ 30 \end{bmatrix} x^4 \in V_C.$$

$$\frac{d^3 p(x)}{dx^3} = \begin{bmatrix} 0 \\ 6 \\ 12 \\ 18 \end{bmatrix} + 4 \begin{bmatrix} 120 \\ 150 \\ 60 \\ 30 \end{bmatrix} x^3$$

$$= \begin{bmatrix} 0 \\ 6 \\ 12 \\ 18 \end{bmatrix} + \begin{bmatrix} 480 \\ 600 \\ 240 \\ 120 \end{bmatrix} x^3 \in V_C.$$

$$\frac{d^4 p(x)}{dx^4} = 3 \begin{bmatrix} 480 \\ 600 \\ 240 \\ 120 \end{bmatrix} x^2 = \begin{bmatrix} 1440 \\ 1800 \\ 720 \\ 360 \end{bmatrix} x^2 \in V_C.$$

Thus we see $V_C$, $V_{m \times n}$, $V_{n \times n}$ and $V_R$ are such that the first derivative and all higher derivatives are in $V_C$, $V_{m \times n}$, $V_{n \times n}$ and $V_R$.

Now we discuss about the integration of matrix coefficient polynomials.

***Example 2.23:*** Let

$$p(x) = \begin{bmatrix} 3 & 0 & 1 \\ 5 & 6 & 0 \\ 1 & 0 & 8 \end{bmatrix} + \begin{bmatrix} 0 & 2 & 1 \\ 6 & 1 & 0 \\ 1 & 2 & 6 \end{bmatrix} x$$

$$+ \begin{bmatrix} 8 & 0 & 0 \\ 0 & 7 & 0 \\ 0 & 0 & 11 \end{bmatrix} x^2 + \begin{bmatrix} 0 & 0 & 2 \\ 0 & 9 & 0 \\ 10 & 0 & 0 \end{bmatrix} x^3.$$

To integrate p(x) . $\int p(x)dx = \begin{bmatrix} 3 & 0 & 1 \\ 5 & 6 & 0 \\ 1 & 0 & 8 \end{bmatrix} x + \frac{1}{2} \begin{bmatrix} 0 & 2 & 1 \\ 6 & 1 & 0 \\ 1 & 2 & 6 \end{bmatrix} x^2 +$

$$1/3 \begin{bmatrix} 8 & 0 & 0 \\ 0 & 7 & 0 \\ 0 & 0 & 11 \end{bmatrix} x^3 + 1/4 \begin{bmatrix} 0 & 0 & 2 \\ 0 & 9 & 0 \\ 10 & 0 & 0 \end{bmatrix} x^4 +$$

$$\begin{bmatrix} a_1 & a_2 & a_3 \\ a_4 & a_5 & a_6 \\ a_7 & a_8 & a_9 \end{bmatrix} \begin{bmatrix} 3 & 0 & 1 \\ 5 & 6 & 0 \\ 1 & 0 & 8 \end{bmatrix} x + \begin{bmatrix} 0 & 1 & 1/2 \\ 3 & 1/2 & 0 \\ 1/2 & 1 & 3 \end{bmatrix} x^2 +$$

$$\begin{bmatrix} 8/3 & 0 & 0 \\ 0 & 7/3 & 0 \\ 0 & 0 & 11/3 \end{bmatrix} x^3 + \begin{bmatrix} 0 & 0 & 1/2 \\ 0 & 9/4 & 0 \\ 5/2 & 0 & 0 \end{bmatrix} x^4 + \begin{bmatrix} a_1 & a_2 & a_3 \\ a_4 & a_5 & a_6 \\ a_7 & a_8 & a_9 \end{bmatrix}.$$

***Example 2.24:*** Let

$$p(x) = (1,2,3,4,5) + (0,1,0,3,-1)x + (5,0,8,1,7)x^2$$
$$+ (1,2,0,4,5)x^3 + (-2, 1,4,3,0)x^4$$

be a row matrix coefficient polynomial.

To integrate p(x), $\int p(x)dx$
= $(1,2,3,4,5)x + \frac{1}{2}(0,1,0,3,-1)x^2 + 1/3(5,0,8,1,7)x^3$
  $+ 1/4(1,2,0,4,5)x^4 + 1/5(-2,1,4,3,0)x^5 + (a_1,a_2,a_3,a_4,a_5)$
  $a_i \in Q; 1 \leq i \leq 5$.

= $(1,2,3,4,5) + (0,1/2,0,3/2,-1/2)x^2 + (5/3,0,8/3,1/3,7/3)x^3$
  $+ (1/4,1/2,0,1,5/4)x^4 + (-2/5, 1/5,4/5,3/5,0)x^5$
  $+ (a_1,a_2,a_3,a_4,a_5)$.

*Example 2.25:* Let

$$p(x) = \begin{bmatrix} 3 \\ 0 \\ 1 \\ 2 \\ 4 \\ 5 \end{bmatrix} + \begin{bmatrix} 0 \\ 1 \\ 2 \\ 0 \\ 4 \\ 8 \end{bmatrix} x + \begin{bmatrix} -1 \\ 0 \\ -9 \\ 8 \\ 7 \\ 0 \end{bmatrix} x^3 + \begin{bmatrix} 7 \\ 8 \\ 9 \\ 10 \\ 3 \\ 7 \end{bmatrix} x^4 + \begin{bmatrix} 8 \\ 2 \\ 4 \\ 5 \\ 5 \\ 10 \end{bmatrix} x^5$$

be a column matrix polynomial.

$$\int p(x)dx = \begin{bmatrix} 3 \\ 0 \\ 1 \\ 2 \\ 4 \\ 5 \end{bmatrix} x + 1/2 \begin{bmatrix} 0 \\ 1 \\ 2 \\ 0 \\ 4 \\ 8 \end{bmatrix} x^2 + 1/4 \begin{bmatrix} -1 \\ 0 \\ -9 \\ 8 \\ 7 \\ 0 \end{bmatrix} x^4$$

$$+ 1/5 \begin{bmatrix} 7 \\ 8 \\ 9 \\ 10 \\ 3 \\ 7 \end{bmatrix} x^5 + 1/6 \begin{bmatrix} 8 \\ 2 \\ 4 \\ 5 \\ 5 \\ 10 \end{bmatrix} x^6 + \begin{bmatrix} a_1 \\ a_2 \\ a_3 \\ a_4 \\ a_5 \\ a_6 \end{bmatrix}$$

$$= \begin{bmatrix} 3 \\ 0 \\ 1 \\ 2 \\ 4 \\ 5 \end{bmatrix} x + \begin{bmatrix} 0 \\ 1/2 \\ 1 \\ 0 \\ 2 \\ 4 \end{bmatrix} x^2 + \begin{bmatrix} -1/4 \\ 0 \\ -9/4 \\ 2 \\ 7/4 \\ 0 \end{bmatrix} x^4 + \begin{bmatrix} 7/5 \\ 8/5 \\ 9/5 \\ 2 \\ 3/5 \\ 7/5 \end{bmatrix} x^5 + \begin{bmatrix} 4/3 \\ 1/3 \\ 2/3 \\ 5/6 \\ 5/6 \\ 5/3 \end{bmatrix} x^6 + \begin{bmatrix} a_1 \\ a_2 \\ a_3 \\ a_4 \\ a_5 \\ a_6 \end{bmatrix}.$$

*Example 2.26:* Let

$$p(x) = \begin{bmatrix} 0 & 2 & 1 & 4 \\ 6 & 0 & 1 & 0 \end{bmatrix} + \begin{bmatrix} 3 & 6 & 2 & 9 \\ 0 & 2 & 1 & 7 \end{bmatrix} x + \begin{bmatrix} 0 & 2 & 4 & 4 \\ 2 & 0 & 1 & 2 \end{bmatrix} x^3$$

$$+ \begin{bmatrix} 2 & 1 & 0 & 0 \\ 0 & 0 & 1 & 2 \end{bmatrix} x^4 + \begin{bmatrix} 0 & 1 & 2 & 0 \\ 6 & 0 & 0 & 3 \end{bmatrix} x^5.$$

We find the integral of p(x).

$$\int p(x)dx = \begin{bmatrix} 0 & 2 & 1 & 4 \\ 6 & 0 & 1 & 0 \end{bmatrix} x + 1/2 \begin{bmatrix} 3 & 6 & 2 & 9 \\ 0 & 2 & 1 & 7 \end{bmatrix} x^2$$

$$+ 1/4 \begin{bmatrix} 0 & 2 & 4 & 4 \\ 2 & 0 & 1 & 2 \end{bmatrix} x^4 + 1/5 \begin{bmatrix} 2 & 1 & 0 & 0 \\ 0 & 0 & 1 & 2 \end{bmatrix} x^5$$

$$+ 1/6 \begin{bmatrix} 0 & 1 & 2 & 0 \\ 6 & 0 & 0 & 3 \end{bmatrix} x^6 \begin{bmatrix} a_1 & a_2 & a_3 & a_4 \\ a_5 & a_6 & a_7 & a_8 \end{bmatrix}$$

$$= \begin{bmatrix} 0 & 2 & 1 & 4 \\ 6 & 0 & 1 & 0 \end{bmatrix} x + \begin{bmatrix} 3/2 & 3 & 1 & 9/2 \\ 0 & 1 & 1/2 & 7/2 \end{bmatrix} x^2$$

$$+ \begin{bmatrix} 0 & 1/2 & 1 & 1 \\ 1/2 & 0 & 1/4 & 1/2 \end{bmatrix} x^4 + \begin{bmatrix} 2/5 & 1/5 & 0 & 0 \\ 0 & 0 & 1/5 & 2/5 \end{bmatrix} x^5$$

$$+ \begin{bmatrix} 0 & 1/6 & 1/3 & 0 \\ 1 & 0 & 0 & 1/2 \end{bmatrix} x^6 + \begin{bmatrix} a_1 & a_2 & a_3 & a_4 \\ a_5 & a_6 & a_7 & a_8 \end{bmatrix}.$$

We see $V_C$, $V_R$, $V_{n \times m}$ and $V_{n \times n}$ under integration is closed, provided the entries of the coefficient matrices take their values from Q or R. If they take the from Z they are not closed under integration only closed under differentiation.

We will illustrate this situation by a few examples.

*Example 2.27:* Let

$p(x) = (3, 8, 4, 0) + (2, 0, 4, 9)x + (1, 2, 1, 1)x^2 + (1, 0, 1, 1)x^3 + (3, 4, 8, 9)x^5$ where the coefficients are $1 \times 4$ row matrices with entries from Z.

We find integral of $p(x)$.
$\int p(x)dx = (3, 8, 4, 0)x + 1/2(2, 0, 4, 9)x^2 + 1/3(1, 2, 1, 1)x^3$
$\qquad + 1/4(1, 0, 1, 1)x^4 + 1/6(3, 4, 8, 9)x^6.$

We see $(1, 0, 2, 9/4)$, $(1/3, 2/3, 1/3, 1/3)$, $(1/4, 0, 1/4, 1/4)$, $(1/2, 2/3, 4/3, 3/2) \notin Z \times Z \times Z \times Z$. Thus we see integral of matrix coefficient polynomials with matrix entries from Z are not closed under intervals that is $\int p(x)dx \notin V_C$ or $V_R$ or $V_{n \times n}$ or $V_{m \times n}$ if the entries are in Z.

*Example 2.28:* Let

$$p(x) = \begin{bmatrix} 2 \\ 3 \\ 4 \\ 0 \end{bmatrix} + \begin{bmatrix} 1 \\ 2 \\ 3 \\ 4 \end{bmatrix} x + \begin{bmatrix} 0 \\ 0 \\ 1 \\ 1 \end{bmatrix} x^2 + \begin{bmatrix} 0 \\ 1 \\ 0 \\ 3 \end{bmatrix} x^3 + \begin{bmatrix} 3 \\ 0 \\ 0 \\ 4 \end{bmatrix} x^4 \text{ where } \begin{bmatrix} a_1 \\ a_2 \\ a_3 \\ a_4 \end{bmatrix}$$

are $4 \times 1$ column matrix with entries from Z; that is $a_i \in Z$; $1 \leq i \leq 4$.

$$\int p(x)\,dx = \begin{bmatrix} 2 \\ 3 \\ 4 \\ 0 \end{bmatrix} x + 1/2 \begin{bmatrix} 1 \\ 2 \\ 3 \\ 4 \end{bmatrix} x^2 + 1/3 \begin{bmatrix} 0 \\ 0 \\ 1 \\ 1 \end{bmatrix} x^3$$

$$+ 1/4 \begin{bmatrix} 0 \\ 1 \\ 0 \\ 3 \end{bmatrix} x^4 + 1/5 \begin{bmatrix} 3 \\ 0 \\ 0 \\ 4 \end{bmatrix} x^5 + \begin{bmatrix} a_1 \\ a_2 \\ a_3 \\ a_4 \end{bmatrix}$$

$$= \begin{bmatrix} 2 \\ 3 \\ 4 \\ 0 \end{bmatrix} x + \begin{bmatrix} 1/2 \\ 1 \\ 3/2 \\ 2 \end{bmatrix} x^2 + \begin{bmatrix} 0 \\ 0 \\ 1/3 \\ 1/3 \end{bmatrix} x^3 + \begin{bmatrix} 0 \\ 1/4 \\ 0 \\ 3/4 \end{bmatrix} x^4$$

$$+ \begin{bmatrix} 3/5 \\ 0 \\ 0 \\ 4/5 \end{bmatrix} x^5 + \begin{bmatrix} a_1 \\ a_2 \\ a_3 \\ a_4 \end{bmatrix}.$$

Clearly these column matrices do not take their entries from Z.

*Example 2.29:* Let

$$p(x) = \begin{bmatrix} 3 & 0 & 1 & 1 \\ 6 & 6 & 0 & 0 \end{bmatrix} + \begin{bmatrix} 1 & 2 & 0 & 0 \\ 0 & 1 & 0 & 2 \end{bmatrix} x + \begin{bmatrix} 3 & 0 & 0 & 1 \\ 0 & 2 & 2 & 0 \end{bmatrix} x^2$$

$$+ \begin{bmatrix} 2 & 1 & 1 & 0 \\ 0 & 2 & 0 & 1 \end{bmatrix} x^3 + \begin{bmatrix} 1 & 0 & 1 & 0 \\ 0 & 1 & 0 & 1 \end{bmatrix} x^4$$

where the coefficient of these are $2 \times 4$ matrices and they take their values from Z.

$$\int p(x)dx = \begin{bmatrix} 3 & 0 & 1 & 1 \\ 6 & 6 & 0 & 0 \end{bmatrix} x + 1/2 \begin{bmatrix} 1 & 2 & 0 & 0 \\ 0 & 1 & 0 & 2 \end{bmatrix} x^2$$

$$+ 1/3 \begin{bmatrix} 3 & 0 & 0 & 1 \\ 0 & 2 & 2 & 0 \end{bmatrix} x^3 + 1/4 \begin{bmatrix} 2 & 1 & 1 & 0 \\ 0 & 2 & 0 & 1 \end{bmatrix} x^4$$

$$+ 1/5 \begin{bmatrix} 1 & 0 & 1 & 0 \\ 0 & 1 & 0 & 1 \end{bmatrix} x^5 + \begin{bmatrix} a_1 & a_2 & a_3 & a_4 \\ a_5 & a_6 & a_7 & a_8 \end{bmatrix} \quad a_i \in Z; \ 1 \le i \le 8.$$

$$\int p(x)dx = \begin{bmatrix} 3 & 0 & 1 & 1 \\ 6 & 6 & 0 & 0 \end{bmatrix} x + \begin{bmatrix} 1/2 & 1 & 0 & 0 \\ 0 & 1/2 & 0 & 1 \end{bmatrix} x^2$$

$$+ \begin{bmatrix} 1 & 0 & 0 & 1/3 \\ 0 & 2/3 & 2/3 & 0 \end{bmatrix} x^3 + \begin{bmatrix} 1/2 & 1/4 & 1/4 & 0 \\ 0 & 1/2 & 0 & 1/4 \end{bmatrix} x^4$$

$$+ \begin{bmatrix} 1/5 & 0 & 1/5 & 0 \\ 0 & 1/5 & 0 & 1/5 \end{bmatrix} x^5 + \begin{bmatrix} a_1 & a_2 & a_3 & a_4 \\ a_5 & a_6 & a_7 & a_8 \end{bmatrix}.$$

In view of this we have the following theorems.

**THEOREM 2.1:** *Let $V_R$ (or $V_C$ or $V_{n \times m}$ or $V_{m \times m}$) be the matrix coefficient polynomials with matrix entries from C or Z or R or Q. The derivatives of every polynomial in $V_R$ (or $V_C$ or $V_{n \times m}$ or $V_{m \times m}$) is in $V_R$ (or $V_C$ or $V_{n \times m}$ or $V_{m \times m}$).*

The proof is simple and hence is left as an exercise to the reader.

**THEOREM 2.2:** *Let $V_R$ (or $V_C$ or $V_{n \times m}$ or $V_{m \times m}$) be the matrix coefficient polynomial with matrix entries from Z. The integrals of every matrix coefficient polynomial need not be in $V_R$.*

**COROLLARY 1:** *If in theorem, Z is replaced by Q or R or C then every integral of the matrix coefficient polynomial is in $V_R$ (or $V_C$ or $V_{n \times m}$ or $V_{m \times m}$).*

Now we find or show some polynomial identities true in case of matrix coefficient polynomials.

Consider $(1, 1, 1, 1, 1)x^2 - (4, 9, 16, 25, 81) = (0)$ in $V_R = \left\{ \sum_{i=0}^{\infty} a_i x^i \;\middle|\; a = (x_1, x_2, x_3, x_4, x_5) \text{ with } x_i \in Z \text{ or } Q \text{ or } C \text{ or } R; 1 \leq i \leq 5 \right\}$. Given $(1, 1, 1, 1, 1)x^2 - (4, 9, 16, 25, 81) = (0)$.

$((1, 1, 1, 1, 1)x - (2, 3, 4, 5, 9)) \times ((1, 1, 1, 1, 1)x + (2, 3, 4, 5, 9)) = (0)$

Thus $x = (2,3,4,5,9)$ or $-(2,3,4,5,9)$.

Take the matrix coefficient polynomial
$$(1,1,1)x^3 - (27,8,125) = (0)$$

We can factorize $(1,1,1)x^3 - (27,8,125) = 0$ as
$[(1,1,1)x - (3,2,5)] [(1,1,1)x^2 + (3,2,5)x + (9,4,25)]$.

Take $(1,1,1,1)x^4 - (16,81,625,16) = (0)$

We can factorize this polynomial as $[(1,1,1,1)x^2 + (4,9,25,4)] [(1,1,1,1)x^2 - (4,9,25,4)] = (0,0,0,0)$

$x^2 = -(4,9,25,4)$

and $x^2 = (4,9,25,4)$, we see now $x^2 = (4,9,25,4)$ can be yet solved as $x = \pm (2,3,5,2)$, we see however $x^2 = -(4,9,25,4)$ gives a imaginary value for x. If $V_R$ is defined over R or Z or Q we see the solution does not exist; that is the equation is not linearly solvable over R or Z or Q but linearly solvable over $V_C$.

Now we see yet another equation
$p(x) = (1,1,1,1)x^2 + (4,4,4,4)x + (4,4,4,4) = (0)$ where $p(x)$ is a matrix coefficient polynomial in the variable x over Z.

$$x = \frac{-(4,4,4,4) \pm \sqrt{(4,4,4,4)^2 - 4 \times (1,1,1,1)(4,4,4,4)}}{(2,2,2,2)}$$

$$= \frac{-(4,4,4,4) \pm \sqrt{(0)}}{(2,2,2,2)}$$

$$= \frac{-(4,4,4,4)}{(2,2,2,2)} = -(2,2,2,2).$$

Thus $p(x)$ has coincident roots.

Consider $(1,1,1)x^3 - (6,3,9)x^2 + (12,3,27)x + (8,1,27)$

$= p(x) = (0,0,0)$ be a matrix coefficient polynomial.

To find the roots of $p(x)$.

$p(x) = (1,1,1)x^3 - 3(2,1,3)x^2 + 3(4,1,9)x - (8,1,27)$

$= ((1,1,1)x - (2,1,3))^3.$

Thus $x = (2,1,3), (2,1,3)$ and $(2,1,3)$.

Now $p(2,1,3) = (1,1,1)(2,1,3)^3 - 3(2,1,3)(2,1,3)^2 + 3(2,1,3)^2 (2,1,3) - (2,1,3)^3$

$= (0,0,0).$

We can also find equation with matrix coefficient polynomials as follows:

Consider

$$\left[\begin{pmatrix}1 & 0\\0 & 1\end{pmatrix}x - \begin{pmatrix}2 & 1\\0 & 2\end{pmatrix}\right]\left[\begin{pmatrix}1 & 0\\0 & 1\end{pmatrix}x + \begin{pmatrix}3 & 7\\0 & 1\end{pmatrix}\right] = p(x).$$

Clearly $p\begin{pmatrix}2 & 1\\0 & 2\end{pmatrix} = (0)$ and $p\begin{pmatrix}3 & 7\\0 & 1\end{pmatrix} = (0)$.

We can consider any product of linear polynomial with matrix coefficients. However we see it is difficult to solve equations in the matrix coefficients as even solving equations in usual polynomials is not an easy problem.

Now having seen the properties of matrix coefficients polynomials we now proceed onto discuss other properties associated with it.

**Chapter Three**

# ALGEBRAIC STRUCTURES USING MATRIX COEFFICIENT POLYNOMIALS

In this chapter we introduce several types of algebraic structures on these matrix coefficient polynomials and study them.

Throughout this chapter $V_R$ denotes the collection of all row matrix coefficient polynomials. $V_R = \left\{ \sum_{i=0}^{\infty} a_i x^i \mid a_i = (y_1, \ldots, y_n) \right.$ where $y_i \in R$ (or Q or C or Z); $1 \leq i \leq n$ and x an indeterminate$\}$.

$V_C$ denotes the collection of all column matrix coefficient polynomials; that is $V_C = \left\{ \sum_{i=0}^{\infty} a_i x^i \mid a_i = \begin{bmatrix} x_1 \\ x_2 \\ \vdots \\ x_m \end{bmatrix} ; x_j \in R \text{ (or Q or C or Z) } 1 \leq j \leq m \right\}$.

Now $V_{n \times m} = \left\{ \sum_{i=0}^{\infty} a_i x^i \,\middle|\, a_k = \begin{pmatrix} a_{11} & \cdots & a_{1m} \\ a_{21} & \cdots & a_{2m} \\ \vdots & \vdots & \vdots \\ a_{n1} & \cdots & a_{nm} \end{pmatrix} a_{ij} \in R \text{ (or Q or } \right.$

Z or C); $1 \le i \le n$; $1 \le j \le m$} denotes the collection of all $n \times m$ matrix coefficient polynomial.

Finally $V_{n \times n} = \left\{ \sum_{i=0}^{\infty} a_i x^i \,\middle|\, a_k = \begin{pmatrix} a_{11} & a_{12} & \cdots & a_{1n} \\ a_{21} & a_{22} & \cdots & a_{2n} \\ \vdots & \vdots & & \vdots \\ a_{n1} & a_{n2} & \cdots & a_{nn} \end{pmatrix} \right.$ with $a_{ij}$

$\in$ R (or Q or Z or C); $1 \le i, j \le n$} denotes the collection of all $n \times n$ matrix coefficient polynomial.

We give algebraic structures on them.

**THEOREM 3.1:** *$V_R$, $V_C$, $V_{n \times m}$ and $V_{n \times n}$ ($m \ne n$) are groups under addition.*

**THEOREM 3.2:** *$V_R$ and $V_{n \times n}$ are semigroups (monoid) under multiplication.*

**THEOREM 3.3:** *$V_R$ and $V_{n \times n}$ are rings*

*(i) $V_R$ is a commutative ring.*
*(ii) $V_{n \times n}$ is a non commutative ring.*

The proof of all these theorems are simple and hence left as an exercise to the reader.

**THEOREM 3.4:** *Both $V_R$ and $V_{n \times n}$ have zero divisors.*

**THEOREM 3.5:** *Both $V_R$ and $V_{n \times n}$ have no idempotents which are not constant matrix coefficient polynomials.*

We give examples of zero divisors.

***Example 3.1:*** Let $V_R[x] = \left\{ \sum_{i=0}^{\infty} a_i x^i \mid a_i = (x_1, x_2, x_3, x_4, x_5); x_j \in Q \text{ or } R \text{ or } Z, 1 \leq j \leq 5 \right\}$ be a matrix coefficient polynomial ring.

Take $p(x) = (3,2,0,0,0) + (6,3,0,0,0)x + (7,0,0,0,0)x^2 + (8,1,0,0,0)x^4$ and

$q(x) = (0,0,1,2,3) + (0,0,0,4,2)x^2 + (0,0,0,1,4)x^3 + (0,0,0,3,4)x^4 + (0,0,0,5,2)x^7$

be elements in $V_R$.

$p(x) q(x) = (0,0,0,0,0)$.

Thus $V_R$ has zero divisors.

Consider
$a(x) = (5,0,0,0,2) + (3,0,0,0,0)x + (0,0,0,0,7)x^2 + (2,0,0,0,-1)x^3 + (6,0,0,0,0)x^5$ and

$b(x) = (0,1,2,3,0) + (0,0,1,2,0)x + (0,1,0,0,0)x^4 + (0,1,0,7,0)x^3 + (0,2,0,4,0)x^8$ in $V_R$.

We see $(a(x)) \times (b(x)) = (0,0,0,0,0)$. We see if $q(x)$ is not a constant polynomial certainly $(q(x))^2 \neq q(x)$ for if $\deg q(x) = n$ then $\deg ((q(x). q(x)) = n^2$.

We show that $V_R$ and $V_{n \times n}$ have several non trivial ideals.

***Example 3.2:*** Let $V_R$ be a ring. Consider the ideal generated by $p(x) = (2,3,1,5,7,8)x^3 + (4,2,0,1,5,7)$ in $V_R$. Clearly $I = \langle p(x) \rangle$ is a two sided ideal. Since $V_R$ is commutative every ideal is two sided. Infact $V_R$ has infinite number of ideals.

*Example 3.3:* Let $V_{n\times n}$ be a ring.

Take $p(x) = \begin{pmatrix} 3 & 1 \\ 6 & 2 \end{pmatrix} x^3 + \begin{pmatrix} 2 & 1 \\ 5 & 7 \end{pmatrix} x^2 + 1$ be in $V_{n\times n}$.

Clearly $\langle p(x) \rangle$ generates a two sided ideal.

But $p(x) V_{n\times n}$ generates only one sided ideal. Similarly $(V_{n\times n})(p(x))$ is not a two sided ideal. Thus $V_{n\times n}$ has infinite number of right ideals which are not left ideal and two sided ideals. Further $V_{n\times n}$ has left ideals which are not right ideals.

*Example 3.4:* Let

$$V_{4\times 4}[x] = \left\{ \sum_{i=0}^{\infty} a_i x^i \;\middle|\; a_i = \begin{pmatrix} x_1 & x_2 \\ x_3 & x_4 \end{pmatrix} \text{ where } x_j \in Q; 1 \leq j \leq 4 \right\}$$

be the matrix coefficient polynomial ring. Let

$$P = \left\{ \sum_{i=0}^{\infty} a_i x^i \;\middle|\; a_i = \begin{pmatrix} x_1 & x_2 \\ x_3 & x_4 \end{pmatrix} \text{ where } x_j \in Z; 1 \leq j \leq 4 \right\} \subseteq V_{4\times 4};$$

P is only a subring of $V_{4\times 4}$ and is not an ideal of $V_{4\times 4}$.

**THEOREM 3.6:** *Let $V_R$ and $V_{n\times n}$ be matrix coefficient polynomial rings. Both $V_R$ and $V_{n\times n}$ have subrings which are not ideals. We see if $p(x) \in V_R$ or $V_{n\times n}$; degree of p(x) as in case of usual polynomials is the highest power of x which has non zero coefficient.*

Consider
$p(x) = (2,3,4) + (0, -1, 2)x + (7,2,5)x^3 + (0,1,0)x^7 \in V_R$.

The degree of p(x) is even.

$$p(x) = \begin{pmatrix} 3 & 1 & 2 \\ 0 & 1 & 5 \\ 0 & 0 & 1 \end{pmatrix} + \begin{pmatrix} 7 & 2 & 1 \\ 0 & 5 & 7 \\ 6 & 1 & 2 \end{pmatrix} x^2 + \begin{pmatrix} 2 & 0 & 1 \\ 0 & 7 & 4 \\ 0 & 1 & 0 \end{pmatrix} x^4$$

$$+ \begin{pmatrix} 2 & 1 & 5 \\ 6 & 7 & 8 \\ 0 & 1 & 2 \end{pmatrix} x^8.$$

The degree of p(x) is 8.

Now in case of usual polynomials if their coefficients are from a field then every polynomial p(x) can be made monic.

However the same does not hold good in case of both $V_R$ and $V_{n \times n}$.

Consider
$p(x) = (0,3,0,0)x^4 + (1,2,3,4)x^3 + (2,0,0,1)x + (1,2,0,5)$ in $V_R$. Clearly p(x) cannot be made into a monic matrix coefficient polynomial for (0,3,0,0) has no inverse with respect to multiplication.

Let $q(x) = (5,7,8,-4)x^5 + (1,2,3,0)x^3 + (7,0,1,5)x + (8,9,0,2)$ be in $V_R$. Now q(x) can be made as a monic matrix coefficient polynomial. For multiply q(x) by

$t = (1/5, 1/7, 1/8, -1/4)$. Now $tq(x) = (1,1,1,1)x^5 + (1/5, 2/7, 3/8, 0)x^3 + (7/5, 0, 1/8, -5/4)x + (8/5, 9/7, 0, -2/4)$ is a monic matrix coefficient polynomial of degree five.

Let $p(x) = \begin{pmatrix} 3 & 0 \\ 1 & 0 \end{pmatrix} x^7 + \begin{pmatrix} 2 & 1 \\ 5 & 7 \end{pmatrix} x^3 + \begin{pmatrix} 8 & 1 \\ 0 & 5 \end{pmatrix} x^2 + \begin{pmatrix} 18 & 7 \\ 0 & 2 \end{pmatrix} x + \begin{pmatrix} 1 & 2 \\ 3 & 4 \end{pmatrix}$ be a matrix coefficient polynomial in $V_{2 \times 2}$. We see the coefficient of the highest power of p(x) is $\begin{pmatrix} 3 & 0 \\ 1 & 0 \end{pmatrix}$.

Clearly $\begin{pmatrix} 3 & 0 \\ 1 & 0 \end{pmatrix}$ has no inverse or the matrix $\begin{pmatrix} 3 & 0 \\ 1 & 0 \end{pmatrix}$ is non invertible.

Hence p(x) cannot be made into a monic matrix coefficient polynomial in $V_{2\times 2}$. Consider

$$p(x) = \begin{pmatrix} 7 & 0 \\ 0 & 8 \end{pmatrix} x^5 + \begin{pmatrix} 1 & 8 \\ 7 & 5 \end{pmatrix} x^4$$

$$+ \begin{pmatrix} 0 & 1 \\ 2 & 0 \end{pmatrix} x^3 + \begin{pmatrix} 0 & 1 \\ 1 & 0 \end{pmatrix} x^2 + \begin{pmatrix} 1 & 0 \\ 2 & 5 \end{pmatrix} \text{ in } V_{2\times 2}.$$

We see p(x) can be made into a monic polynomial.

$A = \begin{pmatrix} 1/7 & 0 \\ 0 & 1/8 \end{pmatrix}$ is such that

$$\begin{pmatrix} 1/7 & 0 \\ 0 & 1/8 \end{pmatrix} \begin{pmatrix} 7 & 0 \\ 0 & 8 \end{pmatrix} = \begin{pmatrix} 1 & 0 \\ 0 & 1 \end{pmatrix}.$$

Thus

$$\begin{pmatrix} 1/7 & 0 \\ 0 & 1/8 \end{pmatrix} p(x) = \begin{pmatrix} 1 & 0 \\ 0 & 1 \end{pmatrix} x^5 + \begin{pmatrix} 1/7 & 0 \\ 0 & 1/8 \end{pmatrix} \begin{pmatrix} 1 & 8 \\ 7 & 5 \end{pmatrix} x^4$$

$$+ \begin{pmatrix} 1/7 & 0 \\ 0 & 1/8 \end{pmatrix} \begin{pmatrix} 0 & 1 \\ 2 & 0 \end{pmatrix} x^3 + \begin{pmatrix} 1/7 & 0 \\ 0 & 1/8 \end{pmatrix} \begin{pmatrix} 0 & 1 \\ 1 & 0 \end{pmatrix} x^2$$

$$+ \begin{pmatrix} 1/7 & 0 \\ 0 & 1/8 \end{pmatrix} \begin{pmatrix} 1 & 0 \\ 2 & 5 \end{pmatrix}$$

$$= \begin{pmatrix} 0 & 1 \\ 1 & 0 \end{pmatrix} x^5 + \begin{pmatrix} 1/7 & 8/7 \\ 7/8 & 5/8 \end{pmatrix} x^4 + \begin{pmatrix} 0 & 1/7 \\ 1/4 & 0 \end{pmatrix} x^3$$

$$+ \begin{pmatrix} 0 & 1/7 \\ 1/8 & 0 \end{pmatrix} x^2 + \begin{pmatrix} 1/7 & 0 \\ 1/4 & 5/8 \end{pmatrix}$$

has been made into a monic polynomial.

We have shown only some of the matrix coefficient polynomials can be made and not all matrix coefficient polynomials as the collection of row matrices or collection of n×n matrices are not field just a ring with zero divisors.

Thus we have seen some of the properties of matrix coefficient polynomials. Unlike the number system which are not zero divisors, we cannot extend, all the results as these matrix coefficients can also be zero divisors.

Thus we can say a matrix coefficient polynomial $p(x) \in V_R$ (or $V_{n \times n}$) divides another matrix coefficient polynomial $q(x) \in V_R$ (or $V_{n \times n}$) if $q(x) = p(x) b(x)$ where $\deg(b(x)) < \deg q(x)$ and $\deg p(x) < \deg q(x)$.

We illustrate this by some examples. Suppose

$p(x) = ((3,2,1) + (7,-1,9)x) ((1,1,2) + (1,1,1)x) ((9,2,1) - (2,5,1)x)$ and

$q(x) = ((7,-1,9)x + (3,2,1)) ((1,1,1)x + (1,1,2)) ((9,2,1) - (2,5,1)x)) ((2,4,6)x^2 + (3,1,2)x + (1,3,6))$ are in $V_R$.

It is easily verified $p(x)/q(x)$ and $\deg(p(x)) = 3$ and $\deg q(x) = 5$.

However it is very difficult to derive all results in case of matrix coefficient polynomials; we have to define the concept of prime row matrix.

Suppose $X = (a_1, a_2, \ldots, a_n)$ is a row matrix with $a_i \in Z$, $1 \leq i \leq n$; we say X is a prime row vector or row matrix if each $a_i$ is a prime and none of the $a_i$ is zero. Thus (3,5,11,13), (7,5,2,19, 23,31) and (11,23,29,43,41,53,59,47,7,11) are prime row matrices.

We say or define the row matrix $(a_1, \ldots, a_n)$ divides the row matrix $(b_1, b_2, \ldots, b_n)$ if none of the $a_i$'s are zero for $i=1,2,\ldots,n$ and $a_i/b_i$ for every i, $1 \leq i \leq n$. That is we say $(a_1, \ldots, a_n) / (b_1, b_2, \ldots, b_n)$ if $(b_1/a_1, \ldots, b_n/a_n) = (c_1, \ldots, c_n)$ and $c_i \in Z$; $1 \leq i \leq n$ ($a_i \neq 0$; $i = 1, 2, \ldots, n$).

We will illustrate this situation by some examples.

Let $(5,7,2,8) = x$ and $y = (10,14,8,8)$ we say x/y and y/x = (10/5, 14/7, 8/2, 8/8) = (2,2,4,1).

Now if $x = (0,2,3,5,7,8)$ and $y = (5,4,6,10,21,24)$ then $x \not{/} y$ or y/x is not defined.

So when matrix coefficient polynomials are dealt with it is very very difficult to define division in $V_R$.

Clearly if $x = (a_1, \ldots, a_n)$ with $a_i \neq 0$ and $a_i$ primes for all i, then we see there does not exist any $y = (b_1, \ldots, b_n)$ with $b_i \neq 0$ and $b_i \neq 1$ dividing x, ($1 \leq i \leq n$). Thus the only divisors of $x = (a_1, \ldots, a_n)$ are $y = (1,1,\ldots, 1)$ and $y = (a_1, a_2, \ldots, a_n)$ only. Since we face a lot of problems in dealing with matrix multiplication and however we only multiply the two row matrices of same order $x = (a_1, a_2, \ldots, a_n)$ with $y = (b_1, b_2, \ldots, b_n)$ as x.y = $(a_1, a_2, \ldots, a_n) (b_1, b_2, \ldots, b_n)$ as x.y = $(a_1b_1, a_2b_2, \ldots, a_nb_n)$ we wish to extend this sort of multiplication for all matrices only criteria being that they should be of same order.

We call such multiplication or product of matrices of same order as natural multiplication of matrices. Thus we define natural multiplication or product of two $n \times 1$ column matrices as follows; if $x = \begin{bmatrix} a_1 \\ a_2 \\ \vdots \\ a_n \end{bmatrix}$ and $y = \begin{bmatrix} b_1 \\ b_2 \\ \vdots \\ b_n \end{bmatrix}$ then the natural product of x with y denoted by $x \times_n y = \begin{bmatrix} a_1 \\ a_2 \\ \vdots \\ a_n \end{bmatrix} \times_n \begin{bmatrix} b_1 \\ b_2 \\ \vdots \\ b_n \end{bmatrix} = \begin{bmatrix} a_1 b_1 \\ a_2 b_2 \\ \vdots \\ a_n b_n \end{bmatrix}$.

This product is defined as natural product of two n×1 column matrices and the natural product operation is denoted by $\times_n$.

*Example 3.5:* Let $x = \begin{bmatrix} 7 \\ 2 \\ 0 \\ 1 \\ 5 \end{bmatrix}$ and $y = \begin{bmatrix} 1 \\ 3 \\ 5 \\ 2 \\ 7 \end{bmatrix}$.

Now the natural product of x with y is $x \times_n y = \begin{bmatrix} 7 \\ 2 \\ 0 \\ 1 \\ 5 \end{bmatrix} \times_n \begin{bmatrix} 1 \\ 3 \\ 5 \\ 2 \\ 7 \end{bmatrix}$

$= \begin{bmatrix} 7.1 \\ 2.3 \\ 0.5 \\ 1.2 \\ 5.7 \end{bmatrix} = \begin{bmatrix} 7 \\ 6 \\ 0 \\ 2 \\ 35 \end{bmatrix}$.

We see the natural product is both associative and commutative.

Now if $x = \begin{bmatrix} 1 \\ 1 \\ \vdots \\ 1 \end{bmatrix}$ and $y = \begin{bmatrix} a_1 \\ a_2 \\ \vdots \\ a_n \end{bmatrix}$ be any two n × 1 column matrices when $x \times_n y = y \times_n x = y$.

Thus $x = \begin{bmatrix} 1 \\ 1 \\ \vdots \\ 1 \end{bmatrix}$ acts as the natural product identity. We see infact any n × 1 collection of column vectors is a semigroup under natural multiplication or natural product and is a monoid and is a commutative monoid.

**THEOREM 3.7:** *Let*

$$V = \left\{ \begin{bmatrix} a_1 \\ a_2 \\ \vdots \\ a_n \end{bmatrix} \,\middle|\, a_i \in Q \text{ (or Z or R)}; \, 1 \leq i \leq n \right\}$$

*be the collection of all n × 1 column matrices. V is a commutative semigroup under natural product (or multiplication) of column matrices.*

Proof is direct and hence is left as an exercise to the reader.

*Example 3.6:* Let

$$x = \begin{bmatrix} 1 \\ 2 \\ 3 \\ 0 \\ 0 \\ 0 \end{bmatrix} \text{ and } y = \begin{bmatrix} 0 \\ 0 \\ 0 \\ 0 \\ 1 \\ 2 \end{bmatrix}$$

be 6 × 1 column matrices.

We see $x \times_n y = \begin{bmatrix} 1 \\ 2 \\ 3 \\ 0 \\ 0 \\ 0 \end{bmatrix} \times_n \begin{bmatrix} 0 \\ 0 \\ 0 \\ 0 \\ 1 \\ 2 \end{bmatrix} = \begin{bmatrix} 0 \\ 0 \\ 0 \\ 0 \\ 0 \\ 0 \end{bmatrix}$.

Thus x is a zero divisor. Inview of this we have the following result.

**THEOREM 3.8:** *Let*

$$V = \left\{ \begin{bmatrix} a_1 \\ a_2 \\ \vdots \\ a_n \end{bmatrix} \middle| \; a_i \in Z \; (or \; Q \; or \; R); \; 1 \leq i \leq n \right\}$$

*be the semigroup under natural multiplication $\times_n$. V has zero divisors.*

This proof is also very simple.

*Example 3.7:* Let

$$V = \left\{ \begin{bmatrix} a_1 \\ a_2 \\ \vdots \\ a_n \end{bmatrix} \middle| \; a_i \in Z; \; 1 \leq i \leq 6 \right\}$$

be a semigroup under natural product.

Take

$$W = \left\{ \begin{bmatrix} a_1 \\ a_2 \\ \vdots \\ a_6 \end{bmatrix} \middle| a_i \in 3Z; 1 \leq i \leq 6 \right\} \subseteq V;$$

W is a subsemigroup of V. Infact W is an ideal of the semigroup V. Thus we have several ideals for V.

*Example 3.8:* Let

$$V = \left\{ \begin{bmatrix} a_1 \\ a_2 \\ \vdots \\ a_{10} \end{bmatrix} \middle| a_i \in Q; 1 \leq i \leq 10 \right\}$$

be the semigroup under natural product.

Consider

$$W = \left\{ \begin{bmatrix} a_1 \\ a_2 \\ \vdots \\ a_{10} \end{bmatrix} \middle| a_i \in Z; 1 \leq i \leq 10 \right\} \subseteq V;$$

W is only a subsemigroup of V and is not an ideal of V.

Take

$$S = \left\{ \begin{bmatrix} a_1 \\ a_2 \\ 0 \\ \vdots \\ 0 \end{bmatrix} \middle| a_i \in Q; 1 \leq i \leq 10 \right\} \subseteq V;$$

S is a subsemigroup of V under usual product. Also S is an ideal of V.

From this example we see a subsemigroup in general is not an ideal.

Inview of this we give the following result the proof of which is simple.

**THEOREM 3.9:** *Let*

$$V = \left\{ \begin{bmatrix} a_1 \\ a_2 \\ \vdots \\ a_n \end{bmatrix} \middle| a_i \in Q \text{ (or R)}; 1 \leq i \leq n \right\}$$

*be a semigroup under natural product. V has subsemigroups which are not ideals. However every ideal is a subsemigroup.*

Proof is left as an exercise for the reader.

Now we have the concept of Smarandache semigroups. We will illustrate this situation by an example.

*Example 3.9:* Let

$$V = \left\{ \begin{bmatrix} a_1 \\ a_2 \\ \vdots \\ a_8 \end{bmatrix} \middle| a_i \in Q; 1 \leq i \leq 8 \right\}$$

be the semigroup under natural multiplication.

Consider

$$M = \left\{ \begin{bmatrix} a_1 \\ a_2 \\ \vdots \\ a_8 \end{bmatrix} \middle| a_i \in Q \setminus \{0\}; 1 \leq i \leq 8 \right\} \subseteq V;$$

M is a subring as well as a group under natural product. Further we see M is not an ideal of V. Thus V is a Smarandache semigroup.

Inview of this we can easily prove the following theorem.

**THEOREM 3.10:** *Let*

$$V = \left\{ \begin{bmatrix} m_1 \\ m_2 \\ \vdots \\ m_n \end{bmatrix} \middle| \; m_i \in Q \; (or \; R); \; 1 \leq i \leq n \right\}$$

*be a semigroup under natural product. V is a Smarandache semigroup.*

*Proof:* For take

$$M = \left\{ \begin{bmatrix} a_1 \\ a_2 \\ \vdots \\ a_m \end{bmatrix} \middle| \; a_i \in Q \setminus \{0\} \; or \; (a_i \in R \setminus \{0\}); \; 1 \leq i \leq m \right\} \subseteq V,$$

M is a group under natural product as every element in M is invertible, hence the theorem.

Now we proceed onto give an example or two.

*Example 3.10:* Let

$$M = \left\{ \begin{bmatrix} a_1 \\ a_2 \\ a_3 \end{bmatrix} \middle| \; a_i \in Z; \; 1 \leq i \leq 3 \right\}$$

be a semigroup under natural product.

Consider the set

$$P = \left\{ \begin{bmatrix} 1 \\ 1 \\ 1 \end{bmatrix}, \begin{bmatrix} -1 \\ -1 \\ -1 \end{bmatrix}, \begin{bmatrix} 1 \\ -1 \\ 1 \end{bmatrix}, \begin{bmatrix} -1 \\ 1 \\ -1 \end{bmatrix}, \begin{bmatrix} -1 \\ -1 \\ 1 \end{bmatrix}, \begin{bmatrix} 1 \\ 1 \\ -1 \end{bmatrix}, \begin{bmatrix} -1 \\ 1 \\ 1 \end{bmatrix}, \begin{bmatrix} 1 \\ -1 \\ -1 \end{bmatrix} \right\} \subseteq M$$

is a group under product. Thus Z is a Smarandache semigroup.

Infact $B = \left\{ \begin{bmatrix} 1 \\ 1 \\ 1 \end{bmatrix}, \begin{bmatrix} -1 \\ -1 \\ -1 \end{bmatrix} \right\} \subseteq M$ is also a group. Thus P is a Smarandache semigroup as $B \subseteq P$.

Now we wish to prove the following theorem.

**THEOREM 3.11:** *Let*

$$M = \left\{ \begin{bmatrix} a_1 \\ a_2 \\ \vdots \\ a_n \end{bmatrix} \middle| a_i \in Z \text{ (or } Q \text{ or } R\text{); } 1 \leq i \leq n \right\}$$

*be a semigroup under natural product. If M has a Smarandache subsemigroup then M is a Smarandache semigroup. However even if M is a Smarandache semigroup, every subsemigroup of M need not be a Smarandache subsemigroup.*

*Proof:* Suppose we have a proper subsemigroup under natural product for M say W. W is a Smarandache subsemigroup of M; then W has a proper subset X such that X is a group under natural product; $X \subseteq W$.

Now $X \subseteq W \subseteq M$; that is $X \subseteq M$ so M is a Smarandache semigroup. Hence the claim.

We prove the other part of the theorem by an example.

Consider
$$Y = \left\{ \begin{bmatrix} x_1 \\ x_2 \\ \vdots \\ x_n \end{bmatrix} \middle| x_i \in Z; 1 \leq i \leq n \right\}$$

be a semigroup under natural product.

Y is a Smarandache semigroup as

$$P = \left\{ \begin{bmatrix} 1 \\ 1 \\ 1 \\ 1 \\ 1 \\ 1 \\ 1 \\ 1 \end{bmatrix}, \begin{bmatrix} -1 \\ -1 \\ -1 \\ -1 \\ -1 \\ -1 \\ -1 \\ -1 \end{bmatrix} \right\} \subseteq Y$$

is a group under natural multiplication. Hence Y is a S-semigroup.

Take

$$W = \left\{ \begin{bmatrix} a_1 \\ a_2 \\ \vdots \\ a_7 \end{bmatrix} \middle| a_i \in 3Z; \ 1 \le i \le 7 \right\} \subseteq Y;$$

W is only a subsemigroup of Y and is not a Smarandache subsemigroup of Y. Hence even if Y is a S-semigroup. Y has subsemigroups whch are not Smarandache subsemigroup. Hence the theorem.

Now we have seen ideals and subsemigroups and S-subsemigroups about column matrix semigroups under natural product.

Now we define natural product on m × n matrix semigroups (m ≠ n).

**DEFINITION 3.1:** *Let*

$$M = \left\{ \begin{bmatrix} a_{11} & a_{12} & \ldots & a_{1n} \\ a_{21} & a_{22} & \ldots & a_{2n} \\ \vdots & \vdots & & \vdots \\ a_{m1} & a_{m2} & \ldots & a_{mn} \end{bmatrix} \middle| a_{ij} \in Z \text{ (or } Q \text{ or } R\text{)}; \right.$$
$$\left. 1 \le i \le m \text{ and } 1 \le j \le n \right\}$$

*be the collection of all $m \times n$ ($m \ne n$) matrices. M under natural multiplication / product $\times_n$ is a semigroup.*

*If* $X = \begin{bmatrix} a_{11} & a_{12} & \ldots & a_{1n} \\ a_{21} & a_{22} & \ldots & a_{2n} \\ \vdots & \vdots & & \vdots \\ a_{m1} & a_{m2} & \ldots & a_{mn} \end{bmatrix}$ *and* $Y = \begin{bmatrix} b_{11} & b_{12} & \ldots & b_{1n} \\ b_{21} & b_{22} & \ldots & b_{2n} \\ \vdots & \vdots & & \vdots \\ b_{m1} & b_{m2} & \ldots & b_{mn} \end{bmatrix}$

*be any two $m \times n$ matrices in M.*

$$\text{We define } X \times_m Y = \begin{bmatrix} a_{11}b_{11} & a_{12}b_{12} & \ldots & a_{1n}b_{1n} \\ a_{21}b_{21} & a_{22}b_{22} & \ldots & a_{2n}b_{2n} \\ \vdots & \vdots & & \vdots \\ a_{m1}b_{m1} & a_{m2}b_{m2} & \ldots & a_{mn}b_{mn} \end{bmatrix}.$$

*Clearly $X \times_m Y$ is in M. $(M, \times_n)$ is defined as the semigroup under natural product.*

We give examples of them.

*Example 3.11:* Let

$$X = \begin{bmatrix} 2 & 1 & 0 & 5 & 1 \\ 0 & 3 & 1 & 2 & 5 \\ -1 & 4 & 3 & 0 & 1 \end{bmatrix} \text{ and } Y = \begin{bmatrix} 3 & 2 & 0 & 1 & 3 \\ 4 & 0 & 1 & 5 & 7 \\ 0 & 1 & 2 & 0 & 5 \end{bmatrix}$$

be any two 3 × 5 matrices. We find the natural product of X with Y.

$$X \times_n Y = \begin{bmatrix} 2 & 1 & 0 & 5 & 1 \\ 0 & 3 & 1 & 2 & 5 \\ -1 & 4 & 3 & 0 & 1 \end{bmatrix} \times_n \begin{bmatrix} 3 & 2 & 0 & 1 & 3 \\ 4 & 0 & 1 & 5 & 7 \\ 0 & 1 & 2 & 0 & 5 \end{bmatrix}$$

$$= \begin{bmatrix} 6 & 2 & 0 & 5 & 3 \\ 0 & 0 & 1 & 10 & 35 \\ 0 & 4 & 6 & 0 & 5 \end{bmatrix}.$$

*Example 3.12:* Let

$$S = \left\{ \begin{bmatrix} a_1 & a_2 & a_3 \\ a_4 & a_5 & a_6 \\ \vdots & \vdots & \vdots \\ a_{28} & a_{29} & a_{30} \end{bmatrix} \middle| a_i \in Z; 1 \le i \le 30 \right\}$$

be the semigroup under natural product. S is a commutative semigroup with identity. S has infinite number of ideals and subsemigroups which are not ideals.

*Example 3.13:* Let

$$S = \left\{ \begin{bmatrix} a_1 & a_2 \\ a_3 & a_4 \\ \vdots & \vdots \\ a_{11} & a_{12} \end{bmatrix} \middle| a_i \in Q; 1 \le i \le 12 \right\}$$

be the semigroup under natural product.

Take

$$I = \left\{ \begin{bmatrix} a_1 & a_2 \\ 0 & 0 \\ 0 & 0 \\ \vdots & \vdots \\ 0 & 0 \\ a_3 & a_4 \end{bmatrix} \middle| a_i \in Q; 1 \leq i \leq 4 \right\}.$$

It is easily verified I is an ideal of P under natural product $\times_n$.

Consider the subsemigroup

$$S = \left\{ \begin{bmatrix} a_1 & a_2 \\ a_3 & a_4 \\ 0 & 0 \\ \vdots & \vdots \\ 0 & 0 \\ a_5 & a_6 \end{bmatrix} \middle| a_i \in Q; 1 \leq i \leq 6 \right\} \subseteq P$$

under natural product. Clearly S is an ideal of P.

Suppose

$$T = \left\{ \begin{bmatrix} a_1 & a_2 \\ a_3 & a_4 \\ \vdots & \vdots \\ a_{11} & a_{12} \end{bmatrix} \middle| a_i \in Z; 1 \leq i \leq 12 \right\} \subseteq P,$$

T is only a subsemigroup of P under natural product and is not an ideal of P.

We can as in case of usual semigroups define in case of these semigroups under natural product the concept of Smarandache-ideals, Smarandache zero divisors and so on. By our natural product we are able to define some form of product on column matrices and rectangular matrices. Now we proceed onto define natural product on usual square matrices.

Let $A = (a_{ij})_{n \times n}$ and $B = (b_{ij})_{n \times n}$ be square matrices; $a_{ij}, b_{ij} \in Z$ (or Q or R); $1 \leq i, j \leq n$. We define the natural product

$$A \times_n B \text{ as } A \times_n B = (a_{ij})_{n \times n} \; (b_{ij})_{n \times n}$$
$$= (a_{ij} \, b_{ij})_{n \times n}$$
$$= (c_{ij})_{n \times n}.$$

We will illustrate this by few examples.

*Example 3.14:* Let

$$A = \begin{pmatrix} 6 & 1 & 2 \\ 0 & 3 & 4 \\ 2 & 1 & 0 \end{pmatrix} \text{ and } B = \begin{pmatrix} 3 & 0 & 1 \\ 2 & 1 & 0 \\ 0 & 1 & 2 \end{pmatrix}$$

be two $3 \times 3$ matrices. To find the natural product of A with B.

$$A \times_n B = \begin{pmatrix} 6 & 1 & 2 \\ 0 & 3 & 4 \\ 2 & 1 & 0 \end{pmatrix} \begin{pmatrix} 3 & 0 & 1 \\ 2 & 1 & 0 \\ 0 & 1 & 2 \end{pmatrix} = \begin{pmatrix} 18 & 0 & 2 \\ 0 & 3 & 0 \\ 0 & 1 & 0 \end{pmatrix}.$$

Now the usual matrix product of A with B is

$$A.B = \begin{pmatrix} 6 & 1 & 2 \\ 0 & 3 & 4 \\ 2 & 1 & 0 \end{pmatrix} \begin{pmatrix} 3 & 0 & 1 \\ 2 & 1 & 0 \\ 0 & 1 & 2 \end{pmatrix}$$

$$= \begin{pmatrix} 20 & 3 & 10 \\ 6 & 7 & 8 \\ 8 & 1 & 2 \end{pmatrix}.$$

We see A.B ≠ A ×$_n$ B in general. Further we see the operation '.' the usual matrix multiplication is non commutative where as the natural product ×$_n$ is commutative.

We just consider the following examples.

*Example 3.15:* Let

$$M = \begin{bmatrix} 3 & 4 \\ 2 & 0 \end{bmatrix} \text{ and } N = \begin{bmatrix} 1 & 2 \\ 0 & 1 \end{bmatrix}$$

be any two 2 × 2 matrices.

$$M.N = \begin{bmatrix} 3 & 4 \\ 2 & 0 \end{bmatrix} \begin{bmatrix} 1 & 2 \\ 0 & 1 \end{bmatrix} = \begin{bmatrix} 3 & 10 \\ 2 & 4 \end{bmatrix}$$

$$\text{and } N.M = \begin{bmatrix} 1 & 2 \\ 0 & 1 \end{bmatrix} \begin{bmatrix} 3 & 4 \\ 2 & 0 \end{bmatrix} = \begin{bmatrix} 7 & 4 \\ 2 & 0 \end{bmatrix}.$$

We see M.N ≠ N.M.

$$\text{However } M \times_n N = \begin{bmatrix} 3 & 4 \\ 2 & 0 \end{bmatrix} \times_n \begin{bmatrix} 1 & 2 \\ 0 & 1 \end{bmatrix}$$

$$= \begin{bmatrix} 3 & 8 \\ 0 & 0 \end{bmatrix} \text{ and}$$

$$N \times_n M = \begin{bmatrix} 1 & 2 \\ 0 & 1 \end{bmatrix} \times \begin{bmatrix} 3 & 4 \\ 2 & 0 \end{bmatrix} = \begin{bmatrix} 3 & 8 \\ 0 & 0 \end{bmatrix}.$$

Thus N ×$_n$ M = M ×$_n$ N.

*Example 3.16:* Let

$$M = \begin{bmatrix} 7 & 0 & 0 & 0 \\ 0 & 8 & 0 & 0 \\ 0 & 0 & 2 & 0 \\ 0 & 0 & 0 & 4 \end{bmatrix} \text{ and } N = \begin{bmatrix} 1 & 0 & 0 & 0 \\ 0 & 2 & 0 & 0 \\ 0 & 0 & 3 & 0 \\ 0 & 0 & 0 & 4 \end{bmatrix}.$$

We find $M.N = \begin{bmatrix} 7 & 0 & 0 & 0 \\ 0 & 8 & 0 & 0 \\ 0 & 0 & 2 & 0 \\ 0 & 0 & 0 & 4 \end{bmatrix} \begin{bmatrix} 1 & 0 & 0 & 0 \\ 0 & 2 & 0 & 0 \\ 0 & 0 & 3 & 0 \\ 0 & 0 & 0 & 4 \end{bmatrix}$

$$= \begin{bmatrix} 7 & 0 & 0 & 0 \\ 0 & 16 & 0 & 0 \\ 0 & 0 & 6 & 0 \\ 0 & 0 & 0 & 16 \end{bmatrix}.$$

Also

$$N.M = \begin{bmatrix} 1 & 0 & 0 & 0 \\ 0 & 2 & 0 & 0 \\ 0 & 0 & 3 & 0 \\ 0 & 0 & 0 & 4 \end{bmatrix} \begin{bmatrix} 7 & 0 & 0 & 0 \\ 0 & 8 & 0 & 0 \\ 0 & 0 & 2 & 0 \\ 0 & 0 & 0 & 4 \end{bmatrix} = \begin{bmatrix} 7 & 0 & 0 & 0 \\ 0 & 16 & 0 & 0 \\ 0 & 0 & 6 & 0 \\ 0 & 0 & 0 & 16 \end{bmatrix}.$$

Now consider $M \times_n N = \begin{bmatrix} 7 & 0 & 0 & 0 \\ 0 & 8 & 0 & 0 \\ 0 & 0 & 2 & 0 \\ 0 & 0 & 0 & 4 \end{bmatrix} \times_n \begin{bmatrix} 1 & 0 & 0 & 0 \\ 0 & 2 & 0 & 0 \\ 0 & 0 & 3 & 0 \\ 0 & 0 & 0 & 4 \end{bmatrix}$

$$= \begin{bmatrix} 7 & 0 & 0 & 0 \\ 0 & 16 & 0 & 0 \\ 0 & 0 & 6 & 0 \\ 0 & 0 & 0 & 16 \end{bmatrix}.$$ We see $M.N = M \times_n N$.

In view of this we have the following theorem.

**THEOREM 3.12:** *Let*

$$M = \left\{ \begin{bmatrix} a_1 & 0 & 0 & 0 & \ldots & 0 \\ 0 & a_2 & 0 & 0 & \ldots & 0 \\ 0 & 0 & a_3 & 0 & \ldots & 0 \\ \vdots & \vdots & \vdots & \vdots & & \vdots \\ 0 & 0 & 0 & 0 & \ldots & a_n \end{bmatrix} \middle| a_i \in Q \text{ (or } Z \text{ or } R \text{ or } C\text{)};\right.$$

$$1 \leq i \leq n \}$$

*be the collection of all $n \times n$ diagonal matrices. M is a semigroup under natural product and M is also a semigroup under usual product of matrices and both the operations are identical on M.*

*Proof:* Let

$$A = \begin{bmatrix} a_1 & 0 & 0 & 0 & \ldots & 0 \\ 0 & a_2 & 0 & 0 & \ldots & 0 \\ 0 & 0 & a_3 & 0 & \ldots & 0 \\ 0 & 0 & 0 & a_4 & & 0 \\ \vdots & \vdots & \vdots & \vdots & \ldots & \vdots \\ 0 & 0 & 0 & 0 & \ldots & a_n \end{bmatrix}$$

and $\begin{bmatrix} b_1 & 0 & 0 & 0 & \ldots & 0 \\ 0 & b_2 & 0 & 0 & \ldots & 0 \\ 0 & 0 & b_3 & 0 & \ldots & 0 \\ 0 & 0 & 0 & b_4 & & 0 \\ \vdots & \vdots & \vdots & \vdots & \ldots & \vdots \\ 0 & 0 & 0 & 0 & \ldots & b_n \end{bmatrix}$ be two matrices from M.

Now let us consider the natural product of

$$A \times_n B = \begin{bmatrix} a_1 b_1 & 0 & 0 & 0 & \ldots & 0 \\ 0 & a_2 b_2 & 0 & 0 & \ldots & 0 \\ 0 & 0 & a_3 b_3 & 0 & \ldots & 0 \\ 0 & 0 & 0 & a_4 b_4 & & 0 \\ \vdots & \vdots & \vdots & \vdots & \ldots & \vdots \\ 0 & 0 & 0 & 0 & \ldots & a_n b_n \end{bmatrix}.$$

Consider the matrix product;

$$A.B = \begin{bmatrix} a_1 b_1 & 0 & 0 & 0 & \ldots & 0 \\ 0 & a_2 b_2 & 0 & 0 & \ldots & 0 \\ 0 & 0 & a_3 b_3 & 0 & \ldots & 0 \\ 0 & 0 & 0 & a_4 b_4 & & 0 \\ \vdots & \vdots & \vdots & \vdots & \ldots & \vdots \\ 0 & 0 & 0 & 0 & \ldots & a_n b_n \end{bmatrix}.$$

It is easily verified $A.B = A \times_n B$. Thus both the operations are identical as diagonal matrices.

**THEOREM 3.13:** *Let*

*$M = \{(a_{ij})_{n \times n} \mid a_{ij} \in R$ (or Q or Z or C); $1 \leq i, j \leq n\}$*
*be the collection of all $n \times n$ matrices, M is a semigroup under natural product and M is a semigroup under matrix*

*multiplication. Both the operations on M are distinct in general.*

The proof is direct and hence left as an exercise to the reader.

**THEOREM 3.14:** *Let*
$$M = \{(a_{ij}) \mid a_{ij} \in Z \text{ (or } Q \text{ or } R \text{ or } C\text{); } 1 \leq i, j \leq n\}$$
*be a semigroup under natural product. M is a Smarandache semigroup.*

*Proof:* Let P = $\{(a_{ij}) \mid a_{ij} \in Z \setminus \{0\}$, (R] $\{0\}$ or Q] $\{0\}$ or C \ $\{0\}$) $1 \leq i \leq n\} \subseteq M$ be a group under natural multiplication. So M is a S-semigroup.

It is pertinent to mention here that these semigroups have ideals subsemigroups, zero divisors and idempotents and their Smarandache analogue.

Now we proceed onto give more structures using this product.

**DEFINITION 3.2:** *Let*

$$M = \left\{ \begin{bmatrix} a_1 \\ a_2 \\ \vdots \\ a_m \end{bmatrix} \middle| a_i \in Q \text{ (or } Z \text{ or } R \text{ or } C\text{); } 1 \leq i \leq m \right\}$$

*be the collection of all m $\times$ 1 column matrices. M is a ring under usual matrix addition and natural product $\times_n$.*

*Example 3.17:* Let

$$M = \left\{ \begin{bmatrix} a_1 \\ a_2 \\ a_3 \\ a_4 \end{bmatrix} \middle| a_i \in R \text{ (or Q or Z)}; \ 1 \leq i \leq 4 \right\}$$

be a ring under + and $\times_n$. The reader can easily verify that $A \times_n (B+C) = A \times_n B + A \times_n C$ where A, B and C are $n \times 1$ column matrices.

$$\text{Consider } A = \begin{bmatrix} a_1 \\ a_2 \\ \vdots \\ a_n \end{bmatrix}, \ B = \begin{bmatrix} b_1 \\ b_2 \\ \vdots \\ b_n \end{bmatrix} \text{ and } C = \begin{bmatrix} c_1 \\ c_2 \\ \vdots \\ c_n \end{bmatrix};$$

$$A \times_n (B+C) = \begin{bmatrix} a_1 \\ a_2 \\ \vdots \\ a_n \end{bmatrix} \times_n \left( \begin{bmatrix} b_1 \\ b_2 \\ \vdots \\ b_n \end{bmatrix} + \begin{bmatrix} c_1 \\ c_2 \\ \vdots \\ c_n \end{bmatrix} \right)$$

$$= \begin{bmatrix} a_1 \\ a_2 \\ \vdots \\ a_n \end{bmatrix} \times_n \begin{bmatrix} b_1 + c_1 \\ b_2 + c_2 \\ \vdots \\ b_n + c_n \end{bmatrix} = \begin{bmatrix} a_1(b_1+c_1) \\ a_2(b_2+c_2) \\ \vdots \\ a_n(b_n+c_n) \end{bmatrix}$$

$$= \begin{bmatrix} a_1 b_1 + a_1 c_1 \\ a_2 b_2 + a_2 c_2 \\ \vdots \\ a_n b_n + a_n c_n \end{bmatrix} = \begin{bmatrix} a_1 b_1 \\ a_2 b_2 \\ \vdots \\ a_n b_n \end{bmatrix} + \begin{bmatrix} a_1 c_1 \\ a_2 c_2 \\ \vdots \\ a_n c_n \end{bmatrix}.$$

Now consider $A \times_n B + A \times_n C$

$$= \begin{bmatrix} a_1 \\ a_2 \\ \vdots \\ a_n \end{bmatrix} \times_n \begin{bmatrix} b_1 \\ b_2 \\ \vdots \\ b_n \end{bmatrix} + \begin{bmatrix} a_1 \\ a_2 \\ \vdots \\ a_n \end{bmatrix} \times_n \begin{bmatrix} c_1 \\ c_2 \\ \vdots \\ c_n \end{bmatrix}$$

$$= \begin{bmatrix} a_1 b_1 \\ a_2 b_2 \\ \vdots \\ a_n b_n \end{bmatrix} + \begin{bmatrix} a_1 c_1 \\ a_2 c_2 \\ \vdots \\ a_n c_n \end{bmatrix} = \begin{bmatrix} a_1 b_1 + a_1 c_1 \\ a_2 b_2 + a_2 c_2 \\ \vdots \\ a_n b_n + a_n c_n \end{bmatrix}.$$

Thus we see $\times_n$ distributes over addition. Now consider the collection of all $m \times n$ matrices ($m \neq n$) with entries taken from Z or Q or C or R. We see this collection also under matrix addition and natural product is a ring. Let

$$M = \{(a_{ij})_{m \times n} \mid m \neq n;\ a_{ij} \in R\ (\text{or Z or Q or C});$$
$$1 \leq i \leq m \text{ and } 1 \leq j \leq n\};$$

M is a ring infact a commutative ring.

However M is not a ring under matrix addition and matrix product.

*Example 3.18:* Let

$$M = \left\{ \begin{bmatrix} a_1 & a_2 \\ a_3 & a_4 \\ a_5 & a_6 \\ a_7 & a_8 \end{bmatrix} \middle| a_i \in Q\ (\text{or Z or R or C});\ 1 \leq i \leq 8 \right\}$$

be a ring under matrix addition and natural product.

M is a commutative ring with unit $\begin{bmatrix} 1 & 1 \\ 1 & 1 \\ 1 & 1 \\ 1 & 1 \end{bmatrix}$.

M has units, zero divisors, subrings and ideals.

Take $a = \begin{bmatrix} 3 & 4 \\ 5 & 8 \\ 1 & 9 \\ 4 & 7 \end{bmatrix}$ and $b = \begin{bmatrix} 1/3 & 1/4 \\ 1/5 & 1/8 \\ 1 & 1/9 \\ 1/4 & 1/7 \end{bmatrix}$;

clearly $ab = ba = \begin{bmatrix} 1 & 1 \\ 1 & 1 \\ 1 & 1 \\ 1 & 1 \end{bmatrix}$.

Consider $a = \begin{bmatrix} 0 & 0 \\ a_1 & a_2 \\ a_3 & a_4 \\ 0 & 0 \end{bmatrix}$ and $b = \begin{bmatrix} a_1 & a_2 \\ 0 & 0 \\ 0 & 0 \\ a_3 & a_4 \end{bmatrix} \in M$.

Clearly $ab = \begin{bmatrix} 0 & 0 \\ 0 & 0 \\ 0 & 0 \\ 0 & 0 \end{bmatrix}$ is a zero divisor in M.

Take
$$P = \left\{ \begin{bmatrix} a_1 & 0 \\ a_2 & 0 \\ a_3 & 0 \\ a_4 & 0 \end{bmatrix} \middle| a_i \in Q; \ 1 \leq i \leq 4 \right\} \subseteq M;$$

P is an ideal of M.

Consider
$$T = \left\{ \begin{bmatrix} a_1 & a_2 \\ a_3 & a_4 \\ a_5 & a_6 \\ a_7 & a_8 \end{bmatrix} \middle| a_i \in Z; \ 1 \leq i \leq 8 \right\} \subseteq M;$$

clearly T is only a subring of M and is not an ideal of M. Thus M has subrings which are not ideals. We can find several subrings which are not ideals.

*Example 3.19:* Let

$$N = \left\{ \begin{bmatrix} a_1 & a_2 & a_3 \\ a_4 & a_5 & a_6 \end{bmatrix} \middle| a_i \in Z; \ 1 \leq i \leq 6 \right\}$$

be the ring under matrix addition and natural product. We see M has no units.

$\begin{bmatrix} 1 & 1 & 1 \\ 1 & 1 & 1 \end{bmatrix}$ is the identity with the natural product $\times_n$ in N.

Consider
$$P = \left\{ \begin{bmatrix} b_1 & b_2 & b_3 \\ b_4 & b_5 & b_6 \end{bmatrix} \middle| b_i \in 3Z; \ 1 \leq i \leq 6 \right\} \subseteq N$$

is an ideal of N.

*Example 3.20:* Let

$$M = \left\{ \begin{bmatrix} a_1 & a_2 & a_3 \\ a_4 & a_5 & a_6 \\ \vdots & \vdots & \vdots \\ a_{31} & a_{32} & a_{33} \end{bmatrix} \middle| a_i \in Q;\ 1 \le i \le 33 \right\}$$

be a ring under matrix addition and natural product. M is a S-ring.

For consider

$$P = \left\{ \begin{bmatrix} b_1 & b_2 & b_3 \\ b_4 & b_5 & b_6 \\ \vdots & \vdots & \vdots \\ b_{31} & b_{32} & b_{33} \end{bmatrix} \middle| b_i \in 3Z;\ 1 \le i \le 33 \right\} \subseteq M;$$

P is not an ideal of M. M has units, zero divisors, subrings and ideals.

Take

$$W = \left\{ \begin{bmatrix} a_1 & a_2 & a_3 \\ 0 & 0 & 0 \\ \vdots & \vdots & \vdots \\ 0 & 0 & 0 \\ a_4 & a_5 & a_6 \end{bmatrix} \middle| a_i \in Q;\ 1 \le i \le 6 \right\} \subseteq M,$$

W is an ideal of M.

Consider

$$S = \left\{ \begin{bmatrix} a_1 & a_2 & a_3 \\ 0 & 0 & 0 \\ \vdots & \vdots & \vdots \\ 0 & 0 & 0 \end{bmatrix} \middle| a_i \in Z; \ 1 \leq i \leq 3 \right\} \subseteq M$$

is only a subring and not an ideal.

*Example 3.21:* Let

$$M = \left\{ \begin{bmatrix} a_1 & a_2 \\ a_3 & a_4 \\ \vdots & \vdots \\ a_{21} & a_{22} \end{bmatrix} \middle| a_i \in Q; \ 1 \leq i \leq 22 \right\}$$

be a ring M is commutative ring with $\begin{bmatrix} 1 & 1 \\ 1 & 1 \\ \vdots & \vdots \\ 1 & 1 \end{bmatrix}$ as unit with respect to natural multiplication. M is not an integral domain. M has zero divisors and every element M is torsion free.

For consider $x = \begin{bmatrix} a_1 & a_2 \\ a_3 & a_4 \\ \vdots & \vdots \\ a_{21} & a_{22} \end{bmatrix} \in M.$

$$x^2 = \begin{bmatrix} a_1^2 & a_2^2 \\ a_3^2 & a_4^2 \\ \vdots & \vdots \\ a_{21}^2 & a_{22}^2 \end{bmatrix} \text{ and so on. } x^n = \begin{bmatrix} a_1^n & a_2^n \\ a_3^n & a_4^n \\ \vdots & \vdots \\ a_{21}^n & a_{22}^n \end{bmatrix}.$$

Thus every $x \in M$ is such that $x^n \neq [1]$ for any positive n.

*Example 3.22:* Let

$$P = \left\{ \begin{bmatrix} a_1 & a_2 & a_3 \\ a_4 & a_5 & a_6 \\ a_7 & a_8 & a_9 \end{bmatrix} \middle| a_i \in Q \text{ (or Z or R or C)}; \ 1 \leq i \leq 9 \right\}$$

be a commutative ring with unit under natural product.

$$M = \left\{ \begin{bmatrix} a_1 & a_2 & a_3 \\ 0 & a_4 & a_5 \\ 0 & 0 & a_6 \end{bmatrix} \middle| a_i \in Z; \ 1 \leq i \leq 6 \right\} \subseteq P;$$

M is a subring but M has no unit. M is not an ideal. M is of infinite order.

*Example 3.23:* Let

$$M = \left\{ \begin{bmatrix} a_1 & a_2 & a_3 & a_4 \\ 0 & a_5 & a_6 & a_7 \\ 0 & 0 & a_8 & a_9 \\ 0 & 0 & 0 & a_{10} \end{bmatrix} \middle| a_i \in Z; \ 1 \leq i \leq 10 \right\}$$

be a ring M is a commutative ring with no identity.

$$P = \left\{ \begin{bmatrix} a_1 & 0 & 0 & 0 \\ 0 & a_2 & a_3 & 0 \\ 0 & 0 & 0 & a_4 \\ 0 & 0 & 0 & 0 \end{bmatrix} \middle| a_i \in Z; \ 1 \leq i \leq 4 \right\} \subseteq M;$$

P is a subring and an ideal of M.

*Example 3.24:* Let

$$P = \left\{ \begin{bmatrix} a_1 & a_2 & a_3 \\ a_4 & a_5 & a_6 \\ \vdots & \vdots & \vdots \\ a_{31} & a_{32} & a_{33} \end{bmatrix} \middle| a_i \in Z;\ 1 \leq i \leq 33 \right\}$$

be a ring P has no units. P has ideals. P has subrings which are not ideals. P has no idempotents or nilpotents.

Every element x in P is such that for no $n \in Z^+$,

$$x^n = \begin{bmatrix} 1 & 1 & 1 \\ 1 & 1 & 1 \\ \vdots & \vdots & \vdots \\ 1 & 1 & 1 \end{bmatrix}.$$

**THEOREM 3.15:** *Let*

$$M = \left\{ \begin{bmatrix} a_{11} & \cdots & a_{1n} \\ a_{21} & \cdots & a_{2n} \\ \vdots & \vdots & \vdots \\ a_{m1} & \cdots & a_{mn} \end{bmatrix} \middle| a_{ij} \in Q;\ 1 \leq i \leq m;\ 1 \leq j \leq n \right\}$$

*be a ring M is a S-ring.*

*Proof:* Consider

$$P = \left\{ \begin{bmatrix} b & 0 & \cdots & 0 \\ 0 & 0 & \cdots & 0 \\ \vdots & \vdots & \vdots & \vdots \\ 0 & 0 & \cdots & 0 \end{bmatrix} \middle| b \in Q \right\} \in M;$$

P is a field. So M is a S-ring.

**COROLLARY 2:** *Every matrix ring under natural product is a S-ring.*

If Q is replaced by Z in the theorem and corollary then the matrix ring is not a S-ring.

*Example 3.25:* Let

$$S = \left\{ \begin{bmatrix} a_1 & a_2 & a_3 \\ a_4 & a_5 & a_6 \\ a_7 & a_8 & a_9 \\ a_{10} & a_{11} & a_{12} \end{bmatrix} \middle| a_i \in Z;\ 1 \le i \le 12 \right\}$$

be a ring, S is a S-ring.

*Example 3.26:* Let

$$M = \left\{ \begin{bmatrix} a_1 & a_2 \\ a_3 & a_4 \\ a_5 & a_6 \\ a_7 & a_8 \end{bmatrix} \middle| a_i \in Q;\ 1 \le i \le 8 \right\}$$

be a ring. M is a S-ring. For M has 8 subfields given by

$$F_1 = \left\{ \begin{bmatrix} a_1 & 0 \\ 0 & 0 \\ 0 & 0 \\ 0 & 0 \end{bmatrix} \middle| a_1 \in Q \right\} \subseteq M \text{ is a field.}$$

$$F_2 = \left\{ \begin{bmatrix} 0 & a_1 \\ 0 & 0 \\ 0 & 0 \\ 0 & 0 \end{bmatrix} \middle| a_1 \in Q \right\} \subseteq M \text{ is a field.}$$

$$F_3 = \left\{ \begin{bmatrix} 0 & 0 \\ a_1 & 0 \\ 0 & 0 \\ 0 & 0 \end{bmatrix} \middle| a_1 \in Q \right\} \subseteq M \text{ is a field.}$$

$$F_4 = \left\{ \begin{bmatrix} 0 & 0 \\ 0 & a_1 \\ 0 & 0 \\ 0 & 0 \end{bmatrix} \middle| a_1 \in Q \right\} \subseteq M \text{ is a field.}$$

$$F_5 = \left\{ \begin{bmatrix} 0 & 0 \\ 0 & 0 \\ a_1 & 0 \\ 0 & 0 \end{bmatrix} \middle| a_1 \in Q \right\} \subseteq M \text{ is a field.}$$

$$F_6 = \left\{ \begin{bmatrix} 0 & 0 \\ 0 & 0 \\ 0 & a_1 \\ 0 & 0 \end{bmatrix} \middle| a_1 \in Q \right\} \subseteq M \text{ is a field.}$$

$$F_7 = \left\{ \begin{bmatrix} 0 & 0 \\ 0 & 0 \\ 0 & 0 \\ a_1 & 0 \end{bmatrix} \middle| a_1 \in Q \right\} \subseteq M \text{ is a field and}$$

$$F_8 = \left\{ \begin{bmatrix} 0 & 0 \\ 0 & 0 \\ 0 & 0 \\ 0 & a_1 \end{bmatrix} \middle| a_1 \in Q \right\} \subseteq M \text{ is a field.}$$

Thus M has only 8 fields.

$$N = \left\{ \begin{bmatrix} a_1 & 0 \\ 0 & a_2 \\ 0 & 0 \\ 0 & 0 \end{bmatrix} \middle| a_1, a_2 \in Q \right\} \subseteq M;$$

N is a subring and an ideal and not a field. Thus M has only 8 fields. M is a S-ring. N also is a S-subring. However all subrings of M are not S-subrings.

For consider

$$S = \left\{ \begin{bmatrix} a_1 & a_2 \\ a_3 & a_4 \\ a_5 & a_6 \\ a_7 & a_8 \end{bmatrix} \middle| a_i \in Z; \ 1 \leq i \leq 8 \right\} \subseteq M;$$

S is only a subring and clearly S is not a S-subring. Infact M has infinite number of subrings which are not S-subrings.

*Example 3.27:* Let

$$W = \left\{ \begin{bmatrix} a_1 & a_2 & a_3 \\ a_4 & a_5 & a_6 \\ a_7 & a_8 & a_9 \end{bmatrix} \middle| a_i \in Q; \ 1 \leq i \leq 9 \right\}$$

be a ring under usual multiplication of matrices. W is a non commutative ring. However W is also a S-ring. W is a commutative ring under natural matrix multiplication, (W, +, $\times_n$) is also a S-ring. Both have zero divisors.

It is interesting to recall that by the natural matrix multiplication we are in a position to extend all the properties of reals into these matrix rings, (R, +, $\times_n$); the only difference being that these rings have zero divisors. So as blocks they do not loose any of the properties over which they are defined.

Further we can get compatability of natural product for both column and rectangular matrices. Hence we see these matrices under natural product can serve better purpose for they almost behave like the real numbers or complex number or rationals or integers on which they are built. Now we give more algebraic structure on them. Consider the set of all row matrices $M = \{(x_1, \ldots, x_n) \mid x_i \in R^+ \cup \{0\}$ (or $Q^+ \cup \{0\}$ or $Z^+ \cup \{0\}$); $1 \leq i \leq n\}$, M under + is a commutative semigroup with $(0, 0, \ldots, 0)$ as its additive identity.

M under $\times_n$ is also semigroup. Thus $(M, +, \times_n)$ is a semiring. We see this semiring is a commutative semiring with zero divisors.

Suppose

$S = \{(x_1, \ldots, x_n) \mid x_i \in R^+ \cup \{0\}$ (or $Z^+$ or $Q^+$); $1 \leq i \leq n\}$.

Now $\{S \cup \{(0, 0, \ldots, 0)\} = T, +, \times_n\}$ is a semifield.

It is easily verified T has no zero divisors and that T is a strict semiring for $a = (x_1, x_2, \ldots, x_n)$ and $b = (y_1, y_2, \ldots, y_n)$ is such that $x+y = 0$ implies $a = (0) = b = (0, 0, \ldots, 0)$. Now we will give examples of them before we proceed onto define and describe more properties.

*Example 3.28:* Let

$M = \{(a_1, a_2, a_3)$ where $a_i \in Z^+ \cup \{0\}$; $1 \leq i \leq 3\}$; $(M, +, \times_n)$

is a semiring. M is not a semifield as $a = (3, 0, 4)$ and $b = (0, 7, 0)$ in M are such that $a.b = (3, 0, 4)(0, 7, 0) = (0, 0, 0)$. However M is a strict commutative semiring which is not a semifield.

*Example 3.29:* Let

$$T = \left\{ \begin{bmatrix} a_1 \\ a_2 \\ a_3 \\ a_4 \\ a_5 \\ a_6 \end{bmatrix} \middle| a_i \in Q^+ \cup \{0\}; \ 1 \le i \le 6 \right\}$$

be a semiring under $+$ and $\times_n$.

We see if $x = \begin{bmatrix} 0 \\ 3 \\ 1 \\ 0 \\ 2 \\ 5 \end{bmatrix}$ and $y = \begin{bmatrix} 1 \\ 0 \\ 0 \\ 2 \\ 0 \\ 0 \end{bmatrix}$ are in T then

$$x \times_n y = \begin{bmatrix} 0 \\ 3 \\ 1 \\ 0 \\ 2 \\ 5 \end{bmatrix} \times_n \begin{bmatrix} 1 \\ 0 \\ 0 \\ 2 \\ 0 \\ 0 \end{bmatrix} = \begin{bmatrix} 0 \\ 0 \\ 0 \\ 0 \\ 0 \\ 0 \end{bmatrix}.$$

Thus T is only a commutative strict semiring and is not a semifield.

*Example 3.30:* Let

$$M = \left\{ \begin{bmatrix} a_1 & a_2 & a_3 \\ a_4 & a_5 & a_6 \\ \vdots & \vdots & \vdots \\ a_{13} & a_{14} & a_{15} \end{bmatrix} \middle| a_i \in R^+ \cup \{0\}; \ 1 \leq i \leq 15 \right\}$$

be a semiring under + and $\times_n$. Clearly M is commutative and is a strict semiring. However M does contain zero divisor, for if T

$$= \begin{bmatrix} a_1 & a_2 & a_3 \\ 0 & 0 & 0 \\ 0 & 0 & 0 \\ 0 & 0 & 0 \\ a_4 & a_5 & a_6 \end{bmatrix} \text{ and } N = \begin{bmatrix} 0 & 0 & 0 \\ a_1 & a_2 & a_3 \\ a_4 & a_5 & a_6 \\ a_7 & a_8 & a_9 \\ 0 & 0 & 0 \end{bmatrix} \text{ with } a_i \in R^+ \cup \{0\} \text{ are}$$

in M then

$$T \times_n N = \begin{bmatrix} a_1 & a_2 & a_3 \\ 0 & 0 & 0 \\ 0 & 0 & 0 \\ 0 & 0 & 0 \\ a_4 & a_5 & a_6 \end{bmatrix} \times_n \begin{bmatrix} 0 & 0 & 0 \\ a_1 & a_2 & a_3 \\ a_4 & a_5 & a_6 \\ a_7 & a_8 & a_9 \\ 0 & 0 & 0 \end{bmatrix} = \begin{bmatrix} 0 & 0 & 0 \\ 0 & 0 & 0 \\ 0 & 0 & 0 \\ 0 & 0 & 0 \\ 0 & 0 & 0 \end{bmatrix}.$$

Thus M is not a semifield.

*Example 3.31:* Let

$$J = \left\{ \begin{bmatrix} a_1 & a_2 & a_3 \\ a_4 & a_5 & a_6 \\ a_7 & a_8 & a_9 \end{bmatrix} \middle| a_i \in Q^+ \cup \{0\}; \ 1 \leq i \leq 9 \right\}$$

be a semiring under + and $\times_n$. J is a commutative strict semiring. However J is not a semifield.

For take $a = \begin{bmatrix} 0 & 0 & a_1 \\ 0 & a_2 & a_3 \\ a_4 & a_5 & a_6 \end{bmatrix}$ and $b = \begin{bmatrix} b_1 & b_2 & 0 \\ b_3 & 0 & 0 \\ 0 & 0 & 0 \end{bmatrix}$ in J,

we see $a.b = \begin{bmatrix} 0 & 0 & a_1 \\ 0 & a_2 & a_3 \\ a_4 & a_5 & a_6 \end{bmatrix} \begin{bmatrix} b_1 & b_2 & 0 \\ b_3 & 0 & 0 \\ 0 & 0 & 0 \end{bmatrix} = \begin{bmatrix} 0 & 0 & 0 \\ 0 & 0 & 0 \\ 0 & 0 & 0 \end{bmatrix}$.

Thus J is only a strict commutative semiring and is not a semifield, we show how we can build semifields.

First we will illustrate this situation by some examples.

*Example 3.32:* Let

$M = \{(0,0,0,0), (x_1, x_2, x_3, x_4) \mid x_i \in Q^+; 1 \le i \le 4\}$; $(M, +, \times_n)$

be a semifield. For we see $(M, +)$ is a commutative semigroup with additive identity $(0,0,0,0)$.

Further $(M, \times_n)$ is a commutative semigroup with $(1,1,1,1)$ as its multiplicative identity.

Also M is a strict semiring for $(a,b,c,d) + (x,y,z,t)$

$= (a + x, b + y, c + z, t + d)$

$= (0,0,0,0)$ if and only if each of a,b,c,d,x,y,z and t is zero.

Also for any $x = (a_1, a_2, a_3, a_4)$ and $y = (b_1, b_2, b_3, b_4)$ in M. We see $x.y = (a_1, a_2, a_3, a_4) \ (b_1, b_2, b_3, b_4) = (a_1 b_1, a_2 b_2, a_3 b_3, a_4 b_4)$ where $a_i b_i$ are in $Q^+$; $1 \le i \le 4$ so $x.y \ne (0,0,0,0)$. Thus $(M, \times_n)$ is a semifield. Thus we can get many semifields.

*Example 3.33:* Let

$$M = \left\{ \begin{bmatrix} a_1 \\ a_2 \\ \vdots \\ a_9 \\ a_{10} \end{bmatrix} \text{ where } a_i \in R^+, 1 \leq i \leq 10 \right\} \text{ and}$$

$$P = M \cup \left\{ \begin{bmatrix} 0 \\ 0 \\ \vdots \\ 0 \\ 0 \end{bmatrix} \right\}; (P, +, \times_n) \text{ is a semifield.}$$

*Example 3.34:* Let

$$S = \left\{ \begin{bmatrix} a_1 & a_2 & a_3 & a_4 \\ a_5 & a_6 & a_7 & a_8 \\ a_9 & a_{10} & a_{11} & a_{12} \\ a_{13} & a_{14} & a_{15} & a_{16} \\ a_{17} & a_{18} & a_{19} & a_{20} \end{bmatrix} \text{ where } a_i \in R^+, 1 \leq i \leq 20 \right\}$$

$$\text{and } P = S \cup \left\{ \begin{bmatrix} 0 & 0 & 0 & 0 \\ 0 & 0 & 0 & 0 \\ 0 & 0 & 0 & 0 \\ 0 & 0 & 0 & 0 \\ 0 & 0 & 0 & 0 \end{bmatrix} \right\}; (P, +, \times_n) \text{ is a semifield.}$$

*Example 3.35:* Let

$$T = \left\{ \begin{bmatrix} a_1 & a_2 \\ a_3 & a_4 \end{bmatrix} \middle| a_i \in Z^+, 1 \leq i \leq 4 \right\} \text{ and }$$

$$P = T \cup \left\{ \begin{bmatrix} 0 & 0 \\ 0 & 0 \end{bmatrix} \right\}; (P, +, \times_n) \text{ is a semifield.}$$

We see by defining natural product on matrices we get infinite number of semifields apart from $R^+ \cup \{0\}$, $Q^+ \cup \{0\}$ and $Z^+ \cup \{0\}$. We proceed onto give examples of Smarandache semirings. Recall a semiring S is a Smarandache semiring if S contains a proper subset T such that T under the operations of S is a semifield.

*Example 3.36:* Let
$$M = \{(a_1, a_2, \ldots, a_{10}) \mid a_i \in Q^+ \cup \{0\}, 1 \leq i \leq 10\}$$
be a semiring under $+$ and $\times_n$. Take

$$T = \{(0,a,0,\ldots,0) \mid a \in Q^+ \cup \{0\}\} \subseteq M;$$
T is a subsemiring of M. T is strict and T has no zero divisors, so T is a semifield under $+$ and $\times_n$. Hence M is a Smarandache semiring.

*Example 3.37:* Let
$$T = \left\{ \begin{bmatrix} a_1 \\ a_2 \\ a_3 \\ a_4 \\ a_5 \\ a_6 \\ a_7 \\ a_8 \end{bmatrix} \text{ where } a_i \in Z^+ \cup \{0\}, 1 \leq i \leq 8 \right\}$$
be a semiring under $+$ and $\times_n$.

Consider

$$P = \left\{ \begin{bmatrix} 0 \\ 0 \\ \vdots \\ 0 \\ a \end{bmatrix} \middle| a \in Z^+ \cup \{0\} \right\} \subseteq T.$$

P is a subsemiring of T which is strict and has no zero divisors. Thus T is a Smarandache semiring.

*Example 3.38:* Let

$$V = \left\{ \begin{bmatrix} a_1 & a_2 & a_3 & a_4 \\ a_5 & a_6 & a_7 & a_8 \\ a_9 & a_{10} & a_{11} & a_{12} \end{bmatrix} \text{ where } a_i \in Z^+ \cup \{0\}, 1 \leq i \leq 12 \right\}$$

be a semiring under + and $\times_n$.

V is a Smarandache semiring as

$$P = \left\{ \begin{bmatrix} a & 0 & 0 & 0 \\ 0 & 0 & 0 & 0 \\ 0 & 0 & 0 & 0 \end{bmatrix} \middle| a \in Z^+ \cup \{0\} \right\} \subseteq V;$$

is a semifield.

*Example 3.39:* Let

$$M = \left\{ \begin{bmatrix} a_1 & a_2 & a_3 & a_4 & a_5 \\ a_6 & a_7 & a_8 & a_9 & a_{10} \\ a_{11} & a_{12} & a_{13} & a_{14} & a_{15} \\ a_{16} & a_{17} & a_{18} & a_{19} & a_{20} \\ a_{21} & a_{22} & a_{23} & a_{24} & a_{25} \end{bmatrix} \right.$$

where $a_i \in R^+ \cup \{0\}$, $1 \le i \le 25\}$

be a semiring under $+$ and $\times_n$.

Consider

$$S = \left\{ \begin{bmatrix} 0 & 0 & 0 & 0 & 0 \\ 0 & 0 & a_1 & 0 & 0 \\ 0 & 0 & 0 & 0 & 0 \\ 0 & 0 & 0 & 0 & 0 \\ 0 & 0 & 0 & 0 & 0 \end{bmatrix} \middle| a \in Z^+ \cup \{0\} \right\} \subseteq M$$

is a semiring as well as a semifield under $+$ and $\times_n$. Hence M is a Smarandache semiring.

We can now define subsemirings and Smarandache subsemirings. These definitions are a matter of routine and hence left as an exercise to the reader. We however provide some examples of them.

*Example 3.40:* Let

$$P = \left\{ \begin{bmatrix} a_1 & a_2 & a_3 & a_4 \\ a_5 & a_6 & a_7 & a_8 \\ a_9 & a_{10} & a_{11} & a_{12} \end{bmatrix} \text{ where } a_i \in Z^+ \cup \{0\}, 1 \le i \le 12 \right\}$$

be a semiring under $+$ and $\times_n$.

Consider

$$X = \left\{ \begin{bmatrix} a_1 & a_2 & a_3 & a_4 \\ 0 & 0 & 0 & 0 \\ 0 & 0 & 0 & 0 \end{bmatrix} \middle| a \in 3Z^+ \cup \{0\}; 1 \leq i \leq 4 \right\} \subseteq P;$$

x under + and $\times_n$ is a subsemiring of P. However x is not a Smarandache subsemiring.

But we see P is a Smarandache semiring for

$$V = \left\{ \begin{bmatrix} d & 0 & 0 & 0 \\ 0 & 0 & 0 & 0 \\ 0 & 0 & 0 & 0 \end{bmatrix} \middle| d \in Z^+ \cup \{0\} \right\} \subseteq P$$

is a semiring as well as a semifield under + and $\times_n$.

*Example 3.41:* Let

$$P = \left\{ \begin{bmatrix} a_1 & a_2 & a_3 & a_4 \\ a_5 & a_6 & a_7 & a_8 \\ a_9 & a_{10} & a_{11} & a_{12} \\ a_{13} & a_{14} & a_{15} & a_{16} \end{bmatrix} \text{ where } a_i \in Q^+ \cup \{0\}, 1 \leq i \leq 16 \right\}$$

be a semiring under + and $\times_n$.

P is a Smarandache semiring for take

$$W = \left\{ \begin{bmatrix} a & 0 & 0 & 0 \\ 0 & 0 & 0 & 0 \\ 0 & 0 & 0 & 0 \\ 0 & 0 & 0 & 0 \end{bmatrix} \middle| a \in Z^+ \cup \{0\} \right\} \subseteq P$$

is a semifield under $+$ and $\times_n$. Hence P is a Smarandache semiring. However P has infinitely many subsemirings which are not Smarandache subsemirings.

Consider

$$P_n = \left\{ \begin{bmatrix} a_1 & a_2 & a_3 & 0 \\ 0 & 0 & 0 & 0 \\ 0 & 0 & 0 & 0 \\ 0 & 0 & 0 & 0 \end{bmatrix} \middle| a \in nZ; 1 \leq i \leq 3; n \geq 2 \right\} \subseteq P;$$

$P_n$ is a subsemiring of P but is not a Smarandache subsemiring of P. Thus we see in general all subsemirings of a Smarandache semiring need not be a Smarandache subsemiring. But if S be a semiring which has a Smarandache subsemiring then S is also a Smarandache semiring.

Inview of this we have the following theorem.

**THEOREM 3.16:** *Let S, be a semiring of $n \times m$ matrices with entries from $R^+ \cup \{0\}$ (or $Q^+ \cup \{0\}$ or $Z^+ \cup \{0\}$). If S has a subsemiring which is a Smarandache subsemiring then S is a Smarandache semiring. However if S is a Smarandache semiring then in general a subsemiring of S need not be a Smarandache subsemiring.*

*Proof:* Let S be a semiring and $P \subseteq S$ be a proper subsemiring of S, which is a Smarandache subsemiring of S. Since P is a Smarandache subsemiring of S we have a proper subset $X \subseteq P$ ($X \neq \phi$ and $X \neq P$) such that X is a semifield. Now we see $P \subseteq S$ and $X \subseteq P$ so $X \subseteq P \subseteq S$ that is X is a proper subset of S and X is a semifield, so S is a Smarandache semiring.

To show every subsemiring of a Smarandache semiring need not be a Smarandache subsemiring, we give an example.

Consider

$$P = \left\{ \begin{bmatrix} a_1 & a_2 & a_3 & a_4 \\ a_5 & a_6 & a_7 & a_8 \end{bmatrix} \text{ where } a_i \in Z^+ \cup \{0\}, 1 \leq i \leq 8 \right\}$$

be a semiring under + and $\times_n$.

If is easily verified P is a Smarandache semiring as

$$X = \left\{ \begin{pmatrix} a & 0 & 0 & 0 \\ 0 & 0 & 0 & 0 \end{pmatrix} \bigg| a \in Z^+ \cup \{0\} \right\} \subseteq P;$$

is a semifield under + and $\times_n$; so P is a Smarandache semiring.

Consider a subsemiring

$$T = \left\{ \begin{bmatrix} a_1 & a_2 & a_3 & a_4 \\ a_5 & a_6 & a_7 & a_8 \end{bmatrix} \text{ where } a_i \in 5Z^+ \cup \{0\}, 1 \leq i \leq 8 \right\} \subseteq P;$$

clearly T is a subsemiring of P; however T is not a Smarandache subsemiring of P, but we know P is a Smarandache semiring. Hence the result.

We will show the existence of zero divisors in semirings.

*Example 3.42:* Let

$$P = \left\{ \begin{bmatrix} a_1 & a_2 & a_3 \\ a_4 & a_5 & a_6 \\ a_7 & a_8 & a_9 \\ a_{10} & a_{11} & a_{12} \\ a_{13} & a_{14} & a_{15} \end{bmatrix} \text{ where } a_i \in Z^+ \cup \{0\}, 1 \le i \le 15 \right\}$$

be a semiring and $+$ and $\times_n$.

To show M has zero divisors.

$$\text{Consider } x = \begin{bmatrix} 0 & 0 & 0 \\ a_1 & a_2 & a_3 \\ 0 & 0 & 0 \\ a_4 & a_5 & a_6 \\ 0 & 0 & 0 \end{bmatrix} \text{ and } y = \begin{bmatrix} a_1 & a_2 & a_3 \\ 0 & 0 & 0 \\ a_4 & a_5 & a_6 \\ 0 & 0 & 0 \\ a_7 & a_8 & a_9 \end{bmatrix} \text{ in M.}$$

We see $x \times_n y = \begin{bmatrix} 0 & 0 & 0 \\ 0 & 0 & 0 \\ 0 & 0 & 0 \\ 0 & 0 & 0 \\ 0 & 0 & 0 \end{bmatrix}$.

Thus M has zero divisors. Infact M has infinitely many zero divisors.

For take $a = \begin{bmatrix} a_1 & a_2 & 0 \\ 0 & 0 & 0 \\ 0 & 0 & 0 \\ 0 & 0 & 0 \\ 0 & 0 & 0 \end{bmatrix}$ where $a_1, a_2 \in 10Z^+$ and

$$b = \begin{bmatrix} 0 & 0 & 0 \\ b_1 & b_2 & 0 \\ 0 & 0 & 0 \\ 0 & 0 & 0 \\ 0 & 0 & 0 \end{bmatrix} \text{ where } b_1, b_2 \in 3Z^+ \text{ in M,}$$

we see $a.b = \begin{bmatrix} 0 & 0 & 0 \\ 0 & 0 & 0 \\ 0 & 0 & 0 \\ 0 & 0 & 0 \\ 0 & 0 & 0 \end{bmatrix}$.

Thus we can get any number of zero divisors in M.

M has no idempotents other than elements of the form

$$x = \begin{bmatrix} 1 & 1 & 1 \\ 0 & 0 & 0 \\ 1 & 1 & 1 \\ 0 & 0 & 0 \\ 0 & 0 & 0 \end{bmatrix} \in M, \text{ we see } x^2 = x \text{ or}$$

$$y = \begin{bmatrix} 0 & 0 & 0 \\ 1 & 0 & 0 \\ 0 & 1 & 0 \\ 0 & 0 & 1 \\ 0 & 1 & 1 \end{bmatrix} \in M \text{ is such that } y^2 = y \text{ and so on.}$$

Inview of this we have the following nice theorem.

**THEOREM 3.17:** *Let $S = \{(a_{ij})_{m \times n} \mid a_{ij} \in Z^+ \cup \{0\}$ (or $Q^+ \cup \{0\}$ or $R^+ \cup \{0\}$); $1 \leq i \leq m$; $1 \leq j \leq n\}$ be a semiring under + and $\times_n$. All elements in S of the form $T = \{(a_{ij})_{m \times n} \mid a_{ij} \in \{0, 1\}\} \subseteq S$ are collection of idempotents.*

The proof is direct hence left as an exercise to the reader.

We call all these idempotents only as trivial or {1, 0} generated idempotents; apart from this these matrix semiring with natural product do not contain any other idempotents.

*Example 3.43:* Let M = {(a, b) | a, b ∈ $Z^+ \cup \{0\}$} be a semiring under + and $\times_n$. The only trivial idempotents of M are (0, 0), (0, 1), (1, 0) and (1, 1).

*Example 3.44:* Let P = $\left\{ \begin{bmatrix} a & b \\ c & d \end{bmatrix} \middle| a, b, c, d \in Z^+ \cup \{0\} \right\}$ be a semiring under + and $\times_n$.

The trivial idempotents of P are

$$\left\{ \begin{bmatrix} 0 & 0 \\ 0 & 0 \end{bmatrix}, \begin{bmatrix} 1 & 0 \\ 0 & 0 \end{bmatrix}, \begin{bmatrix} 0 & 1 \\ 0 & 0 \end{bmatrix}, \begin{bmatrix} 0 & 0 \\ 1 & 0 \end{bmatrix}, \begin{bmatrix} 0 & 0 \\ 0 & 1 \end{bmatrix}, \begin{bmatrix} 1 & 0 \\ 1 & 0 \end{bmatrix}, \begin{bmatrix} 1 & 1 \\ 0 & 0 \end{bmatrix}, \right.$$

$$\begin{bmatrix} 0 & 0 \\ 1 & 1 \end{bmatrix}, \begin{bmatrix} 0 & 1 \\ 0 & 1 \end{bmatrix}, \begin{bmatrix} 1 & 1 \\ 1 & 1 \end{bmatrix}, \begin{bmatrix} 1 & 0 \\ 1 & 1 \end{bmatrix},$$

$$\left. \begin{bmatrix} 0 & 1 \\ 1 & 1 \end{bmatrix}, \begin{bmatrix} 1 & 1 \\ 1 & 0 \end{bmatrix}, \begin{bmatrix} 1 & 1 \\ 0 & 1 \end{bmatrix}, \begin{bmatrix} 1 & 0 \\ 0 & 1 \end{bmatrix}, \begin{bmatrix} 0 & 1 \\ 1 & 0 \end{bmatrix} \right\} = I \subseteq P$$

we see I under natural product $\times_n$ is a semigroup. However I is not closed under +.

**THEOREM 3.18:** *Let*

$$M = \{(a_{ij}) \mid a_{ij} \in Z^+ \cup \{0\}; 1 \leq i \leq m; 1 \leq j \leq n\}$$

*be the collection of m × n matrices. The collection of all trivial idempotents forms a semigroup under $\times_n$ and the number of such trivial idempotents is $2^{m \times n}$.*

The proof involves only simple number theoretic techniques, hence left as an exercise to the reader.

*Example 3.45:* Let

$$M = \left\{ \begin{bmatrix} a_1 & a_2 & a_3 & a_4 \\ a_5 & a_6 & a_7 & a_8 \end{bmatrix} \text{ where } a_i \in Z^+ \cup \{0\}, 1 \leq i \leq 8 \right\}$$

be a semiring under + and $\times_n$.

$$I = \left\{ \begin{bmatrix} 0 & 0 & 0 & 0 \\ 0 & 0 & 0 & 0 \end{bmatrix}, \begin{bmatrix} 1 & 0 & 0 & 0 \\ 0 & 0 & 0 & 0 \end{bmatrix}, \begin{bmatrix} 0 & 1 & 0 & 0 \\ 0 & 0 & 0 & 0 \end{bmatrix}, \right.$$

$$\begin{bmatrix} 0 & 0 & 1 & 0 \\ 0 & 0 & 0 & 0 \end{bmatrix}, \begin{bmatrix} 0 & 0 & 0 & 1 \\ 0 & 0 & 0 & 0 \end{bmatrix},$$

$$\begin{bmatrix} 0 & 0 & 0 & 0 \\ 1 & 0 & 0 & 0 \end{bmatrix}, \begin{bmatrix} 0 & 0 & 0 & 0 \\ 0 & 1 & 0 & 0 \end{bmatrix}, \begin{bmatrix} 0 & 0 & 0 & 0 \\ 0 & 0 & 1 & 0 \end{bmatrix}, \begin{bmatrix} 0 & 0 & 0 & 0 \\ 0 & 0 & 0 & 1 \end{bmatrix},$$

$$\begin{bmatrix} 1 & 1 & 0 & 0 \\ 0 & 0 & 0 & 0 \end{bmatrix}, \begin{bmatrix} 0 & 1 & 1 & 0 \\ 0 & 0 & 0 & 0 \end{bmatrix}, \begin{bmatrix} 0 & 0 & 1 & 1 \\ 0 & 0 & 0 & 0 \end{bmatrix}, \begin{bmatrix} 1 & 0 & 0 & 0 \\ 1 & 0 & 0 & 0 \end{bmatrix}, ...,$$

$$\left. \begin{bmatrix} 1 & 1 & 1 & 1 \\ 1 & 1 & 1 & 0 \end{bmatrix}, \begin{bmatrix} 1 & 1 & 1 & 1 \\ 1 & 1 & 1 & 1 \end{bmatrix} \right\} \subseteq M;$$

is the collection of all trivial idempotents. Clearly they form a semigroup under product. However I is not closed under addition. Further the number of elements in I is $2^8$.

Further the semigroup I has zero divisors and $\begin{bmatrix} 1 & 1 & 1 & 1 \\ 1 & 1 & 1 & 1 \end{bmatrix}$ acts as the multiplicative identity.

$$x = \begin{bmatrix} 1 & 1 & 0 & 0 \\ 0 & 0 & 0 & 0 \end{bmatrix} \text{ and } y = \begin{bmatrix} 0 & 0 & 1 & 1 \\ 1 & 1 & 1 & 1 \end{bmatrix} \text{ is such that}$$

$$x \times_n y = \begin{bmatrix} 0 & 0 & 0 & 0 \\ 0 & 0 & 0 & 0 \end{bmatrix}.$$

Each element in $I \setminus \left\{ \begin{bmatrix} 1 & 1 & 1 & 1 \\ 1 & 1 & 1 & 1 \end{bmatrix} \right\}$ can generate in ideal of the semigroup.

For consider $x = \begin{bmatrix} 1 & 1 & 1 & 1 \\ 0 & 0 & 0 & 0 \end{bmatrix} \in I;$

$\langle x \rangle =$

$$\left\{ \begin{bmatrix} 0 & 0 & 0 & 0 \\ 0 & 0 & 0 & 0 \end{bmatrix}, \begin{bmatrix} 1 & 1 & 1 & 1 \\ 0 & 0 & 0 & 0 \end{bmatrix}, \begin{bmatrix} 0 & 1 & 0 & 1 \\ 0 & 0 & 0 & 0 \end{bmatrix}, \begin{bmatrix} 0 & 0 & 0 & 1 \\ 0 & 0 & 0 & 0 \end{bmatrix}, \right.$$

$$\begin{bmatrix} 1 & 0 & 0 & 0 \\ 0 & 0 & 0 & 0 \end{bmatrix}, \begin{bmatrix} 0 & 0 & 1 & 0 \\ 0 & 0 & 0 & 0 \end{bmatrix}, \begin{bmatrix} 0 & 0 & 1 & 1 \\ 0 & 0 & 0 & 0 \end{bmatrix}, \begin{bmatrix} 0 & 1 & 0 & 0 \\ 0 & 0 & 0 & 0 \end{bmatrix},$$

$$\begin{bmatrix} 1 & 1 & 0 & 0 \\ 0 & 0 & 0 & 0 \end{bmatrix}, \begin{bmatrix} 0 & 1 & 1 & 0 \\ 0 & 0 & 0 & 0 \end{bmatrix}, \begin{bmatrix} 1 & 0 & 1 & 0 \\ 0 & 0 & 0 & 0 \end{bmatrix}, \begin{bmatrix} 1 & 0 & 0 & 1 \\ 0 & 0 & 0 & 0 \end{bmatrix},$$

$$\left. \begin{bmatrix} 1 & 1 & 1 & 0 \\ 0 & 0 & 0 & 0 \end{bmatrix}, \begin{bmatrix} 0 & 1 & 1 & 1 \\ 0 & 0 & 0 & 0 \end{bmatrix}, \begin{bmatrix} 1 & 1 & 0 & 1 \\ 0 & 0 & 0 & 0 \end{bmatrix}, \begin{bmatrix} 1 & 0 & 1 & 1 \\ 0 & 0 & 0 & 0 \end{bmatrix} \right\} \subseteq I$$

is an ideal of the semigroup and order of the ideal generated by $\langle x \rangle$ is 16.

$$J = \left\langle x = \begin{bmatrix} 1 & 1 & 1 & 0 \\ 1 & 1 & 1 & 0 \end{bmatrix} \right\rangle$$

$$= \left\{ \begin{bmatrix} 0 & 0 & 0 & 0 \\ 0 & 0 & 0 & 0 \end{bmatrix}, \begin{bmatrix} 1 & 0 & 0 & 0 \\ 0 & 0 & 0 & 0 \end{bmatrix}, \begin{bmatrix} 0 & 0 & 0 & 0 \\ 1 & 0 & 0 & 0 \end{bmatrix}, \begin{bmatrix} 0 & 0 & 1 & 0 \\ 0 & 0 & 0 & 0 \end{bmatrix}, \right.$$

$$\begin{bmatrix} 0 & 0 & 0 & 0 \\ 0 & 0 & 1 & 0 \end{bmatrix}, \begin{bmatrix} 0 & 1 & 0 & 0 \\ 0 & 0 & 0 & 0 \end{bmatrix}, \begin{bmatrix} 0 & 0 & 0 & 0 \\ 0 & 1 & 0 & 0 \end{bmatrix}, \begin{bmatrix} 1 & 1 & 0 & 0 \\ 0 & 0 & 0 & 0 \end{bmatrix},$$

$$\begin{bmatrix} 0 & 1 & 1 & 0 \\ 0 & 0 & 0 & 0 \end{bmatrix}, \begin{bmatrix} 1 & 0 & 1 & 0 \\ 0 & 0 & 0 & 0 \end{bmatrix}, \begin{bmatrix} 0 & 0 & 0 & 0 \\ 1 & 1 & 0 & 0 \end{bmatrix}, \begin{bmatrix} 0 & 0 & 0 & 0 \\ 0 & 1 & 1 & 0 \end{bmatrix},$$

$$\begin{bmatrix} 0 & 0 & 0 & 0 \\ 1 & 0 & 1 & 0 \end{bmatrix}, \begin{bmatrix} 1 & 0 & 0 & 0 \\ 1 & 0 & 0 & 0 \end{bmatrix}, \begin{bmatrix} 0 & 1 & 0 & 0 \\ 0 & 1 & 0 & 0 \end{bmatrix}, \begin{bmatrix} 0 & 0 & 1 & 0 \\ 0 & 0 & 1 & 0 \end{bmatrix},$$

$$\begin{bmatrix} 0 & 0 & 1 & 0 \\ 0 & 0 & 1 & 0 \end{bmatrix}, \begin{bmatrix} 1 & 0 & 0 & 0 \\ 0 & 1 & 0 & 0 \end{bmatrix}, \begin{bmatrix} 0 & 1 & 0 & 0 \\ 1 & 0 & 0 & 0 \end{bmatrix}, \begin{bmatrix} 0 & 1 & 0 & 0 \\ 0 & 0 & 1 & 0 \end{bmatrix},$$

$$\left. \begin{bmatrix} 1 & 0 & 0 & 0 \\ 0 & 0 & 1 & 0 \end{bmatrix}, \begin{bmatrix} 0 & 0 & 1 & 0 \\ 1 & 0 & 0 & 0 \end{bmatrix}, \begin{bmatrix} 0 & 0 & 1 & 0 \\ 0 & 1 & 0 & 0 \end{bmatrix} \right.$$

and so on}.

We see order J is $2^6$. Thus every singleton other than $\{0\}$ and identity generate an ideal in the trivial idempotent semigroup.

Infact $\{0\}$ generates $\{0\}$ the trivial zero ideal and $\begin{bmatrix} 1 & 1 & 1 & 1 \\ 1 & 1 & 1 & 1 \end{bmatrix}$ generates the totality of the semigroup.

Now having seen the collection of trivial idempotents we proceed onto define other properties. Using the semifields we can build semivector spaces. More properties like zero divisors and Smarandache zero divisors are left as an exercise to the reader.

**Chapter Four**

# NATURAL PRODUCT ON MATRICES

In this chapter we construct semivector space over semifields and vector spaces over fields using these collection of matrices under natural product.

**DEFINITION 4.1:** *Let V be the collection of all $n \times m$ matrices with entries from Q (or R) or C. (V, +) is an abelian group. V is a vector space over Q (or R) according as V takes its entries from Q (or R). If V takes its entries from Q; V is not a vector space over R however if V takes its entries from R, V is a vector space over Q as well as vector spaces over R. We see all vector spaces V ($m \neq n$) are also linear algebras for using the natural product we get the linear algebra.*

We first illustrate this situation by some examples.

***Example 4.1:*** Let $V = \{(x_1, \ldots, x_5) \mid x_i \in Q; 1 \leq i \leq 5\}$ be a vector space over Q. Infact V is a linear algebra over Q. Clearly dimension of V is five and a basis for V is

$\{(1,0,0,0,0)\ (0,1,0,0,0)\ (0,0,1,0,0),\ (0,0,0,0,1),\ (0,0,0,1,0)\}$.

*Example 4.2:* Let

$$M = \left\{ \begin{bmatrix} a_1 & a_2 \\ a_3 & a_4 \\ a_5 & a_6 \\ a_7 & a_8 \\ a_9 & a_{10} \end{bmatrix} \middle| a_i \in R; 1 \leq i \leq 10 \right\}$$

be a vector space over Q. Clearly M is of infinite dimension.

*Example 4.3:* Let

$$P = \left\{ \begin{bmatrix} a_1 & a_2 \\ a_3 & a_4 \end{bmatrix} \middle| a_i \in R; 1 \leq i \leq 4 \right\}$$

be a vector space over R. Clearly dimension of P over R is four.

*Example 4.4:* Let

$$M = \left\{ \begin{bmatrix} a_1 \\ a_2 \\ \vdots \\ a_{20} \end{bmatrix} \middle| a_i \in Q; 1 \leq i \leq 20 \right\}$$

be a vector space over R. Clearly M is not a vector space over R. Clearly dimension of M over Q is 20.

*Example 4.5:* Let

$$T = \left\{ \begin{bmatrix} a_1 & a_2 & a_3 \\ a_4 & a_5 & a_6 \\ a_7 & a_8 & a_9 \end{bmatrix} \middle| a_i \in Q; 1 \leq i \leq 9 \right\}$$

be a vector space of dimension nine over Q. We see T is a linear algebra over Q under natural product as well as under the usual matrix product.

The concept of subspace is a matter of routine and hence is left as an exercise to the reader.

However we give examples of them.

*Example 4.6:* Let

$$M = \left\{ \begin{bmatrix} a_1 & a_2 \\ a_3 & a_4 \\ a_5 & a_6 \end{bmatrix} \middle| a_i \in Q; 1 \le i \le 6 \right\}$$

be a vector space over Q.

Consider

$$T = \left\{ \begin{bmatrix} a_1 & 0 \\ 0 & a_2 \\ a_3 & 0 \end{bmatrix} \middle| a_i \in Q; 1 \le i \le 3 \right\} \subseteq M;$$

it is easily verified T is a subspace of M over Q.

Consider

$$P = \left\{ \begin{bmatrix} a_1 & 0 \\ a_2 & 0 \\ a_3 & 0 \end{bmatrix} \middle| a_i \in Q; 1 \le i \le 3 \right\} \subseteq M;$$

P is also a subspace of M over Q. Now we consider $P \cap T$ (the intersection) of these two subspaces,

$$P \cap T = \left\{ \begin{bmatrix} a_1 & 0 \\ 0 & 0 \\ a_3 & 0 \end{bmatrix} \middle| a_i \in Q; i=1,3 \right\} \subseteq M;$$

$P \cap T$ is also a subspace of M over Q.

*Example 4.7:* Let

$$P = \left\{ \begin{bmatrix} a_1 & a_2 & a_3 \\ a_4 & a_5 & a_6 \\ a_7 & a_8 & a_9 \end{bmatrix} \middle| a_i \in Q; 1 \leq i \leq 9 \right\}$$

be a vector space Q.

Consider

$$M_1 = \left\{ \begin{bmatrix} a_1 & a_2 & 0 \\ 0 & 0 & 0 \\ 0 & 0 & a_3 \end{bmatrix} \middle| a_i \in Q; 1 \leq i \leq 3 \right\} \subseteq P,$$

$M_1$ is a subspace of P over Q.

Consider

$$M_2 = \left\{ \begin{bmatrix} 0 & 0 & a_1 \\ 0 & a_2 & 0 \\ 0 & 0 & 0 \end{bmatrix} \middle| a_1, a_2 \in Q \right\} \subseteq P$$

is also a subspace of P over Q.

However we see $M_1 \cap M_2 = \begin{bmatrix} 0 & 0 & 0 \\ 0 & 0 & 0 \\ 0 & 0 & 0 \end{bmatrix}$.

Consider

$$M_3 = \left\{ \begin{bmatrix} 0 & 0 & 0 \\ a_1 & 0 & a_2 \\ 0 & a_3 & 0 \end{bmatrix} \middle| a_i \in Q; 1 \leq i \leq 3 \right\} \subseteq P,$$

$M_3$ is a subspace of P over Q.

Take

$$M_4 = \left\{ \begin{bmatrix} 0 & 0 & 0 \\ 0 & 0 & 0 \\ a_1 & 0 & 0 \end{bmatrix} \middle| a_i \in Q \right\} \subseteq P$$

is also a subspace of P over Q.

We see $P = M_1 + M_2 + M_3 + M_4$ and $M_i \cap M_j = \begin{bmatrix} 0 & 0 & 0 \\ 0 & 0 & 0 \\ 0 & 0 & 0 \end{bmatrix}$

if $i \neq j$. Thus we can write P as a direct sum of subspaces of P.

*Example 4.8:* Let

$$P = \left\{ \begin{bmatrix} a_1 \\ a_2 \\ \vdots \\ a_{12} \end{bmatrix} \middle| a_i \in Q; 1 \leq i \leq 12 \right\}$$

be a vector space over Q.

Consider

$$X_1 = \left\{ \begin{bmatrix} a_1 \\ a_2 \\ 0 \\ \vdots \\ 0 \end{bmatrix} \middle| a_1, a_2 \in Q \right\} \subseteq P,$$

$X_1$ is a subspace of P over Q.

$$X_2 = \left\{ \begin{bmatrix} 0 \\ a_1 \\ a_2 \\ a_3 \\ 0 \\ \vdots \\ 0 \end{bmatrix} \middle| a_i \in Q; 1 \leq i \leq 3 \right\} \subseteq P$$

is again a subspace of P over Q.

Take

$$X_3 = \left\{ \begin{bmatrix} 0 \\ 0 \\ a_1 \\ a_2 \\ a_3 \\ a_4 \\ 0 \\ 0 \\ 0 \\ 0 \\ 0 \end{bmatrix} \middle| a_i \in Q; 1 \leq i \leq 4 \right\} \subseteq P$$

is again a subspace of P over Q.

Consider

$$X_4 = \left\{ \begin{bmatrix} 0 \\ 0 \\ 0 \\ 0 \\ 0 \\ a_1 \\ a_2 \\ a_3 \\ a_4 \\ a_5 \end{bmatrix} \middle| a_i \in Q; 1 \leq i \leq 5 \right\} \subseteq P$$

is a subspace of P over Q.

$$\text{We see } X_i \cap X_j \neq \begin{bmatrix} 0 \\ 0 \\ 0 \\ \vdots \\ 0 \\ 0 \end{bmatrix} \text{ if } i \neq j.$$

Thus P is not a direct sum. However we see

$P \subseteq X_1 + X_2 + X_3 + X_4$, thus we say P is only a pseudo direct sum of subspaces of P over Q.

Thus we have seen examples of direct sum and pseudo direct sum of subspaces. Interested reader can supply with more examples of them. Our main motivation is to define Smarandache strong vector spaces.

It is important to mention that usual matrix vector space over the fields Q or R are not that interesting except for the fact if V is set of $n \times 1$ column matrices then V is a vector space over Q or R but V is never a linear algebra under matrix multiplication, however V is a linear algebra under the natural matrix product $\times_n$. This is the vital difference and importance of defining natural product $\times_n$ of matrices of same order.

Now we define special strong Smarandache vector space.

**DEFINITION 4.2:** *Let*

$$M = \left\{ \begin{bmatrix} a_1 \\ \vdots \\ a_n \end{bmatrix} \middle| a_i \in Q; 1 \leq i \leq n \right\}.$$

*We define M as a natural Smarandache special field of characteristic zero under usual addition of matrices and the natural product $\times_n$. Thus $\{M, +, \times_n\}$ is natural Smarandache special field.*

We give an example or two.

*Example 4.9:* Let

$$V = \left\{ \begin{bmatrix} a_1 \\ a_2 \\ \vdots \\ a_7 \end{bmatrix} \middle| a_i \in Q; 1 \leq i \leq 7 \right\}$$

is a natural Smarandache special field of characteristic zero.

Consider

$$x = \begin{bmatrix} 1 \\ 7 \\ -9 \\ 2 \\ -7 \\ 4 \end{bmatrix} \in V \text{ then } x^{-1} = \begin{bmatrix} 1 \\ 1/7 \\ -1/9 \\ 1/2 \\ -1/7 \\ 1/4 \end{bmatrix} \in V \text{ and } x.x^{-1} = \begin{bmatrix} 1 \\ 1 \\ 1 \\ 1 \\ 1 \\ 1 \\ 1 \end{bmatrix}.$$

Thus $\begin{bmatrix} 1 \\ 1 \\ 1 \\ 1 \\ 1 \\ 1 \\ 1 \end{bmatrix}$ acts as the multiplicative identity.

*Example 4.10:* Consider the collection of $7 \times 1$ column matrices V with entries from Q.

We see $\left\{ \begin{bmatrix} a \\ 0 \\ \vdots \\ 0 \end{bmatrix} \;\middle|\; a \in Q \right\} = P$ is a proper subset of V which is a field hence natural S-special field.

*Example 4.11:* Let

$$M = \left\{ \begin{bmatrix} a_1 \\ a_2 \\ a_3 \\ a_4 \\ a_5 \end{bmatrix} \;\middle|\; a_i \in Q;\ 1 \leq i \leq 5 \right\}$$

be a natural S-special field of characteristic zero. M is a column matrix of natural Smarandache special field.

Now if we consider
$S = \{(x_1, x_2, \ldots, x_n) \mid x_i \in Q \text{ (or R) } 1 \leq i \leq n\}$.
S under usual addition of row matrices and natural matrix product $\times_n$ is a natural Smarandache special field called the special rational natural Smarandache special field of characteristic zero.

***Example 4.12:*** Let V = {($x_1$, $x_2$, $x_3$, $x_4$, $x_5$) | $x_i \in$ R, $1 \le i \le 4$} be the special real natural Smarandache special field of column matrices of characteristic zero.

***Example 4.13:*** Let V = {($x_1$, $x_2$) | $x_i \in$ R, $1 \le i \le 2$} be special row matrix natural Smarandache special field of characteristic zero.

All these fields are non prime natural Smarandache special fields for they have several natural S-special subfields.

Now we can define the natural Smarandache special field of m × n matrices (m ≠ n).

Let

$$V = \left\{ \begin{bmatrix} a_{11} & a_{12} & \ldots & a_{1n} \\ a_{21} & a_{22} & \ldots & a_{2n} \\ \vdots & \vdots & & \vdots \\ a_{m1} & a_{m2} & \ldots & a_{mn} \end{bmatrix} \middle| a_{ij} \in R, 1 \le i \le m, 1 \le j \le n \right\}$$

V is the special m × n matrix of natural Smarandache special field of characteristic zero.

***Example 4.14:*** Let

$$M = \left\{ \begin{bmatrix} a_1 & a_2 & a_3 & a_4 \\ a_5 & a_6 & a_7 & a_8 \\ a_9 & a_{10} & a_{11} & a_{12} \\ a_{13} & a_{14} & a_{15} & a_{16} \\ a_{17} & a_{18} & a_{19} & a_{20} \end{bmatrix} \middle| a_{ij} \in R, 1 \le i \le 5, 1 \le j \le 4 \right\}$$

be the special 5 × 4 matrix of natural Smarandache special field of characteristic zero.

*Example 4.15:* Let

$$M = \left\{ \begin{bmatrix} a_1 & a_2 \\ a_3 & a_4 \\ \vdots & \vdots \\ a_{21} & a_{22} \end{bmatrix} \middle| a_i \in R, 1 \le i \le 22 \right\}$$

be the $11 \times 2$ matrix of natural Smarandache special field of characteristic zero.

*Example 4.16:* Let

$$M = \left\{ \begin{bmatrix} a_1 & a_2 & \ldots & a_{16} \\ a_{17} & a_{18} & \ldots & a_{32} \end{bmatrix} \middle| a_i \in R, 1 \le i \le 32 \right\}$$

be the $2 \times 16$ matrix of natural Smarandache special field.

Now having seen natural S-special fields of $m \times n$ matrices ($m \ne n$). We now proceed onto define the notion of natural special Smarandache field of square matrices.

Let

$$P = \left\{ \begin{bmatrix} a_{11} & a_{12} & \ldots & a_{1n} \\ a_{21} & a_{22} & \ldots & a_{2n} \\ \vdots & \vdots & & \vdots \\ a_{n1} & a_{n2} & \ldots & a_{nn} \end{bmatrix} \middle| a_{ij} \in R, 1 \le i, j \le n \right\}$$

be a square matrix natural Smarandache special field of characteristic zero.

We give examples of them.

*Example 4.17:* Let

$$M = \left\{ \begin{bmatrix} a_1 & a_2 & a_3 \\ a_4 & a_5 & a_6 \\ a_7 & a_8 & a_9 \end{bmatrix} \middle| a_i \in R, 1 \leq i \leq 9 \right\}$$

be the $3 \times 3$ square matrix of natural special Smarandache field.

*Example 4.18:* Let

$$M = \left\{ \begin{bmatrix} a_1 & a_2 & a_3 & a_4 & a_5 \\ a_6 & a_7 & a_8 & a_9 & a_{10} \\ \vdots & \vdots & \vdots & \vdots & \vdots \\ a_{21} & a_{22} & a_{23} & a_{24} & a_{25} \end{bmatrix} \middle| a_i \in R, 1 \leq i \leq 25 \right\}$$

be the $5 \times 5$ square matrix of natural special Smarandache field of characteristic zero.

Now having seen and defined the concept of matrix natural special Smarandache field we are in a position to define natural Smarandache special strong matrix vector spaces.

**DEFINITION 4.3:** *Let*

$$V = \{(x_1, x_2, \ldots, x_n) \mid x_i \in Q \text{ (or } R\text{)}, 1 \leq i \leq n\}$$
*be an additive abelian group.*

$$F_R = \{(a_1, \ldots, a_n) \mid a_i \in Q, 1 \leq i \leq n\}$$
*be the natural special row matrix Smarandache special field. We see for every $x = (a_1, \ldots, a_n) \in F_R$ and $v = (x_1, \ldots, x_n) \in V$. $xv = vx \in V$; Further $(x+y)v = xv + yv$ and $x(v+u) = xv + xu$ for all $x, y \in F_R$ and $v, u \in V$.*

*Finally $(1, 1, \ldots, 1)v = v \in V$ for $(1, 1, \ldots, 1)$ multiplicative identity under the natural product $\times_n$. Thus V is a Smarandache vector space over $F_R$ known as the Smarandache special strong*

*row vector space over the natural row matrix Smarandache special field $F_R$.*

First we proceed onto give a few examples of them.

***Example 4.19:*** Let $M = \{(x_1, x_2, x_3, x_4) \mid x_i \in Q; 1 \leq i \leq 4\}$ be a Smarandache special strong row vector space over the natural Smarandache special row matrix field

$$F_R = \{(x_1, x_2, x_3, x_4) \mid x_i \in Q; 1 \leq i \leq 4\}.$$

***Example 4.20:*** Let $P = \{(x_1, x_2, x_3, \ldots, x_{10}) \mid x_i \in R; 1 \leq i \leq 10\}$ be a Smarandache special strong row vector space over the natural special row matrix Smarandache field

$$F_R = \{(x_1, x_2, \ldots, x_{10}) \mid x_i \in Q; 1 \leq i \leq 10\}.$$

***Example 4.21:*** Let $T = \{(x_1, x_2, x_3, \ldots, x_7) \mid x_i \in R; 1 \leq i \leq 10\}$ be a Smarandache special strong row vector space over the natural special row matrix Smarandache field

$$F_R = \{(x_1, x_2, \ldots, x_7) \mid x_i \in Q; 1 \leq i \leq 10\}.$$

Now we proceed onto define natural S-special strong column matrix vector space over the special column matrix natural S-special field $F_C$.

**DEFINITION 4.4:** *Let*

$$V = \left\{ \begin{bmatrix} x_1 \\ x_2 \\ \vdots \\ x_n \end{bmatrix} \middle| x_i \in Q \text{ (or } R\text{)},\ 1 \leq i \leq n \right\}$$

*be an addition abelian group. Let*

$$F_C = \left\{ \begin{bmatrix} a_1 \\ a_2 \\ \vdots \\ a_n \end{bmatrix} \middle| a_i \in Q \text{ (or R)}, 1 \leq i \leq n \right\}$$

*be the special column matrix natural S-field. Clearly V is a S-vector space over the natural Smarandache special field $F_C$, we define V as a S-special strong column matrix vector space over the special column matrix natural S-field $F_C$.*

We will illustrate this situation by some examples.

*Example 4.22:* Let

$$V = \left\{ \begin{bmatrix} x_1 \\ x_2 \\ \vdots \\ x_{10} \end{bmatrix} \middle| x_i \in Q, 1 \leq i \leq 10 \right\}$$

is a S-special strong column matrix vector space over the special column matrix natural S-field

$$F_C = \left\{ \begin{bmatrix} x_1 \\ x_2 \\ \vdots \\ x_{10} \end{bmatrix} \middle| x_i \in Q \text{ (or R)}, 1 \leq i \leq 10 \right\}.$$

*Example 4.23:* Let

$$V = \left\{ \begin{bmatrix} a_1 \\ a_2 \\ a_3 \\ a_4 \\ a_5 \end{bmatrix} \middle| a_i \in R;\ 1 \leq i \leq 5 \right\}$$

be a S-special strong column matrix vector space over the special column matrix natural S-field

$$F_C = \left\{ \begin{bmatrix} x_1 \\ x_2 \\ x_3 \\ x_4 \\ x_5 \end{bmatrix} \middle| x_i \in R;\ 1 \leq i \leq 5 \right\}.$$

*Example 4.24:* Let

$$V = \left\{ \begin{bmatrix} x_1 \\ x_2 \end{bmatrix} \middle| x_i \in R;\ 1 \leq i \leq 2 \right\}$$

be a S-special strong column matrix vector space over the special column matrix natural S-field

$$F_C = \left\{ \begin{bmatrix} x_1 \\ x_2 \end{bmatrix} \middle| x_1, x_2 \in Q \right\}.$$

Now we proceed onto define the notion of Smarandache special strong $m \times n$ ($m \neq n$) matrix vector space over the special $m \times n$ matrix natural Smarandache special field $F_{m \times n}$ ($m \neq n$).

Let

$$M = \left\{ \begin{bmatrix} a_{11} & a_{12} & \cdots & a_{1n} \\ a_{21} & a_{22} & \cdots & a_{2n} \\ \vdots & \vdots & & \vdots \\ a_{m1} & a_{m2} & \cdots & a_{mn} \end{bmatrix} \middle| a_{ij} \in Q \text{ (or R)}, \right.$$

$$1 \le i \le m, 1 \le j \le n; m \ne n \}$$

be a group under matrix addition. Define

$$F_{m \times n} (m \ne n) = \left\{ \begin{bmatrix} a_{11} & a_{12} & \cdots & a_{1n} \\ a_{21} & a_{22} & \cdots & a_{2n} \\ \vdots & \vdots & & \vdots \\ a_{m1} & a_{m2} & \cdots & a_{mn} \end{bmatrix} \right.$$

$$a_{ij} \in Q \text{ (or R)}, 1 \le i \le m, 1 \le j \le n \}$$

to be special m × n matrix natural S-special field. Now we see M is a vector space over $F_{m \times n}$ called the S-special strong m × n matrix vector space over the special m × n matrix natural S-special field $F_{m \times n}$.

We will illustrate this situation by an example or two.

*Example 4.25:* Let

$$P = \left\{ \begin{bmatrix} a_1 & a_2 & a_3 \\ a_4 & a_5 & a_6 \\ a_7 & a_8 & a_9 \\ a_{10} & a_{11} & a_{12} \\ a_{13} & a_{14} & a_{15} \\ a_{16} & a_{17} & a_{18} \end{bmatrix} \middle| a_i \in Q; 1 \le i \le 18 \right\}$$

be a S-special strong 6 × 3 matrix vector space over the special 6 × 3 matrix natural special S-field

$$F_{6\times3} = \left\{ \begin{bmatrix} a_1 & a_2 & a_3 \\ a_4 & a_5 & a_6 \\ \vdots & \vdots & \vdots \\ a_{16} & a_{17} & a_{18} \end{bmatrix} \middle| \; a_i \in Q; \; 1 \leq i \leq 18 \right\}.$$

*Example 4.26:* Let

$$V = \left\{ \begin{bmatrix} a_1 & a_2 & a_3 & a_4 & a_5 & a_6 & a_7 & a_8 \\ a_9 & a_{10} & a_{11} & a_{12} & a_{13} & a_{14} & a_{15} & a_{16} \\ a_{17} & a_{18} & a_{19} & a_{20} & a_{21} & a_{22} & a_{23} & a_{24} \end{bmatrix} \middle| \; a_i \in R; \right.$$

$$1 \leq i \leq 24\}$$

be a S-special strong $3 \times 8$ matrix vector space over the special $3 \times 8$ natural special matrix S-field

$$F_{3\times8} = \left\{ \begin{bmatrix} a_1 & a_2 & \ldots & a_8 \\ a_9 & a_{10} & \ldots & a_{16} \\ a_{17} & a_{18} & \ldots & a_{24} \end{bmatrix} \middle| \; a_i \in Q; \; 1 \leq i \leq 24 \right\}.$$

*Example 4.27:* Let

$$V = \left\{ \begin{bmatrix} a_1 & a_2 \\ a_3 & a_4 \\ a_5 & a_6 \\ a_7 & a_8 \end{bmatrix} \middle| \; a_i \in R; \; 1 \leq i \leq 8 \right\}$$

be a S-special strong $4 \times 2$ matrix vector space over the special $4 \times 2$ matrix natural special S-field

$$F_{4\times 2} = \left\{ \begin{bmatrix} a_1 & a_2 \\ a_3 & a_4 \\ a_5 & a_6 \\ a_7 & a_8 \end{bmatrix} \middle| \; a_i \in R; \; 1 \leq i \leq 8 \right\}.$$

Now finally we define the S-special strong square matrix vector space over the special square matrix natural special S-field $F_{n\times n}$.

$$F_{n\times n} = \left\{ \begin{bmatrix} a_{11} & \cdots & a_{1n} \\ a_{21} & \cdots & a_{2n} \\ \vdots & & \vdots \\ a_{n1} & \cdots & a_{nn} \end{bmatrix} \middle| \; a_{ij} \in R \text{ (or Q)}; \; 1 \leq i, j \leq n \right\}.$$

We just define this structure.

$$\text{Let } M = \left\{ \begin{bmatrix} a_{11} & a_{12} & \cdots & a_{1n} \\ a_{21} & a_{22} & \cdots & a_{2n} \\ \vdots & \vdots & & \vdots \\ a_{n1} & a_{n2} & \cdots & a_{nn} \end{bmatrix} \middle| \; a_{ij} \in Q \text{ (or R)}, \; 1 \leq i, j \leq n \right\}$$

be the group under addition of square matrices.

Let

$$F_{n\times n} = \left\{ \begin{bmatrix} a_{11} & a_{12} & \cdots & a_{1n} \\ a_{21} & a_{22} & \cdots & a_{2n} \\ \vdots & \vdots & & \vdots \\ a_{n1} & a_{n2} & \cdots & a_{nn} \end{bmatrix} \middle| \; a_{ij} \in R, \; 1 \leq i, j \leq n \right\}$$

be the S-special square matrix field M is a vector space over $F_{n\times n}$. M is defined as the S-strong special square n×n matrix vector space over the special square matrix natural special S-field $F_{n\times n}$.

We will illustrate this situation by some simple examples.

*Example 4.28:* Let

$$M = \left\{ \begin{bmatrix} a_1 & a_2 & a_3 \\ a_4 & a_5 & a_6 \\ a_7 & a_8 & a_9 \end{bmatrix} \middle| a_i \in Q; 1 \leq i \leq 9 \right\}$$

be a special strong square $3 \times 3$ matrix S-vector space over the special $3 \times 3$ square matrix natural special S-field.

$$F_{3\times3} = \left\{ \begin{bmatrix} a_1 & a_2 & a_3 \\ a_4 & a_5 & a_6 \\ a_7 & a_8 & a_9 \end{bmatrix} \middle| a_i \in Q; 1 \leq i \leq 9 \right\}.$$

*Example 4.29:* Let

$$V = \left\{ \begin{bmatrix} a_1 & a_2 & a_3 & a_4 & a_5 \\ a_6 & a_7 & a_8 & a_9 & a_{10} \\ \vdots & \vdots & \vdots & \vdots & \vdots \\ a_{21} & a_{22} & a_{23} & a_{24} & a_{25} \end{bmatrix} \middle| a_i \in R, 1 \leq i \leq 25 \right\}$$

be a S-special strong $5 \times 5$ square matrix vector space over the special $5 \times 5$ natural special Smarandache field.

$$F_{5\times5} = \left\{ \begin{bmatrix} a_1 & a_2 & a_3 & a_4 & a_5 \\ a_6 & a_7 & a_8 & a_9 & a_{10} \\ \vdots & \vdots & \vdots & \vdots & \vdots \\ a_{21} & a_{22} & a_{23} & a_{24} & a_{25} \end{bmatrix} \middle| a_i \in R, 1 \leq i \leq 25 \right\}.$$

*Example 4.30:* Let

$$A = \left\{ \begin{bmatrix} a_1 & a_2 \\ a_3 & a_4 \end{bmatrix} \middle| a_i \in R; 1 \le i \le 4 \right\}$$

be a S-special strong 2×2 square matrix vector space over the special 2×2 square matrix natural special S-field

$$F_{2\times 2} = \left\{ \begin{bmatrix} a_1 & a_2 \\ a_3 & a_4 \end{bmatrix} \middle| a_i \in Q; 1 \le i \le 4 \right\}.$$

Now seen various types of S-special vector spaces; we now proceed onto define S-subspaces over natural special S-fields.

**DEFINITION 4.5:** *Let V be a S-strong special row matrix (or column matrix or m × n matrix (m≠n) or square matrix) vector space over the special row matrix natural S-field $F_R$ (or $F_C$ or $F_{n\times n}$ (n ≠ m) or $F_{n\times n}$).*

*Consider W ⊆ V (W a proper subset of V); if W itself is a S-strong special row matrix (or column matrix or m × n matrix (m ≠ n) or square matrix) S-vector space over $F_R$ (or $F_C$ or $F_{m\times n}$ or $F_{n\times n}$) then we define W to be a S-special strong row matrix (column matrix or m × n matrix or square matrix) vector space of V over $F_R$ (or $F_C$ or $F_{m\times n}$ or $F_{n\times n}$).*

We will illustrate this situation by some simple examples.

*Example 4.31:* Let V = {($a_1$, $a_2$, $a_3$) | $a_i$ ∈ Q; 1 ≤ i ≤ 3} be a S-special row matrix vector space over

$F_R$ = {($x_1$, $x_2$, $x_3$) | $x_i$ ∈ Q, 1 ≤ i ≤ 3} the natural special S-row field. Consider M = {($a_1$, 0, 0) | $a_1$ ∈ Q} ⊆ V; M is a S-special strong row matrix vector subspace of V over $F_R$.

*Example 4.32:* Let

$$V = \left\{ \begin{bmatrix} a_1 \\ a_2 \\ \vdots \\ a_6 \end{bmatrix} \middle| a_i \in R; 1 \leq i \leq 6 \right\}$$

be a S-special strong vector space over the S-field;

$$F_C = \left\{ \begin{bmatrix} a_1 \\ a_2 \\ \vdots \\ a_6 \end{bmatrix} \middle| a_i \in Q; 1 \leq i \leq 6 \right\}.$$

Consider

$$M = \left\{ \begin{bmatrix} a_1 \\ 0 \\ a_2 \\ 0 \\ a_3 \\ 0 \end{bmatrix} \middle| a_i \in R; 1 \leq i \leq 3 \right\} \subseteq V;$$

M is a Smarandache special strong vector subspace of V over the S-field $F_C$.

Take

$$P = \left\{ \begin{bmatrix} a_1 \\ a_2 \\ 0 \\ 0 \\ 0 \\ 0 \end{bmatrix} \middle| a_1, a_2 \in R \right\} \subseteq V;$$

P is a Smarandache special strong vector subspace of V over the S-field $F_C$.

*Example 4.33:* Let

$$M = \left\{ \begin{bmatrix} a_1 & a_2 & a_3 & a_4 \\ a_5 & a_6 & a_7 & a_8 \\ a_9 & a_{10} & a_{11} & a_{12} \end{bmatrix} \middle| a_i \in Q; 1 \leq i \leq 12 \right\}$$

be a S-special strong vector space over the S-field

$$F_{3\times 4} = \left\{ \begin{bmatrix} a_1 & a_2 & a_3 & a_4 \\ a_5 & a_6 & a_7 & a_8 \\ a_9 & a_{10} & a_{11} & a_{12} \end{bmatrix} \middle| a_i \in Q; 1 \leq i \leq 12 \right\} \subseteq M,$$

P is a S-special strong vector subspace of M over the S-field $F_{3\times 4}$.

*Example 4.34:* Let $V = \{(a_1, a_2, \ldots, a_9) \mid a_i \in R; 1 \leq i \leq 9\}$ be a Smarandache special strong vector space over the S-field,
$$F_R = \{(a_1, a_2, \ldots, a_9) \mid a_i \in Q; 1 \leq i \leq 9\}.$$

Consider $M_1 = \{(a_1, 0, a_2, 0, \ldots, 0) \mid a_1, a_2 \in R\} \subseteq V$, $M_1$ is a S-special strong vector subspace of V over the S-field $F_R$.

Consider

$M_2 = \{(0, a_1, 0, a_2, 0, \ldots, 0) \mid a_i \in R; 1 \leq i \leq 2\} \subseteq V$, $M_2$ is a again a S-special strong vector subspace of V over $F_R$.

$M_3 = \{(0, 0, 0, 0, a_1, a_2, 0, 0, 0) \mid a_1, a_2 \in R\} \subseteq V$, $M_3$ is a again a S-special strong vector subspace of V over $F_R$.

$M_4 = \{(0, 0, 0, 0, 0, 0, a_1, a_2, a_3) \mid a_i \in R; 1 \leq i \leq 3\} \subseteq V$, $M_4$ is again a S-special strong vector subspace of V over $F_R$.

It is easily verified $M_i \cap M_j = (0, 0, 0, \ldots, 0)$ if $i \neq j$, $1 \leq i, j \leq 4$ and $V = M_1 + M_2 + M_3 + M_4$. Thus V is the direct sum of S-strong vector subspaces.

*Example 4.35:* Let

$$M = \left\{ \begin{bmatrix} a_1 \\ a_2 \\ \vdots \\ a_{11} \\ a_{12} \end{bmatrix} \middle| a_i \in Q; 1 \leq i \leq 12 \right\}$$

be a S-special strong vector space over the S-field,

$$F_C = \left\{ \begin{bmatrix} a_1 \\ a_2 \\ \vdots \\ a_{11} \\ a_{12} \end{bmatrix} \middle| a_i \in Q; 1 \leq i \leq 12 \right\}.$$

Clearly dimension of M is also 12. Consider the following S-special strong vector subspaces.

$$P_1 = \left\{ \begin{bmatrix} a_1 \\ 0 \\ \vdots \\ 0 \\ a_2 \end{bmatrix} \middle| a_i \in Q; 1 \leq i \leq 2 \right\} \subseteq M,$$

a S-strong vector subspace of M.

$$P_2 = \left\{ \begin{bmatrix} 0 \\ a_1 \\ 0 \\ \vdots \\ 0 \\ a_2 \\ 0 \end{bmatrix} \middle| a_i \in Q; 1 \leq i \leq 2 \right\} \subseteq M,$$

a S-strong vector subspace of M over $F_C$.

$$P_3 = \left\{ \begin{bmatrix} 0 \\ 0 \\ a_1 \\ a_2 \\ a_3 \\ a_4 \\ 0 \\ \vdots \\ 0 \end{bmatrix} \middle| a_i \in Q; 1 \leq i \leq 4 \right\} \subseteq M,$$

is a S-strong special vector subspace of V over $F_C$.

$$P_4 = \left\{ \begin{bmatrix} 0 \\ 0 \\ 0 \\ 0 \\ 0 \\ 0 \\ a_1 \\ a_2 \\ 0 \\ 0 \\ 0 \\ 0 \end{bmatrix} \middle| a_i \in Q; 1 \leq i \leq 2 \right\} \subseteq M,$$

is again a S-strong special vector subspace of V over $F_C$ and

$$P_5 = \left\{ \begin{bmatrix} 0 \\ 0 \\ 0 \\ 0 \\ 0 \\ 0 \\ 0 \\ 0 \\ a_1 \\ a_2 \\ 0 \\ 0 \end{bmatrix} \middle| a_i \in Q; 1 \leq i \leq 2 \right\} \subseteq M$$

is again a S-strong special vector subspace of V over $F_C$.

We see $P_i \cap P_j = \begin{bmatrix} 0 \\ 0 \\ \vdots \\ 0 \\ 0 \end{bmatrix}$ if $i \neq j$; $1 \leq i, j \leq 5$ and

$V = P_1 + P_2 + P_3 + P_4 + P_5$.

Thus V is a direct sum of S-special strong vector subspaces.

*Example 4.36:* Let

$$M = \left\{ \begin{bmatrix} a_1 & a_2 & a_3 \\ a_4 & a_5 & a_6 \\ a_7 & a_8 & a_9 \\ a_{10} & a_{11} & a_{12} \end{bmatrix} \middle| a_i \in R; 1 \leq i \leq 12 \right\}$$

be a S-strong special vector space over the S-field,

$$F_{4 \times 3} = \left\{ \begin{bmatrix} a_1 & a_2 & a_3 \\ a_4 & a_5 & a_6 \\ a_7 & a_8 & a_9 \\ a_{10} & a_{11} & a_{12} \end{bmatrix} \middle| a_i \in R; 1 \leq i \leq 12 \right\}.$$

Let $B_1 = \left\{ \begin{bmatrix} a_1 & a_2 & 0 \\ 0 & 0 & 0 \\ 0 & 0 & 0 \\ 0 & a_3 & a_4 \end{bmatrix} \middle| a_i \in R; 1 \leq i \leq 4 \right\} \subseteq M$

be a S-strong special subvector space of M over $F_{4 \times 3}$.

$$B_2 = \left\{ \begin{bmatrix} a_1 & 0 & a_2 \\ 0 & a_3 & a_4 \\ 0 & 0 & 0 \\ 0 & 0 & 0 \end{bmatrix} \middle| a_i \in R;\ 1 \le i \le 4 \right\} \subseteq M$$

be a S-strong special vector subspace of M over $F_{4 \times 3}$.

$$B_3 = \left\{ \begin{bmatrix} a_1 & 0 & 0 \\ a_2 & 0 & 0 \\ a_3 & 0 & 0 \\ 0 & 0 & 0 \end{bmatrix} \middle| a_i \in R;\ 1 \le i \le 3 \right\} \subseteq M$$

is a S-strong special vector subspace of M over $F_{4 \times 3}$.

Finally

$$B_4 = \left\{ \begin{bmatrix} a_1 & 0 & 0 \\ 0 & 0 & 0 \\ 0 & a_2 & a_3 \\ a_4 & 0 & 0 \end{bmatrix} \middle| a_i \in R;\ 1 \le i \le 4 \right\} \subseteq M$$

is again a S-strong special vector subspace of V over $F_{4 \times 3}$.

We see $B_i \cap B_j \ne \begin{bmatrix} 0 & 0 & 0 \\ 0 & 0 & 0 \\ 0 & 0 & 0 \\ 0 & 0 & 0 \end{bmatrix}$; even if $i \ne j$, $1 \le i, j \le 4$.

But $V \subseteq B_1 + B_2 + B_3 + B_4$; so we define V to be the pseudo direct sum of the S-special strong vector subspaces of V. Now we have seen examples of direct sum and pseudo direct sum of S-special strong vector subspaces of a S-special strong vector space over a S-field.

The main use of this structure will be found as and when this sort of study becomes familiar and in due course of time they may find applications in all places where the result is not a real number (or a rational number) but an array of numbers.

We can define orthogonal vectors of S-special strong matrix vector spaces also.

First we see how orthogonal vector matrices are defined when they are defined over R or Q or C.

Let $V_R = \{(a_1, \ldots, a_n) \mid a_i \in Q; 1 \leq i \leq n\}$ be a row matrix vector space defined over the field Q.

We define for any $x = (a_1, \ldots, a_n)$ and $y = (b_1, \ldots, b_n)$ in $V_R$, x is perpendicular to y if $x \times_n y = (0)$.

Thus if $V_R = \{(x_1, x_2, x_3, x_4, x_5) \mid x_i \in R; 1 \leq i \leq 5\}$ be a row matrix vector space defined over Q and if $x = (0, 4, -5, 0, 7)$ and $y = (1, 0, 0, 8, 0)$ are in $V_R$. We see $x \times_n y = (0)$ so x is orthogonal with y.

$$V_C = \left\{ \begin{bmatrix} a_1 \\ \vdots \\ a_n \end{bmatrix} \middle| a_i \in Q, \text{(or R)}; 1 \leq i \leq n \right\} \text{ be the vector space}$$

of column matrices over Q (or R) respectively.

We say two elements $x = \begin{bmatrix} x_1 \\ x_2 \\ \vdots \\ x_n \end{bmatrix}$ and $y = \begin{bmatrix} y_1 \\ y_2 \\ \vdots \\ y_n \end{bmatrix}$ in $V_C$ are

orthogonal if $x \times_n y = \begin{bmatrix} x_1 \\ x_2 \\ \vdots \\ x_n \end{bmatrix} \times_n \begin{bmatrix} y_1 \\ y_2 \\ \vdots \\ y_n \end{bmatrix} = (0)$.

For instance if $x = \begin{bmatrix} 2 \\ -1 \\ 0 \\ 0 \\ 0 \\ 7 \end{bmatrix}$ and $y = \begin{bmatrix} 0 \\ 0 \\ 1 \\ 2 \\ 3 \\ 0 \end{bmatrix}$ then

$$x \times_n y = \begin{bmatrix} 2 \\ -1 \\ 0 \\ 0 \\ 0 \\ 7 \end{bmatrix} \times_n \begin{bmatrix} 0 \\ 0 \\ 1 \\ 2 \\ 3 \\ 0 \end{bmatrix} = \begin{bmatrix} 0 \\ 0 \\ 0 \\ 0 \\ 0 \\ 0 \end{bmatrix}.$$

Now we can define, unlike in other matrix vector space in case of these vector spaces $V_{m \times n}$ ($m \neq n$) and $V_{n \times n}$ the orthogonality under natural product. This is a special feature enjoyed only by vector spaces on which natural product can be defined.

We just illustrate this situation by some examples.

*Example 4.37:* Let

$$V_{5 \times 3} = \left\{ \begin{bmatrix} a_1 & a_2 & a_3 \\ a_4 & a_5 & a_6 \\ a_7 & a_8 & a_9 \\ a_{10} & a_{11} & a_{12} \\ a_{13} & a_{14} & a_{15} \end{bmatrix} \middle| a_i \in R; 1 \leq i \leq 15 \right\}$$

be a $5 \times 3$ matrix linear algebra (vector space) over Q.

Now let

$$x = \begin{bmatrix} 3 & 2 & 0 \\ 0 & 1 & 5 \\ 1 & 1 & 0 \\ 2 & 0 & 7 \\ 0 & 1 & 8 \end{bmatrix} \text{ and } y = \begin{bmatrix} 0 & 0 & 7 \\ 9 & 0 & 0 \\ 0 & 0 & 9 \\ 0 & 8 & 0 \\ 4 & 0 & 0 \end{bmatrix} \text{ be in } V_{5 \times 3},$$

we see $x \times_n y = \begin{bmatrix} 0 & 0 & 0 \\ 0 & 0 & 0 \\ 0 & 0 & 0 \\ 0 & 0 & 0 \\ 0 & 0 & 0 \end{bmatrix}$ thus we say x is orthogonal to y under natural product in $V_{5 \times 3}$.

It is pertinent to mention here that we can have several y's in $V_{5 \times 3}$ such that for a given x in $V_{5 \times 3}$. $x \times_n y = (0)$.

Now we see all elements in $V_{5 \times 3}$ are orthogonal under natural product to the zero $5 \times 3$ matrix $\begin{bmatrix} 0 & 0 & 0 \\ 0 & 0 & 0 \\ 0 & 0 & 0 \\ 0 & 0 & 0 \\ 0 & 0 & 0 \end{bmatrix}$.

*Example 4.38:* Let

$$V_{2 \times 6} = \left\{ \begin{bmatrix} a_1 & a_2 & a_3 & a_4 & a_5 & a_6 \\ a_7 & a_8 & a_9 & a_{10} & a_{11} & a_{12} \end{bmatrix} \middle| a_i \in Q; 1 \leq i \leq 12 \right\}$$

be a vector space over Q. Under natural product defined on $V_{2 \times 6}$ we can get orthogonal elements.

Take $x = \begin{bmatrix} a_1 & a_3 & a_5 & 0 & 0 & 0 \\ a_2 & a_4 & a_6 & 0 & 0 & 0 \end{bmatrix}$ and

$$y = \begin{bmatrix} 0 & 0 & 0 & a_1 & 0 & a_3 \\ 0 & 0 & 0 & a_2 & 0 & a_4 \end{bmatrix} \text{ in } V_{2 \times 6}. \ a_i \in Q; \ 1 \leq i \leq 6.$$

We see $x \times_n y = \begin{bmatrix} 0 & 0 & 0 & 0 & 0 & 0 \\ 0 & 0 & 0 & 0 & 0 & 0 \end{bmatrix}$.

Thus x is orthogonal with y. Infact we have several such y's which are orthogonal with x.

*Example 4.39:* Let

$$V = \left\{ \begin{bmatrix} a_1 & a_2 \\ a_3 & a_4 \end{bmatrix} \middle| a_i \in R; \ 1 \leq i \leq 4 \right\}$$

be a $2 \times 2$ matrix vector space over the field R. $\times_n$ be the natural product on V.

We define two matrices in V to be orthogonal if

$$x \times_n y = \begin{bmatrix} 0 & 0 \\ 0 & 0 \end{bmatrix} \text{ for } y, x \in V. \text{ We see } x = \begin{bmatrix} a_1 & 0 \\ 0 & a_2 \end{bmatrix},$$

then $y = \begin{bmatrix} 0 & b_1 \\ 0 & 0 \end{bmatrix}$ is orthogonal with x.

Also $\begin{bmatrix} 0 & 0 \\ b_1 & 0 \end{bmatrix} = a$ is such that $x \times_n a = (0)$.

Further $\begin{bmatrix} 0 & b_1 \\ b_2 & 0 \end{bmatrix}$ = b, is orthogonal with x under natural product as $x \times_n b = \begin{bmatrix} 0 & 0 \\ 0 & 0 \end{bmatrix}$.

Now $x^\perp = \left\{ \begin{pmatrix} 0 & 0 \\ 0 & 0 \end{pmatrix}, \begin{pmatrix} 0 & b_1 \\ b_2 & 0 \end{pmatrix}, \begin{pmatrix} 0 & b_1 \\ 0 & 0 \end{pmatrix}, \begin{pmatrix} 0 & 0 \\ b_1 & 0 \end{pmatrix} \right\}$. Clearly $x^\perp$ is additively closed and also $\times_n$ product also is closed; infact $x^\perp$ is a proper subspace of V defined as a subspace perpendicular with x.

Consider $x = \begin{bmatrix} 0 & a \\ b & 0 \end{bmatrix}$ in V; now the elements perpendicular with x are $\left\{ \begin{pmatrix} 0 & 0 \\ 0 & 0 \end{pmatrix}, \begin{pmatrix} t & 0 \\ 0 & 0 \end{pmatrix}, \begin{pmatrix} 0 & 0 \\ 0 & u \end{pmatrix}, \begin{pmatrix} t & 0 \\ 0 & u \end{pmatrix} \right\}$.

We see this is also a subspace of V.

Infact if

$$B = \left\{ \begin{bmatrix} 0 & 0 \\ 0 & 0 \end{bmatrix}, \begin{bmatrix} t & 0 \\ 0 & 0 \end{bmatrix}, \begin{bmatrix} 0 & 0 \\ 0 & u \end{bmatrix}, \begin{bmatrix} x & 0 \\ 0 & y \end{bmatrix} \right\} \subseteq V$$

and

$$C = \left\{ \begin{bmatrix} 0 & 0 \\ 0 & 0 \end{bmatrix}, \begin{bmatrix} 0 & a \\ 0 & 0 \end{bmatrix}, \begin{bmatrix} 0 & 0 \\ b & 0 \end{bmatrix}, \begin{bmatrix} 0 & a \\ b & 0 \end{bmatrix} \right\} \subseteq V$$

are such that they are orthogonal subspaces under product.

For $B \times_n C = \left\{ \begin{bmatrix} 0 & 0 \\ 0 & 0 \end{bmatrix} \right\}$.

*Example 4.40:* Let

$$M = \left\{ \begin{bmatrix} a_1 & a_2 & a_3 \\ a_4 & a_5 & a_6 \\ a_7 & a_8 & a_9 \end{bmatrix} \middle| a_i \in Q; 1 \leq i \leq 9 \right\}$$

be a $3 \times 3$ vector space over the field Q. Consider an element,

$$x = \begin{bmatrix} a & b & c \\ 0 & 0 & 0 \\ 0 & 0 & d \end{bmatrix} \text{ in M.}$$

The elements perpendicular to x be denoted by
$$x^\perp = B =$$

$$\left\{ \begin{bmatrix} 0 & 0 & 0 \\ 0 & 0 & 0 \\ 0 & 0 & 0 \end{bmatrix}, \begin{bmatrix} 0 & 0 & 0 \\ 0 & 0 & a_1 \\ 0 & 0 & 0 \end{bmatrix}, \begin{bmatrix} 0 & 0 & 0 \\ a_2 & 0 & 0 \\ 0 & 0 & 0 \end{bmatrix}, \begin{bmatrix} 0 & 0 & 0 \\ 0 & b & 0 \\ 0 & 0 & 0 \end{bmatrix}, \right.$$

$$\begin{bmatrix} 0 & 0 & 0 \\ 0 & 0 & 0 \\ d & 0 & 0 \end{bmatrix}, \begin{bmatrix} 0 & 0 & 0 \\ 0 & 0 & 0 \\ 0 & e & 0 \end{bmatrix}, \begin{bmatrix} 0 & 0 & 0 \\ a & b & 0 \\ 0 & 0 & 0 \end{bmatrix}, \begin{bmatrix} 0 & 0 & 0 \\ a & 0 & b \\ 0 & 0 & 0 \end{bmatrix}, \begin{bmatrix} 0 & 0 & 0 \\ 0 & a & b \\ 0 & 0 & 0 \end{bmatrix},$$

$$\begin{bmatrix} 0 & 0 & 0 \\ 0 & 0 & a \\ b & 0 & 0 \end{bmatrix}, \begin{bmatrix} 0 & 0 & 0 \\ 0 & 0 & a \\ 0 & b & 0 \end{bmatrix}, \begin{bmatrix} 0 & 0 & 0 \\ a & 0 & 0 \\ b & 0 & 0 \end{bmatrix}, \begin{bmatrix} 0 & 0 & 0 \\ a & 0 & 0 \\ 0 & b & 0 \end{bmatrix}, \begin{bmatrix} 0 & 0 & 0 \\ 0 & a & 0 \\ b & 0 & 0 \end{bmatrix},$$

$$\left. \begin{bmatrix} 0 & 0 & 0 \\ 0 & a & 0 \\ 0 & b & 0 \end{bmatrix}, \begin{bmatrix} 0 & 0 & 0 \\ 0 & 0 & 0 \\ a & b & 0 \end{bmatrix}, \begin{bmatrix} 0 & 0 & 0 \\ a & b & c \\ 0 & 0 & 0 \end{bmatrix}, \begin{bmatrix} 0 & 0 & 0 \\ a & b & 0 \\ c & 0 & 0 \end{bmatrix}, \begin{bmatrix} 0 & 0 & 0 \\ a & b & 0 \\ 0 & c & 0 \end{bmatrix}, \right.$$

$$\left[\begin{array}{ccc}0&0&0\\0&a&b\\c&0&0\end{array}\right], \left[\begin{array}{ccc}0&0&0\\0&a&b\\0&c&0\end{array}\right], \left[\begin{array}{ccc}0&0&0\\a&0&b\\c&0&0\end{array}\right], \left[\begin{array}{ccc}0&0&0\\a&0&b\\0&c&0\end{array}\right], \left[\begin{array}{ccc}0&0&0\\a&0&0\\b&c&0\end{array}\right],$$

$$\left[\begin{array}{ccc}0&0&0\\0&a&0\\b&c&0\end{array}\right], \left[\begin{array}{ccc}0&0&0\\0&0&a\\b&c&0\end{array}\right], \left[\begin{array}{ccc}0&0&0\\a&b&c\\d&b&0\end{array}\right], \left[\begin{array}{ccc}0&0&0\\a&b&c\\d&0&0\end{array}\right], \left[\begin{array}{ccc}0&0&0\\a&b&c\\0&d&0\end{array}\right],$$

$$\left.\left[\begin{array}{ccc}0&0&0\\0&a&b\\c&d&0\end{array}\right], \left[\begin{array}{ccc}0&0&0\\a&0&b\\d&c&0\end{array}\right], \left[\begin{array}{ccc}0&0&0\\a&b&0\\d&c&0\end{array}\right]\right\} \subseteq M$$

is again a subspace orthogonal with x.

Inview of this we have the following theorem.

**THEOREM 4.1:** *Let V be a matrix vector space over a field F. Let $0 \neq x \in V$. The elements orthogonal to x under natural product $\times_n$ is a subspace of V over F.*

***Proof:*** Let $0 \neq x \in V$. Suppose $x^{\perp} = B = \{y \in V \mid y \times_n x = 0\}$, to show B is a subspace of V.

Clearly $B \subseteq V$, by the very definition of orthogonal element of x. Further $(0) \in V$ is such that $x \times_n (0) = (0)$ so $(0)$ is orthogonal with x. Now let $z, y \in B$ be orthogonal with x under $\times_n$. To show $y + z$ is orthogonal to x. Given $x \times_n y = (0)$ and $x \times_n z = (0)$. Now consider

$$\begin{aligned}x \times_n (y+z) &= x \times_n y + x \times_n z\\ &= (0) + (0)\\ &= (0);\end{aligned}$$

so y + z ∈ B. Also if y ∈ B is such that x ×$_n$ y = 0 then x ×$_n$ (–y) = 0 so if y ∈ B then –y ∈ B. Finally let a ∈ F and y ∈ B to show ay ∈ B. Consider x ×$_n$ ay = a (x ×$_n$ y) = a.0 = 0.

Thus B ⊆ V is a vector subspace of V.

Now we can define orthogonality of two elements in matrix vector spaces under natural product, we now derive various properties associated with orthogonal natural product.

*Example 4.41:* Let

$$W = \left\{ \begin{bmatrix} a & b & c \\ d & e & f \\ g & h & i \end{bmatrix} \middle| a, b, c, d, e, f, g, h, i \in Q \right\}$$

be a vector space over Q.

Consider $x^\perp$ = B =

$$\left\{ \begin{bmatrix} 0 & 0 & 0 \\ 0 & 0 & 0 \\ 0 & 0 & 0 \end{bmatrix}, \begin{bmatrix} 0 & 0 & 0 \\ 0 & 0 & 0 \\ a & 0 & 0 \end{bmatrix}, \begin{bmatrix} 0 & 0 & 0 \\ 0 & 0 & 0 \\ 0 & b & 0 \end{bmatrix}, \begin{bmatrix} 0 & 0 & 0 \\ 0 & 0 & 0 \\ 0 & 0 & c \end{bmatrix}, \right.$$

$$\left. \begin{bmatrix} 0 & 0 & 0 \\ 0 & 0 & 0 \\ a & b & 0 \end{bmatrix}, \begin{bmatrix} 0 & 0 & 0 \\ 0 & 0 & 0 \\ a & 0 & b \end{bmatrix}, \begin{bmatrix} 0 & 0 & 0 \\ 0 & 0 & 0 \\ 0 & a & b \end{bmatrix}, \begin{bmatrix} 0 & 0 & 0 \\ 0 & 0 & 0 \\ a & b & c \end{bmatrix} \right\} \subseteq W$$

is a subspace of V. Consider the complementary space of B.

$$B^\perp = \left\{ \begin{bmatrix} 0 & 0 & 0 \\ 0 & 0 & 0 \\ 0 & 0 & 0 \end{bmatrix}, \begin{bmatrix} a & 0 & 0 \\ 0 & 0 & 0 \\ 0 & 0 & 0 \end{bmatrix}, \begin{bmatrix} 0 & b & 0 \\ 0 & 0 & 0 \\ 0 & 0 & 0 \end{bmatrix}, \begin{bmatrix} 0 & 0 & c \\ 0 & 0 & 0 \\ 0 & 0 & c \end{bmatrix}, \right.$$

$$\begin{bmatrix} 0 & 0 & 0 \\ 0 & 0 & 0 \\ 0 & 0 & 0 \end{bmatrix}, \begin{bmatrix} 0 & 0 & 0 \\ d & 0 & 0 \\ 0 & 0 & 0 \end{bmatrix}, \begin{bmatrix} 0 & 0 & 0 \\ 0 & e & 0 \\ 0 & 0 & 0 \end{bmatrix}, \begin{bmatrix} 0 & 0 & 0 \\ 0 & 0 & f \\ 0 & 0 & 0 \end{bmatrix}, \begin{bmatrix} a & b & 0 \\ 0 & 0 & 0 \\ 0 & 0 & 0 \end{bmatrix},$$

$$\begin{bmatrix} 0 & a & b \\ 0 & 0 & 0 \\ 0 & 0 & 0 \end{bmatrix}, \begin{bmatrix} a & 0 & b \\ 0 & 0 & 0 \\ 0 & 0 & 0 \end{bmatrix}, \begin{bmatrix} a & 0 & 0 \\ b & 0 & 0 \\ 0 & 0 & 0 \end{bmatrix}, \begin{bmatrix} 0 & 0 & 0 \\ a & b & 0 \\ 0 & 0 & 0 \end{bmatrix}, \begin{bmatrix} 0 & 0 & 0 \\ a & 0 & c \\ 0 & 0 & 0 \end{bmatrix},$$

$$\begin{bmatrix} 0 & 0 & 0 \\ 0 & a & b \\ 0 & 0 & 0 \end{bmatrix}, \begin{bmatrix} a & 0 & 0 \\ 0 & b & 0 \\ 0 & 0 & 0 \end{bmatrix}, \begin{bmatrix} a & 0 & 0 \\ 0 & 0 & b \\ 0 & 0 & 0 \end{bmatrix}, \begin{bmatrix} 0 & a & 0 \\ b & 0 & 0 \\ 0 & 0 & 0 \end{bmatrix}, \begin{bmatrix} 0 & a & 0 \\ 0 & b & 0 \\ 0 & 0 & 0 \end{bmatrix},$$

$$\begin{bmatrix} 0 & a & 0 \\ 0 & 0 & b \\ 0 & 0 & 0 \end{bmatrix}, \begin{bmatrix} 0 & 0 & a \\ b & 0 & 0 \\ 0 & 0 & 0 \end{bmatrix}, \begin{bmatrix} 0 & 0 & a \\ 0 & b & 0 \\ 0 & 0 & 0 \end{bmatrix}, \begin{bmatrix} 0 & 0 & a \\ 0 & 0 & b \\ 0 & 0 & 0 \end{bmatrix}, \begin{bmatrix} a & b & c \\ 0 & 0 & 0 \\ 0 & 0 & 0 \end{bmatrix},$$

$$\begin{bmatrix} 0 & 0 & 0 \\ a & b & c \\ 0 & 0 & 0 \end{bmatrix}, \begin{bmatrix} a & 0 & 0 \\ 0 & b & c \\ 0 & 0 & 0 \end{bmatrix}, \begin{bmatrix} a & 0 & 0 \\ b & c & 0 \\ 0 & 0 & 0 \end{bmatrix}, \begin{bmatrix} a & 0 & b \\ 0 & c & 0 \\ 0 & 0 & 0 \end{bmatrix}, \begin{bmatrix} a & 0 & 0 \\ b & 0 & c \\ 0 & 0 & 0 \end{bmatrix},$$

$$\begin{bmatrix} a & 0 & b \\ c & 0 & 0 \\ 0 & 0 & 0 \end{bmatrix}, \begin{bmatrix} a & 0 & b \\ 0 & 0 & c \\ 0 & 0 & 0 \end{bmatrix}, \begin{bmatrix} a & 0 & 0 \\ b & c & 0 \\ 0 & 0 & 0 \end{bmatrix}, \begin{bmatrix} a & 0 & 0 \\ 0 & b & c \\ 0 & 0 & 0 \end{bmatrix}, \begin{bmatrix} a & 0 & 0 \\ b & 0 & c \\ 0 & 0 & 0 \end{bmatrix},$$

$$\begin{bmatrix} 0 & a & 0 \\ b & c & 0 \\ 0 & 0 & 0 \end{bmatrix}, \begin{bmatrix} 0 & a & 0 \\ 0 & b & c \\ 0 & 0 & 0 \end{bmatrix}, \begin{bmatrix} 0 & a & 0 \\ b & 0 & c \\ 0 & 0 & 0 \end{bmatrix}, \begin{bmatrix} 0 & 0 & a \\ b & c & 0 \\ 0 & 0 & 0 \end{bmatrix}, \begin{bmatrix} 0 & 0 & a \\ b & 0 & c \\ 0 & 0 & 0 \end{bmatrix},$$

$$\begin{bmatrix} 0 & 0 & a \\ 0 & b & c \\ 0 & 0 & 0 \end{bmatrix}, \begin{bmatrix} 0 & a & b \\ c & 0 & 0 \\ 0 & 0 & 0 \end{bmatrix}, \begin{bmatrix} 0 & a & b \\ 0 & c & 0 \\ 0 & 0 & 0 \end{bmatrix}, \begin{bmatrix} 0 & a & b \\ 0 & 0 & c \\ 0 & 0 & 0 \end{bmatrix}, \begin{bmatrix} a & b & c \\ 0 & 0 & 0 \\ 0 & 0 & 0 \end{bmatrix},$$

$$\begin{bmatrix} a & b & c \\ 0 & d & 0 \\ 0 & 0 & 0 \end{bmatrix}, \begin{bmatrix} a & b & c \\ 0 & 0 & d \\ 0 & 0 & 0 \end{bmatrix}, \begin{bmatrix} 0 & a & b \\ d & 0 & c \\ 0 & 0 & 0 \end{bmatrix}, \begin{bmatrix} 0 & a & b \\ d & c & 0 \\ 0 & 0 & 0 \end{bmatrix}, \begin{bmatrix} 0 & a & b \\ 0 & d & c \\ 0 & 0 & 0 \end{bmatrix},$$

$$\begin{bmatrix} a & 0 & 0 \\ b & c & d \\ 0 & 0 & 0 \end{bmatrix}, \begin{bmatrix} 0 & a & 0 \\ b & c & d \\ 0 & 0 & 0 \end{bmatrix}, \begin{bmatrix} 0 & 0 & a \\ b & c & d \\ 0 & 0 & 0 \end{bmatrix}, \begin{bmatrix} a & b & 0 \\ c & d & 0 \\ 0 & 0 & 0 \end{bmatrix}, \begin{bmatrix} a & 0 & d \\ b & 0 & c \\ 0 & 0 & 0 \end{bmatrix},$$

$$\begin{bmatrix} a & 0 & b \\ 0 & c & d \\ 0 & 0 & 0 \end{bmatrix}, \begin{bmatrix} a & 0 & b \\ c & d & 0 \\ 0 & 0 & 0 \end{bmatrix}, \begin{bmatrix} a & b & 0 \\ c & 0 & d \\ 0 & 0 & 0 \end{bmatrix}, \begin{bmatrix} a & b & c \\ 0 & d & e \\ 0 & 0 & 0 \end{bmatrix}, \begin{bmatrix} a & b & c \\ d & e & 0 \\ 0 & 0 & 0 \end{bmatrix},$$

$$\left. \begin{bmatrix} 0 & a & b \\ c & d & e \\ 0 & 0 & 0 \end{bmatrix}, \begin{bmatrix} a & b & 0 \\ c & d & e \\ 0 & 0 & 0 \end{bmatrix}, \begin{bmatrix} a & 0 & b \\ d & e & c \\ 0 & 0 & 0 \end{bmatrix}, \begin{bmatrix} a & b & c \\ d & e & f \\ 0 & 0 & 0 \end{bmatrix} \right\}$$

We see $B^\perp \oplus B = W$.

Thus we have the following theorem.

**THEOREM 4.2:** *Let V be a matrix vector space over a field F. Suppose $0 \neq x \in V$ and let W be the subspace of V perpendicular to V then the complement of W denoted by $W^\perp$ is such that for every $a \in W$ and $b \in W^\perp$  $a \times_n b = (0)$. Further V $= W^\perp \oplus W$.*

The proof is simple hence left as an exercise to the reader.

**COROLLARY 4.1:** *Let V be a matrix vector space over the field F. {0} ∈ V; the space perpendicular to V under natural product $\times_n$ is V, that is $\{0\}^\perp = V$.*

**COROLLARY 4.2:** *Let V be a matrix vector space over the field F. The space perpendicular to V under natural product is {0} that is $\{V\}^\perp = \{0\}$.*

*Example 4.42:* Let

$$V = \left\{ \begin{bmatrix} a_1 & a_2 & a_3 \\ a_4 & a_5 & a_6 \end{bmatrix} \middle| a_i \in Q; 1 \leq i \leq 6 \right\}$$

be a vector space over Q.

Now consider $x = \begin{bmatrix} 0 & a_1 & a_2 \\ 0 & a_3 & a_4 \end{bmatrix}$ be the element of V. The vectors perpendicular or orthogonal to x are given by

$$\left\{ \begin{bmatrix} a & 0 & 0 \\ 0 & 0 & 0 \end{bmatrix}, \begin{bmatrix} 0 & 0 & 0 \\ a & 0 & 0 \end{bmatrix}, \begin{bmatrix} a & 0 & 0 \\ b & 0 & 0 \end{bmatrix}, \begin{bmatrix} 0 & 0 & 0 \\ 0 & 0 & 0 \end{bmatrix} \right\} = B.$$

Now take $y = \begin{bmatrix} a & 0 & 0 \\ b & 0 & 0 \end{bmatrix} \in B$, the vectors perpendicular to y are $y^\perp =$

$$\left\{ \begin{bmatrix} 0 & 0 & 0 \\ 0 & 0 & 0 \end{bmatrix}, \begin{bmatrix} 0 & a & 0 \\ 0 & 0 & 0 \end{bmatrix}, \begin{bmatrix} 0 & 0 & a_1 \\ 0 & 0 & 0 \end{bmatrix}, \begin{bmatrix} 0 & 0 & 0 \\ 0 & a_2 & 0 \end{bmatrix}, \begin{bmatrix} 0 & 0 & 0 \\ 0 & 0 & a_3 \end{bmatrix}, \right.$$

$$\begin{bmatrix} 0 & a_1 & a_2 \\ 0 & 0 & 0 \end{bmatrix}, \begin{bmatrix} 0 & 0 & 0 \\ 0 & a_1 & a_2 \end{bmatrix}, \begin{bmatrix} 0 & a_1 & 0 \\ 0 & a_2 & 0 \end{bmatrix}, \begin{bmatrix} 0 & 0 & x \\ 0 & 0 & y \end{bmatrix}, \begin{bmatrix} 0 & 0 & a \\ 0 & b & c \end{bmatrix},$$

$$\left. \begin{bmatrix} 0 & a & 0 \\ 0 & 0 & b \end{bmatrix}, \begin{bmatrix} 0 & 0 & c \\ 0 & d & 0 \end{bmatrix}, \begin{bmatrix} 0 & a & b \\ 0 & c & 0 \end{bmatrix}, \begin{bmatrix} 0 & a & b \\ 0 & 0 & b \end{bmatrix}, \begin{bmatrix} 0 & a & 0 \\ 0 & b & c \end{bmatrix}, \begin{bmatrix} 0 & a & b \\ 0 & c & d \end{bmatrix} \right\}$$

We see $x \in y^\perp$ under natural product. Now $\langle x^\perp \rangle = \{y \in V \mid x \times_n y = (0)\}$ and $\langle y^\perp \rangle = \{x \in V \mid x \times_n y = \{0\}\}$ are not only subspaces of V but are such that $\langle x^\perp \rangle \cup \langle y^\perp \rangle = V$ and $\langle x^\perp \rangle \cap \langle y^\perp \rangle = \{(0)\}$.

Further $x = \begin{bmatrix} 0 & a & b \\ 0 & c & d \end{bmatrix}$ is in $\langle y^\perp \rangle$.

Suppose $z^\perp = \begin{bmatrix} a & 0 & 0 \\ 0 & 0 & 0 \end{bmatrix} \in B$ is taken

$$z^\perp = \{m \in V \mid m \times_n z = (0)\}.$$

$$= \left\{ \begin{bmatrix} 0 & 0 & 0 \\ 0 & 0 & 0 \end{bmatrix}, \begin{bmatrix} 0 & a & 0 \\ 0 & 0 & 0 \end{bmatrix}, \begin{bmatrix} 0 & 0 & b \\ 0 & 0 & 0 \end{bmatrix}, \begin{bmatrix} 0 & 0 & 0 \\ x & 0 & 0 \end{bmatrix}, \begin{bmatrix} 0 & 0 & 0 \\ 0 & y & 0 \end{bmatrix}, \right.$$

$$\begin{bmatrix} 0 & 0 & 0 \\ 0 & 0 & t \end{bmatrix}, \begin{bmatrix} 0 & a & b \\ 0 & 0 & 0 \end{bmatrix}, \begin{bmatrix} 0 & b & 0 \\ a & 0 & 0 \end{bmatrix}, \begin{bmatrix} 0 & a & 0 \\ 0 & b & 0 \end{bmatrix}, \begin{bmatrix} 0 & a & 0 \\ 0 & 0 & b \end{bmatrix},$$

$$\begin{bmatrix} 0 & 0 & a \\ b & 0 & 0 \end{bmatrix}, \begin{bmatrix} 0 & 0 & a \\ 0 & b & 0 \end{bmatrix}, \begin{bmatrix} 0 & 0 & a \\ 0 & 0 & b \end{bmatrix}, \begin{bmatrix} 0 & 0 & 0 \\ a & b & 0 \end{bmatrix}, \begin{bmatrix} 0 & 0 & 0 \\ a & 0 & b \end{bmatrix},$$

$$\begin{bmatrix} 0 & 0 & 0 \\ 0 & a & b \end{bmatrix}, \begin{bmatrix} 0 & a & b \\ c & 0 & 0 \end{bmatrix}, \begin{bmatrix} 0 & a & b \\ 0 & c & 0 \end{bmatrix}, \begin{bmatrix} 0 & a & b \\ 0 & 0 & c \end{bmatrix}, \begin{bmatrix} 0 & a & 0 \\ b & c & 0 \end{bmatrix},$$

$$\begin{bmatrix} 0 & a & 0 \\ b & 0 & c \end{bmatrix}, \begin{bmatrix} 0 & a & 0 \\ 0 & b & c \end{bmatrix}, \begin{bmatrix} 0 & 0 & a \\ b & c & 0 \end{bmatrix}, \begin{bmatrix} 0 & 0 & a \\ b & 0 & c \end{bmatrix}, \begin{bmatrix} 0 & 0 & a \\ 0 & b & c \end{bmatrix},$$

$$\left.\begin{bmatrix} 0 & 0 & 0 \\ a & b & c \end{bmatrix}, \begin{bmatrix} 0 & a & b \\ 0 & c & d \end{bmatrix}, \begin{bmatrix} 0 & a & b \\ c & 0 & d \end{bmatrix}, \begin{bmatrix} 0 & a & b \\ c & d & 0 \end{bmatrix}, \begin{bmatrix} 0 & 0 & a \\ b & c & d \end{bmatrix},\right.$$

$$\left.\begin{bmatrix} 0 & a & 0 \\ b & c & d \end{bmatrix}, \begin{bmatrix} 0 & a & b \\ d & c & x \end{bmatrix}\right\} = T.$$

We see though $z \in B$ still $\langle z^\perp \rangle \neq \langle y^\perp \rangle$.

Further every element in T is not perpendicular to x under the natural product $\times_n$.

For consider $m = \begin{bmatrix} 0 & a & 0 \\ b & c & d \end{bmatrix}$ in T and

$$x \times_n m = \begin{bmatrix} 0 & a & b \\ 0 & c & d \end{bmatrix} \times_n \begin{bmatrix} 0 & a & 0 \\ b & c & d \end{bmatrix}$$

$$= \begin{bmatrix} 0 & x_1 & 0 \\ 0 & y_1 & z_1 \end{bmatrix} \neq \begin{bmatrix} 0 & 0 & 0 \\ 0 & 0 & 0 \end{bmatrix}.$$

Thus we can say for any $x \in V$ we have one and only one y in V such that x is the complement of y with respect to natural product $\times_n$.

We say complement, if $x^\perp$ generates the space B and $y^\perp$ generates another space say C.

Thus for $x = \begin{bmatrix} 0 & a & b \\ 0 & c & d \end{bmatrix}$ we say, $y = \begin{bmatrix} m & 0 & 0 \\ n & 0 & 0 \end{bmatrix}$

to be the main complement; a, b, c, d, m, n $\in Q \setminus \{0\}$. We say y is also the main complement of x with respect to natural product $\times_n$.

**THEOREM 4.3:** *Let V be a matrix vector space over the field Q (or R). Let $\times_n$ be the natural product defined on V. If for x in V, y is the main complement of V and vice versa, then $\langle x^\perp \rangle + \langle y^\perp \rangle = V$ and $\langle x^\perp \rangle \cap \langle y^\perp \rangle = (0)$.*

*(1) However for no other element t in $\langle y^\perp \rangle$ t will be the main complement of x.*

*(2) Also no element in $\langle x^\perp \rangle$ will be the main complement of y only x will be the main complement of y.*

The proof is left as an exercise. However we illustrate this situation by some example.

*Example 4.43:* Let

$$M = \left\{ \begin{bmatrix} a & b \\ c & d \end{bmatrix} \middle| a, b, c, d \in Q \right\}$$

be the vector space of $4 \times 4$ matrices over the field $F = Q$.

Take $p = \begin{bmatrix} x & 0 \\ y & 0 \end{bmatrix} \in M$, now the complements of p under natural product in M are $\left\{ \begin{bmatrix} 0 & 0 \\ 0 & 0 \end{bmatrix}, \begin{bmatrix} 0 & a \\ 0 & b \end{bmatrix}, \begin{bmatrix} 0 & a \\ 0 & 0 \end{bmatrix}, \right.$

$\left. \begin{bmatrix} 0 & 0 \\ 0 & b \end{bmatrix} \middle| a, b \in Q \right\} = T$.

The main complement of p under natural product $\times_n$ is

$$q = \begin{bmatrix} 0 & a \\ 0 & b \end{bmatrix} \in T.$$

Now the complements of q under natural product $\times_n$ in M are $\left\{\begin{bmatrix} 0 & 0 \\ 0 & 0 \end{bmatrix}, \begin{bmatrix} a & 0 \\ b & 0 \end{bmatrix}, \begin{bmatrix} a & 0 \\ 0 & 0 \end{bmatrix}, \begin{bmatrix} 0 & 0 \\ a & 0 \end{bmatrix} \middle| a, b \in Q\right\}$.

We see $V + T = M$ and $V \cap T = \begin{bmatrix} 0 & 0 \\ 0 & 0 \end{bmatrix}$.

Now $x = \begin{bmatrix} 0 & a \\ 0 & 0 \end{bmatrix}$ is in T. We find the elements orthogonal to x in M under the natural product $\times_n$.

$$\langle x^\perp \rangle = \left\langle \begin{bmatrix} 0 & a \\ 0 & 0 \end{bmatrix} \right\rangle^\perp =$$

$$\left\{ \begin{bmatrix} 0 & 0 \\ 0 & 0 \end{bmatrix}, \begin{bmatrix} a & 0 \\ 0 & 0 \end{bmatrix}, \begin{bmatrix} 0 & 0 \\ b & 0 \end{bmatrix}, \begin{bmatrix} 0 & 0 \\ 0 & c \end{bmatrix}, \right.$$

$$\left. \begin{bmatrix} a & 0 \\ b & 0 \end{bmatrix}, \begin{bmatrix} a & 0 \\ 0 & c \end{bmatrix}, \begin{bmatrix} 0 & 0 \\ b & c \end{bmatrix}, \begin{bmatrix} a & 0 \\ b & c \end{bmatrix} \right\}.$$

The main complement of x in M is $\begin{bmatrix} a & 0 \\ b & c \end{bmatrix}$ others are just complements.

Consider

$$\begin{bmatrix} a & 0 \\ b & 0 \end{bmatrix} \in \left\langle \begin{bmatrix} a & 0 \\ b & 0 \end{bmatrix} \right\rangle^\perp = \left\{ \begin{bmatrix} 0 & 0 \\ 0 & 0 \end{bmatrix}, \begin{bmatrix} 0 & 0 \\ 0 & a \end{bmatrix}, \begin{bmatrix} 0 & b \\ 0 & 0 \end{bmatrix}, \begin{bmatrix} 0 & a \\ 0 & b \end{bmatrix} \right\}.$$

The main complement of $\begin{bmatrix} a & 0 \\ b & 0 \end{bmatrix}$ is $\begin{bmatrix} 0 & a \\ 0 & b \end{bmatrix}$; other elements being just complements.

*Example 4.44*: Let

$$V = \left\{ \begin{bmatrix} a_1 \\ a_2 \\ \vdots \\ a_8 \end{bmatrix} \middle| a_i \in Q; 1 \le i \le 8 \right\}$$

be a vector space of column matrices over the field Q.

Consider the element $a = \begin{bmatrix} a_1 \\ a_2 \\ 0 \\ 0 \\ 0 \\ 0 \\ 0 \\ a_3 \end{bmatrix}$ in V.

To find complements or elements orthogonal to

$$\langle a \rangle^\perp = \left\{ \begin{bmatrix} 0 \\ 0 \\ 0 \\ a \\ \vdots \\ 0 \end{bmatrix}, \begin{bmatrix} 0 \\ 0 \\ 0 \\ a \\ 0 \\ \vdots \\ 0 \end{bmatrix}, \begin{bmatrix} 0 \\ 0 \\ 0 \\ 0 \\ a \\ 0 \\ \vdots \\ 0 \end{bmatrix}, \begin{bmatrix} 0 \\ 0 \\ 0 \\ 0 \\ 0 \\ a \\ 0 \\ 0 \end{bmatrix}, \begin{bmatrix} 0 \\ 0 \\ 0 \\ 0 \\ 0 \\ 0 \\ a \\ 0 \end{bmatrix}, \begin{bmatrix} 0 \\ 0 \\ 0 \\ 0 \\ 0 \\ a \\ 0 \\ 0 \end{bmatrix}, \begin{bmatrix} 0 \\ 0 \\ a_1 \\ a_2 \\ 0 \\ 0 \\ 0 \\ 0 \end{bmatrix}, \begin{bmatrix} 0 \\ 0 \\ a_1 \\ 0 \\ a_2 \\ 0 \\ 0 \\ 0 \end{bmatrix}, \right.$$

$$\left[\begin{array}{c}0\\0\\a_1\\0\\0\\a_2\\0\\0\end{array}\right],\left[\begin{array}{c}0\\0\\a_1\\0\\0\\0\\a_2\\0\end{array}\right],\left[\begin{array}{c}0\\0\\0\\0\\a_1\\a_2\\0\\0\end{array}\right],\left[\begin{array}{c}0\\0\\0\\a_1\\a_2\\0\\0\\0\end{array}\right],\left[\begin{array}{c}0\\0\\0\\a_1\\0\\a_2\\0\\0\end{array}\right],\left[\begin{array}{c}0\\0\\0\\a_1\\0\\0\\a_2\\0\end{array}\right],\left[\begin{array}{c}0\\0\\0\\0\\a_1\\0\\a_2\\0\end{array}\right],\left[\begin{array}{c}0\\0\\0\\0\\0\\a_1\\a_2\\0\end{array}\right],$$

$$\left[\begin{array}{c}0\\0\\a_1\\a_2\\a_3\\0\\0\\0\end{array}\right],\left[\begin{array}{c}0\\0\\a_1\\a_2\\0\\a_3\\0\\0\end{array}\right],\left[\begin{array}{c}0\\0\\a_1\\a_2\\0\\0\\a_3\\0\end{array}\right],\left[\begin{array}{c}0\\0\\a_1\\0\\a_2\\a_3\\0\\0\end{array}\right],\left[\begin{array}{c}0\\0\\a_1\\0\\0\\a_2\\a_3\\0\end{array}\right],\left[\begin{array}{c}0\\0\\a_1\\0\\a_2\\0\\a_3\\0\end{array}\right],\left[\begin{array}{c}0\\0\\0\\a_1\\a_2\\a_3\\0\\0\end{array}\right],\left[\begin{array}{c}0\\0\\0\\a_1\\0\\a_2\\a_3\\0\end{array}\right],$$

$$\left.\left[\begin{array}{c}0\\0\\0\\a_1\\a_2\\0\\a_3\\0\end{array}\right],\left[\begin{array}{c}0\\0\\0\\0\\a_1\\a_2\\a_3\\0\end{array}\right],\left[\begin{array}{c}0\\0\\a_1\\a_2\\a_3\\a_4\\0\\0\end{array}\right],\left[\begin{array}{c}0\\0\\a_1\\a_2\\a_3\\0\\a_4\\0\end{array}\right],\left[\begin{array}{c}0\\0\\a_1\\0\\a_2\\a_3\\a_4\\0\end{array}\right],\left[\begin{array}{c}0\\0\\0\\a_1\\a_2\\a_3\\a_4\\0\end{array}\right],\left[\begin{array}{c}0\\0\\a_1\\a_2\\a_3\\a_4\\a_5\\0\end{array}\right]\right\}.$$

The main complement of $\begin{bmatrix} a_1 \\ a_2 \\ 0 \\ 0 \\ 0 \\ 0 \\ 0 \\ a_3 \end{bmatrix}$ is $\begin{bmatrix} 0 \\ 0 \\ a_1 \\ a_2 \\ a_3 \\ a_4 \\ a_5 \\ 0 \end{bmatrix}$.

Certainly this type of study will be a boon to algebraic coding theory as in case of algebraic coding theory we mainly use only matrices which are m × n (m ≠ n) as parity check matrix and generator matrix.

Now we define other related properties of these matrices with natural product on them. Suppose S is a subset of a matrix vector space defined on R or Q we can define

$$S^\perp = \{x \in V \mid x \times_n s = (0) \text{ for every } s \in S\}.$$

We will illustrate this situation by some simple example.

*Example 4.45:* Let

$$V = \left\{ \begin{bmatrix} a_1 & a_2 \\ a_3 & a_4 \\ a_5 & a_6 \\ a_7 & a_8 \end{bmatrix} \middle| a_i \in Q; 1 \leq i \leq 8 \right\}$$

be a vector space of 4 × 2 matrices defined over the field V. V is a linear algebra under the natural product.

Consider

$$S = \left\{ \begin{bmatrix} a & b \\ 0 & 0 \\ 0 & 0 \\ 0 & 0 \end{bmatrix}, \begin{bmatrix} 0 & 0 \\ 0 & 0 \\ 0 & 0 \\ c & d \end{bmatrix} \middle| a, b, c, d \in Q \right\} \subseteq V.$$

To find $S^\perp$. $S^\perp =$

$$\left\{ \begin{bmatrix} 0 & 0 \\ 0 & 0 \\ 0 & 0 \\ 0 & 0 \end{bmatrix}, \begin{bmatrix} 0 & 0 \\ a & b \\ d & e \\ 0 & 0 \end{bmatrix}, \begin{bmatrix} 0 & 0 \\ a & 0 \\ 0 & 0 \\ 0 & 0 \end{bmatrix}, \begin{bmatrix} 0 & 0 \\ 0 & a \\ 0 & 0 \\ 0 & 0 \end{bmatrix}, \right.$$

$$\begin{bmatrix} 0 & 0 \\ 0 & 0 \\ b & 0 \\ 0 & 0 \end{bmatrix}, \begin{bmatrix} 0 & 0 \\ 0 & 0 \\ 0 & d \\ 0 & 0 \end{bmatrix}, \begin{bmatrix} 0 & 0 \\ 0 & 0 \\ a & b \\ 0 & 0 \end{bmatrix}, \begin{bmatrix} 0 & 0 \\ a & b \\ 0 & 0 \\ 0 & 0 \end{bmatrix},$$

$$\begin{bmatrix} 0 & 0 \\ a & 0 \\ 0 & b \\ 0 & 0 \end{bmatrix}, \begin{bmatrix} 0 & 0 \\ a & 0 \\ b & 0 \\ 0 & 0 \end{bmatrix}, \begin{bmatrix} 0 & 0 \\ 0 & a \\ 0 & b \\ 0 & 0 \end{bmatrix}, \begin{bmatrix} 0 & 0 \\ 0 & b \\ a & 0 \\ 0 & 0 \end{bmatrix},$$

$$\left. \begin{bmatrix} 0 & 0 \\ a & b \\ c & 0 \\ 0 & 0 \end{bmatrix}, \begin{bmatrix} 0 & 0 \\ a & b \\ 0 & c \\ 0 & 0 \end{bmatrix}, \begin{bmatrix} 0 & 0 \\ 0 & a \\ b & c \\ 0 & 0 \end{bmatrix}, \begin{bmatrix} 0 & 0 \\ a & 0 \\ b & d \\ 0 & 0 \end{bmatrix} \right\} \subseteq V$$

is the orthogonal complement $S^\perp$.

We see $S^\perp$ is a subspace of V. Further the main complement of $x = \begin{bmatrix} a & b \\ 0 & 0 \\ 0 & 0 \\ 0 & 0 \end{bmatrix}$ and $\begin{bmatrix} 0 & 0 \\ 0 & 0 \\ 0 & 0 \\ b & c \end{bmatrix} = y$ are not in $S^\perp$ for the main complement of x is $\begin{bmatrix} 0 & 0 \\ a & b \\ c & d \\ e & f \end{bmatrix}$ and that of y is $\begin{bmatrix} a & b \\ c & d \\ e & f \\ 0 & 0 \end{bmatrix}$ but the main complement of x is not orthogonal with y; for $y \times_n x^\perp$.

$$= \begin{bmatrix} 0 & 0 \\ 0 & 0 \\ 0 & 0 \\ b & c \end{bmatrix} \begin{bmatrix} 0 & 0 \\ a & b \\ c & d \\ e & f \end{bmatrix} = \begin{bmatrix} 0 & 0 \\ 0 & 0 \\ 0 & 0 \\ be & cf \end{bmatrix} \neq \begin{bmatrix} 0 & 0 \\ 0 & 0 \\ 0 & 0 \\ 0 & 0 \end{bmatrix}.$$

Similarly the main complement of $\begin{bmatrix} 0 & 0 \\ 0 & 0 \\ 0 & 0 \\ a & b \end{bmatrix}$, that is $\begin{bmatrix} a & b \\ c & d \\ e & f \\ 0 & 0 \end{bmatrix}$ is not orthogonal with $\begin{bmatrix} a & b \\ 0 & 0 \\ 0 & 0 \\ 0 & 0 \end{bmatrix}$ under natural product as

$$\begin{bmatrix} a & b \\ 0 & 0 \\ 0 & 0 \\ 0 & 0 \end{bmatrix} \times_n \begin{bmatrix} a_1 & b_1 \\ c & d \\ e & f \\ 0 & 0 \end{bmatrix} = \begin{bmatrix} a_1 a & bb_1 \\ 0 & 0 \\ 0 & 0 \\ 0 & 0 \end{bmatrix} \neq \begin{bmatrix} 0 & 0 \\ 0 & 0 \\ 0 & 0 \\ 0 & 0 \end{bmatrix}.$$

We can define as in case of usual vector spaces define linear transformation for matrix vector spaces. However it is meaning less to define linear transformation in case of S-special strong vector spaces defined over the S-field. However in that case only linear operators can be defined. The definition properties etc in case of the former vector space is a matter of routine and we see no difference with usual spaces. However in case of latter S-special strong vector spaces over a S-field we can define only linear operators.

We just illustrate this situation by an example.

*Example 4.46:* Let

$$V = \left\{ \begin{bmatrix} a_1 & a_2 \\ a_3 & a_4 \\ a_5 & a_6 \\ a_7 & a_8 \end{bmatrix} \middle| a_i \in R; 1 \leq i \leq 8 \right\}$$

be a S-special super vector space over the S-field.

$$F_{4\times 2} = \left\{ \begin{bmatrix} a & b \\ c & d \\ e & f \\ g & h \end{bmatrix} \middle| a, b, c, d, e, f, g, h \in Q \right\}.$$

Now we define a map $\eta : V \to V$ as

$$\eta\left( \begin{bmatrix} a_1 & a_2 \\ a_3 & a_4 \\ a_5 & a_6 \\ a_7 & a_8 \end{bmatrix} \right) = \begin{pmatrix} a_1 & 0 \\ 0 & a_2 \\ a_3 & 0 \\ 0 & a_4 \end{pmatrix}.$$

It is easily verified η is a linear operator on V.

$$\text{We can find kernel } \eta = \left\{ \begin{bmatrix} a & b \\ c & d \\ e & f \\ g & h \end{bmatrix} \in V \;\middle|\; \eta\begin{pmatrix} a & b \\ c & d \\ e & f \\ g & h \end{pmatrix} = (0) \right\}$$

$$= \left\{ \begin{bmatrix} 0 & a \\ b & 0 \\ 0 & c \\ d & 0 \end{bmatrix} \;\middle|\; a, b, c, d \in Q \right\}$$

is a subspace of V. Now interested reader can work with linear operators on S-special strong vector spaces over the S-field.

Thus all matrix vector spaces are linear algebras under the natural product. Now as in case of usual vector spaces we can for the case of matrix vector spaces also define the notion of linear functional. But in case of S-strong special matrix vector spaces we can not define only Smarandache linear functional as matter of routine as it needs more modifications and changes.

Now we have discussed some of properties about vector spaces. We now define n - row matrix vector space over a field F.

**DEFINITION 4.6:** *Let*

$$V = \{(a_1, \ldots, a_m) \mid a_i = (x_1, \ldots, x_n);\; x_j \in Q;\; 1 \leq i \leq m,\; 1 \leq j \leq n\};$$

*V is a vector space over Q defined as the n-row matrix structured vector space over Q.*

We will illustrate this situation by some examples.

***Example 4.47:*** Let $V = \{(a_1, a_2, a_3, a_4) \mid a_j = (x_1, x_2, x_3, x_4, x_5)$ where $x_i \in Q$, $1 \leq i \leq 5$ and $1 \leq j \leq 4\}$ be a 5-row matrix structured vector space over Q.

We will just show how addition and scalar multiplication is performed on V.

Suppose $x = ((3, 0, 2, 4, 5), (0, 0, 0, 1, 2), (1, 1, 1, 3, 0), (2, 0, 1, 0, 5)) \in V$ and $a = 7$ then $7x = ((21, 0, 14, 28, 35), (0, 0, 0, 7, 14), (7, 7, 7, 21, 0), (14, 0, 7, 0, 35)) \in V$.

Let $y = ((4, 0, 1, 1, 1), (0, 1, 0, 1, 2), (1, 0, 1, 1, 1), (2, 0, 0, 0, 1)) \in V$ then $x + y = ((7, 0, 3, 5, 6), (0, 1, 0, 2, 4), (2, 1, 2, 4, 1), (4, 0, 1, 0, 6)) \in V$. We see V is a row matrix structured vector space over Q.

***Example 4.48:*** Let

$$P = \{(a_1, a_2, a_3) \mid a_i = (x_1, x_2, \ldots, x_{15})\ x_j \in Q;\ 1 \leq i \leq 3;\ 1 \leq j \leq 15\}$$

a row matrix structured vector space over Q. We can define row matrix structured subvector space as in case of usual vector space. On P we can always define the natural product hence P is always a row matrix structured linear algebra under the natural product $\times_n$.

***Example 4.49:*** Let

$V = \{(a_1, \ldots, a_{10}) \mid a_i = (x_1, x_2, x_3);\ x_j \in Q;\ 1 \leq j \leq 3;\ 1 \leq i \leq 10\}$

be a row matrix structured vector space over Q. Consider $H = \{(a_1, a_2, a_3, 0, 0, 0, 0, 0, 0, 0) \mid a_i = (x_1, x_2, x_3);\ x_j \in Q;\ 1 \leq i \leq 3\} \subseteq V$; H is a row matrix structured vector subspace of V over Q. Also

$$P = \{(a_1, a_2, \ldots, a_{10}) \mid a_i = (x_1, 0, x_2)$$
$$\text{with } x_1, x_2 \in Q;\ 1 \leq i \leq 10\} \subseteq V$$

is a row matrix structured vector subspace of V over Q. The reader can see the difference between the subspaces P and H.

Let $V = \{(x_1, \ldots, x_n) \mid x_i \in R^+ \cup \{0\}, 1 \leq i \leq n\}$ be a semigroup under addition. V is a semivector space over the semifield $R^+ \cup \{0\}$ or $Q^+ \cup \{0\}$ or $Z^+ \cup \{0\}$.

Likewise $M = \left\{ \begin{bmatrix} x_1 \\ x_2 \\ \vdots \\ x_m \end{bmatrix} \mid x_i \in Q^+ \cup \{0\}, 1 \leq i \leq m \right\}$ be a

semigroup under addition. M is a semivector space over the semifield $Z^+ \cup \{0\}$ or $Q^+ \cup \{0\}$. M is not a semivector space over $R^+ \cup \{0\}$.

$$P = \left\{ \begin{bmatrix} a_{11} & a_{12} & \ldots & a_{1n} \\ a_{21} & a_{22} & \ldots & a_{2n} \\ \vdots & \vdots & & \vdots \\ a_{m1} & a_{m2} & \ldots & a_{mn} \end{bmatrix} \mid a_{ij} \in Z^+ \cup \{0\}, 1 \leq i \leq m, 1 \leq j \leq n \right\}$$

is a semivector space over the semifield $Z^+ \cup \{0\}$.

$$T = \left\{ \begin{bmatrix} a_{11} & a_{12} & \ldots & a_{1n} \\ a_{21} & a_{22} & \ldots & a_{2n} \\ \vdots & \vdots & & \vdots \\ a_{n1} & a_{n2} & \ldots & a_{nn} \end{bmatrix} \mid a_{ij} \in Q^+ \cup \{0\}, 1 \leq i, j \leq n \right\}$$

is a semivector space over the semifield $Q^+ \cup \{0\}$. M, P T and V are also semilinear algebras over the respective semifields under the natural product $\times_n$.

We will illustrate these situations by some examples.

*Example 4.50:* Let
$$V = \{(x_1, x_2, x_3, x_4, x_5, x_6) \mid x_i \in 3Z^+ \cup \{0\}; 1 \leq i \leq 6\}$$

be the semivector space over the semifield $S = Z^+ \cup \{0\}$.

V is also a semilinear algebra over the semifield S.

*Example 4.51:* Let
$$V_1 = \{(x_1, x_2, x_3, x_4, x_5, x_6) \text{ where } x_i \in R^+ \cup \{0\}; 1 \leq i \leq 6\}$$

be the semivector space over the semifield $S = Z^+ \cup \{0\}$.

It is interesting to compare V and $V_1$ for $V_1$ is finite dimensional where as $V_2$ is of infinite dimension.

*Example 4.52*: Let

$$V = \left\{ \begin{bmatrix} x_1 \\ x_2 \\ x_3 \\ x_4 \\ x_5 \\ x_6 \\ x_7 \\ x_8 \end{bmatrix} \middle| x_i \in Z^+ \cup \{0\}, 1 \leq i \leq 8 \right\}$$

be a semivector space over the semifield $S = Z^+ \cup \{0\}$. V is a semilinear algebra over S.

Dimension of V over S is eight.

*Example 4.53*: Let

$$M = \left\{ \begin{bmatrix} x_1 \\ x_2 \\ x_3 \\ x_4 \\ x_5 \\ x_6 \\ x_7 \\ x_8 \end{bmatrix} \middle| \, x_i \in Q^+ \cup \{0\}, 1 \leq i \leq 8 \right\}$$

be a semivector space over the semifield $S = Z^+ \cup \{0\}$. Clearly dimension of M over S is infinite.

*Example 4.54*: Let

$$M = \left\{ \begin{bmatrix} a_1 & a_2 & a_3 & a_4 \\ a_5 & a_6 & a_7 & a_8 \\ a_9 & a_{10} & a_{11} & a_{12} \\ a_{13} & a_{14} & a_{15} & a_{16} \\ a_{17} & a_{18} & a_{19} & a_{20} \\ a_{21} & a_{22} & a_{23} & a_{24} \\ a_{25} & a_{26} & a_{27} & a_{28} \\ a_{29} & a_{30} & a_{31} & a_{32} \\ a_{33} & a_{34} & a_{35} & a_{36} \\ a_{37} & a_{38} & a_{39} & a_{40} \end{bmatrix} \middle| \, a_i \in Z^+ \cup \{0\}, 1 \leq i \leq 40 \right\}$$

be a semivector space over the semifield $S = Z^+ \cup \{0\}$. Clearly V is not defined over the semifield $T = Q^+ \cup \{0\}$ or $R^+ \cup \{0\}$. Dimension of V over S in 40.

*Example 4.55*: Let

$$M = \left\{ \begin{bmatrix} a_1 & a_2 & a_3 \\ a_4 & a_5 & a_6 \\ a_7 & a_8 & a_9 \end{bmatrix} \,\middle|\, a_i \in Q^+ \cup \{0\}, 1 \leq i \leq 9 \right\}$$

be a semivector space over the semifield $S = Z^+ \cup \{0\} = S$.

Take

$$P = \left\{ \begin{bmatrix} a_1 & a_2 & a_3 \\ a_4 & 0 & a_5 \\ a_6 & 0 & 0 \end{bmatrix} \,\middle|\, a_i \in Q^+ \cup \{0\}, 1 \leq i \leq 9 \right\},$$

P is a semivector subspace of V over $S = Z^+ \cup \{0\}$.

*Example 4.56*: Let

$$M = \left\{ \begin{bmatrix} a_1 & a_2 & a_3 \\ a_4 & a_5 & a_6 \\ a_7 & a_8 & a_9 \\ a_{10} & a_{11} & a_{12} \\ a_{13} & a_{14} & a_{15} \\ a_{16} & a_{17} & a_{18} \end{bmatrix} \,\middle|\, a_i \in Q^+ \cup \{0\}, 1 \leq i \leq 18 \right\}$$

be a semivector space over the semifield $S = Q^+ \cup \{0\}$.

$$M_1 = \left\{ \begin{bmatrix} a_1 & a_2 & a_3 \\ 0 & 0 & 0 \\ \vdots & \vdots & \vdots \\ 0 & 0 & 0 \end{bmatrix} \,\middle|\, a_i \in Q^+ \cup \{0\}, 1 \leq i \leq 3 \right\} \subseteq M$$

is a semivector subspace of M over $S = Q^+ \cup \{0\}$, the semivector space.

$$M_2 = \left\{ \begin{bmatrix} 0 & 0 & 0 \\ a_1 & a_2 & a_3 \\ a_4 & a_5 & a_6 \\ 0 & 0 & 0 \\ 0 & 0 & 0 \\ 0 & 0 & 0 \end{bmatrix} \middle| a_i \in Q^+ \cup \{0\}, 1 \le i \le 6 \right\} \subseteq M$$

is a semivector subspace of M over $S = Q^+ \cup \{0\}$.

$$M_3 = \left\{ \begin{bmatrix} 0 & 0 & 0 \\ 0 & 0 & 0 \\ 0 & 0 & 0 \\ a_1 & a_2 & a_3 \\ a_4 & a_5 & a_6 \\ 0 & 0 & 0 \end{bmatrix} \middle| a_i \in Q^+ \cup \{0\}, 1 \le i \le 6 \right\} \subseteq M$$

is a semivector subspace of M over $S = Q^+ \cup \{0\}$ and

$$M_4 = \left\{ \begin{bmatrix} 0 & 0 & 0 \\ 0 & 0 & 0 \\ 0 & 0 & 0 \\ 0 & 0 & 0 \\ a_1 & a_2 & a_3 \\ a_4 & a_5 & a_6 \end{bmatrix} \middle| a_i \in Q^+ \cup \{0\}, 1 \le i \le 6 \right\} \subseteq M$$

is the semivector subspace of M over S.

We see $M = M_1 + M_2 + M_3 + M_4$ and $M_i \cap M_j = \begin{bmatrix} 0 & 0 & 0 \\ 0 & 0 & 0 \\ 0 & 0 & 0 \\ 0 & 0 & 0 \\ 0 & 0 & 0 \\ 0 & 0 & 0 \end{bmatrix}$

if $i \ne j$, $1 \le i, j \le 4$.

Thus M is the direct sum of semivector subspaces of M over S.

*Example 4.57:* Let

$$V = \left\{ \begin{bmatrix} a_1 & a_2 & a_3 & a_4 & a_5 \\ a_6 & a_7 & a_8 & a_9 & a_{10} \\ a_{11} & a_{12} & a_{13} & a_{14} & a_{15} \end{bmatrix} \middle| a_i \in Z^+ \cup \{0\}, 1 \leq i \leq 15 \right\}$$

be a semivector space over the semifield $S = Z^+ \cup \{0\}$.

Consider

$$P_1 = \left\{ \begin{bmatrix} a_1 & 0 & 0 & 0 & 0 \\ a_2 & 0 & 0 & 0 & 0 \\ a_3 & 0 & 0 & 0 & 0 \end{bmatrix} \middle| a_i \in Z^+ \cup \{0\}, 1 \leq i \leq 3 \right\} \subseteq V$$

be a semivector subspace of V over S.

Let

$$P_2 = \left\{ \begin{bmatrix} a_4 & a_1 & a_2 & 0 & 0 \\ 0 & 0 & a_3 & 0 & 0 \\ 0 & 0 & 0 & 0 & 0 \end{bmatrix} \middle| a_i \in Z^+ \cup \{0\}, 1 \leq i \leq 4 \right\} \subseteq V$$

be a semivector subspace of V over S.

Consider

$$P_3 = \left\{ \begin{bmatrix} a_1 & 0 & a_2 & 0 & 0 \\ 0 & a_4 & a_3 & 0 & 0 \\ 0 & 0 & 0 & 0 & 0 \end{bmatrix} \middle| a_i \in Z^+ \cup \{0\}, 1 \leq i \leq 4 \right\} \subseteq V,$$

a semivector subspace of V over S.

Further

$$P_4 = \left\{ \begin{bmatrix} a_1 & 0 & a_2 & 0 & 0 \\ 0 & 0 & 0 & 0 & 0 \\ a_3 & 0 & a_4 & 0 & 0 \end{bmatrix} \middle| a_i \in Z^+ \cup \{0\}, 1 \leq i \leq 4 \right\} \subseteq V,$$

is a semivector subspace of V over S.

$$P_5 = \left\{ \begin{bmatrix} a_1 & 0 & 0 & a_3 & a_4 \\ 0 & 0 & 0 & a_2 & 0 \\ 0 & a_5 & 0 & 0 & 0 \end{bmatrix} \middle| a_i \in Z^+ \cup \{0\}, 1 \leq i \leq 5 \right\} \subseteq V$$

is a semivector subspace of V over S.

$$P_6 = \left\{ \begin{bmatrix} a_1 & 0 & 0 & 0 & 0 \\ 0 & 0 & 0 & 0 & a_5 \\ 0 & 0 & a_2 & a_3 & a_4 \end{bmatrix} \middle| a_i \in Z^+ \cup \{0\}, 1 \leq i \leq 5 \right\} \subseteq V,$$

is a semivector subspace of V over S.

We see $P_i \cap P_j \neq \begin{bmatrix} 0 & 0 & 0 & 0 & 0 \\ 0 & 0 & 0 & 0 & 0 \\ 0 & 0 & 0 & 0 & 0 \end{bmatrix}$ if $i \neq j$, $1 \leq i, j \leq 6$.

We see $V \subseteq P_1 + P_2 + P_3 + P_4 + P_5 + P_6$; thus V is a pseudo direct sum of semivector subspaces.

Now we have seen examples of direct sum and pseudo direct sum of semivector subspaces over the semifield $S = Z^+ \cup \{0\}$. Now we can as in case of other semivector spaces define linear transformation and linear operator.

We give examples of semivector space with polynomial matrix coefficient elements.

*Example 4.58:* Let

$$V = \left\{ \sum_{i=0}^{\infty} a_i x^i \,\middle|\, a_i = (x_1, x_2, x_3, x_4) \mid x_j \in Z^+ \cup \{0\}; 1 \leq j \leq 4 \right\}$$

be a semivector space of infinite dimension over $S = Z^+ \cup \{0\}$.

Clearly V is also a semilinear algebra over S under the natural product $\times_n$.

*Example 4.59:* Let

$$M = \left\{ \sum_{i=0}^{\infty} a_i x^i \,\middle|\, a_i = \begin{bmatrix} x_1 \\ x_2 \\ x_3 \\ \vdots \\ x_{10} \end{bmatrix} \text{ where } x_j \in Z^+ \cup \{0\}; 1 \leq j \leq 10 \right\}$$

be a semivector space of infinite dimension over $S = Z^+ \cup \{0\}$.

*Example 4.60:* Let

$$P = \left\{ \sum_{i=0}^{\infty} d_i x^i \,\middle|\, d_i = \begin{bmatrix} a_1 & a_2 & a_3 & a_4 \\ a_5 & a_6 & a_7 & a_8 \\ a_9 & a_{10} & a_{11} & a_{12} \\ a_{13} & a_{14} & a_{15} & a_{16} \end{bmatrix} \right.$$
$$\left. \text{where } a_i \in Q^+ \cup \{0\}; 1 \leq i \leq 16 \right\}$$

be a semivector space of infinite dimension over $S = Q^+ \cup \{0\}$.

*Example 4.61:* Let

$$M = \left\{ \sum_{i=0}^{\infty} a_i x^i \,\middle|\, a_i = \begin{pmatrix} d_1 & d_2 & \ldots & d_{10} \\ d_{11} & d_{12} & \ldots & d_{20} \\ d_{21} & d_{22} & \ldots & d_{30} \end{pmatrix} \right.$$
$$\left. \text{where } d_i \in Z^+ \cup \{0\}; 1 \leq i \leq 30 \right\}$$

be a semivector space of infinite dimension over $S = Z^+ \cup \{0\}$.

*Example 4.62:* Let

$$S = \left\{ \sum_{i=0}^{\infty} a_i x^i \;\middle|\; a_i = \begin{bmatrix} d_1 & d_2 & \ldots & d_{10} \\ d_{11} & d_{12} & \ldots & d_{20} \\ \vdots & \vdots & & \vdots \\ d_{81} & d_{82} & \ldots & d_{90} \end{bmatrix} \right.$$

with $d_i \in Z^+ \cup \{0\}$; $1 \le i \le 9\}$

be a semivector space of infinite dimension over $S = Z^+ \cup \{0\}$.

*Example 4.63:* Let

$$T = \left\{ \sum_{i=0}^{\infty} a_i x^i \;\middle|\; a_j = \begin{bmatrix} d_1 & d_2 & \ldots & d_n \\ d_{n+1} & d_{n+2} & \ldots & d_{2n} \\ d_{2n+1} & d_{2n+2} & \ldots & d_{3n} \\ \vdots & \vdots & & \vdots \\ d_{7n+1} & d_{7n+2} & \ldots & d_{8n} \end{bmatrix} \right.$$

with $d_i \in Z^+ \cup \{0\}$; $1 \le i \le 8n\}$

be a semivector space over the semifield of infinite dimension over $Z^+ \cup \{0\}$.

Now having seen such examples we now proceed onto define other properties of these semivector spaces.

All these are also semilinear algebras over the semifields.

Now if the natural product is defined we can speak of complements of the semivector subspaces (semilinear algebras).

We will illustrate this situation by some simple examples.

*Example 4.64:* Let

$$V = \left\{ \begin{bmatrix} a_1 \\ a_2 \\ a_3 \\ a_4 \\ a_5 \end{bmatrix} \text{ with } a_i \in Z^+ \cup \{0\}, 1 \leq i \leq 5 \right\}$$

be a semivector space over the semifield $S = Z^+ \cup \{0\}$.

Consider

$$M_1 = \left\{ \begin{bmatrix} 0 \\ 0 \\ a_1 \\ a_2 \\ 0 \end{bmatrix} \text{ with } a_i \in Z^+ \cup \{0\}, 1 \leq i \leq 2 \right\} \subseteq V$$

be a semivector subspace of V over S.

$$M_2 = \left\{ \begin{bmatrix} a_1 \\ a_2 \\ 0 \\ 0 \\ a_3 \end{bmatrix} \middle| a_i \in Z^+ \cup \{0\}, 1 \leq i \leq 3 \right\} \subseteq V,$$

be a semivector subspace of V over S.  Clearly every $x \in M_1$ and $y \in M_2$ are such that $x \times_n y = (0)$.

We see $V = M_1 \oplus M_2$; $M_1 \cap M_2 = (0)$.

*Example 4.65:* Let $V = \{(a_1, a_2, \ldots, a_{10}) \mid a_i \in Z^+ \cup \{0\}, 1 \leq i \leq 10\}$ be a semivector space over the semifield $S = Z^+ \cup \{0\}$.

Consider $P_1 = \{(0, 0, 0, a_1, 0, 0, 0, 0, 0, a_7) \mid a_i \in Z^+ \cup \{0\}, 1 \le i \le 7\} \subseteq V$ be a semivector subspace of V over S.

$P_2 = \{(a_1, a_2, 0, 0, \ldots, 0) \mid a_1, a_2 \in Z^+ \cup \{0\}, 1 \le i \le 7\} \subseteq V$ be a semivector subspace of V over S.

$P_3 = \{(0, 0, a_1, 0, a_2, a_3, 0, 0, 0, 0) \mid a_i \in Z^+ \cup \{0\}, 1 \le i \le 3\} \subseteq V$ be a semivector subspace of V over S.

$P_4 = \{(0, 0, 0, 0, 0, 0, a_1, a_2, a_3, 0) \mid a_i \in Z^+ \cup \{0\}, 1 \le i \le 3\} \subseteq V$ is again a semivector subspace of V over S.

We see every vector in $P_i$ is orthogonal with every other vector in P if $i \ne j$; $1 \le i, j \le 4$.

Further $V = P_1 + P_2 + P_3 + P_4$ and $P_i \cap P_j = (0)$ if $i \ne j$.

*Example 4.66:* Let

$$M = \left\{ \begin{bmatrix} a_1 & a_2 & a_3 \\ a_4 & a_5 & a_6 \\ a_7 & a_8 & a_9 \\ a_{10} & a_{11} & a_{12} \\ a_{13} & a_{14} & a_{15} \end{bmatrix} \middle| a_i \in Z^+ \cup \{0\}, 1 \le i \le 15 \right\}$$

be a semivector space over the semifield $S = Z^+ \cup \{0\}$.

Take

$$P_1 = \left\{ \begin{bmatrix} 0 & 0 & 0 \\ a_1 & a_2 & a_3 \\ 0 & 0 & 0 \\ a_4 & a_5 & a_6 \\ 0 & 0 & 0 \end{bmatrix} \middle| a_i \in Z^+ \cup \{0\}, 1 \le i \le 6 \right\} \subseteq M,$$

is a semivector subspace of M over $S = Z^+ \cup \{0\}$.

$$P_2 = \left\{ \begin{bmatrix} a_1 & a_2 & a_3 \\ 0 & 0 & 0 \\ a_4 & a_5 & a_6 \\ 0 & 0 & 0 \\ a_7 & a_8 & a_9 \end{bmatrix} \middle| a_i \in Z^+ \cup \{0\}, 1 \le i \le 9 \right\} \subseteq M,$$

is a semivector subspace of M over $S = Z^+ \cup \{0\}$.

We see for every $x \in P_1$ we have $x \times_n y = (0)$ for every $y \in P_2$.

Thus $M = M_1 + M_2$ and $P_1 \cap P_2 = (0)$. We say the space $P_1$ is orthogonal with the space $P_2$ of M.

However if

$$P_3 = \left\{ \begin{bmatrix} 0 & 0 & 0 \\ 0 & 0 & 0 \\ 0 & 0 & 0 \\ 0 & 0 & 0 \\ a_1 & a_2 & a_3 \end{bmatrix} \middle| a_i \in Z^+ \cup \{0\}, 1 \le i \le 3 \right\} \subseteq M,$$

we see $P_1$ and $P_3$ are such that for every $x \in P_1$ we have $x \times_n y = (0)$ for every $y \in P_3$; however we do not call $P_3$ the complementary space of $P_1$ as $M \ne P_1 + P_3$.

*Example 4.67:* Let

$$P = \left\{ \begin{bmatrix} a_1 & a_2 & a_3 & a_4 \\ a_5 & a_6 & a_7 & a_8 \\ a_9 & a_{10} & a_{11} & a_{12} \\ a_{13} & a_{14} & a_{15} & a_{16} \end{bmatrix} \middle| a_i \in Q^+ \cup \{0\}, 1 \le i \le 16 \right\}$$

be a semivector space over the semifield $Z^+ \cup \{0\} = S$.

Consider

$$M_1 = \left\{ \begin{bmatrix} a_1 & a_2 & 0 & 0 \\ a_3 & a_4 & 0 & 0 \\ 0 & 0 & 0 & 0 \\ 0 & 0 & 0 & 0 \end{bmatrix} \middle| a_i \in Q^+ \cup \{0\}, 1 \le i \le 4 \right\} \subseteq P;$$

$M_1$ is a semivector subspace of P over S.

$$M_2 = \left\{ \begin{bmatrix} 0 & 0 & a_1 & a_2 \\ 0 & 0 & a_3 & a_4 \\ 0 & 0 & 0 & 0 \\ 0 & 0 & 0 & 0 \end{bmatrix} \middle| a_i \in Q^+ \cup \{0\}, 1 \le i \le 4 \right\} \subseteq P;$$

$M_2$ is a semivector subspace of P over S.

Consider

$$M_3 = \left\{ \begin{bmatrix} 0 & 0 & 0 & 0 \\ 0 & 0 & 0 & 0 \\ a_1 & a_2 & 0 & 0 \\ a_3 & a_4 & 0 & 0 \end{bmatrix} \middle| a_i \in Q^+ \cup \{0\}, 1 \le i \le 4 \right\} \subseteq P;$$

$M_3$ is a semivector subspace of P over S.

Now

$$M_4 = \left\{ \begin{bmatrix} 0 & 0 & 0 & 0 \\ 0 & 0 & 0 & 0 \\ 0 & 0 & a_1 & a_2 \\ 0 & 0 & a_3 & a_4 \end{bmatrix} \middle| a_i \in Q^+ \cup \{0\}, 1 \le i \le 4 \right\} \subseteq P$$

is a semivector subspace of P over S.

We see $P = M_1 + M_2 + M_3 + M_4$, $M_i \cap M_j = (0)$ if $i \ne j$; $1 \le i, j \le 4$.

*Example 4.68:* Let

$$V = \left\{ \begin{bmatrix} a_1 & a_2 & a_3 \\ a_4 & a_5 & a_6 \\ a_7 & a_8 & a_9 \\ a_{10} & a_{11} & a_{12} \\ a_{13} & a_{14} & a_{15} \end{bmatrix} \middle| a_i \in Q^+ \cup \{0\}, 1 \le i \le 15 \right\}$$

be a semivector space over the semifield $S = \{0\} \cup Z^+$.

So we can define as in case of other spaces complements in case of semivector space of polynomials with matrix coefficients also. We will only illustrate this situation by some examples.

*Example 4.69:* Let

$$V = \left\{ \sum_{i=0}^{\infty} a_i x^i \middle| a_i = (x_1, x_2, x_3, x_4, x_5, x_6) \mid \right.$$
$$\left. x_j \in Z^+ \cup \{0\}; 1 \le j \le 6 \right\}$$

be a semivector space over the semifield $S = Z^+ \cup \{0\}$.

Consider

$$M = \left\{ \sum_{i=0}^{\infty} a_i x^i \middle| a_i = (0, 0, 0, x_1, x_2, x_3) \right.$$
$$\left. \text{with } x_j \in Z^+ \cup \{0\}; 1 \le j \le 3 \right\} \subseteq V,$$

M is a semivector subspace of V over S.

Take

$$N = \left\{ \sum_{i=0}^{\infty} a_i x^i \,\middle|\, a_i = (x_1, x_2, x_3, 0,0,0,0) \right.$$

with $x_j \in Z^+ \cup \{0\}; 1 \le j \le 3\} \subseteq V$,

N is a semivector subspace of V over S. We see M+N = V and infact M is the orthogonal complement of N and vice versa. That is $M^\perp = N$ and $N^\perp = M$ and $M \cap N = (0)$.

***Example 4.70:*** Let

$$V = \left\{ \sum_{i=0}^{\infty} a_i x^i \,\middle|\, a_i = \begin{bmatrix} d_1 & d_2 & d_3 \\ d_4 & d_5 & d_6 \\ d_7 & d_8 & d_9 \\ d_{10} & d_{11} & d_{12} \\ d_{13} & d_{14} & d_{15} \end{bmatrix} ; d_j \in Z^+ \cup \{0\}; 1 \le j \le 15 \right\}$$

be a semivector space define over the semifield $S = Z^+ \cup \{0\}$.

Now

$$W_1 = \left\{ \sum_{i=0}^{\infty} a_i x^i \,\middle|\, x_i = \begin{bmatrix} x_1 & x_2 & x_3 \\ 0 & 0 & 0 \\ x_4 & x_5 & x_6 \\ 0 & 0 & 0 \\ x_7 & x_8 & x_9 \end{bmatrix} \right.$$

where $x_j \in Z^+ \cup \{0\}; 1 \le j \le 9\} \subseteq V$

is a semivector subspace of V over S.

We see the complement of W;

$$W_1^\perp = \left\{ \sum_{i=0}^{\infty} a_i x^i \,\middle|\, a_i = \begin{bmatrix} 0 & 0 & 0 \\ y_1 & y_2 & y_3 \\ 0 & 0 & 0 \\ y_4 & y_5 & y_6 \\ 0 & 0 & 0 \end{bmatrix} \right.$$

where $y_j \in Z^+ \cup \{0\}$; $1 \leq j \leq 6\} \subseteq V$.

We see $W_1^\perp + W_1 = V$ and $W_1 \cap W_1^\perp = (0)$.

Suppose

$$M_1 = \left\{ \sum_{i=0}^{\infty} a_i x^i \,\middle|\, x_i = \begin{bmatrix} 0 & 0 & 0 \\ 0 & 0 & 0 \\ 0 & 0 & 0 \\ x_1 & x_2 & x_3 \\ 0 & 0 & 0 \end{bmatrix} \right. \text{with } x_j \in Z^+ \cup \{0\};$$

$1 \leq j \leq 3\} \subseteq V$

be another semivector subspace of V; we see $M_1$ is not the orthogonal complement of $W_1$ but however $W_1 \cap M_1 = (0)$. But $W_1 + M_1 \subseteq V$. Hence we can have orthogonal semivector subspaces but they do not serve as the orthogonal complement of $W_1$.

*Example 4.71:* Let

$$V = \left\{ \sum_{i=0}^{\infty} a_i x^i \,\middle|\, a_i = \begin{bmatrix} d_1 & d_2 & d_3 \\ d_4 & d_5 & d_6 \\ d_7 & d_8 & d_9 \end{bmatrix} ; d_j \in Z^+ \cup \{0\}; 1 \leq j \leq 9 \right\}$$

be a semivector space over the semifield $S = Z^+ \cup \{0\}$.

Consider

$$P_1 = \left\{ \sum_{i=0}^{\infty} a_i x^i \,\middle|\, a_i = \begin{bmatrix} d_1 & d_2 & 0 \\ d_3 & 0 & 0 \\ 0 & 0 & 0 \end{bmatrix} \text{with } d_j \in Z^+ \cup \{0\}; \right.$$

$$\left. 1 \leq j \leq 3 \right\} \subseteq V$$

a semivector subspace of V over S. We see the complement of

$$P_1 \text{ is } P_2 = \left\{ \sum_{i=0}^{\infty} a_i x^i \,\middle|\, a_i = \begin{bmatrix} 0 & 0 & d_1 \\ 0 & d_2 & d_3 \\ d_4 & d_5 & d_6 \end{bmatrix} \text{with } d_j \in Z^+ \cup \{0\}; \right.$$

$$\left. 1 \leq j \leq 6 \right\} \subseteq V.$$

We see $P_1 + P_2 = V$ and $P_1 \cap P_2 = \{0\}$. We call $P_2$ the orthogonal complement of $P_1$ and vice versa.

However if

$$N = \left\{ \sum_{i=0}^{\infty} a_i x^i \,\middle|\, a_i = \begin{bmatrix} 0 & 0 & d_1 \\ 0 & 0 & d_2 \\ 0 & 0 & d_3 \end{bmatrix} \right.$$

with $d_j \in Z^+ \cup \{0\}$; $1 \leq j \leq 3 \} \subseteq V$;

N is also a semivector subspace of V over S and N is orthogonal with $P_1$ however N is not the orthogonal complement of $P_1$ as $P_1 + N \neq V$ only properly contained in V.

Thus a given semivector subspace can have more than one orthogonal semivector subspace but only one orthogonal complement.

For

$$T = \left\{ \sum_{i=0}^{\infty} a_i x^i \,\middle|\, a_i = \begin{bmatrix} 0 & 0 & 0 \\ 0 & d_1 & 0 \\ d_2 & d_3 & 0 \end{bmatrix} \right.$$

with $d_j \in Z^+ \cup \{0\}$; $1 \leq j \leq 3\} \subseteq V$

is such that T is a semivector subspace of V and T is also orthogonal with $P_1$ but is not the orthogonal complement of $P_1$ as $T + P_1 \subset V$. Thus we see there is a difference between a semivector subspace orthogonal with a semivector subspace and an orthogonal complement of a semivector subspace.

Now having seen examples of complement semivector subspace and orthogonal complement of a semivector subspace we now proceed onto give one or two examples of pseudo direct sum of semivector subspaces.

*Example 4.72:* Let

$$V = \left\{ \sum_{i=0}^{\infty} a_i x^i \,\middle|\, a_i = (x_1, x_2, \ldots, x_8) \right.$$

where $x_j \in Z^+ \cup \{0\}$; $1 \leq j \leq 8\}$

be a semivector space over the semifield $S = Z^+ \cup \{0\}$.

Take

$$W_1 = \left\{ \sum_{i=0}^{\infty} a_i x^i \,\middle|\, a_i = (0, 0, x_1, x_2, x_3, 0, 0, 0) \right.$$

where $x_j \in Z^+ \cup \{0\}$; $1 \leq j \leq 3\} \subseteq V$

be a semivector subspace of V over S.

Consider

$$W_2 = \left\{ \sum_{i=0}^{\infty} a_i x^i \,\middle|\, a_i = (x_1, x_2, x_3, 0,0,0,0,0) \right.$$

with $x_j \in Z^+ \cup \{0\}$; $1 \leq j \leq 3\} \subseteq V$;

another semivector subspace of V.

Take

$$W_3 = \left\{ \sum_{i=0}^{\infty} a_i x^i \,\middle|\, a_i = (0, 0, d_1, 0, 0, d_2, 0,0) \right.$$

with $d_j \in Z^+ \cup \{0\}$; $1 \leq j \leq 2\} \subseteq V$

another semivector subspace of V over S.

Finally let

$$W_4 = \left\{ \sum_{i=0}^{\infty} a_i x^i \,\middle|\, a_i = (x_1, 0, 0, x_2, 0,0, x_3, x_4) \right.$$

with $x_j \in Z^+ \cup \{0\}$; $1 \leq j \leq 4\} \subseteq V$,

a semivector subspace of V over S. We see $V \subseteq W_1 + W_2 + W_3 + W_4$ and $W_i \cap W_j \neq (0)$ if $i \neq j$. Thus V is a pseudo direct sum of semivector subspace of V over the semifields.

*Example 4.73:* Let

$$V = \left\{ \sum_{i=0}^{\infty} a_i x^i \;\middle|\; a_i = \begin{bmatrix} x_1 \\ x_2 \\ x_3 \\ x_4 \\ x_5 \\ x_6 \\ x_7 \\ x_8 \end{bmatrix} \text{ where } x_j \in Q^+ \cup \{0\}; \; 1 \le j \le 8 \right\}$$

be a semivector space over the semifield $S = Q^+ \cup \{0\}$.

Take

$$W_1 = \left\{ \sum_{i=0}^{\infty} a_i x^i \;\middle|\; a_i = \begin{bmatrix} x_1 \\ x_2 \\ x_3 \\ 0 \\ \vdots \\ 0 \end{bmatrix} \text{ where } x_j \in Q^+ \cup \{0\}; \; 1 \le j \le 3 \right\} \subseteq V$$

be a semivector subspace of V over S.

Consider

$$W_2 = \left\{ \sum_{i=0}^{\infty} a_i x^i \;\middle|\; a_i = \begin{bmatrix} 0 \\ 0 \\ x_1 \\ x_2 \\ x_3 \\ 0 \\ 0 \\ 0 \end{bmatrix} \text{ where } x_j \in Q^+ \cup \{0\}; \; 1 \le j \le 3 \right\} \subseteq V$$

be a semivector subspace of V over S.

$$W_3 = \left\{ \sum_{i=0}^{\infty} a_i x^i \,\middle|\, a_i = \begin{bmatrix} 0 \\ 0 \\ 0 \\ 0 \\ x_1 \\ x_2 \\ x_3 \\ 0 \end{bmatrix} \text{ where } x_j \in Q^+ \cup \{0\};\ 1 \leq j \leq 3 \right\} \subseteq V$$

be a semivector subspace of V over S and

$$W_4 = \left\{ \sum_{i=0}^{\infty} a_i x^i \,\middle|\, a_i = \begin{bmatrix} 0 \\ 0 \\ 0 \\ 0 \\ 0 \\ x_1 \\ x_2 \\ x_3 \end{bmatrix} \text{ where } x_j \in Q^+ \cup \{0\};\ 1 \leq j \leq 3 \right\} \subseteq V$$

be a semivector subspace of V over S, the semifield.

We see $W_i \cap W_j = (0)$; $i \neq j$. But $V \subseteq W_1 + W_2 + W_3 + W_4$; $1 \leq i, j \leq 4$. Thus V is the pseudo direct sum of semivector subspaces of V over S.

*Example 4.74:* Let

$$V = \left\{ \sum_{i=0}^{\infty} a_i x^i \,\middle|\, a_i = \begin{bmatrix} d_1 & d_2 & d_3 & d_4 \\ d_5 & d_6 & d_7 & d_8 \\ d_9 & d_{10} & d_{11} & d_{12} \\ d_{13} & d_{14} & d_{15} & d_{16} \end{bmatrix} \right.$$

$$\left. \text{where } d_j \in Z^+ \cup \{0\};\ 1 \leq j \leq 16 \right\}$$

be a semivector space over the semifield $S = Z^+ \cup \{0\}$.

Take

$$M_1 = \left\{\sum_{i=0}^{\infty} a_i x^i \;\middle|\; a_i = \begin{bmatrix} d_1 & d_2 & d_3 & 0 \\ 0 & 0 & 0 & 0 \\ 0 & 0 & 0 & 0 \\ 0 & 0 & 0 & 0 \end{bmatrix}\right.$$

where $d_1, d_2, d_3 \in Z^+ \cup \{0\}\} \subseteq V$,

$$M_2 = \left\{\sum_{i=0}^{\infty} a_i x^i \;\middle|\; a_i = \begin{bmatrix} 0 & 0 & d_1 & d_2 \\ d_3 & 0 & 0 & 0 \\ 0 & 0 & 0 & 0 \\ 0 & 0 & 0 & 0 \end{bmatrix}\right.$$

where $d_1, d_2, d_3 \in Z^+ \cup \{0\}\} \subseteq V$,

$$M_3 = \left\{\sum_{i=0}^{\infty} a_i x^i \;\middle|\; a_i = \begin{bmatrix} 0 & 0 & 0 & 0 \\ d_1 & d_2 & d_3 & 0 \\ 0 & 0 & 0 & 0 \\ 0 & 0 & 0 & 0 \end{bmatrix}\right.$$

where $d_1, d_2, d_3 \in Z^+ \cup \{0\}\} \subseteq V$,

$$M_4 = \left\{\sum_{i=0}^{\infty} a_i x^i \;\middle|\; a_i = \begin{bmatrix} 0 & 0 & 0 & 0 \\ 0 & 0 & d_1 & d_2 \\ d_3 & 0 & 0 & 0 \\ 0 & 0 & 0 & 0 \end{bmatrix}\right.$$

where $d_1, d_2, d_3 \in Z^+ \cup \{0\}\} \subseteq V$,

$$M_5 = \left\{ \sum_{i=0}^{\infty} a_i x^i \middle| a_i = \begin{bmatrix} 0 & 0 & 0 & 0 \\ 0 & 0 & 0 & 0 \\ d_1 & d_2 & d_3 & 0 \\ 0 & 0 & 0 & 0 \end{bmatrix} \right.$$

where $d_1, d_2, d_3 \in Z^+ \cup \{0\}\} \subseteq V$,

$$M_6 = \left\{ \sum_{i=0}^{\infty} a_i x^i \middle| a_i = \begin{bmatrix} 0 & 0 & 0 & 0 \\ 0 & 0 & 0 & 0 \\ 0 & d_1 & d_2 & d_3 \\ 0 & 0 & 0 & 0 \end{bmatrix} \right.$$

where $d_1, d_2, d_3 \in Z^+ \cup \{0\}\} \subseteq V$,

$$M_7 = \left\{ \sum_{i=0}^{\infty} a_i x^i \middle| a_i = \begin{bmatrix} 0 & 0 & 0 & 0 \\ 0 & 0 & 0 & 0 \\ 0 & 0 & 0 & d_1 \\ d_2 & d_3 & 0 & 0 \end{bmatrix} \right.$$

where $d_1, d_2, d_3 \in Z^+ \cup \{0\}\} \subseteq V$,

and

$$M_8 = \left\{ \sum_{i=0}^{\infty} a_i x^i \middle| a_i = \begin{bmatrix} 0 & 0 & 0 & 0 \\ 0 & 0 & 0 & 0 \\ 0 & 0 & 0 & 0 \\ 0 & d_1 & d_2 & d_3 \end{bmatrix} \right.$$

where $d_1, d_2, d_3 \in Z^+ \cup \{0\}\} \subseteq V$

be semivector subspaces of V over the semifield S.

Clearly $V \subseteq M_1 + M_2 + \ldots + M_8$, $M_i \cap M_j \neq (0)$ if $i \neq j$, $1 \leq j, i \leq 8$.

*Example 4.75:* Let

$$V = \left\{ \sum_{i=0}^{\infty} a_i x^i \,\middle|\, a_i = \begin{bmatrix} d_1 & d_2 & d_3 & d_4 & d_5 & d_6 \\ d_7 & d_8 & d_9 & d_{10} & d_{11} & d_{12} \\ d_{13} & d_{14} & d_{15} & d_{16} & d_{17} & d_{18} \\ d_{19} & d_{20} & d_{21} & d_{22} & d_{23} & d_{24} \end{bmatrix} \right.$$

where $d_j \in Z^+ \cup \{0\}$; $1 \leq j \leq 24\}$
be a semivector space over the semifield $S = Z^+ \cup \{0\}$.

Consider

$$P_1 = \left\{ \sum_{i=0}^{\infty} a_i x^i \,\middle|\, a_i = \begin{bmatrix} d_1 & d_2 & d_3 & d_4 & 0 & 0 \\ 0 & 0 & 0 & 0 & 0 & 0 \\ 0 & 0 & 0 & 0 & 0 & 0 \\ 0 & 0 & 0 & 0 & 0 & 0 \end{bmatrix} \right.$$

where $d_1, d_2, d_3, d_4 \in Z^+ \cup \{0\}\} \subseteq V$,

$$P_2 = \left\{ \sum_{i=0}^{\infty} a_i x^i \,\middle|\, a_i = \begin{bmatrix} 0 & 0 & 0 & d_1 & d_2 & d_3 \\ d_4 & 0 & 0 & 0 & 0 & 0 \\ 0 & 0 & 0 & 0 & 0 & 0 \\ 0 & 0 & 0 & 0 & 0 & 0 \end{bmatrix} \right.$$

with $d_1, d_2, d_3, d_4 \in Z^+ \cup \{0\}\} \subseteq V$,

$$P_3 = \left\{ \sum_{i=0}^{\infty} a_i x^i \,\middle|\, a_i = \begin{bmatrix} 0 & 0 & 0 & 0 & 0 & 0 \\ d_1 & d_2 & d_3 & d_4 & 0 & 0 \\ 0 & 0 & 0 & 0 & 0 & 0 \\ 0 & 0 & 0 & 0 & 0 & 0 \end{bmatrix} \right.$$

with $d_1, d_2, d_3, d_4 \in Z^+ \cup \{0\}\} \subseteq V$,

$$P_4 = \left\{ \sum_{i=0}^{\infty} a_i x^i \,\middle|\, a_i = \begin{bmatrix} 0 & 0 & 0 & 0 & 0 & 0 \\ 0 & 0 & 0 & d_1 & d_2 & d_3 \\ d_4 & 0 & 0 & 0 & 0 & 0 \\ 0 & 0 & 0 & 0 & 0 & 0 \end{bmatrix} \right.$$

with $d_1, d_2, d_3, d_4 \in Z^+ \cup \{0\}\} \subseteq V$,

$$P_5 = \left\{ \sum_{i=0}^{\infty} a_i x^i \,\middle|\, a_i = \begin{bmatrix} 0 & 0 & 0 & 0 & 0 & 0 \\ 0 & 0 & 0 & 0 & 0 & 0 \\ d_1 & d_2 & d_3 & d_4 & 0 & 0 \\ 0 & 0 & 0 & 0 & 0 & 0 \end{bmatrix} \right.$$

with $d_1, d_2, d_3, d_4 \in Z^+ \cup \{0\}\} \subseteq V$,

$$P_6 = \left\{ \sum_{i=0}^{\infty} a_i x^i \,\middle|\, a_i = \begin{bmatrix} 0 & 0 & 0 & 0 & 0 & 0 \\ 0 & 0 & 0 & 0 & 0 & 0 \\ 0 & 0 & 0 & d_1 & d_2 & d_3 \\ d_4 & 0 & 0 & 0 & 0 & 0 \end{bmatrix} \right.;$$

$d_1, d_2, d_3, d_4 \in Z^+ \cup \{0\}\} \subseteq V$,

$$P_7 = \left\{ \sum_{i=0}^{\infty} a_i x^i \,\middle|\, a_i = \begin{bmatrix} 0 & 0 & 0 & 0 & 0 & 0 \\ 0 & 0 & 0 & 0 & 0 & 0 \\ 0 & 0 & 0 & 0 & 0 & 0 \\ d_1 & d_2 & d_3 & d_4 & 0 & 0 \end{bmatrix} \right.$$

with $d_1, d_2, d_3, d_4 \in Z^+ \cup \{0\}\} \subseteq V$

and

$$P_8 = \left\{ \sum_{i=0}^{\infty} a_i x^i \;\middle|\; a_i = \begin{bmatrix} 0 & 0 & 0 & 0 & 0 & 0 \\ 0 & 0 & 0 & 0 & 0 & 0 \\ 0 & 0 & 0 & 0 & 0 & 0 \\ 0 & d_1 & d_2 & d_3 & d_4 & d_5 \end{bmatrix} \right.$$

with $d_1$, $d_2$, $d_3$, $d_4$, $d_5 \in Z^+ \cup \{0\}\} \subseteq V$

be semivector subspaces of V over $S = Z^+ \cup \{0\}$. We see $V \subseteq P_1 + P_2 + P_3 + P_4 + P_5 + P_6 + P_7 + P_8$; and $P_i \cap P_j \neq \{0\}$ if $i \neq j$, $1 \leq i, j \leq 8$. Thus V is only a pseudo direct sum of semivector subspaces.

**Chapter Five**

# NATURAL PRODUCT ON SUPERMATRICES

In this chapter we define the new notion of natural product in supermatrices. Products in supermatrices are very different from usual product on matrices and product on super matrices.

Throughout this chapter

$F_R^S = \{(a_1\ a_2\ a_3\ |\ a_4\ a_5\ |\ \ldots\ |\ a_{n-1}\ a_n)\ |\ a_i \in Q \text{ or } R \text{ or } Z\}$ that is collection of $1 \times n$ super row matrices with same type of partition in it.

$$F_C^S = \left\{ \begin{bmatrix} a_1 \\ a_2 \\ \hline a_3 \\ a_4 \\ a_5 \\ \vdots \\ \hline a_m \end{bmatrix} \;\middle|\; a_i \in Z \text{ or } Q \text{ or } R;\ 1 \leq i \leq m \right\}$$

denotes the collection of all m × 1 super column matrices with same type of partition on it.

$$F^S_{m \times n} \ (m \neq n) = \left\{ \begin{bmatrix} \begin{array}{ccc|c} a_{11} & a_{12} & \cdots & a_{1n} \\ a_{21} & a_{22} & \cdots & a_{2n} \\ \vdots & \vdots & \vdots & \vdots \\ \hline a_{m1} & a_{m2} & \cdots & a_{mn} \end{array} \end{bmatrix} \ \middle| \ a_{ij} \in Q \text{ or } R \text{ or } Z; \right.$$

$$1 \leq i \leq m;\ 1 \leq j \leq n \}$$

denotes the collection of all m × n super matrices with same type of partition on it.

$$F^S_{n \times n} = \left\{ \begin{bmatrix} \begin{array}{ccc|c} a_{11} & a_{12} & \cdots & a_{1n} \\ a_{21} & a_{22} & \cdots & a_{2n} \\ \vdots & \vdots & \vdots & \vdots \\ \hline a_{n1} & a_{n2} & \cdots & a_{nn} \end{array} \end{bmatrix} \ \middle| \ a_{ij} \in Q \text{ or } Z \text{ or } R; \right.$$

$$1 \leq i, j \leq n \}$$

denotes the collection of n × n super matrices with same type of partition on it.

We will first illustrate this situation before we proceed onto give any form of algebraic structure on them.

*Example 5.1:* Let

$M = \{(x_1\ x_2\ |\ x_3\ |\ x_4\ x_5)$ where $x_i \in Z$ or $Q$ or $R;\ 1 \leq i \leq 5\}$

be the collection of 1 × 5 row super matrices with same type of partition on it $M = F^S_R$.

*Example 5.2:* Let

$$F_R^S = \{(x_1 \mid x_2 \ x_3 \mid x_4 \ x_5 \ x_6 \mid x_7 \ x_8 \ x_9 \ x_{10} \mid x_{11} \ x_{12}) \mid x_i \in R; \ 1 \leq i \leq 10\}$$

be again a collection of $1 \times 10$ super row matrices with same type of partition on it.

*Example 5.3:* Let

$$P = \{(x_1 \ x_2 \mid x_3 \ x_4 \ x_5) \text{ where } x_i \in Q \text{ or } R \text{ or } Z; \ 1 \leq i \leq 5\}$$

be again a collection of $1 \times 5$ super row super matrices of same type.

Now we will see examples of column super matrices of same type.

*Example 5.4:* Let

$$F_C^S = \left\{ \begin{bmatrix} x_1 \\ x_2 \\ \hline x_3 \\ x_4 \\ x_5 \\ \hline x_6 \\ \hline x_7 \end{bmatrix} \; \middle| \; x_i \in R; \ 1 \leq i \leq 7 \right\}$$

be the $7 \times 1$ column super matrices of same type.

*Example 5.5:* Let

$$F_C^S = \left\{ \begin{bmatrix} x_1 \\ x_2 \\ x_3 \\ \hline x_4 \\ x_5 \\ \hline x_6 \\ \hline x_7 \\ x_8 \\ x_9 \end{bmatrix} \middle| \; x_i \in R; \; 1 \le i \le 9 \right\}$$

be the $9 \times 1$ column super matrix of same type.

*Example 5.6:* Let

$$F_C^S = \left\{ \begin{bmatrix} a_1 \\ a_2 \\ a_3 \\ \hline a_4 \\ \hline a_5 \\ a_6 \\ a_7 \\ a_8 \end{bmatrix} \middle| \; a_i \in Q; \; 1 \le i \le 8 \right\}$$

be the $8 \times 1$ column super matrix of same type.

Now we proceed onto give examples of $F_{m \times n}^S$ $(m \ne n)$.

*Example 5.7:* Let

$$F^S_{3\times 5} = \left\{ \left[ \begin{array}{cc|c|cc} a_1 & a_2 & a_3 & a_4 & a_5 \\ a_6 & a_7 & a_8 & a_9 & a_{10} \\ a_{11} & a_{12} & a_{13} & a_{14} & a_{15} \end{array} \right] \middle| \; a_i \in Q; \; 1 \leq i \leq 15 \right\}$$

be the $3 \times 5$ super matrix of same type.

*Example 5.8:* Let

$$F^S_{7\times 4} = \left\{ \left[ \begin{array}{c|cc|c} a_1 & a_2 & a_3 & a_4 \\ a_5 & a_6 & a_7 & a_8 \\ \vdots & \vdots & \vdots & \vdots \\ \hline a_{25} & a_{26} & a_{27} & a_{28} \end{array} \right] \middle| \; a_i \in Z; \; 1 \leq i \leq 28 \right\}$$

be a $7 \times 4$ super matrix of same type.

*Example 5.9:* Let

$$F^S_{6\times 7} = \left\{ \left[ \begin{array}{cccc|cc} a_1 & a_2 & a_3 & a_4 & a_5 & a_6 \\ a_7 & a_8 & a_9 & a_{10} & a_{11} & a_{12} \\ \vdots & \vdots & \vdots & \vdots & \vdots & \vdots \\ \hline a_{31} & a_{32} & a_{33} & a_{34} & a_{35} & a_{36} \\ a_{37} & a_{38} & a_{39} & a_{40} & a_{41} & a_{42} \end{array} \right] \middle| \; a_i \in Q; \; 1 \leq i \leq 42 \right\}$$

be a $6 \times 7$ super matrix of same type.

Now we proceed onto give examples of $F^S_{n\times n}$ square super matrices of same type.

*Example 5.10:* Let

$$F^S_{4\times 4} = M = \left\{ \left[ \begin{array}{cc|c|cc} a_1 & a_2 & a_3 & a_4 \\ \hline a_5 & a_6 & a_7 & a_8 \\ a_9 & a_{10} & a_{11} & a_{12} \\ a_{13} & a_{14} & a_{15} & a_{16} \end{array} \right] \right. \left. a_i \in Q; 1 \leq i \leq 16 \right\}$$

be a square super matrix of same type.

*Example 5.11:* Let

$$F^S_{4\times 4} = \left\{ \left[ \begin{array}{cc|cc} a_1 & a_2 & a_3 & a_4 \\ a_5 & a_6 & a_7 & a_8 \\ \hline a_9 & a_{10} & a_{11} & a_{12} \\ a_{13} & a_{14} & a_{15} & a_{16} \end{array} \right] \right. \left. a_i \in Q; 1 \leq i \leq 16 \right\}$$

be a square super matrix of same type.

*Example 5.12:* Let

$$F^S_{4\times 4} = \left\{ \left[ \begin{array}{c|c|c|c} a_1 & a_2 & a_3 & a_4 \\ a_5 & a_6 & a_7 & a_8 \\ a_9 & a_{10} & a_{11} & a_{12} \\ a_{13} & a_{14} & a_{15} & a_{16} \end{array} \right] \right. \left. a_i \in Q; 1 \leq i \leq 16 \right\}$$

be a square super matrix of same type.

*Example 5.13:* Let

$$F^S_{3\times 3} = \left\{ \left[ \begin{array}{ccc} a_1 & a_2 & a_3 \\ \hline a_4 & a_5 & a_6 \\ \hline a_7 & a_8 & a_9 \end{array} \right] \right. \left. a_i \in Z; 1 \leq i \leq 9 \right\}$$

be a square super matrix of same order.

Now we can define $F_C^S, F_R^S, F_{m \times n}^S$ (m ≠ n) and $F_{n \times n}^S$ the usual matrix addition and natural product $\times_n$.

Under usual matrix addition $F_C^S, F_R^S, F_{m \times n}^S$ (m ≠ n) and $F_{n \times n}^S$ are abelian (commutative) groups.

How under natural product $F_C^S, F_R^S, F_{m \times n}^S$ (m ≠ n) and $F_{n \times n}^S$ are semigroups with unit.

Suppose

$x = (x_1 \; x_2 \; x_3 \mid x_4 \; x_5 \mid x_6 \; x_7)$ and $y = (y_1 \; y_2 \; y_3 \mid y_4 \; y_5 \mid y_6 \; y_7)$ be two super row matrices of same type

$x + y = (x_1 + y_1, x_2 + y_2, x_3 + y_3 \mid x_4 + y_4, x_5 + y_5 \mid x_6 + y_6, x_7 + y_7)$. Thus $F_R^S$ is closed under '+'.

Likewise if $x = \begin{bmatrix} x_1 \\ x_2 \\ x_3 \\ \hline x_4 \\ x_5 \\ \hline x_6 \\ x_7 \\ x_8 \end{bmatrix}$ and $y = \begin{bmatrix} y_1 \\ y_2 \\ y_3 \\ \hline y_4 \\ y_5 \\ \hline y_6 \\ y_7 \\ y_8 \end{bmatrix}$ are two super column matrices of same type then $x + y = \begin{bmatrix} x_1 + y_1 \\ \hline x_2 + y_2 \\ \hline x_3 + y_3 \\ \hline x_4 + y_4 \\ \hline x_5 + y_5 \\ \hline x_6 + y_6 \\ x_7 + y_7 \\ x_8 + y_8 \end{bmatrix}$.

We see $F_C^S$ is closed under '+' and infact a group under '+'.

Consider $x = \begin{bmatrix} a_1 & a_2 & a_3 & a_4 & a_5 & a_6 & a_7 \\ a_8 & a_9 & a_{10} & a_{11} & a_{12} & a_{13} & a_{14} \\ a_{15} & a_{16} & a_{17} & a_{18} & a_{19} & a_{20} & a_{21} \\ a_{22} & a_{23} & a_{24} & a_{25} & a_{26} & a_{27} & a_{28} \\ a_{29} & a_{30} & a_{31} & a_{32} & a_{33} & a_{34} & a_{35} \end{bmatrix}$ and

$y = \begin{bmatrix} b_1 & b_2 & b_3 & b_4 & b_5 & b_6 & b_7 \\ b_8 & b_9 & b_{10} & b_{11} & b_{12} & b_{13} & b_{14} \\ b_{15} & b_{16} & b_{17} & b_{18} & b_{19} & b_{20} & b_{21} \\ b_{22} & b_{23} & b_{24} & b_{25} & b_{26} & b_{27} & b_{28} \\ b_{29} & b_{30} & b_{31} & b_{32} & b_{33} & b_{34} & b_{35} \end{bmatrix}$

be two $5 \times 7$ super matrices in $F_{5 \times 7}^S$.

Now x + y =

$$\begin{bmatrix} a_1+b_1 & a_2+b_9 & a_3+b_3 & a_4+b_4 & a_5+b_5 & a_6+b_6 & a_7+b_7 \\ a_8+b_8 & a_9+b_9 & a_{10}+b_{10} & a_{11}+b_{11} & a_{12}+b_{12} & a_{13}+b_{13} & a_{14}+b_{14} \\ \hline a_{15}+b_{15} & a_{16}+b_{16} & a_{17}+b_{17} & a_{18}+b_{18} & a_{19}+b_{19} & a_{20}+b_{20} & a_{21}+b_{17} \\ \hline a_{22}+b_{22} & a_{23}+b_{23} & a_{24}+b_{24} & a_{25}+b_{25} & a_{26}+b_{26} & a_{27}+b_{27} & a_{28}+b_{28} \\ a_{29}+b_{29} & a_{30}+b_{30} & a_{31}+b_{31} & a_{32}+b_{32} & a_{33}+b_{33} & a_{34}+b_{34} & a_{35}+b_{35} \end{bmatrix}$$

is in $F^S_{5\times 7}$.

Thus addition can be performed on $F^S_{m\times n}$ ($m \neq n$) and infact $F^S_{m\times n}$ is a group under addition.

Now we give examples of addition of square super matrices $F^S_{n\times n}$.

Let x = $\begin{bmatrix} a_1 & a_2 & a_3 & a_4 & a_5 \\ a_6 & a_7 & a_8 & a_9 & a_{10} \\ a_{11} & a_{12} & a_{13} & a_{14} & a_{15} \\ a_{16} & a_{17} & a_{18} & a_{19} & a_{20} \\ a_{21} & a_{22} & a_{23} & a_{24} & a_{25} \end{bmatrix}$

and y = $\begin{bmatrix} b_1 & b_2 & b_3 & b_4 & b_5 \\ b_6 & b_7 & b_8 & b_9 & b_{10} \\ b_{11} & b_{12} & b_{13} & b_{14} & b_{15} \\ b_{16} & b_{17} & b_{18} & b_{19} & b_{20} \\ b_{21} & b_{22} & b_{23} & b_{24} & b_{25} \end{bmatrix}$

$$x+y = \begin{bmatrix} a_1+b_1 & a_2+b_2 & a_3+b_3 & a_4+b_4 & a_5+b_5 \\ \hline a_6+b_6 & a_7+b_7 & a_8+b_8 & a_9+b_9 & a_{10}+b_{10} \\ a_{11}+b_{11} & a_{12}+b_{12} & a_{13}+b_{13} & a_{14}+b_{14} & a_{15}+b_{15} \\ a_{16}+b_{16} & a_{17}+b_{17} & a_{18}+b_{18} & a_{19}+b_{19} & a_{20}+b_{20} \\ \hline a_{21}+b_{21} & a_{22}+b_{22} & a_{23}+b_{23} & a_{24}+b_{24} & a_{25}+b_{25} \end{bmatrix}$$

$\in F_{5\times 5}^S$.

Infact $F_{5\times 5}^S$ is a group under addition.

Now we proceed onto define natural product $\times_n$ on $F_C^S, F_R^S, F_{n\times m}^S$ ($n \neq m$) and $F_{n\times n}^S$.

Consider $x = (a_1\ a_2\ a_3 \mid a_4\ a_5 \mid a_6\ a_7\ a_8\ a_9)$ and $y = (b_1\ b_2\ b_3 \mid b_4\ b_5 \mid b_6\ b_7\ b_8\ b_9) \in F_R^S$.

$x \times_n y = (a_1b_1\ a_2b_2\ a_3b_3 \mid a_4b_4\ a_5b_5 \mid a_6b_6\ a_7b_7\ a_8b_8\ a_9b_9) \in F_R^S$. $F_R^S$ under product is a semigroup infact $F_R^S$ has zero divisors under natural product $\times_n$.

Suppose $x = (0\ 0\ 0 \mid 2\ 1\ 0 \mid 9\ 2 \mid 3)$ and $y = (3\ 9\ 0 \mid 0\ 0\ 7 \mid 0\ 0 \mid 0)$ be in $F_R^S$. $x \times_n y = (0\ 0\ 0 \mid 0\ 0\ 0 \mid 0\ 0 \mid 0)$. Thus we see $F_R^S$ has zero divisors.

Consider $x = \begin{bmatrix} 2 \\ 0 \\ 3 \\ \hline \frac{1}{2} \\ 0 \\ \hline 0 \\ \hline \frac{4}{5} \end{bmatrix}$ and $y = \begin{bmatrix} 0 \\ 1 \\ 0 \\ \hline 0 \\ 1 \\ \hline \frac{2}{3} \\ \hline 0 \end{bmatrix}$ in $F_C^S$;

we see under the natural product $x \times_n y = \begin{bmatrix} 0 \\ 0 \\ 0 \\ 0 \\ \hline 0 \\ 0 \\ 0 \\ \hline 0 \\ 0 \\ \hline 0 \\ 0 \end{bmatrix}$.

Take $F_{3\times 5}^S$, $F_{3\times 5}^S$ under natural product, they have zero divisors. Infact $F_{3\times 5}^S$ under natural product $\times_n$ is a semigroup.

Consider $x = \begin{bmatrix} 9 & 0 & 2 & 0 & 1 \\ 0 & 1 & 0 & 5 & 0 \\ 1 & 0 & 0 & 2 & 0 \end{bmatrix}$ and

$y = \begin{bmatrix} 0 & 7 & 0 & 8 & 0 \\ 9 & 0 & 2 & 0 & 7 \\ 0 & 7 & 9 & 0 & 2 \end{bmatrix}$ in $F_{3\times 5}^S$;

we see $x \times_n y = \begin{bmatrix} 0 & 0 & 0 & 0 & 0 \\ 0 & 0 & 0 & 0 & 0 \\ 0 & 0 & 0 & 0 & 0 \end{bmatrix}$.

Now $F_{n\times n}^S$ also has zero divisors.

Consider $x = \left[\begin{array}{c|cc|ccc} 7 & 8 & 0 & 9 & 4 & 2 \\ \hline 0 & 1 & 2 & 5 & 7 & 8 \\ 1 & 2 & 3 & 0 & 1 & 0 \\ \hline 5 & 7 & 0 & 9 & 2 & 0 \\ 1 & 2 & 3 & 0 & 2 & 3 \\ 0 & 8 & 7 & 0 & 5 & 4 \end{array}\right]$ and

$y = \left[\begin{array}{c|cc|ccc} 0 & 0 & 9 & 0 & 0 & 0 \\ \hline 7 & 0 & 0 & 0 & 0 & 0 \\ 0 & 0 & 0 & 6 & 0 & 8 \\ \hline 0 & 0 & 6 & 0 & 0 & 2 \\ 0 & 0 & 0 & 6 & 0 & 0 \\ 5 & 0 & 0 & 7 & 0 & 0 \end{array}\right] \in F_{6\times 6}^S$.

$x \times_n y = \left[\begin{array}{c|cc|ccc} 0 & 0 & 0 & 0 & 0 & 0 \\ \hline 0 & 0 & 0 & 0 & 0 & 0 \\ 0 & 0 & 0 & 0 & 0 & 0 \\ \hline 0 & 0 & 0 & 0 & 0 & 0 \\ 0 & 0 & 0 & 0 & 0 & 0 \\ 0 & 0 & 0 & 0 & 0 & 0 \end{array}\right] \in F_{6\times 6}^S$.

Thus $F_{n\times n}^S$ is a semigroup under $\times_n$ and has zero divisors and ideals. We will now give the following theorems the proofs of which are simple.

**THEOREM 5.1**:

$$F_R^S = \{(x_1 \, x_2 \mid ... \mid x_n) \mid x_i \in Q \text{ or } R; \, 1 \leq i \leq n\}$$

*is a group under '+'.*

**THEOREM 5.2:**

$$F_C^S = \left\{ \begin{bmatrix} x_1 \\ x_2 \\ \hline x_3 \\ \vdots \\ \hline \vdots \\ \hline x_{m-1} \\ \hline x_m \end{bmatrix} \mid x_i \in Q \text{ or } R \text{ or } C \text{ or } Z; \ 1 \leq i \leq m \right\}$$

*is the group under '+'.*

**THEOREM 5.3 :**

$$F_{3\times 3}^S \ (m \neq n) = \left\{ \begin{bmatrix} a_{11} & a_{12} & \cdots & a_{1n} \\ a_{21} & a_{22} & \cdots & a_{2n} \\ \vdots & \vdots & & \vdots \\ \hline a_{m1} & a_{m2} & \cdots & a_{mn} \end{bmatrix} \mid a_{ij} \in Q \text{ or } R \text{ or } C \text{ or } Z; \right.$$

$$1 \leq i \leq m; \ 1 \leq j \leq n \}$$

*is a group under '+'.*

**THEOREM 5.4:** *($F_{n\times n}^S$, +) is a group.*

**THEOREM 5.5:** *($F_R^S$, $\times_n$) is a semigroup and has zero divisors, ideals and subsemigroups which are not ideals.*

**THEOREM 5.6**: *($F_C^S$, $\times_n$) is a semigroup with unit and has zero divisors, units, ideals, and subsemigroups.*

**THEOREM 5.7**: *($F_{m\times n}^S (m \neq n)$, $\times_n$) is a semigroup with zero divisors and units.*

**THEOREM 5.8**: $(F_{n\times n}^S, \times_n)$ *is a commutative semigroup and has zero divisors units and ideals.*

We will now give examples of zero divisors units and ideals of $F_{m\times n}^S$ ($m \neq n$), $F_C^S, F_R^S, F_{n\times m}^S$ ($n \neq m$) and $F_{n\times n}^S$.

*Example 5.14:* Let

$$F_R^S = \{(x_1 \mid x_2\ x_3 \mid x_4) \text{ where } x_i \in Z; 1 \leq i \leq 4\}$$

be a commutative semigroup under natural product. $(1 \mid 1\ 1 \mid 1)$ is the unit of $F_R^S$ under $\times_n$.

$P = \{(x_1 \mid x_2\ x_3 \mid x_4) \mid x_i \in 3Z; 1 \leq i \leq 4\} \subseteq F_R^S$ is an ideal of $F_R^S$.

Infact $F_R^S$ has infinite number of ideals under the natural product $\times_n$. Further $F_R^S$ has zero divisors.

Consider $x = (x_1 \mid 0\ 0 \mid x_2) \in F_R^S$, $y = (0 \mid y_1\ y_2 \mid 0)$ in $F_R^S$ is such that $x \times_n y = (0 \mid 0\ 0 \mid 0)$. Also $P = (x_1 \mid 0\ 0 \mid x_2) \mid x_i \in Z, 1 \leq i \leq 2\} \subseteq F_R^S$ is also an ideal.

*Example 5.15:* Let
$$F_R^S = \{(x_1 \mid x_2\ x_3 \mid x_4\ x_5\ x_6) \text{ where } x_i \in Q; 1 \leq i \leq 6\}$$
be a semigroup under $\times_n$.

$S = \{(a_1 \mid a_2\ a_3 \mid a_4\ a_5\ a_6) \mid a_i \in Z, 1 \leq i \leq 6\} \subseteq F_R^S$; $S$ is a subsemigroup of $F_R^S$ under $\times_n$. Clearly $S$ is not an ideal of $F_R^S$.

Consider $P = \{(a_1 \mid a_2\ a_3 \mid 0\ 0\ 0) \mid a_i \in Q, 1 \leq i \leq 3\} \subseteq F_R^S$, is an ideal of $F_R^S$. If Q in P is replaced by Z that is $T = \{(a_1 \mid a_2 \ a_3 \mid 0 \ 0 \ 0) \mid a_i \in Z, 1 \leq i \leq 3\} \subseteq F_R^S$, then T is only a subsemigroup of $F_R^S$ and is not an ideal of $F_R^S$. Thus $F_R^S$ has subsemigroups which are not ideals.

Take $M = \{(a \mid b \ c \mid 0 \ 0 \ 0) \mid a, b, c \in Q\} \subseteq F_R^S$ and
$N = \{(0 \mid 0 \ 0 \mid a \ b \ c) \mid a, b, c \in Q\} \subseteq F_R^S$.

$M \times_n N = (0 \mid 0 \ 0 \mid 0 \ 0 \ 0)$ or $M \cap N = (0 \mid 0 \ 0 \mid 0 \ 0 \ 0)$.

$F_R^S = M + N$.

Suppose $M_1 = \{(a \mid b \ 0 \mid 0 \ 0 \ 0) \mid a, b \in Q\} \subseteq F_R^S$ and

$N_1 = \{(0 \mid 0 \ 0 \mid a \ b \ 0) \mid a, b \in Q\} \subseteq F_R^S$, we see

$M_1 \times_n N_1 = \{(0 \mid 0 \ 0 \mid 0 \ 0 \ 0)\}$. Also $M_1 \cap N_1 = \{(0 \mid 0 \ 0 \mid 0 \ 0 \ 0)\}$ but $N_1 + M_1 \subsetneq F_R^S$; and $N_1 + M_1 \neq F_R^S$. We see that some special properties are enjoyed by M and N that are not true in case $M_1$ and $N_1$.

Now we give an example in case of $(F_C^S, \times_n)$.

*Example 5.16:* Let

$$F_C^S = \left\{ \left[ \begin{array}{c} \underline{a_1} \\ \underline{a_2} \\ \underline{a_3} \\ a_4 \\ \underline{a_5} \\ \underline{a_6} \\ a_7 \end{array} \right] \middle| \ a_i \in Q, 1 \leq i \leq 7 \right\}$$

be a semigroup under natural product $\times_n$.

Consider

$$P = \left\{ \begin{bmatrix} \begin{bmatrix} a_1 \\ \hline a_2 \\ \hline a_3 \\ \hline a_4 \\ \hline 0 \\ \hline 0 \\ \hline 0 \end{bmatrix} \end{bmatrix} \;\middle|\; a_i \in Q, 1 \leq i \leq 4 \right\} \subseteq F_C^S$$

is a subsemigroup of $F_C^S$ and is not an ideal of $F_C^S$.

Take

$$M = \left\{ \begin{bmatrix} \begin{bmatrix} 0 \\ \hline 0 \\ \hline 0 \\ \hline 0 \\ \hline a_1 \\ \hline a_2 \\ \hline a_3 \end{bmatrix} \end{bmatrix} \;\middle|\; a_i \in Z, 1 \leq i \leq 3 \right\} \subseteq F_C^S,$$

M is a subsemigroup of $F_C^S$.

Clearly if $x \in P$ and $y \in M$ then

$$x \times_n y = \left\{ \begin{bmatrix} 0 \\ \hline 0 \\ \hline 0 \\ \hline 0 \\ \hline 0 \\ \hline 0 \\ \hline 0 \end{bmatrix} \right\}.$$

Now take

$$S = \left\{ \begin{bmatrix} \overline{a_1} \\ \overline{a_2} \\ \overline{a_3} \\ \overline{a_4} \\ \overline{a_5} \\ \overline{a_6} \\ \overline{a_7} \end{bmatrix} \middle| a_i \in Z, 1 \leq i \leq 7 \right\} \subseteq F_C^S,$$

S is a subsemigroup of $F_C^S$ but is not an ideal of $F_C^S$.

Suppose

$$J = \left\{ \begin{bmatrix} \overline{a_1} \\ \overline{0} \\ \overline{0} \\ \overline{0} \\ \overline{a_2} \\ \overline{a_3} \\ \overline{a_4} \end{bmatrix} \middle| a_i \in Q, 1 \leq i \leq 4 \right\} \subseteq F_C^S,$$

J is a subsemigroup as well as an ideal of $F_C^S$.

Now we give yet another example.

*Example 5.17:* Let

$$F_C^S = \left\{ \begin{bmatrix} a_1 \\ \hline a_2 \\ a_3 \\ \hline a_4 \\ \hline a_5 \\ a_6 \\ a_7 \end{bmatrix} \;\middle|\; a_i \in Q, 1 \leq i \leq 7 \right\}$$

be a semigroup under natural product $\times_n$.

Take

$$P = \left\{ \begin{bmatrix} 0 \\ \hline a_1 \\ 0 \\ \hline a_2 \\ \hline 0 \\ a_3 \\ 0 \end{bmatrix} \;\middle|\; a_i \in Q, 1 \leq i \leq 3 \right\} \subseteq F_C^S$$

is a subsemigroup as well as an ideal of $F_C^S$.

Take

$$M = \left\{ \begin{bmatrix} 0 \\ \hline a_1 \\ 0 \\ \hline a_2 \\ \hline 0 \\ a_3 \\ 0 \end{bmatrix} \;\middle|\; a_i \in 3Z, 1 \leq i \leq 3 \right\} \subseteq F_C^S$$

is a subsemigroup of $F_C^S$ and is not an ideal of $F_C^S$. $M \subseteq P$; M is a subsemigroup of P also.

Now take

$$T = \left\{ \begin{bmatrix} a_1 \\ \hline 0 \\ \hline a_2 \\ \hline 0 \\ \hline a_3 \\ \hline 0 \\ \hline a_4 \end{bmatrix} \middle| \ a_i \in Q, 1 \le i \le 4 \right\} \subseteq F_C^S,$$

T is a subsemigroup of $F_C^S$ also an ideal of $F_C^S$. We see for every $x \in P$ and for every $y \in T$, $x \times_n y = (0)$.

Now we give examples of zero divisors and ideals in $F_{n \times m}^S$ $(n \ne m)$.

*Example 5.18:* Let

$$F_{5 \times 3}^S = \left\{ \begin{bmatrix} a_1 & a_2 & a_3 & a_4 & a_5 \\ a_6 & a_7 & a_8 & a_9 & a_{10} \\ a_{11} & a_{12} & a_{13} & a_{14} & a_{15} \end{bmatrix} \middle| \ a_i \in Q, 1 \le i \le 15 \right\}$$

be a semigroup of $5 \times 3$ super matrices under natural multiplication.

If $x = \begin{bmatrix} 1 & 2 & 3 & 4 & 5 \\ 9 & 8 & 7 & 6 & 5 \\ 0 & 1 & 2 & 7 & 1 \end{bmatrix}$ and $y = \begin{bmatrix} 0 & 1 & 2 & 3 & 5 \\ 9 & 0 & 1 & 3 & 4 \\ 7 & 2 & 3 & 1 & 2 \end{bmatrix}$

are in $F_{5 \times 3}^S$;

then $x \times_n y = \begin{bmatrix} 0 & 2 & 6 & 12 & 25 \\ 81 & 0 & 7 & 18 & 20 \\ 0 & 2 & 6 & 7 & 2 \end{bmatrix}$.

Now consider $P = \left\{ \begin{bmatrix} a & b & c & d & e \\ 0 & 0 & 0 & 0 & 0 \\ 0 & 0 & 0 & 0 & 0 \end{bmatrix} \middle| a, b, c, d, e \in Q \right\} \subseteq F_{5\times3}^S$; P is an ideal of $F_{5\times3}^S$.

We see for $a = \begin{bmatrix} 0 & 0 & 0 & 0 & 0 \\ x & y & a & b & z \\ 0 & 0 & c & d & m \end{bmatrix} \in F_{5\times3}^S$ is such that

$a \times_n x = \begin{bmatrix} 0 & 0 & 0 & 0 & 0 \\ 0 & 0 & 0 & 0 & 0 \\ 0 & 0 & 0 & 0 & 0 \end{bmatrix}$ for every $x \in P$.

Take

$M = \left\{ \begin{bmatrix} 0 & 0 & e & f & 0 \\ a & b & 0 & 0 & g \\ c & d & 0 & 0 & h \end{bmatrix} \middle| a, b, c, d, e, f, g, h \in Z \right\} \subseteq F_{5\times3}^S$;

clearly M is only a subring of $F_{5\times3}^S$; and is not an ideal of $F_{5\times3}^S$.

Further if

$y = \begin{bmatrix} a & b & 0 & 0 & p \\ 0 & 0 & e & f & 0 \\ 0 & 0 & g & h & 0 \end{bmatrix} \in F_{5\times3}^S$. We see $y \times_n m = (0)$

for every $m \in M$.

*Example 5.19:* Let

$$F_{7\times 3}^{S} = \left\{ \begin{bmatrix} a_1 & a_2 & a_3 \\ a_4 & a_5 & a_6 \\ \hline a_7 & a_8 & a_9 \\ a_{10} & a_{11} & a_{12} \\ a_{13} & a_{14} & a_{15} \\ \hline a_{16} & a_{17} & a_{18} \\ a_{19} & a_{20} & a_{21} \end{bmatrix} \middle| \ a_i \in Q, 1 \leq i \leq 21 \right\}$$

be a $7 \times 3$ super matrix semigroup under natural product.

Consider

$$P = \left\{ \begin{bmatrix} 0 & 0 & 0 \\ 0 & 0 & 0 \\ \hline a_1 & a_2 & a_3 \\ 0 & 0 & 0 \\ 0 & 0 & 0 \\ \hline 0 & 0 & 0 \\ a_4 & a_5 & a_6 \end{bmatrix} \middle| \ a_i \in Z, 1 \leq i \leq 6 \right\} \subseteq F_{7\times 3}^{S};$$

P is a subsemigroup under natural product, how ever P is not an ideal of $F_{7\times 3}^{S}$.

Now take

$$x = \begin{bmatrix} a_1 & a_2 & a_3 \\ a_4 & a_5 & a_6 \\ \hline 0 & 0 & 0 \\ \hline a_7 & a_8 & a_9 \\ a_{10} & a_{11} & a_{12} \\ a_{13} & a_{14} & a_{15} \\ 0 & 0 & 0 \end{bmatrix} \in F_{7\times 3}^S.$$

Clearly $x \times_n p = (0)$ for every $p \in P$.

Thus we have a collection of zero divisors in the semigroup under natural product.

Now consider the set

$$T = \left\{ \begin{bmatrix} a_1 & a_2 & a_3 \\ 0 & 0 & 0 \\ \hline a_4 & a_5 & a_6 \\ a_7 & a_8 & a_9 \\ 0 & 0 & 0 \\ 0 & 0 & 0 \\ \hline a_{10} & a_{11} & a_{12} \end{bmatrix} \middle| a_i \in Q, 1 \le i \le 12 \right\} \subseteq F_{7\times 3}^S;$$

T is an ideal of $F_{7\times 3}^S$.

Further

$$m = \begin{bmatrix} 0 & 0 & 0 \\ a_1 & a_2 & a_3 \\ 0 & 0 & 0 \\ 0 & 0 & 0 \\ a_4 & a_5 & a_6 \\ a_7 & a_8 & a_9 \\ 0 & 0 & 0 \end{bmatrix} \in F_{7\times 3}^S$$

is such that $m \times_n t = (0)$ for every $t \in T$. Thus $F_{7\times 3}^S$ has several zero divisors and has ideals.

*Example 5.20:* Let

$$F_{3\times 7}^S = \left\{ \begin{bmatrix} a_1 & a_2 & a_3 & a_4 & a_5 & a_6 & a_7 \\ a_8 & a_9 & a_{10} & a_{11} & a_{12} & a_{13} & a_{14} \\ a_{15} & a_{16} & a_{17} & a_{18} & a_{19} & a_{20} & a_{21} \end{bmatrix} \right.$$
$$\left. a_i \in Q, 1 \leq i \leq 21 \right\}$$

be $3 \times 7$ matrix semigroup under natural product.

Take

$$P = \left\{ \begin{bmatrix} 0 & a_1 & 0 & a_4 & 0 & a_7 & 0 \\ 0 & a_2 & 0 & a_5 & 0 & a_8 & 0 \\ 0 & a_3 & 0 & a_6 & 0 & a_9 & 0 \end{bmatrix} \; \middle| \; a_i \in Q, 1 \leq i \leq 9 \right\} \subseteq F_{3\times 7}^S;$$

P is an ideal of $F_{3\times 7}^S$.

$$x = \begin{bmatrix} a_1 & 0 & a_4 & 0 & a_7 & 0 & a_{10} \\ a_2 & 0 & a_5 & 0 & a_8 & 0 & a_{11} \\ a_3 & 0 & a_6 & 0 & a_9 & 0 & a_{12} \end{bmatrix} \in F_{3\times 7}^S$$

is such that $x \times_n p = (0)$ for every p in P.

Thus $F^S_{3\times 7}$ has several zero divisors.

Take

$$Y = \left\{ \begin{bmatrix} a_1 & a_2 & a_3 & a_4 & a_5 & a_6 & a_7 \\ a_8 & a_9 & a_{10} & a_{11} & a_{12} & a_{13} & a_{14} \\ a_{15} & a_{16} & a_{17} & a_{18} & a_{19} & a_{20} & a_{21} \end{bmatrix} \middle| a_i \in Q, 1 \le i \le 21 \right\}$$

$\subseteq F^S_{3\times 7}$, Y is only a subsemigroup and not an ideal of $F^S_{3\times 7}$.

*Example 5.21:* Let

$$M = \left\{ \begin{bmatrix} a_1 & a_2 & a_3 & a_4 \\ a_5 & a_6 & a_7 & a_8 \\ a_9 & a_{10} & a_{11} & a_{12} \\ a_{13} & a_{14} & a_{15} & a_{16} \end{bmatrix} \middle| a_i \in Q, 1 \le i \le 16 \right\}$$

be a semigroup under natural product $\times_n$.

Consider

$$P = \left\{ \begin{bmatrix} 0 & a_1 & a_2 & 0 \\ a_3 & 0 & 0 & a_5 \\ a_4 & 0 & 0 & a_6 \\ 0 & a_7 & a_8 & 0 \end{bmatrix} \middle| a_i \in Z, 1 \le i \le 8 \right\} \subseteq F^S_{4\times 4} = M;$$

P is a subsemigroup under $\times_n$. However P is not an ideal of $F^S_{4\times 4}$.

Let

$$X = \left\{ \left[ \begin{array}{c|cc|c} a_1 & 0 & 0 & a_2 \\ \hline 0 & a_5 & a_6 & 0 \\ 0 & a_7 & a_8 & 0 \\ \hline a_3 & 0 & 0 & a_4 \end{array} \right] \,\bigg|\, a_i \in Q,\ 1 \le i \le 8 \right\} \subseteq F_{4\times 4}^S ;$$

X is an ideal of $F_{4\times 4}^S$. Further every $x \in P$ and $m \in X$. $x \times_n m = (0)$. Thus $F_{4\times 4}^S$ has zero divisors and subsemigroups which are not ideals.

Now consider another example.

*Example 5.22*: Let

$$F_{3\times 3}^S = \left\{ \left[ \begin{array}{ccc} a_1 & a_2 & a_3 \\ \hline a_4 & a_5 & a_6 \\ a_7 & a_8 & a_9 \end{array} \right] \,\bigg|\, a_i \in Z,\ 1 \le i \le 9 \right\}$$

be a $3 \times 3$ super matrix semigroup under natural product. It is important to observe $F_{3\times 3}^S$ is not compatible with usual matrix product. Also no type of product on square super matrices can be defined on elements in $F_{3\times 3}^S$.

Take

$$X = \left\{ \left[ \begin{array}{ccc} 0 & 0 & 0 \\ \hline a_1 & a_2 & a_3 \\ 0 & 0 & 0 \end{array} \right] \,\bigg|\, a_i \in Z,\ 1 \le i \le 3 \right\} \subseteq F_{3\times 3}^S ,$$

X is a subsemigroup as well as an ideal of $F_{3\times 3}^S$.

Take

$$M = \left\{ \left[ \begin{array}{ccc} a_1 & a_2 & a_3 \\ \hline 0 & 0 & 0 \\ \hline a_4 & a_5 & a_6 \end{array} \right] \; \middle| \; a_i \in Z, \, 1 \leq i \leq 6 \right\} \subseteq F^S_{3\times 3}.$$

M is a subsemigroup as well as an ideal of $F^S_{3\times 3}$. We see for every $x \in X$ and $m \in M$, $x \times_n m = (0)$.

Now we describe the unit element of $F^S_C, F^S_R, F^S_{m\times n}$ (m ≠ n) and $F^S_{n\times n}$.

In $F^S_C$, $\begin{bmatrix} 1 \\ 1 \\ \hline 1 \\ 1 \\ 1 \\ \hline \vdots \\ \hline 1 \end{bmatrix}$ acts as the supercolumn unit under the natural product $\times_n$.

For $F^S_R$; (1 1 | 1 1 1 | 1 … | 1 1) acts as the super row unit element under the natural product $\times_n$.

For $F^S_{7\times 3}$; $\begin{bmatrix} 1 & 1 & 1 \\ 1 & 1 & 1 \\ \hline 1 & 1 & 1 \\ 1 & 1 & 1 \\ 1 & 1 & 1 \\ \hline 1 & 1 & 1 \\ 1 & 1 & 1 \end{bmatrix}$ acts as the super $7 \times 3$ unit under the natural product $\times_n$.

For $F^S_{4\times 4}$, $\begin{bmatrix} 1 & 1 & 1 & 1 \\ 1 & 1 & 1 & 1 \\ \hline 1 & 1 & 1 & 1 \\ 1 & 1 & 1 & 1 \end{bmatrix}$ acts as the $4 \times 3$ super unit under product.

Take $x = (1 \mid 1\ 1 \mid 1\ 1\ 1 \mid 1\ 1)\ (7 \mid 3\ 2 \mid 5\ 7\ -1) \mid 2\ 0)$
$= (7 \mid 3\ 2 \mid 5\ 7\ -1 \mid 2\ 0)$.

Likewise for $x = \begin{bmatrix} 3 \\ 2 \\ \hline -1 \\ 0 \\ 3 \\ \hline \frac{1}{7} \\ 0 \\ 2 \end{bmatrix}$, $\begin{bmatrix} 1 \\ 1 \\ \hline 1 \\ 1 \\ 1 \\ \hline 1 \\ 1 \\ 1 \end{bmatrix}$ act as the multiplicative super $8 \times 1$ identity for, $\begin{bmatrix} 3 \\ 2 \\ -1 \\ 0 \\ 3 \\ \frac{1}{7} \\ 0 \\ 2 \end{bmatrix} \times_n \begin{bmatrix} 1 \\ 1 \\ 1 \\ 1 \\ 1 \\ 1 \\ 1 \\ 1 \end{bmatrix} = \begin{bmatrix} 3 \\ 2 \\ -1 \\ 0 \\ 3 \\ \frac{1}{7} \\ 0 \\ 2 \end{bmatrix}$.

For $x = \left[\begin{array}{c|ccc|cccc|cc} 3 & 7 & 2 & 5 & 1 & 0 & -1 & 7 & 0 & 8 \\ \hline 0 & 1 & 2 & 3 & 4 & 5 & 6 & 7 & 8 & 9 \\ \hline 0 & 3 & 4 & 0 & 1 & 0 & 7 & 0 & 1 & 1 \\ \hline 4 & 0 & 2 & 1 & 0 & 2 & 0 & 4 & 0 & 0 \end{array}\right]$

$I = \left[\begin{array}{c|ccc|cccc|cc} 1 & 1 & 1 & 1 & 1 & 1 & 1 & 1 & 1 \\ \hline 1 & 1 & 1 & 1 & 1 & 1 & 1 & 1 & 1 \\ \hline 1 & 1 & 1 & 1 & 1 & 1 & 1 & 1 & 1 \\ \hline 1 & 1 & 1 & 1 & 1 & 1 & 1 & 1 & 1 \end{array}\right]$

acts as the super identity under $\times_n$. For $x \times_n I = I \times_n x = x$.

Consider

$I = \left[\begin{array}{c|cc|cc} 1 & 1 & 1 & 1 & 1 \\ \hline 1 & 1 & 1 & 1 & 1 \\ \hline 1 & 1 & 1 & 1 & 1 \\ \hline 1 & 1 & 1 & 1 & 1 \\ \hline 1 & 1 & 1 & 1 & 1 \end{array}\right]$ and $y = \left[\begin{array}{c|cc|cc} 0 & 2 & 3 & 0 & 1 \\ \hline 1 & 0 & 0 & 4 & 0 \\ \hline -1 & -2 & 6 & -1 & 3 \\ \hline 0 & 0 & 5 & 0 & 7 \\ \hline -1 & -4 & -1 & 2 & 0 \end{array}\right]$

are such that $I \times_n y = y \times_n I = y$.

Now having seen how the units look like we now proceed onto see how inverse of an element look under natural product $\times_n$.

Let
$$x = \left[\begin{array}{cc|c} 7 & 3 & -1 \\ 1 & 2 & 9 \\ \hline 8 & 5 & 1 \\ 4 & 7 & 2 \end{array}\right],$$

if x takes its entries either from Q or from R and if no entry in x is zero then alone inverse exists otherwise inverse of x does not exist.

$$\text{Take } y = \left[\begin{array}{cc|c} 1/7 & 1/3 & -1 \\ 1 & 1/2 & 1/9 \\ \hline 1/8 & 1/5 & 1 \\ 1/4 & 1/7 & 1/2 \end{array}\right] \text{ then we see}$$

$$x \times_n y = \left[\begin{array}{cc|c} 1 & 1 & 1 \\ 1 & 1 & 1 \\ \hline 1 & 1 & 1 \\ 1 & 1 & 1 \end{array}\right].$$

Let
$$x = \left[\begin{array}{cc|c} 0 & 0 & 1 \\ 3 & 4 & 5 \\ \hline 8 & 9 & 1 \\ 1 & 0 & 1 \end{array}\right]$$

clearly for this x we do not have a y such that

$$x \times_n y = \begin{bmatrix} 1 & 1 & | & 1 \\ 1 & 1 & | & 1 \\ \hline 1 & 1 & | & 1 \\ 1 & 1 & | & 1 \end{bmatrix}.$$

Consider $x = (1/8 \mid 7\ 5 \mid 3\ 2\ 4\ -1)$ then the inverse for x is $y = (8 \mid 1/7\ 1/5 \mid 1/3\ 1/2\ 1/4\ -1)$ we $x \times_n y = (1 \mid 1\ 1 \mid 1\ 1\ 1\ 1)$.

Consider

$$x = \begin{bmatrix} 8 & 1 \\ \hline 3 & 5 \\ -1 & 1/7 \\ \hline -8 & 4 \\ 1 & -1 \\ 3 & -2 \end{bmatrix} \text{ then } y = \begin{bmatrix} 1/8 & 1 \\ \hline 1/3 & 1/5 \\ -1 & 7 \\ \hline -1/8 & 1/4 \\ 1 & -1 \\ 1/3 & -1/2 \end{bmatrix}$$

is such that

$$x \times_n y = \begin{bmatrix} 1 & 1 \\ \hline 1 & 1 \\ 1 & 1 \\ \hline 1 & 1 \\ 1 & 1 \\ 1 & 1 \end{bmatrix}.$$

Now having seen inverse and unit we just give the statement of a theorem, the proof is left as an exercise to the reader.

**THEOREM 5.9:** *Let $F_C^S$ (or $F_R^S$ or $F_{m \times n}^S$ ($m \neq n$) or $F_{n \times n}^S$) be the super matrix semigroup under natural product. No super matrix other than those super matrices with entries from {1, –1} have inverse if $F_C^S$ (or $F_R^S$ or $F_{m \times n}^S$ ($m \neq n$) or $F_{n \times n}^S$) take its entries from Z.*

**THEOREM 5.10:** *Let $F_C^S$ (or $F_R^S$ or $F_{m \times n}^S$ ($m \neq n$) or $F_{n \times n}^S$) be the super matrix semigroup under natural product, with entries from Q or R. Every super matrix M in which no element of M takes 0 has inverse.*

The proof of this theorem is also left as an exercise to the reader.

Consider $x = (1\ -1\ |\ 1\ 1\ -1\ |\ -1\ -1) \in F_R^S = \{(a_1\ a_2\ |\ a_3\ a_4\ a_5\ |\ a_6\ a_7)\ |\ a_i \in Z,\ 1 \leq i \leq 7\}$; clearly $x = (1\ -1\ |\ 1\ 1\ -1\ |\ -1\ -1)$ acts as its inverse that is $x \times_n x = (1\ 1\ |\ 1\ 1\ 1\ |\ 1\ 1)$.

Consider

$$y = \begin{bmatrix} -1 \\ \hline 1 \\ \hline -1 \\ \hline 1 \\ -1 \\ 1 \\ \hline -1 \\ 1 \\ -1 \end{bmatrix} \in F_C^S = \left\{ \begin{bmatrix} a_1 \\ \hline a_2 \\ \hline a_3 \\ \hline a_4 \\ a_5 \\ a_6 \\ \hline a_7 \\ a_8 \\ a_9 \end{bmatrix} \text{ where } a_i \in Z;\ 1 \leq i \leq 9 \right\}.$$

Now

$$y^2 = \begin{bmatrix} 1 \\ 1 \\ 1 \\ 1 \\ 1 \\ 1 \\ 1 \\ 1 \\ 1 \end{bmatrix};$$

all y whose entries are from Z \ {1, -1} does not have inverse under natural product.

Take

$$y = \begin{bmatrix} -1 & 1 & 1 & -1 & 1 \\ 1 & -1 & 1 & 1 & -1 \\ -1 & 1 & -1 & -1 & 1 \\ 1 & 1 & 1 & 1 & -1 \\ -1 & 1 & -1 & 1 & -1 \end{bmatrix}$$

we see $y^2 = \begin{bmatrix} 1 & 1 & 1 & 1 & 1 \\ 1 & 1 & 1 & 1 & 1 \\ 1 & 1 & 1 & 1 & 1 \\ 1 & 1 & 1 & 1 & 1 \\ 1 & 1 & 1 & 1 & 1 \end{bmatrix}$.

Consider $x = \begin{bmatrix} 0 & 1 & 2 \\ 3 & 0 & -3 \\ 4 & 1 & 2 \end{bmatrix}$ we see $x^{-1}$ does not exist.

Take x = (1 0 | 5 7 2 | 1 5 7 –1 2), $x^{-1}$ does not exist.

$$\text{Consider } y = \begin{bmatrix} -1 \\ 0 \\ 2 \\ \dfrac{3}{-5} \\ 7 \\ 0 \\ \overline{2} \\ \dfrac{1}{4} \end{bmatrix};$$

clearly $y^{-1}$ does not exist.

Consider
$$x = \begin{bmatrix} 7 & -1 & 0 & 2 & 3 & 4 & 0 & 3 & -1 & 0 \\ 0 & 2 & 1 & 0 & 7 & 0 & 1 & 0 & 5 & 1 \end{bmatrix}.$$

Clearly $x^{-1}$ does not exist.

Now we have seen inverse of a super matrix under natural product and the condition under which the inverse exists.

Now we proceed onto discuss the operation '+' on $F_C^S$ or $F_R^S$ or $F_{m \times n}^S$ ($m \neq n$) or $F_{n \times n}^S$, which is stated as the following theorems.

**THEOREM 5.11:** *( $F_C^S$, +) is an additive abelian group of super column matrices.*

**THEOREM 5.12:** *( $F_R^S$, +) is an additive abelian group of super row matrices.*

**THEOREM 5.13:** *( $F_{m \times n}^S$ ($m \neq n$), +) is an additive abelian group of super $m \times n$ ($m \neq n$) matrices.*

**THEOREM 5.14:** *( $F^S_{n \times n}$, +) is an additive abelian group of super square matrices.*

We can define subgroups. All subgroups are normal as these groups are abelian. We will just give some examples.

*Example 5.23:* Let
$$F^S_R = \{(a_1\ a_2\ a_3\ a_4\ |\ a_5\ a_6\ |\ a_7\ a_8\ |\ a_9)\ |\ a_i \in Q;\ 1 \leq i \leq 9\}$$
be an abelian group of super row matrices under addition.

*Example 5.24:* Let

$$F^S_{2\times 3} = \left\{ \begin{bmatrix} a_1 & a_2 & a_3 \\ \hline a_4 & a_5 & a_6 \end{bmatrix} \text{ where } a_i \in Q;\ 1 \leq i \leq 6 \right\}$$

be an additive abelian group of $2 \times 3$ super matrices.

*Example 5.25:* Let

$$F^S_C = \left\{ \begin{bmatrix} a_1 \\ a_2 \\ a_3 \\ a_4 \\ \hline a_5 \\ a_6 \\ a_7 \\ \hline a_8 \\ a_9 \\ \hline a_{10} \\ a_{11} \end{bmatrix} \ \middle|\ a_i \in Q,\ 1 \leq i \leq 11 \right\}$$

be an abelian group of column super matrices.

*Example 5.26:* Let

$$F^S_{4\times 4} = \left\{ \begin{bmatrix} a_1 & a_2 & a_3 & a_4 \\ a_5 & a_6 & a_7 & a_8 \\ a_9 & a_{10} & a_{11} & a_{12} \\ a_{13} & a_{14} & a_{15} & a_{16} \end{bmatrix} \;\middle|\; a_i \in R,\; 1 \le i \le 16 \right\}$$

be an additive abelian group of $4 \times 4$ super matrices.

Now we can define $\{F^S_R,\, +,\, \times_n\}$ as the ring of super row matrices, $(F^S_C,\, +,\, \times_n)$ as the ring of super column matrices, $\{F^S_{m\times n}\; (m \ne n),\, \times_n,\, +\}$ is the ring of super $m \times n$ matrices and $\{F^S_{n\times n},\, \times_n,\, +\}$ be the ring of super $n \times n$ matrices.

We describe properties associated with them.

*Example 5.27:* Let

$$F^S_C = \left\{ \begin{bmatrix} a_1 \\ a_2 \\ \hline a_3 \\ \hline a_4 \\ \hline a_5 \\ a_6 \\ a_7 \\ \hline a_8 \end{bmatrix} \;\middle|\; a_i \in Q,\; 1 \le i \le 8 \right\}$$

be the $8 \times 1$ super column matrix ring under '+' and '$\times_n$'.

*Example 5.28:* Let

$$F^S_R = \{(a_1 \mid a_2 \; a_3 \mid a_4 \; a_5) \text{ where } a_i \in Q,\; 1 \le i \le 5,\, +,\, \times_n\}$$

be the ring of super row matrices.

*Example 5.29:* Let

$$F^S_{3\times 4} = \left\{ \left[ \begin{array}{c|cc|c} a_1 & a_4 & a_7 & a_{10} \\ a_2 & a_5 & a_8 & a_{11} \\ a_3 & a_6 & a_9 & a_{12} \end{array} \right] \middle| a_i \in Q, 1 \leq i \leq 12, +, \times_n \right\}$$

be the ring of $3 \times 4$ supermatrices.

*Example 5.30:* Let

$$F^S_{3\times 4} = \left\{ \left[ \begin{array}{c|ccc} a_1 & a_2 & a_3 & a_4 \\ \hline a_5 & a_6 & a_7 & a_8 \\ \hline a_9 & a_{10} & a_{11} & a_{12} \\ a_{13} & a_{14} & a_{15} & a_{16} \end{array} \right] \middle| a_i \in Z, 1 \leq i \leq 4 \right\}$$

be the ring of square supermatrices.

*Example 5.31:* Let

$$F^S_{9\times 3} = \left\{ \left[ \begin{array}{ccc} a_1 & a_2 & a_3 \\ \hline a_4 & a_5 & a_6 \\ \hline a_7 & a_8 & a_9 \\ \hline a_{10} & a_{11} & a_{12} \\ a_{13} & a_{14} & a_{15} \\ a_{16} & a_{17} & a_{18} \\ \hline a_{19} & a_{20} & a_{21} \\ \hline a_{22} & a_{23} & a_{24} \\ a_{25} & a_{26} & a_{27} \end{array} \right] \middle| a_i \in Q, 1 \leq i \leq 27, +, \times_n \right\}$$

be a ring of column supervectors.

*Example 5.32:* Let

$$F^S_{3\times 10} = \left\{ \begin{bmatrix} a_1 & a_2 & a_3 & a_4 & a_5 & a_6 & a_7 & a_8 & a_9 & a_{10} \\ a_{11} & a_{12} & a_{13} & a_{14} & a_{15} & a_{16} & a_{17} & a_{18} & a_{19} & a_{20} \\ a_{21} & a_{22} & a_{23} & a_{24} & a_{25} & a_{26} & a_{27} & a_{28} & a_{29} & a_{30} \end{bmatrix} \right.$$

$$a_i \in Q,\ 1 \leq i \leq 30,\ +,\ \times_n \}$$

be a ring of super row vectors.

All these rings are commutative have zero divisor and have unit. However we will give examples of ring of super matrices which have no unit.

*Example 5.33:* Let

$$F^S_R = \{(x_1 \mid x_2\ x_3\ x_4 \mid x_5\ x_6 \mid x_7\ x_8\ x_9\ x_{10}) \mid x_i \in 3Z\ ;$$
$$1 \leq i \leq 10,\ +,\ \times_n \}$$

be the ring of super row matrices. Clearly $F^S_R$ does not contain the unit $(1 \mid 1\ 1\ 1 \mid 1\ 1 \mid 1\ 1\ 1\ 1)$.

*Example 5.34:* Let

$$M = \left\{ \begin{bmatrix} a_1 \\ a_2 \\ a_3 \\ a_4 \\ \hline a_5 \\ a_6 \\ \hline a_7 \\ a_8 \\ \hline a_9 \\ a_{10} \end{bmatrix} \middle|\ a_i \in 7Z,\ 1 \leq i \leq 10,\ +,\ \times_n \right\}$$

be the ring of super column matrices. Clearly this ring has no super identity.

*Example 5.35*: Let

$$F_{3\times 4}^S = \left\{ \left[ \begin{array}{ccc|c} a_1 & a_2 & a_3 & a_4 \\ a_5 & a_6 & a_7 & a_8 \\ \hline a_9 & a_{10} & a_{11} & a_{12} \\ a_{13} & a_{14} & a_{15} & a_{16} \end{array} \right] \middle| a_i \in 10Z; 1 \leq i \leq 16, +, \times_n \right\}$$

be a ring of $4 \times 4$ super matrices.

Clearly the super unit $\left[ \begin{array}{ccc|c} 1 & 1 & 1 & 1 \\ 1 & 1 & 1 & 1 \\ \hline 1 & 1 & 1 & 1 \\ 1 & 1 & 1 & 1 \end{array} \right] \notin M$.

*Example 5.36:* Let $F_{3\times 4}^S =$

$$\left\{ \left[ \begin{array}{ccc|ccc|cccc|cc} x_1 & x_2 & x_3 & x_4 & x_5 & x_6 & x_7 & x_8 & x_9 & x_{10} & x_{11} & x_{12} \\ x_{13} & x_{14} & x_{15} & x_{16} & x_{17} & x_{18} & x_{19} & x_{20} & x_{21} & x_{22} & x_{23} & x_{24} \\ x_{25} & x_{26} & x_{27} & x_{28} & x_{29} & x_{30} & x_{31} & x_{32} & x_{33} & x_{34} & x_{35} & x_{36} \end{array} \right] \right.$$

$$\left. x_i \in 5Z; 1 \leq i \leq 36, +, \times \right\}$$

be a ring of super row vector. Clearly P does not contain the super identity

$$I = \left[ \begin{array}{ccc|ccc|cccc|cc} 1 & 1 & 1 & 1 & 1 & 1 & 1 & 1 & 1 & 1 & 1 & 1 \\ 1 & 1 & 1 & 1 & 1 & 1 & 1 & 1 & 1 & 1 & 1 & 1 \\ 1 & 1 & 1 & 1 & 1 & 1 & 1 & 1 & 1 & 1 & 1 & 1 \end{array} \right].$$

*Example 5.37:* Let

$$V = \left\{ \begin{bmatrix} a_1 & a_2 & a_3 & a_4 \\ \hline a_5 & a_6 & a_7 & a_8 \\ \hline a_9 & a_{10} & a_{11} & a_{12} \\ a_{13} & a_{14} & a_{15} & a_{16} \\ a_{17} & a_{18} & a_{19} & a_{20} \\ \hline a_{21} & a_{22} & a_{23} & a_{24} \\ \hline a_{25} & a_{26} & a_{27} & a_{28} \\ \hline a_{29} & a_{30} & a_{31} & a_{32} \end{bmatrix} \middle| a_j \in 15Z; 1 \leq j \leq 32 \right\}$$

be a ring of super column vectors which has no unit.

Now we proceed onto study super matrix structure using $R^+ \cup \{0\}$ or $Q^+ \cup \{0\}$ or $Z^+ \cup \{0\}$.

Let $S_R^+ = \{(x_1 \ x_2 \ x_3 \mid x_4 \ \ldots \mid x_{n-1} \ x_n) \mid x_i \in R^+ \cup \{0\}$ or $Q^+ \cup \{0\}$ or $Z^+ \cup \{0\}\}$ denotes the collection of all super row matrices of same type from $R^+ \cup \{0\}$ or $Q^+ \cup \{0\}$ or $Z^+ \cup \{0\}$. This notation will be used throughout this book.

$$S_C^+ = \left\{ \begin{bmatrix} a_1 \\ \hline a_2 \\ \hline a_3 \\ \vdots \\ \hline a_m \end{bmatrix} \middle| a_j \in Z^+ \cup \{0\} \text{ or } R^+ \cup \{0\} \right.$$

$$\left. \text{or } Q^+ \cup \{0\}; 1 \leq i \leq m \right\}$$

denotes the collection of all column super matrices of same type with entries from $R^+ \cup \{0\}$ or $Q^+ \cup \{0\}$ or $Z^+ \cup \{0\}$.

$$S_{m \times n}^+ \ (m \neq n) = \left\{ \begin{bmatrix} a_{11} & a_{12} & \cdots & a_{1n} \\ a_{21} & a_{22} & \cdots & a_{2n} \\ \vdots & \vdots & & \vdots \\ a_{m1} & a_{m2} & \cdots & a_{mn} \end{bmatrix} \middle| a_{ij} \in Q^+ \cup \{0\} \right.$$

or $Z^+ \cup \{0\}$ or $R^+ \cup \{0\}$; $1 \leq i \leq m$, $1 \leq j \leq n\}$

denotes the collection of all m × n super matrices of same type with entries from $Q^+ \cup \{0\}$ or $Z^+ \cup \{0\}$ or $R^+ \cup \{0\}$.

$$S_{n \times n}^+ = \left\{ \begin{bmatrix} a_{11} & a_{12} & \cdots & a_{1n} \\ a_{21} & a_{22} & \cdots & a_{2n} \\ \vdots & \vdots & & \vdots \\ a_{n1} & a_{n2} & \cdots & a_{nn} \end{bmatrix} \middle| a_{ij} \in Z^+ \cup \{0\} \right.$$

or $Q^+ \cup \{0\}$ or $R^+ \cup \{0\}$; $1 \leq i, j \leq n\}$

denotes the collection of all n × n super matrices of same type with entries from $R^+ \cup \{0\}$ or $Q^+ \cup \{0\}$ or $Z^+ \cup \{0\}$.

We will first illustrate these situations by some examples.

*Example 5.38:* Let

$S_R^+ = \{(x_1 \ x_2 \mid x_3 \ x_4 \ x_5 \mid x_6 \ x_7 \mid x_8) \mid x_i \in Q^+ \cup \{0\}, 1 \leq i \leq 8\}$

be the super row matrices of same type with entries from $Q^+ \cup \{0\}$.

*Example 5.39:* Let $S_R^+ =$

$$\left\{ \begin{bmatrix} a_1 & a_2 & a_3 & a_4 & a_5 & a_6 & a_7 & a_8 & a_9 \\ a_{10} & a_{11} & a_{12} & a_{13} & a_{14} & a_{15} & a_{16} & a_{17} & a_{18} \\ a_{19} & a_{20} & a_{21} & a_{22} & a_{23} & a_{25} & a_{25} & a_{26} & a_{27} \end{bmatrix} \right.$$

$$a_i \in Z^+ \cup \{0\}; 1 \leq i \leq 27\}$$

be the set of all super row vectors of same type with entries from $Z^+ \cup \{0\}$.

*Example 5.40:* Let

$$S_C^+ = \left\{ \begin{bmatrix} a_1 \\ a_2 \\ a_3 \\ a_4 \\ a_5 \\ a_6 \\ a_7 \\ a_8 \\ a_9 \\ a_{10} \\ a_{11} \end{bmatrix} \middle| a_i \in Q^+ \cup \{0\}; 1 \leq i \leq 11 \right\}$$

denote the collection of super column matrices of same type with entries from $Q^+ \cup \{0\}$.

*Example 5.41:* Let

$$S_C^+ = \left\{ \begin{bmatrix} a_1 & a_2 & a_3 \\ a_4 & a_5 & a_6 \\ \hline a_7 & a_8 & a_9 \\ a_{10} & a_{11} & a_{12} \\ a_{13} & a_{14} & a_{15} \\ a_{16} & a_{17} & a_{18} \\ \hline a_{19} & a_{20} & a_{21} \\ a_{22} & a_{23} & a_{24} \\ a_{25} & a_{26} & a_{27} \\ \hline a_{28} & a_{29} & a_{30} \\ a_{31} & a_{32} & a_{33} \\ \hline a_{34} & a_{35} & a_{36} \end{bmatrix} \middle| a_i \in R^+ \cup \{0\}; 1 \le i \le 36 \right\}$$

denote the collection of all super column vectors of same type with entries from $R^+ \cup \{0\}$.

*Example 5.42:* Let

$$S_{3 \times 4}^+ = \left\{ \begin{bmatrix} a_1 & a_2 & a_3 & a_4 \\ a_5 & a_6 & a_7 & a_8 \\ a_9 & a_{10} & a_{11} & a_{12} \end{bmatrix} \middle| a_i \in R^+ \cup \{0\}; 1 \le i \le 12 \right\}$$

be the collection of all $3 \times 4$ super matrices of same type with entries from $R^+ \cup \{0\}$.

*Example 5.43:* Let

$$S_{5\times 5}^{+} = \left\{ \begin{bmatrix} a_1 & a_2 & a_3 & a_4 & a_5 \\ a_6 & a_7 & a_8 & a_9 & a_{10} \\ a_{11} & a_{12} & a_{13} & a_{14} & a_{15} \\ a_{16} & a_{17} & a_{18} & a_{19} & a_{20} \\ a_{21} & a_{22} & a_{23} & a_{24} & a_{25} \end{bmatrix} \middle| a_i \in R^{+} \cup \{0\}; 1 \le i \le 25 \right\}$$

be the collection of $5 \times 5$ super matrices of same type with entries from $R^{+} \cup \{0\}$.

Now we proceed onto give all possible algebraic structures on $F_C^{+}$, $F_R^{+}$ or $F_{m\times n}^{+}$ ($m \ne n$) and $F_{n\times n}^{+}$.

Consider $F_C^{+}$ the collection of all super column matrices of same type with entries from $Q^{+} \cup \{0\}$ or $R^{+} \cup \{0\}$ or $Z^{+} \cup \{0\}$; $S_C^{+}$ is a semigroup under '+' usual addition. Infact it has the additive identity and $(S_C^{+}, +)$ is a commutative semigroup. Likewise $F_R^{+}$ or $F_{m\times n}^{+}$ ($m \ne n$) and $F_{n\times n}^{+}$ are all abelian semigroups with respect to addition. Infact all of them are monoids.

We will illustrate this by some examples.

*Example 5.44:* Let

$$S_C^+ = \left\{ \begin{bmatrix} a_1 \\ a_2 \\ \hline a_3 \\ a_4 \\ \hline a_5 \\ \hline a_6 \\ \hline a_7 \\ a_8 \end{bmatrix} \;\middle|\; a_i \in Z^+ \cup \{0\};\, 1 \leq i \leq 8 \right\}$$

be a commutative semigroup of super column matrices with entries from $Z^+ \cup \{0\}$.

*Example 5.45:* Let

$$S_{m \times n}^+ = \left\{ \begin{bmatrix} a_1 & a_8 & a_{15} & a_{22} \\ a_2 & a_9 & a_{16} & a_{23} \\ a_3 & a_{10} & a_{17} & a_{24} \\ a_4 & a_{11} & a_{18} & a_{25} \\ a_5 & a_{12} & a_{19} & a_{26} \\ a_6 & a_{13} & a_{20} & a_{27} \\ \hline a_7 & a_{14} & a_{21} & a_{28} \end{bmatrix} \;\middle|\; a_i \in Q^+ \cup \{0\};\, 1 \leq i \leq 28 \right\}$$

be the semigroup of super $7 \times 4$ matrices under addition with entries from $Q^+ \cup \{0\}$.

*Example 5.46:* Let

$$S_{4\times 4}^+ = \left\{ \begin{bmatrix} a_1 & a_2 & a_3 & a_4 \\ a_5 & a_6 & a_7 & a_8 \\ a_9 & a_{10} & a_{11} & a_{12} \\ a_{13} & a_{14} & a_{15} & a_{16} \end{bmatrix} \middle| a_i \in R^+ \cup \{0\}; 1 \le i \le 16 \right\}$$

be the semigroup of super $4 \times 4$ super matrices with entries from $R^+ \cup \{0\}$.

*Example 5.47:* Let

$$S_{3\times 9}^+ = \left\{ \begin{bmatrix} a_1 & a_4 & a_7 & a_{10} & a_{13} & a_{16} & a_{19} & a_{22} & a_{25} \\ a_2 & a_5 & a_8 & a_{11} & a_{14} & a_{17} & a_{20} & a_{23} & a_{26} \\ a_3 & a_6 & a_9 & a_{12} & a_{15} & a_{18} & a_{21} & a_{24} & a_{27} \end{bmatrix} \right.$$

$$\left. a_i \in Q^+ \cup \{0\}; 1 \le i \le 27 \right\}$$

be the semigroup of super row vector under addition with elements from $Q^+ \cup \{0\}$.

*Example 5.48:* Let

$$S_{10\times 9}^+ = \left\{ \begin{bmatrix} a_1 & a_2 & a_3 & a_4 \\ a_5 & a_6 & a_7 & a_8 \\ a_9 & a_{10} & a_{11} & a_{12} \\ a_{13} & a_{14} & a_{15} & a_{16} \\ a_{17} & a_{18} & a_{19} & a_{20} \\ a_{21} & a_{22} & a_{23} & a_{24} \\ a_{25} & a_{26} & a_{27} & a_{28} \\ a_{29} & a_{30} & a_{31} & a_{32} \\ a_{33} & a_{34} & a_{35} & a_{36} \\ a_{37} & a_{38} & a_{39} & a_{40} \end{bmatrix} \middle| a_i \in Q^+ \cup \{0\}; 1 \le i \le 40 \right\}$$

be the semigroup of super column vector under addition.

Now we proceed onto define natural product on $S_C^+$, $S_R^+, S_{n\times n}^+$ and $S_{m\times n}^+$ $(m \neq n)$. Clearly $S_C^+$, $S_R^+, S_{n\times m}^+$ $(m \neq n)$ and $S_{m\times m}^+$ are semigroups under the natural product $\times_n$.

Depending on the set from which they take their entries they will be semigroups with multiplicative identity or otherwise.

We will illustrate this situation by some examples.

***Example 5.49:*** Let $S_R^+ = \{(a_1 \mid a_2 \ a_3 \mid a_4 \ a_5 \ a_6 \mid a_7 \ a_8 \ a_9 \ a_{10} \mid a_{11} \mid a_{12}) \mid a_i \in Z^+ \cup \{0\}, 1 \leq i \leq 12\}$ be a semigroup of super row matrices under the natural product $\times_n$.

***Example 5.50:*** Let

$$S_C^+ = \left\{ \begin{bmatrix} x_1 \\ \hline x_2 \\ \hline x_3 \\ x_4 \\ \hline x_5 \\ x_6 \\ \hline x_7 \\ \hline x_8 \\ x_9 \\ \hline x_{10} \\ \hline x_{11} \end{bmatrix} \middle| \ x_i \in Z^+ \cup \{0\}, 1 \leq i \leq 11 \right\}$$

be a semigroup of super column matrices under the natural product $\times_n$. Infact $S_C^+$ has units and zero divisors.

*Example 5.51:* Let

$$S_C^+ = \left\{ \begin{bmatrix} a_1 & a_2 & a_3 & a_4 & a_5 & a_6 & a_7 & a_8 & a_9 & a_{10} \\ a_{11} & a_{12} & a_{13} & a_{14} & a_{15} & a_{16} & a_{17} & a_{18} & a_{19} & a_{20} \end{bmatrix} \right.$$

$$\left. a_i \in R^+ \cup \{0\}, 1 \leq i \leq 20 \right\}$$

be the semigroup of super row vector under the natural product $\times_n$. $S_{2\times 10}^+$ has identity elements, units and zero divisors.

*Example 5.52:* Let

$$S_{8\times 3}^+ = \left\{ \begin{bmatrix} a_1 & a_2 & a_3 \\ a_4 & a_5 & a_6 \\ a_7 & a_8 & a_9 \\ a_{10} & a_{11} & a_{12} \\ a_{13} & a_{14} & a_{15} \\ a_{16} & a_{17} & a_{18} \\ a_{19} & a_{20} & a_{21} \\ a_{22} & a_{23} & a_{24} \end{bmatrix} \middle| a_i \in Z^+ \cup \{0\}, 1 \leq i \leq 24 \right\}$$

be a semigroup of super column vectors under natural product $\times_n$.

Clearly $S_{8\times 3}^+$ has identity $I = \begin{bmatrix} 1 & 1 & 1 \\ 1 & 1 & 1 \\ 1 & 1 & 1 \\ 1 & 1 & 1 \\ 1 & 1 & 1 \\ 1 & 1 & 1 \\ 1 & 1 & 1 \\ 1 & 1 & 1 \end{bmatrix}$ but no element in $S_{8\times 3}^+$ has units but has several zero divisors.

*Example 5.53:* Let

$$S_{2\times 2}^+ = \left\{ \begin{bmatrix} a_1 & a_2 \\ \hline a_3 & a_4 \end{bmatrix} \middle| a_i \in Z^+ \cup \{0\}, 1 \leq i \leq 2 \right\}$$

be the semigroup of super square matrices under natural product $\times_n$. Clearly $\begin{bmatrix} 1 & 1 \\ 1 & 1 \end{bmatrix}$ is the identity element of $S_{2\times 2}^+$, but has no units, that is no element in $S_{2\times 2}^+$ has inverse. Further $S_{2\times 2}^+$ has zero divisors. For take $x = \begin{bmatrix} 0 & 0 \\ \hline a_1 & a_2 \end{bmatrix}$ in $S_{2\times 2}^+$ then $y_1 = \begin{bmatrix} a_1 & a_2 \\ \hline 0 & 0 \end{bmatrix}$ and $y_3 = \begin{bmatrix} 0 & a_1 \\ \hline 0 & 0 \end{bmatrix}$ are all zero divisors in $S_{2\times 2}^+$.

*Example 5.54:* Let

$$S_{8\times 4}^+ = \left\{ \begin{bmatrix} a_1 & a_2 & a_3 & a_4 \\ a_5 & a_6 & a_7 & a_8 \\ \hline a_9 & a_{10} & a_{11} & a_{12} \\ a_{13} & a_{14} & a_{15} & a_{16} \\ \hline a_{17} & a_{18} & a_{19} & a_{20} \\ \hline a_{21} & a_{22} & a_{23} & a_{24} \\ a_{25} & a_{26} & a_{27} & a_{28} \\ a_{29} & a_{30} & a_{31} & a_{32} \end{bmatrix} \middle| a_i \in Q^+ \cup \{0\}, 1 \leq i \leq 32 \right\}$$

be the semigroup of $8 \times 4$ super matrices with entries from $Q^+ \cup \{0\}$.

$S_{8\times 4}^+$ is a semigroup with unit $I = \begin{bmatrix} 1 & 1 & 1 & 1 \\ 1 & 1 & 1 & 1 \\ \hline 1 & 1 & 1 & 1 \\ 1 & 1 & 1 & 1 \\ 1 & 1 & 1 & 1 \\ \hline 1 & 1 & 1 & 1 \\ 1 & 1 & 1 & 1 \\ \hline 1 & 1 & 1 & 1 \end{bmatrix}$ and has zero divisors and inverses. Now we can find ideals, subsemigroups, zero divisors and units in semigroups under the natural product $\times_n$. These will be only illustrated by some examples.

*Example 5.55:* Let

$$S_C^+ = \left\{ \begin{bmatrix} a_1 & a_2 & a_3 \\ \hline a_4 & a_5 & a_6 \\ a_7 & a_8 & a_9 \\ a_{10} & a_{11} & a_{12} \\ a_{13} & a_{14} & a_{15} \\ a_{16} & a_{17} & a_{18} \\ \hline a_{19} & a_{20} & a_{21} \\ a_{22} & a_{23} & a_{24} \\ a_{25} & a_{26} & a_{27} \\ a_{28} & a_{29} & a_{30} \end{bmatrix} \middle| \; a_i \in Z^+ \cup \{0\}, 1 \leq i \leq 30 \right\}$$

be the semigroup of super column vectors.

This has units and zero divisors.

Take P = $\left\{ \begin{bmatrix} \begin{array}{|ccc|} \hline a_1 & a_2 & a_3 \\ a_4 & a_5 & a_6 \\ 0 & 0 & 0 \\ \hline 0 & 0 & 0 \\ 0 & 0 & 0 \\ 0 & 0 & 0 \\ \hline a_7 & a_8 & a_9 \\ a_{10} & a_{11} & a_{12} \\ 0 & 0 & 0 \\ 0 & 0 & 0 \\ \hline \end{array} \end{bmatrix} \;\middle|\; a_i \in Z^+ \cup \{0\}, 1 \leq i \leq 12 \right\} \subseteq S_C^+$

is an ideal of $S_C^+$ under natural product $\times_n$.

Now

M = $\left\{ \begin{bmatrix} \begin{array}{|ccc|} \hline 0 & 0 & 0 \\ a_1 & a_2 & a_3 \\ a_4 & a_5 & a_6 \\ 0 & 0 & 0 \\ \hline 0 & 0 & 0 \\ 0 & 0 & 0 \\ 0 & 0 & 0 \\ \hline a_7 & a_8 & a_9 \\ 0 & 0 & 0 \\ 0 & 0 & 0 \\ 0 & 0 & 0 \\ \hline \end{array} \end{bmatrix} \;\middle|\; a_i \in Z^+ \cup \{0\}, 1 \leq i \leq 9 \right\} \subseteq S_C^+$

is an ideal of $S_C^+$ under natural product.

Take x = $\begin{bmatrix} 0 & 0 & 0 \\ \hline 0 & 0 & 0 \\ 0 & 0 & 0 \\ \hline a_1 & a_2 & a_3 \\ a_4 & a_5 & a_6 \\ a_7 & a_8 & a_9 \\ \hline 0 & 0 & 0 \\ 0 & 0 & 0 \\ 0 & 0 & 0 \\ 0 & 0 & 0 \end{bmatrix}$ and y = $\begin{bmatrix} a_1 & a_2 & a_3 \\ a_4 & a_5 & a_6 \\ a_7 & a_8 & a_9 \\ \hline 0 & 0 & 0 \\ 0 & 0 & 0 \\ 0 & 0 & 0 \\ \hline a_{10} & a_{11} & a_{12} \\ 0 & 0 & 0 \\ 0 & 0 & 0 \\ 0 & 0 & 0 \end{bmatrix}$ in $S_C^+$,

we see x $\times_n$ y = $\begin{bmatrix} 0 & 0 & 0 \\ 0 & 0 & 0 \\ 0 & 0 & 0 \\ \hline 0 & 0 & 0 \\ 0 & 0 & 0 \\ 0 & 0 & 0 \\ \hline 0 & 0 & 0 \\ 0 & 0 & 0 \\ 0 & 0 & 0 \\ 0 & 0 & 0 \end{bmatrix}$. No element in $S_C^+$ has inverse.

***Example 5.56:*** Let $S_R^+$ =

$\left\{ \begin{bmatrix} a_1 & a_2 & a_3 & a_4 & a_5 & a_6 & a_7 & a_8 \\ a_9 & a_{10} & a_{11} & a_{12} & a_{13} & a_{14} & a_{15} & a_{16} \\ a_{17} & a_{18} & a_{19} & a_{20} & a_{21} & a_{22} & a_{23} & a_{24} \\ a_{25} & a_{26} & a_{27} & a_{28} & a_{29} & a_{30} & a_{31} & a_{32} \end{bmatrix} \middle| a_i \in Q^+ \cup \{0\}, \right.$

$1 \le i \le 32\}$

be the semigroup of super row vectors under the natural product $\times_n$.

Take $I = \begin{bmatrix} 1 & 1 & 1 & 1 & 1 & 1 & 1 & 1 \\ 1 & 1 & 1 & 1 & 1 & 1 & 1 & 1 \\ 1 & 1 & 1 & 1 & 1 & 1 & 1 & 1 \\ 1 & 1 & 1 & 1 & 1 & 1 & 1 & 1 \end{bmatrix}$ be the unit in $S_R^+$.

Consider $X = \left\{ \begin{bmatrix} a_1 & 0 & 0 & a_5 & a_9 & a_{13} & 0 & 0 \\ a_2 & 0 & 0 & a_6 & a_{10} & a_{14} & 0 & 0 \\ a_3 & 0 & 0 & a_7 & a_{11} & a_{15} & 0 & 0 \\ a_4 & 0 & 0 & a_8 & a_{12} & a_{16} & 0 & 0 \end{bmatrix} \right.$

$a_i \in 5Z^+ \cup \{0\}, 1 \leq i \leq 16 \} \subseteq S_R^+$;

$X$ is only a subsemigroup under $\times_n$. $X$ has no identity. Further $X$ is not an ideal of $S_R^+$. However $X$ has zero divisors.

Take

$Y = \left\{ \begin{bmatrix} 0 & a_1 & a_2 & 0 & 0 & 0 & a_9 & a_{10} \\ 0 & a_3 & a_4 & 0 & 0 & 0 & a_{11} & a_{12} \\ 0 & a_5 & a_6 & 0 & 0 & 0 & a_{13} & a_{14} \\ 0 & a_7 & a_8 & 0 & 0 & 0 & a_{15} & a_{16} \end{bmatrix} \right.$ $a_i \in Q^+ \cup \{0\},$

$1 \leq i \leq 16 \} \subseteq S_R^+$;

$Y$ is an ideal of $S_R^+$. It is still interesting to note that every element $x$ in $X$ is such that $x \times_n y = (0)$ for every $y \in Y$. However $X + Y \neq S_R^+$.

*Example 5.57:* Let

$$S_{4\times 4}^+ = \left\{ \begin{bmatrix} a_1 & a_2 & a_3 & a_4 & a_5 \\ a_6 & a_7 & a_8 & a_9 & a_{10} \\ a_{11} & a_{12} & a_{13} & a_{14} & a_{15} \\ a_{16} & a_{17} & a_{18} & a_{19} & a_{20} \\ a_{21} & a_{22} & a_{23} & a_{24} & a_{25} \end{bmatrix} \middle| a_i \in Q^+ \cup \{0\},\ 1 \le i \le 32 \right\}$$

be the semigroup of super square matrices under the natural product $\times_n$.

$$I = \begin{bmatrix} 1 & 1 & 1 & 1 & 1 \\ 1 & 1 & 1 & 1 & 1 \\ 1 & 1 & 1 & 1 & 1 \\ 1 & 1 & 1 & 1 & 1 \\ 1 & 1 & 1 & 1 & 1 \end{bmatrix}$$

acts as the identity with respect to the natural product $\times_n$.

Consider

$$P = \left\{ \begin{bmatrix} a_1 & 0 & 0 & a_2 & a_3 \\ 0 & a_4 & a_5 & 0 & 0 \\ 0 & a_6 & a_7 & 0 & 0 \\ a_8 & 0 & 0 & a_{10} & a_{11} \\ a_9 & 0 & 0 & a_{12} & a_{13} \end{bmatrix} \middle| a_i \in Z^+ \cup \{0\}, \right.$$

$$\left. 1 \le i \le 13 \right\} \subseteq S_{4\times 4}^+ .$$

P is only a subsemigroup and not an ideal of $S_{4\times 4}^+$.

Consider

$$M = \left\{ \begin{bmatrix} 0 & a_1 & a_2 & 0 & 0 \\ a_3 & 0 & 0 & a_5 & a_6 \\ a_4 & 0 & 0 & a_7 & a_8 \\ 0 & a_9 & a_{10} & a_{10} & a_{11} \\ 0 & a_{11} & a_{12} & a_{12} & a_{13} \end{bmatrix} \middle| \; a_i \in Q^+ \cup \{0\}, \right.$$

$$\left. 1 \leq i \leq 12 \right\} \subseteq S_{4 \times 4}^+ ;$$

M is an ideal of $S_{4\times 4}^+$. Every element p in P is such that $p \times_n m = (0)$ for every $m \in M$.

Inview of this we have the following theorem.

**THEOREM 5.15:** *Let $S_C^+$ (or $S_R^+$ or $S_{m\times n}^+$ ($m \neq n$) or $S_{n\times n}^+$) be a semigroup under the natural product. Every ideal I in $S_C^+$ (or $S_R^+$ or $S_{n\times m}^+$ ($m \neq n$) or $S_{n\times n}^+$) is a subsemigroup of $S_C^+$ (or $S_R^+$ or $S_{n\times m}^+$ ($m \neq n$) or $S_{n\times n}^+$) but however every subsemigroup of $S_C^+$ (or $S_R^+$ or $S_{n\times m}^+$ ($m \neq n$) or $S_{n\times n}^+$) need not in general be an ideal of $S_C^+$ (or $S_R^+$ or $S_{n\times m}^+$ ($m \neq n$) or $S_{n\times m}^+$).*

The proof is simple and direct hence left as an exercise to the reader.

Now having seen all these we now proceed onto give two binary operations on $S_R^+$ (or $S_C^+$ or $S_{m\times n}^+$ ($m \neq n$) or $S_{n\times n}^+$) so $S_R^+$ (or $S_C^+$ or $S_{m\times n}^+$ ($m \neq n$) or $S_{n\times n}^+$), so that $S_R^+$ is the semiring with respect to addition and natural product.

Consider $(S_C^+, +, \times_n) = P_C^+$, it is easily verified $P_C^+$ is a semiring which is a strict semiring of super column matrices and is not a semifield as it has zero divisors under the product $\times_n$.

Likewise $P_R^+ = \{S_R^+, +, \times_n\}$ is a semiring of super row matrices which is not a semifield, infact a strict commutative semiring.

$P_{m \times n}^+$ (m ≠ n) = $\{S_{m \times n}^+$ (m ≠ n), +, $\times_n\}$ is a strict commutative semiring of m × n super matrices. Finally $P_{n \times n}^+ = \{S_{n \times n}^+, +, \times_n\}$ is a strict commutative semiring of super square matrices.

Now throughout this book $P_C^+$ will denote the semiring of super column matrices, $P_R^+$ will denote the semiring of super row matrices, $P_{m \times n}^+$ (m ≠ n) will denote the semiring of m × n super matrices and $P_{n \times n}^+$ will denote the semiring of square super matrices.

Now having seen the notation we proceed onto give examples of them.

*Example 5.58:* Let

$$P_C^+ = \left\{ \begin{bmatrix} a_1 \\ a_2 \\ a_3 \\ \hline a_4 \\ a_5 \\ \hline a_6 \\ \hline a_7 \\ a_8 \\ \hline a_9 \\ a_{10} \end{bmatrix} \middle| a_i \in Q^+ \cup \{0\}, 1 \leq i \leq 10, +, \times_n \right\}$$

be the semiring of super column matrices; $P_C^+$ is not a semifield as it has zero divisors.

Further $P_C^+$ has subsemirings for take

$$M = \left\{ \begin{bmatrix} a_1 \\ a_2 \\ a_3 \\ a_4 \\ \hline 0 \\ 0 \\ 0 \\ \hline a_5 \\ a_6 \\ \hline 0 \end{bmatrix} \middle| a_i \in Z^+ \cup \{0\}, 1 \leq i \leq 6, +, \times_n \right\} \subseteq P_C^+$$

is a subsemiring of $P_C^+$ and is not an ideal of $P_C^+$.

However $P_C^+$ has ideals for consider

$$N = \left\{ \begin{bmatrix} 0 \\ 0 \\ 0 \\ 0 \\ \hline a_1 \\ a_2 \\ a_3 \\ \hline 0 \\ 0 \\ \hline a_4 \end{bmatrix} \middle| a_i \in Q^+ \cup \{0\}, 1 \leq i \leq 4, +, \times_n \right\} \subseteq P_C^+$$

is an ideal of $P_C^+$.

We see for every $x \in M$ is such that

$$x \times_n y = \begin{bmatrix} 0 \\ 0 \\ 0 \\ 0 \\ \overline{0} \\ 0 \\ 0 \\ \overline{0} \\ 0 \\ \overline{0} \end{bmatrix} \text{ for every } y \in N.$$

Thus $P_C^+$ has infinitely many zero divisors so is not a semifield.

Finally $\begin{bmatrix} 1 \\ 1 \\ 1 \\ \overline{1} \\ 1 \\ \overline{1} \\ 1 \\ \overline{1} \end{bmatrix}$ acts as the identity with respect to $\times_n$.

***Example 5.59:*** Let $P_R^+ = \{(a_1 \mid a_2\ a_3 \mid a_4\ a_5\ a_6 \mid a_7\ a_8\ a_9 \mid a_{10}) \mid a_i \in 3Z^+ \cup \{0\}, 1 \leq i \leq 10, +, \times_n\}$ be a strict semiring of super row matrices. Clearly $(1 \mid 1\ 1 \mid 1\ 1\ 1 \mid 1\ 1\ 1 \mid 1) \notin P_R^+$. We see $P_R^+$ has zero divisors for take x = (4 | 0 3 | 2 0 1 | 0 3 9 | 0) and y = (0 | 7 0 | 0 8 0 | 9 0 0 | 1 0) in $P_R^+$. Clearly x $\times_n$ y = (0 | 0 0 | 0 0 0 | 0 0 0 | 0). So $P_R^+$ is a semiring which is not a semifield and has no identity.

*Example 5.60:* Let

$$P_{3\times 4}^+ = \left\{ \begin{bmatrix} \begin{array}{ccc} a_1 & a_{12} & a_{23} \\ \hline a_2 & a_{13} & a_{24} \\ a_3 & a_{14} & a_{25} \\ a_4 & a_{15} & a_{26} \\ a_5 & a_{16} & a_{27} \\ \hline a_6 & a_{17} & a_{28} \\ a_7 & a_{18} & a_{29} \\ a_8 & a_{19} & a_{30} \\ a_9 & a_{20} & a_{31} \\ \hline a_{10} & a_{21} & a_{32} \\ a_{11} & a_{22} & a_{33} \end{array} \end{bmatrix} \middle| a_i \in Q^+ \cup \{0\}, 1 \leq i \leq 33, +, \times_n \right\}$$

be the semiring of super column vectors. $P_{3\times 11}^+$ has zero divisors, units, ideals and subsemirings which are not ideals. However $P_{3\times 11}^+$ is not a semifield.

*Example 5.61:* Let $P_{12\times 2}^+$ =

$$\left\{ \begin{bmatrix} a_1 & a_2 & a_3 & a_4 & a_5 & a_6 & a_7 & a_8 & a_9 & a_{10} & a_{11} & a_{12} \\ a_{13} & a_{14} & a_{15} & a_{16} & a_{17} & a_{18} & a_{19} & a_{20} & a_{21} & a_{22} & a_{23} & a_{24} \end{bmatrix} \right.$$
$$\left. a_i \in Z^+ \cup \{0\}, 1 \leq i \leq 24, +, \times_n \right\}$$

be a semiring, $P_{12\times 2}^+$ has unit element, however no element in $P_{12\times 2}^+$ has inverse. Further $P_{12\times 2}^+$ has zero divisors so not a semifield. Has ideals. Thus $P^+_{12\times 2}$ is a super row vector semiring.

*Example 5.62:* Let

$$P^+_{4\times 4} = \left\{ \left[\begin{array}{c|ccc} a_1 & a_2 & a_3 & a_4 \\ a_5 & a_6 & a_7 & a_8 \\ a_9 & a_{10} & a_{11} & a_{12} \\ \hline a_{13} & a_{14} & a_{15} & a_{16} \end{array}\right] \;\middle|\; a_i \in R^+ \cup \{0\},\, 1 \le i \le 16,\, +,\, \times_n \right\}$$

be the semiring of square super matrices. $P^+_{4\times 4}$ has zero divisors units, subsemirings which are not ideals and ideals.

Clearly $\begin{bmatrix} 1 & 1 & 1 & 1 \\ 1 & 1 & 1 & 1 \\ 1 & 1 & 1 & 1 \\ 1 & 1 & 1 & 1 \end{bmatrix} = I$ is the unit of $P^+_{4\times 4}$.

$$S = \left\{ \left[\begin{array}{c|ccc} a_1 & 0 & 0 & 0 \\ a_2 & 0 & 0 & 0 \\ a_3 & 0 & 0 & 0 \\ \hline a_4 & a_5 & a_6 & a_7 \end{array}\right] \;\middle|\; a_i \in R^+ \cup \{0\},\, 1 \le i \le 7,\, +,\, \times_n \right\}$$

is a subsemiring which is also an ideal of $P^+_{4\times 4}$.

Consider

$$L = \left\{ \left[\begin{array}{c|ccc} 0 & a_1 & a_2 & a_3 \\ 0 & a_4 & a_5 & a_6 \\ 0 & a_7 & a_8 & a_9 \\ \hline 0 & 0 & 0 & 0 \end{array}\right] \;\middle|\; a_i \in Z^+ \cup \{0\},\, 1 \le i \le 9,\, +,\, \times_n \right\} \subseteq P^+_{4\times 4};$$

L is a subsemiring of $P^+_{4\times 4}$ which is not an ideal of $P^+_{4\times 4}$. Finally for every $x \in S$ and $y \in L$ we have $x \times_n y = 0$ for every $y \in L$.

*Example 5.63:* Let

$$P^+_{6\times 4} = \left\{ \begin{bmatrix} a_1 & a_2 & a_3 & a_4 \\ a_5 & a_6 & a_7 & a_8 \\ a_9 & a_{10} & a_{11} & a_{12} \\ a_{13} & a_{14} & a_{15} & a_{16} \\ a_{17} & a_{18} & a_{19} & a_{20} \\ a_{21} & a_{22} & a_{23} & a_{24} \end{bmatrix} \middle| a_i \in Q^+ \cup \{0\}, 1 \le i \le 24, +, \times_n \right\}$$

be a semiring of $6 \times 4$ super matrices.

$$I = \begin{bmatrix} 1 & 1 & 1 & 1 \\ 1 & 1 & 1 & 1 \\ 1 & 1 & 1 & 1 \\ 1 & 1 & 1 & 1 \\ 1 & 1 & 1 & 1 \\ 1 & 1 & 1 & 1 \end{bmatrix}$$ is the unit of $P^+_{6\times 4}$ under natural product.

Further

$$S = \left\{ \begin{bmatrix} a_1 & 0 & 0 & a_3 \\ a_5 & 0 & 0 & a_4 \\ 0 & a_7 & a_8 & 0 \\ 0 & a_9 & a_{10} & 0 \\ 0 & a_{11} & a_{12} & 0 \\ a_5 & 0 & 0 & a_6 \end{bmatrix} \middle| a_i \in Z^+ \cup \{0\}, 1 \le i \le 12, +, \times_n \right\}$$

$\subseteq P^+_{6\times 4}$ is only a subsemiring of $P^+_{6\times 4}$ and is not an ideal of $P^+_{6\times 4}$.

Now

$$T = \left\{ \begin{bmatrix} 0 & a_9 & a_{10} & 0 \\ 0 & a_{11} & a_{12} & 0 \\ \hline a_1 & 0 & 0 & a_6 \\ a_2 & 0 & 0 & a_7 \\ a_3 & 0 & 0 & a_8 \\ \hline 0 & a_4 & a_5 & 0 \end{bmatrix} \middle| \; a_i \in Q^+ \cup \{0\}, \right.$$

$$\left. 1 \leq i \leq 12, +, \times_n \right\} \subseteq P_{6\times 4}^+$$

is a subsemiring as well as an ideal of $P_{6\times 4}^+$. Now we see for every $x \in S$ we have every $y \in T$ are such that $x \times_n y = (0)$. Thus $P_{6\times 4}^+$ is only a semiring and not a semifield.

Now we proceed onto define semifields of super matrices.

Let $J_R^+ = \{(x_1 \; x_2 \mid x_3 \; x_4 \; x_5 \mid \ldots \mid x_n) \mid x_i \in R^+ \text{ (or } Q^+ \text{ or } Z^+\text{)},$ $1 \leq i \leq n\} \cup \{(0 \; 0 \mid 0 \; 0 \; 0 \mid \ldots \mid 0)\}, \times, +_n\}$ be the semifield of super row matrices.

For $J_R^+$ has no zero divisors with respect to $\times_n$ and $J_R^+$ is a strict commutative semiring.

Now

$$J_C^+ = \left\{ \begin{bmatrix} m_1 \\ m_2 \\ m_3 \\ \hline m_4 \\ \vdots \\ \hline \vdots \\ \hline m_{n-1} \\ m_n \end{bmatrix} \middle| m_i \in Z^+ \text{ (or } Q^+ \text{ or } R^+\text{)}, 1 \le i \le n \right\} \cup \begin{bmatrix} 0 \\ 0 \\ 0 \\ \hline 0 \\ \vdots \\ \hline \vdots \\ \hline 0 \\ 0 \end{bmatrix}, +, \times_n \}$$

is the semifield of super column matrices.

$$J_{n \times m}^+ \ (m \ne n) = \left\{ \begin{bmatrix} a_{11} & a_{12} & \cdots & a_{1m} \\ a_{21} & a_{22} & \cdots & a_{2m} \\ \vdots & \vdots & & \vdots \\ \hline a_{n1} & a_{n2} & \cdots & a_{nm} \end{bmatrix} \middle| a_{ij} \in R^+ \right.$$

(or $Z^+$ or $Q^+$), $1 \le i \le n; 1 \le j \le m\} \cup \begin{bmatrix} 0 & 0 & \cdots & 0 \\ 0 & 0 & \cdots & 0 \\ \vdots & \vdots & & \vdots \\ \hline 0 & 0 & \cdots & 0 \end{bmatrix}, +, \times_n \}$ is

the semifield of $n \times m$ super matrices.

$$J_{n \times n}^+ = \left\{ \begin{bmatrix} a_{11} & a_{12} & a_{13} & \cdots & a_{1n} \\ a_{21} & a_{22} & a_{23} & \cdots & a_{2n} \\ \vdots & \vdots & & & \vdots \\ \hline a_{n1} & a_{n2} & a_{n3} & \cdots & a_{nn} \end{bmatrix} \middle| a_{ij} \in Z^+ \text{ (or } Q^+ \text{ or } R^+\text{)}, \right.$$

$1 \le i, j \le n\}$

$$\cup \begin{bmatrix} 0 & 0 & 0 & \ldots & 0 \\ 0 & 0 & 0 & \ldots & 0 \\ \vdots & \vdots & \vdots & & \vdots \\ 0 & 0 & 0 & \ldots & 0 \end{bmatrix}, \times_n, +\}$$

is the semifield of square super matrices.

Now we just give some examples of them.

*Example 5.64:* Let

$$J_C^+ = \left\{ \begin{bmatrix} a_1 \\ a_2 \\ a_3 \\ a_4 \\ \underline{a_5} \\ \underline{a_6} \\ \underline{a_7} \\ a_8 \\ a_9 \end{bmatrix} \middle| a_i \in Z^+; 1 \leq i \leq 9 \right\} \cup \begin{bmatrix} 0 \\ 0 \\ 0 \\ 0 \\ \underline{0} \\ \underline{0} \\ \underline{0} \\ 0 \\ 0 \end{bmatrix}, \times_n, +\}$$

be the semifield of super column matrices. This has no proper subsemifields.

*Example 5.65:* Let $J_R^+ = \{(a_1\ a_2\ |\ a_3\ a_4\ a_5\ a_6\ a_7\ |\ a_8\ a_9\ a_{10}\ |\ a_{11})\ |\ a_i \in Q^+; 1 \leq i \leq 11\} \cup (0\ 0\ |\ 0\ 0\ 0\ 0\ 0\ |\ 0\ 0\ 0\ |\ 0), \times_n, +\}$ be the semifield of super row matrices.

*Example 5.66:* Let

$$J^+_{5\times 5} = \left\{ \left\{ \left[\begin{array}{c|cccc} a_1 & a_2 & a_3 & a_4 & a_5 \\ a_6 & a_7 & a_8 & a_9 & a_{10} \\ a_{11} & a_{12} & a_{13} & a_{14} & a_{15} \\ \hline a_{16} & a_{17} & a_{18} & a_{19} & a_{20} \\ \hline a_{21} & a_{22} & a_{23} & a_{24} & a_{25} \end{array}\right] \right| a_i \in R^+; 1 \le i \le 25 \right\}$$

$$\cup \left[\begin{array}{c|cccc} 0 & 0 & 0 & 0 & 0 \\ 0 & 0 & 0 & 0 & 0 \\ 0 & 0 & 0 & 0 & 0 \\ \hline 0 & 0 & 0 & 0 & 0 \\ \hline 0 & 0 & 0 & 0 & 0 \end{array}\right], \times_n, + \right\}$$

is the semifield of super square matrices. This semifield has subsemifields.

*Example 5.67:* Let

$$J^+_{8\times 4} = \left\{ \left\{ \left[\begin{array}{c|ccc} a_1 & a_2 & a_3 & a_4 \\ a_5 & a_6 & a_7 & a_8 \\ \hline a_9 & a_{10} & a_{11} & a_{12} \\ \hline a_{13} & a_{14} & a_{15} & a_{16} \\ \hline a_{17} & a_{18} & a_{19} & a_{20} \\ \hline a_{21} & a_{22} & a_{23} & a_{24} \\ \hline a_{25} & a_{26} & a_{27} & a_{28} \\ \hline a_{29} & a_{30} & a_{31} & a_{32} \end{array}\right] \right| a_i \in Q^+; 1 \le i \le 32 \right\}$$

$$\cup \begin{bmatrix} 0 & 0 & 0 & 0 \\ 0 & 0 & 0 & 0 \\ \hline 0 & 0 & 0 & 0 \\ 0 & 0 & 0 & 0 \\ 0 & 0 & 0 & 0 \\ \hline 0 & 0 & 0 & 0 \\ 0 & 0 & 0 & 0 \\ 0 & 0 & 0 & 0 \end{bmatrix}, \times_n, +\}$$

is a semifield of super $8 \times 4$ matrices. This semifield has subsemifields.

*Example 5.68:*

$$J^+_{9\times 3} = \left\{ \left[ \begin{array}{cc|ccc|c|ccc} a_1 & a_2 & a_3 & a_4 & a_5 & a_6 & a_7 & a_8 & a_9 \\ a_{10} & a_{11} & a_{12} & a_{13} & a_{14} & a_{15} & a_{16} & a_{17} & a_{18} \\ a_{19} & a_{20} & a_{21} & a_{22} & a_{23} & a_{24} & a_{25} & a_{26} & a_{27} \end{array} \right] \right.$$

$$a_i \in R^+; 1 \leq i \leq 27\}$$

$$\cup \begin{bmatrix} 0 & 0 & 0 & 0 & 0 & 0 & 0 & 0 & 0 \\ 0 & 0 & 0 & 0 & 0 & 0 & 0 & 0 & 0 \\ 0 & 0 & 0 & 0 & 0 & 0 & 0 & 0 & 0 \end{bmatrix}, +, \times_n\}$$

is the semifield of super row vectors. This has subsemifields.

*Example 5.69:* Let

$$J^+ = \left\{ \left[ \begin{array}{cc} a_1 & a_2 \\ a_3 & a_4 \\ \hline a_5 & a_6 \\ a_7 & a_8 \\ \hline a_9 & a_{10} \\ a_{11} & a_{12} \\ \hline a_{13} & a_{14} \\ \hline a_{15} & a_{16} \\ \hline a_{17} & a_{18} \\ a_{19} & a_{20} \end{array} \right] \middle| a_i \in Z^+; 1 \leq i \leq 20 \right\} \cup \left[ \begin{array}{cc} 0 & 0 \\ 0 & 0 \\ \hline 0 & 0 \\ 0 & 0 \\ \hline 0 & 0 \\ 0 & 0 \\ \hline 0 & 0 \\ \hline 0 & 0 \\ \hline 0 & 0 \\ 0 & 0 \end{array} \right], +, \times_n \right\}$$

be the semifield of super column vectors. This has no subsemifields but has subsemirings.

**Chapter Six**

# SUPERMATRIX LINEAR ALGEBRAS

In this chapter we introduce the notion of super matrix vector space (linear algebra) and super matrix, polynomials with super matrix coefficients. Several properties enjoyed by them are defined, described and discussed.

Let $V = \{(x_1 \ x_2 \mid x_3 \ \ldots \mid x_{m-1} \ x_n) \mid x_i \in Q \text{ (or R)}; 1 \leq i \leq n\}$ be the collection of super row vectors of same type under addition. V is a vector space over Q (or R). Now we can define the natural product on V so that V is a super linear algebra of super row matrices or linear algebras of super row matrices or super row matrices linear algebras.

We will illustrate this by some simple examples.

***Example 6.1:*** Let $V = \{(x_1 \ x_2 \ x_3 \mid x_4 \mid x_5 \ x_6) \mid x_i \in R; 1 \leq i \leq 6\}$ be a super vector space over the field $F = R$.

V is also a super linear algebra under the natural product $\times_n$.

For $x = (x_1\ x_2\ x_3 \mid x_4 \mid x_5\ x_6)$ and $y = (y_1\ y_2\ y_3 \mid y_4 \mid y_5\ y_6)$ we have $x \times_n y = (x_1y_1\ x_2y_2\ x_3y_3 \mid x_4y_4 \mid x_5y_5\ x_6y_6)$.

*Example 6.2:* Let

$$V = \{(x_1\ x_2 \mid x_3 \mid x_4\ x_5\ x_6 \mid x_7\ x_8 \mid x_9) \mid x_i \in Q;\ 1 \le i \le 9\}$$

be a linear algebra of super row matrices over the field Q, with natural product $\times_n$.

Consider $P_1 = \{(x_1\ x_2 \mid 0 \mid 0\ 0\ 0 \mid 0\ 0 \mid 0) \mid x_1, x_2 \in Q\} \subseteq V$; $P_2 = \{(0\ 0 \mid a_1 \mid 0\ 0\ 0 \mid a_2\ a_3 \mid 0) \mid a_1, a_2, a_3 \in Q\} \subseteq V$ and $P_3 = \{(0\ 0 \mid 0 \mid a_1\ a_2\ a_3 \mid 0\ 0 \mid a_4) \mid a_i \in Q;\ 1 \le i \le 4\} \subseteq V$ be linear subalgebras of V over the field Q.

Clearly $V = P_1 + P_2 + P_3$ and $P_i \cap P_j = (0\ 0 \mid 0 \mid 0\ 0\ 0 \mid 0\ 0 \mid 0)$ if $i \ne j$; $1 \le i, j \le 3$, so V is a direct sum of $P_1$, $P_2$ and $P_3$.

*Example 6.3:* Let

$$V = \{(x_1 \mid x_2\ x_3 \mid x_4\ x_5 \mid x_6\ x_7\ x_8\ x_9\ x_{10}) \mid x_i \in Q;\ 1 \le i \le 10\}$$

be a super row matrix linear algebra over the field Q. Consider $M_1 = \{(a_1 \mid 0 \mid 0 \mid 0\ 0 \mid 0\ 0\ 0\ a_2\ a_3) \mid a_1, a_2, a_3 \in Q\} \subseteq V$, $M_2 = \{(0 \mid a_1 \mid 0 \mid 0\ 0 \mid a_2\ 0\ 0\ 0\ a_3) \mid a_1, a_2, a_3 \in Q\} \subseteq V$, $M_3 = \{(0 \mid 0 \mid a_1 \mid a_2\ 0 \mid 0\ 0\ 0\ 0\ a_3) \mid a_1, a_2, a_3 \in Q\} \subseteq V$, $M_4 = \{(0 \mid 0 \mid 0 \mid 0\ a_1 \mid 0\ a_2\ 0\ 0\ a_3) \mid a_1, a_2, a_3 \in Q\} \subseteq V$ and $M_5 = \{(0 \mid 0 \mid 0 \mid 0\ 0 \mid 0\ 0\ a_1\ a_2\ a_3) \mid a_1, a_2, a_3 \in Q\} \subseteq V$ be a super sublinear algebras of V over Q.

Clearly $M_i \cap M_j \ne (0 \mid 0 \mid 0 \mid 0\ 0 \mid 0\ 0\ 0\ 0\ a)$ if $i \ne j$ for $1 \le i, j \le 5$. Thus $V \subseteq M_1 + M_2 + M_3 + M_4 + M_5$; so V is not a direct sum only a pseudo direct sum.

*Example 6.4:* Let

$$V = \{(x_1\ x_2\ x_3 \mid x_4 \mid x_5 \mid x_6\ x_7\ x_8) \mid x_i \in Q;\ 1 \le i \le 8\}$$

be a linear algebra of super row matrices.

Consider $X = (x_1\ x_2\ x_3\ |\ 0\ |\ 0\ |\ 0\ 0\ x_4)$ where $x_i \in Q$; $1 \le i \le 4\} \subseteq V$ and $Y = \{(0\ 0\ 0\ |\ x_1\ |\ x_2\ |\ x_3\ x_4\ 0)\ |\ x_i \in Q;\ 1 \le i \le 4\} \subseteq V$ be two linear subalgebras of super row matrices, we see for every $x \in X$ and $y \in Y$, $x \times_n y = (0\ 0\ 0\ |\ 0\ |\ 0\ |\ 0\ 0\ 0)$. Thus we see $X^\perp = Y$ and $Y^\perp = X$. Further $V = X+Y$ and $X \cap Y = (0\ 0\ 0\ |\ 0\ |\ 0\ |\ 0\ 0\ 0)$.

We have seen examples of super linear algebras of super row matrices.

*Example 6.5:* Let

$$V = \left\{ \begin{bmatrix} a_1 \\ a_2 \\ \hline a_3 \\ a_4 \\ a_5 \\ \hline a_6 \\ a_7 \\ \hline a_8 \end{bmatrix} \middle| a_i \in Q;\ 1 \le i \le 8 \right\}$$

be a super column linear algebra over the field Q. Dimension of V over Q is eight. Consider

$$x = \begin{bmatrix} 0 \\ 0 \\ \hline a_1 \\ a_2 \\ a_3 \\ \hline 0 \\ 0 \\ \hline 0 \end{bmatrix} \text{ in V then } y = \begin{bmatrix} x \\ y \\ \hline 0 \\ 0 \\ 0 \\ \hline a_1 \\ a_2 \\ \hline a_3 \end{bmatrix}$$

in V is such that

$$x \times_n y = \begin{bmatrix} 0 \\ 0 \\ \hline 0 \\ 0 \\ 0 \\ \hline 0 \\ 0 \\ \hline 0 \end{bmatrix}.$$

Now we can find sublinear algebras of V.

***Example 6.6:*** Let

$$V = \left\{ \begin{bmatrix} a_1 \\ \hline a_2 \\ a_3 \\ \hline a_4 \\ a_5 \\ a_6 \\ \hline a_7 \\ a_8 \\ a_9 \\ \hline a_{10} \\ a_{11} \end{bmatrix} \middle| a_i \in Q; 1 \le i \le 11 \right\}$$

be a linear algebra of super column matrices under the natural product $\times_n$.

Consider

$$M_1 = \left\{ \begin{bmatrix} a_1 \\ a_2 \\ 0 \\ \hline 0 \\ 0 \\ 0 \\ \hline 0 \\ 0 \\ 0 \\ 0 \\ \hline 0 \end{bmatrix} \middle| \ a_1, a_2 \in Q \right\} \subseteq V,$$

$$M_2 = \left\{ \begin{bmatrix} 0 \\ \hline 0 \\ a_1 \\ \hline a_2 \\ a_3 \\ 0 \\ \hline 0 \\ 0 \\ 0 \\ 0 \\ \hline 0 \end{bmatrix} \middle| \ a_1, a_2, a_3 \in Q \right\} \subseteq V,$$

$$M_3 = \left\{ \begin{bmatrix} 0 \\ \hline 0 \\ \hline 0 \\ \hline 0 \\ 0 \\ a_1 \\ a_2 \\ 0 \\ 0 \\ \hline 0 \\ \hline 0 \end{bmatrix} \middle| \; a_1, a_2 \in Q \right\} \subseteq V$$

and

$$M_4 = \left\{ \begin{bmatrix} 0 \\ \hline 0 \\ \hline 0 \\ \hline 0 \\ 0 \\ 0 \\ \hline 0 \\ a_1 \\ a_2 \\ \hline a_3 \\ \hline a_4 \end{bmatrix} \middle| \; a_i \in Q; \; 1 \leq i \leq 4 \right\} \subseteq V$$

be linear subalgebras of super column matrices.

Clearly $M_i \cap M_j = \begin{bmatrix} 0 \\ \hline 0 \\ \hline 0 \\ \hline 0 \\ \hline 0 \\ 0 \\ \hline 0 \\ \hline 0 \\ 0 \\ \hline 0 \\ \hline 0 \end{bmatrix}$ if $i \neq j$, $1 \leq i, j \leq 4$.

Also $V = M_1 + M_2 + M_3 + M_4$; thus V is a direct sum of linear subalgebras of V over Q.

We see for every $x \in M_1$ every $y \in M_2$ is such that $x \times_n y = (0)$. Likewise for every $x \in M_1$, every $z \in M_3$ is such that $x \times_n z = (0)$ and for every $x \in M_1$, every $t \in M_4$ is such that $x \times_n t = (0)$.

Hence we can say for every $x \in M_j$ every element in $M_i$ ($i \neq j$) ($i=1$ or 2 or 3 or 4) is orthogonal with x; however $m_j^\perp$ is not $M_i$ for $i \neq j$, $i = 1$ or 2 or 3 or 4.

Thus we see $m_1^\perp = (M_2 + M_3 + M_4)$; similarly for $M_2 = M_1 + M_3 + M_4$ and so on. These subspaces are not complements of each other.

*Example 6.7:* Let

$$P = \left\{ \begin{bmatrix} x_1 \\ \hline x_2 \\ x_3 \\ \hline x_4 \\ x_5 \\ x_6 \\ \hline x_7 \end{bmatrix} \middle| x_i \in Q; 1 \leq i \leq 7 \right\}$$

be a super linear algebra of super column matrices under the natural product $\times_n$.

Consider

$$M_1 = \left\{ \begin{bmatrix} 0 \\ \hline x_1 \\ x_2 \\ \hline 0 \\ 0 \\ 0 \\ \hline x_3 \end{bmatrix} \middle| x_i \in Q; 1 \leq i \leq 7 \right\} \subseteq P$$

and

$$M_2 = \left\{ \begin{bmatrix} x_1 \\ \hline 0 \\ 0 \\ \hline x_2 \\ x_3 \\ x_4 \\ \hline 0 \end{bmatrix} \middle| x_i \in Q; 1 \leq i \leq 4 \right\} \subseteq P$$

be two super linear subalgebras of P over Q.

We see $M_1 \cap M_2 = \begin{bmatrix} 0 \\ \overline{0} \\ 0 \\ \overline{0} \\ 0 \\ 0 \\ \overline{0} \end{bmatrix}$ and $P = M_1 + M_2$.

Further the complementary subspace of $M_1$ is $M_2$ and vice versa. We see every element in $M_1$ is orthogonal with every element in $M_2$ under orthogonal product.

*Example 6.8:* Let

$$V = \left\{ \begin{bmatrix} a_1 & a_2 & a_3 & a_4 \\ a_5 & a_6 & a_7 & a_8 \\ a_9 & a_{10} & a_{11} & a_{12} \\ \hline a_{13} & a_{14} & a_{15} & a_{16} \end{bmatrix} \middle| a_i \in Q; 1 \leq i \leq 16 \right\}$$

be a super linear algebra of square super matrices under natural product $\times_n$.

We see for

$$x = \begin{bmatrix} 0 & a_1 & a_2 & a_3 \\ 0 & 0 & 0 & 0 \\ 0 & a_4 & a_5 & a_6 \\ \hline 0 & 0 & 0 & 0 \end{bmatrix} \text{ in V, } y = \begin{bmatrix} x_1 & 0 & 0 & 0 \\ x_2 & x_5 & x_6 & x_7 \\ x_3 & 0 & 0 & 0 \\ \hline x_4 & x_5 & x_6 & x_7 \end{bmatrix}$$

in V is such that

$$x \times_n y = \begin{bmatrix} 0 & 0 & 0 & 0 \\ 0 & 0 & 0 & 0 \\ 0 & 0 & 0 & 0 \\ \hline 0 & 0 & 0 & 0 \end{bmatrix} \text{ in V}$$

$$\text{take } y_1 = \begin{bmatrix} 0 & a_1 & a_2 & a_3 \\ 0 & 0 & 0 & 0 \\ 0 & 0 & 0 & 0 \\ \hline 0 & 0 & 0 & 0 \end{bmatrix} \text{ in V}$$

$$\text{we see } x \times_n y_1 = \begin{bmatrix} 0 & 0 & 0 & 0 \\ 0 & 0 & 0 & 0 \\ 0 & 0 & 0 & 0 \\ \hline 0 & 0 & 0 & 0 \end{bmatrix}.$$

Thus we call $y_1$ as a partial complement, only y is the real or total complement of x.

We can have more than one partial complement but one and only one total complement.

*Example 6.9:* Let

$$M = \left\{ \begin{bmatrix} a_1 & a_2 & a_3 & a_4 & a_5 \\ a_6 & a_7 & a_8 & a_9 & a_{10} \\ \hline a_{11} & a_{12} & a_{13} & a_{14} & a_{15} \\ a_{16} & a_{17} & a_{18} & a_{19} & a_{20} \\ a_{21} & a_{22} & a_{23} & a_{24} & a_{25} \end{bmatrix} \middle| a_i \in Q; 1 \leq i \leq 25 \right\}$$

be a super linear of square super matrices under natural product $\times_n$. We can have the following subspaces of M.

Consider

$$P_1 = \left\{ \begin{bmatrix} a_1 & 0 & 0 & a_3 & a_4 \\ a_2 & 0 & 0 & a_5 & a_6 \\ \hline 0 & 0 & 0 & 0 & 0 \\ 0 & 0 & 0 & 0 & 0 \\ 0 & 0 & 0 & 0 & 0 \end{bmatrix} \middle| a_i \in Q; 1 \le i \le 6 \right\} \subseteq M,$$

$$P_2 = \left\{ \begin{bmatrix} 0 & 0 & 0 & a_4 & 0 \\ 0 & 0 & 0 & a_5 & 0 \\ \hline a_1 & 0 & 0 & 0 & a_6 \\ a_2 & 0 & 0 & 0 & a_7 \\ a_3 & 0 & 0 & 0 & a_8 \end{bmatrix} \middle| a_i \in Q; 1 \le i \le 8 \right\} \subseteq M,$$

$$P_3 = \left\{ \begin{bmatrix} 0 & a_1 & a_2 & a_5 & 0 \\ 0 & a_3 & a_4 & a_6 & 0 \\ \hline 0 & 0 & 0 & 0 & 0 \\ 0 & 0 & 0 & 0 & 0 \\ 0 & 0 & 0 & 0 & 0 \end{bmatrix} \middle| a_i \in Q; 1 \le i \le 6 \right\} \subseteq M,$$

$$P_4 = \left\{ \begin{bmatrix} 0 & 0 & 0 & a_5 & 0 \\ 0 & 0 & 0 & a_6 & 0 \\ \hline 0 & a_1 & 0 & a_6 & 0 \\ 0 & a_2 & 0 & a_7 & 0 \\ 0 & a_3 & 0 & a_8 & 0 \end{bmatrix} \middle| a_i \in Q; 1 \le i \le 8 \right\} \subseteq M \text{ and}$$

$$P_5 = \left\{ \begin{bmatrix} 0 & 0 & 0 & a_4 & 0 \\ 0 & 0 & 0 & a_5 & 0 \\ \hline 0 & 0 & a_1 & 0 & 0 \\ 0 & 0 & a_2 & 0 & 0 \\ 0 & 0 & a_3 & 0 & 0 \end{bmatrix} \,\middle|\, a_i \in Q; 1 \leq i \leq 5 \right\} \subseteq M$$

are super linear subalgebras of super square matrices over Q.

We see $V \subseteq P_1 + P_2 + P_3 + P_4 + P_5$ and

$$P_i \cap P_j = \begin{bmatrix} 0 & 0 & 0 & 0 & 0 \\ 0 & 0 & 0 & 0 & 0 \\ \hline 0 & 0 & 0 & 0 & 0 \\ 0 & 0 & 0 & 0 & 0 \\ 0 & 0 & 0 & 0 & 0 \end{bmatrix};$$

if $i \neq j$; $1 \leq i, j \leq 5$.

Thus V is only a pseudo direct sum of linear subalgebras and is not a direct sum of linear subalgebras.

*Example 6.10:* Let

$$V = \left\{ \begin{bmatrix} a_1 & a_2 & a_3 \\ \hline a_4 & a_5 & a_6 \\ a_7 & a_8 & a_9 \end{bmatrix} \,\middle|\, a_i \in Q; 1 \leq i \leq 9 \right\}$$

be a super linear algebra of square super matrices over Q.

Consider

$$P_1 = \left\{ \begin{bmatrix} 0 & 0 & 0 \\ \hline a_1 & a_2 & 0 \\ a_3 & a_4 & 0 \end{bmatrix} \,\middle|\, a_i \in Q; 1 \leq i \leq 4 \right\} \subseteq V$$

and

$$P_2 = \left\{ \begin{array}{|ccc|} \hline a_1 & a_2 & a_3 \\ \hline 0 & 0 & a_4 \\ 0 & 0 & a_5 \\ \hline \end{array} \right| a_i \in Q; 1 \leq i \leq 5 \right\} \subseteq V$$

be super linear subalgebras of V over Q. Clearly the space orthogonal with $P_1$ is $P_2$ and vice versa. No other space of V can be orthogonal (complement) of $P_1$ in V.

Further $V = P_1 + P_2$ and $P_1 \cap P_2 = \begin{bmatrix} 0 & 0 & 0 \\ 0 & 0 & 0 \\ 0 & 0 & 0 \end{bmatrix}$.

*Example 6.11:* Let

$$M = \left\{ \begin{array}{|ccc|} \hline a_1 & a_2 & a_3 \\ \hline a_4 & a_5 & a_6 \\ \hline a_7 & a_8 & a_9 \\ a_{10} & a_{11} & a_{12} \\ a_{13} & a_{14} & a_{15} \\ a_{16} & a_{17} & a_{18} \\ a_{19} & a_{20} & a_{21} \\ \hline \end{array} \right| a_i \in Q; 1 \leq i \leq 21 \right\}$$

be a super linear algebra of super column vector under natural product $\times_n$. Consider

$$P_1 = \left\{ \begin{array}{|ccc|} \hline a_1 & a_2 & a_3 \\ \hline 0 & 0 & 0 \\ \hline 0 & 0 & 0 \\ 0 & 0 & 0 \\ 0 & 0 & 0 \\ 0 & 0 & 0 \\ 0 & 0 & 0 \\ \hline \end{array} \right| a_i \in Q; 1 \leq i \leq 3 \right\} \subseteq M,$$

$$P_2 = \left\{ \begin{bmatrix} 0 & 0 & 0 \\ \hline a_1 & a_2 & a_3 \\ a_4 & a_5 & a_6 \\ \hline 0 & 0 & 0 \\ 0 & 0 & 0 \\ 0 & 0 & 0 \\ 0 & 0 & 0 \end{bmatrix} \middle| a_i \in Q; 1 \leq i \leq 6 \right\} \subseteq M,$$

$$P_3 = \left\{ \begin{bmatrix} 0 & 0 & 0 \\ \hline 0 & 0 & 0 \\ 0 & 0 & 0 \\ \hline 0 & 0 & 0 \\ 0 & 0 & 0 \\ a_1 & a_2 & a_3 \\ a_4 & a_5 & a_6 \end{bmatrix} \middle| a_i \in Q; 1 \leq i \leq 6 \right\} \subseteq M$$

and

$$P_4 = \left\{ \begin{bmatrix} 0 & 0 & 0 \\ \hline 0 & 0 & 0 \\ 0 & 0 & 0 \\ \hline a_1 & a_2 & a_3 \\ a_4 & a_5 & a_6 \\ 0 & 0 & 0 \\ 0 & 0 & 0 \end{bmatrix} \middle| a_i \in Q; 1 \leq i \leq 6 \right\} \subseteq M$$

be super linear subalgebras of super column vectors over the field Q.

We see $P_1 + P_2 + P_3 + P_4 = V$ and $P_i \cap P_j = \begin{bmatrix} 0 & 0 & 0 \\ 0 & 0 & 0 \\ 0 & 0 & 0 \\ \hline 0 & 0 & 0 \\ 0 & 0 & 0 \\ 0 & 0 & 0 \\ 0 & 0 & 0 \end{bmatrix}$, $i \neq j$;

$1 \leq i, j \leq 4$.

*Example 6.12:* Let

$$M = \left\{ \begin{bmatrix} a_1 & a_2 & a_3 & a_4 \\ a_5 & a_6 & a_7 & a_8 \\ \hline a_9 & a_{10} & a_{11} & a_{12} \\ a_{13} & a_{14} & a_{15} & a_{16} \\ a_{17} & a_{18} & a_{19} & a_{20} \\ \hline a_{21} & a_{22} & a_{23} & a_{24} \\ a_{25} & a_{26} & a_{27} & a_{28} \\ a_{29} & a_{30} & a_{31} & a_{32} \\ a_{33} & a_{34} & a_{35} & a_{36} \end{bmatrix} \middle| a_i \in Q; 1 \leq i \leq 36 \right\}$$

be a super linear algebra of super column vectors. Take

$$P_1 = \left\{ \begin{bmatrix} a_1 & a_2 & a_3 & a_4 \\ 0 & 0 & 0 & 0 \\ \hline a_5 & a_6 & a_7 & a_8 \\ 0 & 0 & 0 & 0 \\ 0 & 0 & 0 & 0 \\ \hline 0 & 0 & 0 & 0 \\ 0 & 0 & 0 & 0 \\ 0 & 0 & 0 & 0 \\ a_9 & a_{10} & a_{11} & a_{12} \end{bmatrix} \middle| \, a_i \in Q; 1 \leq i \leq 12 \right\} \subseteq M;$$

$$P_2 = \left\{ \begin{bmatrix} 0 & 0 & 0 & 0 \\ a_1 & a_2 & a_3 & a_4 \\ \hline 0 & 0 & 0 & 0 \\ a_5 & a_6 & a_7 & a_8 \\ 0 & 0 & 0 & 0 \\ \hline 0 & 0 & 0 & 0 \\ 0 & 0 & 0 & 0 \\ 0 & 0 & 0 & 0 \\ a_9 & a_{10} & a_{11} & a_{12} \end{bmatrix} \middle| \, a_i \in Q; 1 \leq i \leq 12 \right\} \subseteq M,$$

$$P_3 = \left\{ \begin{bmatrix} 0 & 0 & 0 & 0 \\ 0 & 0 & 0 & 0 \\ \hline 0 & 0 & 0 & 0 \\ 0 & 0 & 0 & 0 \\ a_1 & a_2 & a_3 & a_4 \\ \hline a_5 & a_6 & a_7 & a_8 \\ 0 & 0 & 0 & 0 \\ 0 & 0 & 0 & 0 \\ a_9 & a_{10} & a_{11} & a_{12} \end{bmatrix} \middle| \, a_i \in Q; 1 \leq i \leq 12 \right\} \subseteq M$$

and

$$P_4 = \left\{ \begin{bmatrix} 0 & 0 & 0 & 0 \\ 0 & 0 & 0 & 0 \\ \hline 0 & 0 & 0 & 0 \\ 0 & 0 & 0 & 0 \\ 0 & 0 & 0 & 0 \\ \hline 0 & 0 & 0 & 0 \\ a_1 & a_2 & a_3 & a_4 \\ a_5 & a_6 & a_7 & a_8 \\ a_9 & a_{10} & a_{11} & a_{12} \end{bmatrix} \middle| a_i \in Q; 1 \le i \le 12 \right\} \subseteq M$$

be super linear subalgebras M over the field Q.

$$\text{Clearly } P_i \cap P_j \ne \begin{bmatrix} 0 & 0 & 0 & 0 \\ 0 & 0 & 0 & 0 \\ \hline 0 & 0 & 0 & 0 \\ 0 & 0 & 0 & 0 \\ 0 & 0 & 0 & 0 \\ \hline 0 & 0 & 0 & 0 \\ 0 & 0 & 0 & 0 \\ 0 & 0 & 0 & 0 \\ 0 & 0 & 0 & 0 \end{bmatrix}; \text{ if } i \ne j; 1 \le i, j \le 4.$$

Further $V \subseteq P_1 + P_2 + P_3 + P_4$; thus V is only a pseudo direct sum of super linear subalgebras of super column vectors over Q.

*Example 6.13:* Let M =

$$\left\{ \begin{bmatrix} a_1 & a_5 & a_9 & a_{13} & a_{17} & a_{21} & a_{25} & a_{26} & a_{27} & a_{28} & a_{29} \\ a_2 & a_6 & a_{10} & a_{14} & a_{18} & a_{22} & a_{30} & a_{31} & a_{32} & a_{33} & a_{34} \\ a_3 & a_7 & a_{11} & a_{15} & a_{19} & a_{23} & a_{35} & a_{36} & a_{37} & a_{38} & a_{39} \\ a_4 & a_8 & a_{12} & a_{16} & a_{20} & a_{24} & a_{40} & a_{41} & a_{42} & a_{43} & a_{44} \end{bmatrix} \middle| \right.$$
$$\left. a_i \in Q; 1 \le i \le 44 \right\}$$

be the super linear algebra of super row vectors over Q.
Consider

$$P_1 = \left\{ \begin{bmatrix} a_1 & a_5 & 0 & 0 & 0 & 0 & 0 & 0 & 0 & 0 & 0 \\ a_2 & a_6 & 0 & 0 & 0 & 0 & 0 & 0 & 0 & 0 & 0 \\ a_3 & a_7 & 0 & 0 & 0 & 0 & 0 & 0 & 0 & 0 & 0 \\ a_4 & a_8 & 0 & 0 & 0 & 0 & 0 & 0 & 0 & 0 & 0 \end{bmatrix} \right.$$

$a_i \in Q; 1 \leq i \leq 8\} \subseteq M$,

$$P_2 = \left\{ \begin{bmatrix} 0 & 0 & a_1 & a_5 & 0 & 0 & 0 & 0 & 0 & 0 & 0 \\ 0 & 0 & a_2 & a_6 & 0 & 0 & 0 & 0 & 0 & 0 & 0 \\ 0 & 0 & a_3 & a_7 & 0 & 0 & 0 & 0 & 0 & 0 & 0 \\ 0 & 0 & a_4 & a_8 & 0 & 0 & 0 & 0 & 0 & 0 & 0 \end{bmatrix} \right.$$

$a_i \in Q; 1 \leq i \leq 8\} \subseteq M$,

$$P_3 = \left\{ \begin{bmatrix} 0 & 0 & 0 & 0 & a_1 & a_5 & 0 & 0 & 0 & 0 & 0 \\ 0 & 0 & 0 & 0 & a_2 & a_6 & 0 & 0 & 0 & 0 & 0 \\ 0 & 0 & 0 & 0 & a_3 & a_7 & 0 & 0 & 0 & 0 & 0 \\ 0 & 0 & 0 & 0 & a_4 & a_8 & 0 & 0 & 0 & 0 & 0 \end{bmatrix} \right.$$

$a_i \in Q; 1 \leq i \leq 8\} \subseteq M$,

$$P_4 = \left\{ \begin{bmatrix} 0 & 0 & 0 & 0 & 0 & 0 & a_1 & a_2 & 0 & 0 & 0 \\ 0 & 0 & 0 & 0 & 0 & 0 & a_3 & a_4 & 0 & 0 & 0 \\ 0 & 0 & 0 & 0 & 0 & 0 & a_5 & a_6 & 0 & 0 & 0 \\ 0 & 0 & 0 & 0 & 0 & 0 & a_7 & a_8 & 0 & 0 & 0 \end{bmatrix} \right| a_i \in Q;$$

$1 \leq i \leq 8\} \subseteq M$

and

$$P_5 = \left\{ \begin{bmatrix} 0 & 0 & 0 & 0 & 0 & 0 & 0 & 0 & a_1 & a_2 & a_3 \\ 0 & 0 & 0 & 0 & 0 & 0 & 0 & 0 & a_4 & a_5 & a_6 \\ 0 & 0 & 0 & 0 & 0 & 0 & 0 & 0 & a_7 & a_8 & a_9 \\ 0 & 0 & 0 & 0 & 0 & 0 & 0 & 0 & a_{10} & a_{11} & a_{12} \end{bmatrix} \middle| a_i \in Q; \right.$$

$$1 \leq i \leq 12\} \subseteq M,$$

be super linear subalgebras of V of super row vectors over Q.

Clearly $P_i \cap P_j = \begin{bmatrix} 0 & 0 & 0 & 0 & 0 & 0 & 0 & 0 & 0 & 0 & 0 \\ 0 & 0 & 0 & 0 & 0 & 0 & 0 & 0 & 0 & 0 & 0 \\ 0 & 0 & 0 & 0 & 0 & 0 & 0 & 0 & 0 & 0 & 0 \\ 0 & 0 & 0 & 0 & 0 & 0 & 0 & 0 & 0 & 0 & 0 \end{bmatrix}$ if

$i \neq j$ and $1 \leq i, j \leq 5$.

$M = P_1 + P_2 + P_3 + P_4 + P_5$ so M is the direct sum of super linear subalgebras of super row vectors over the field Q.

*Example 6.14:* Let V =

$$\left\{ \begin{bmatrix} a_1 & a_2 & a_3 & a_4 & a_5 & a_6 & a_7 & a_8 & a_9 & a_{10} \\ a_{11} & a_{12} & a_{13} & a_{14} & a_{15} & a_{16} & a_{17} & a_{18} & a_{19} & a_{20} \end{bmatrix} \middle| a_i \in Q; \right.$$

$$1 \leq i \leq 20\}$$

be a super linear algebra of super row vectors over Q.

Consider

$$H_1 = \left\{ \begin{bmatrix} a_1 & a_3 & 0 & 0 & 0 & 0 & 0 & 0 & 0 \\ a_2 & a_4 & 0 & 0 & 0 & 0 & 0 & 0 & 0 \end{bmatrix} \middle| a_i \in Q; \right.$$

$$1 \leq i \leq 4\} \subseteq V,$$

$$H_2 = \left\{ \begin{bmatrix} 0 & a_2 & a_3 & 0 & 0 & 0 & 0 & 0 & 0 \\ 0 & a_1 & a_4 & 0 & 0 & 0 & 0 & 0 & 0 \end{bmatrix} \right\} \subseteq V,$$

$$H_3 = \left\{ \begin{bmatrix} 0 & a_1 & 0 & a_3 & 0 & 0 & 0 & 0 & 0 \\ 0 & a_2 & 0 & a_4 & 0 & 0 & 0 & 0 & 0 \end{bmatrix} \right\} \subseteq V,$$

$$H_4 = \left\{ \begin{bmatrix} 0 & a_1 & 0 & 0 & a_3 & 0 & 0 & 0 & 0 \\ 0 & a_2 & 0 & 0 & a_4 & 0 & 0 & 0 & 0 \end{bmatrix} \right\} \subseteq V,$$

$$H_5 = \left\{ \begin{bmatrix} 0 & a_1 & 0 & 0 & 0 & a_3 & 0 & 0 & 0 \\ 0 & a_2 & 0 & 0 & 0 & a_4 & 0 & 0 & 0 \end{bmatrix} \right\} \subseteq V,$$

$$H_6 = \left\{ \begin{bmatrix} 0 & a_1 & 0 & 0 & 0 & 0 & a_3 & 0 & 0 \\ 0 & a_2 & 0 & 0 & 0 & 0 & a_4 & 0 & 0 \end{bmatrix} \right\} \subseteq V,$$

$$H_7 = \left\{ \begin{bmatrix} 0 & a_1 & 0 & 0 & 0 & 0 & 0 & a_3 & 0 \\ 0 & a_2 & 0 & 0 & 0 & 0 & 0 & a_4 & 0 \end{bmatrix} \right\} \subseteq V,$$

$$H_8 = \left\{ \begin{bmatrix} 0 & a_1 & 0 & 0 & 0 & 0 & 0 & 0 & a_3 \\ 0 & a_2 & 0 & 0 & 0 & 0 & 0 & 0 & a_4 \end{bmatrix} \; \middle| \; a_i \in Q; 1 \leq i \leq 4 \right\} \subseteq V$$

are super linear subalgebras of super row vector over the field Q.

Clearly

$$H_i \cap H_j \neq \left\{ \begin{bmatrix} 0 & a_1 & 0 & 0 & 0 & 0 & 0 & 0 & 0 \\ 0 & a_2 & 0 & 0 & 0 & 0 & 0 & 0 & 0 \end{bmatrix} \; \middle| \; a_i \in Q; 1 \leq i \leq 4 \right\}$$

if $i \neq j$, $1 \leq i, j \leq 8$. We see $V \subseteq H_1 + H_2 + H_3 + H_4 + H_5 + H_6 + H_7 + H_8$ is only a pseudo direct sum of the sublinear algebras of V over Q.

*Example 6.15:* Let

$$V = \left\{ \begin{bmatrix} a_1 & a_2 & a_3 & a_4 & a_5 \\ a_6 & a_7 & a_8 & a_9 & a_{10} \\ \hline a_{11} & a_{12} & a_{13} & a_{14} & a_{15} \\ a_{16} & a_{17} & a_{18} & a_{19} & a_{20} \\ \hline a_{21} & a_{22} & a_{23} & a_{24} & a_{25} \\ \hline a_{26} & a_{27} & a_{28} & a_{29} & a_{30} \end{bmatrix} \middle| a_i \in Q; 1 \le i \le 30 \right\}$$

be a super linear algebra of $6 \times 5$ super matrices over the field Q under the natural product $\times_n$.

Consider

$$M = \left\{ \begin{bmatrix} a_1 & 0 & 0 & 0 & 0 \\ a_2 & 0 & 0 & 0 & 0 \\ \hline 0 & 0 & 0 & 0 & 0 \\ 0 & 0 & 0 & 0 & 0 \\ \hline 0 & 0 & 0 & 0 & 0 \\ \hline 0 & a_3 & a_4 & 0 & 0 \end{bmatrix} \middle| a_i \in Q; 1 \le i \le 4 \right\} \subseteq V,$$

$$M_2 = \left\{ \begin{bmatrix} 0 & 0 & 0 & 0 & 0 \\ 0 & 0 & 0 & 0 & 0 \\ \hline a_1 & 0 & 0 & 0 & 0 \\ a_2 & 0 & 0 & 0 & 0 \\ \hline a_3 & 0 & 0 & 0 & 0 \\ \hline a_4 & 0 & 0 & 0 & 0 \end{bmatrix} \middle| a_i \in Q; 1 \le i \le 4 \right\} \subseteq V,$$

$$M_3 = \left\{ \begin{bmatrix} 0 & a_1 & a_2 & 0 & 0 \\ 0 & a_3 & a_4 & 0 & 0 \\ \hline 0 & 0 & 0 & 0 & 0 \\ 0 & 0 & 0 & 0 & 0 \\ 0 & 0 & 0 & 0 & 0 \\ \hline 0 & 0 & 0 & 0 & 0 \end{bmatrix} \middle| a_i \in Q; 1 \le i \le 4 \right\} \subseteq V,$$

$$M_4 = \left\{ \begin{bmatrix} 0 & 0 & 0 & 0 & 0 \\ 0 & 0 & 0 & 0 & 0 \\ \hline 0 & a_1 & a_2 & 0 & 0 \\ 0 & a_3 & a_4 & 0 & 0 \\ 0 & a_5 & a_6 & 0 & 0 \\ \hline 0 & 0 & 0 & 0 & 0 \end{bmatrix} \middle| a_i \in Q; 1 \le i \le 4 \right\} \subseteq V,$$

$$M_5 = \left\{ \begin{bmatrix} 0 & 0 & 0 & a_1 & a_2 \\ 0 & 0 & 0 & a_3 & a_4 \\ \hline 0 & 0 & 0 & 0 & 0 \\ 0 & 0 & 0 & 0 & 0 \\ 0 & 0 & 0 & 0 & 0 \\ \hline 0 & 0 & 0 & 0 & 0 \end{bmatrix} \middle| a_i \in Q; 1 \le i \le 4 \right\} \subseteq V,$$

$$M_6 = \left\{ \begin{bmatrix} 0 & 0 & 0 & 0 & 0 \\ 0 & 0 & 0 & 0 & 0 \\ \hline 0 & 0 & 0 & 0 & 0 \\ 0 & 0 & 0 & 0 & 0 \\ 0 & 0 & 0 & 0 & 0 \\ \hline 0 & 0 & 0 & a_1 & a_2 \end{bmatrix} \middle| a_i \in Q; 1 \le i \le 4 \right\} \subseteq V,$$

and

$$M_7 = \left\{ \begin{bmatrix} 0 & 0 & 0 & 0 & 0 \\ 0 & 0 & 0 & 0 & 0 \\ \hline 0 & 0 & 0 & a_1 & a_2 \\ 0 & 0 & 0 & a_3 & a_4 \\ 0 & 0 & 0 & a_5 & a_6 \\ \hline 0 & 0 & 0 & 0 & 0 \end{bmatrix} \,\middle|\, a_i \in Q;\ 1 \le i \le 4 \right\} \subseteq V$$

be super linear subalgebra of super matrices over Q under the natural product $\times_n$.

Clearly $M_i \cap M_j = \begin{bmatrix} 0 & 0 & 0 & 0 & 0 \\ 0 & 0 & 0 & 0 & 0 \\ \hline 0 & 0 & 0 & 0 & 0 \\ 0 & 0 & 0 & 0 & 0 \\ 0 & 0 & 0 & 0 & 0 \\ \hline 0 & 0 & 0 & 0 & 0 \end{bmatrix}$ if $i \ne j$, $1 \le i, j \le 7$.

We see $V = M_1 + M_2 + M_3 + M_4 + M_5 + M_6 + M_7$, that is V is the direct sum of sublinear algebras of V.

Now we proceed onto give examples of linear transformation and linear operators on super linear algebra of super matrices with natural product $\times_n$.

*Example 6.16:* Let

$$M = \left\{ \begin{bmatrix} a_1 & a_2 & a_3 \\ \hline a_4 & a_5 & a_6 \\ \hline a_7 & a_8 & a_9 \end{bmatrix} \,\middle|\, a_i \in Q;\ 1 \le i \le 9 \right\}$$

be a super linear algebra of square super matrices over the field Q under the natural product $\times_n$.

Let

$$P = \left\{ \begin{bmatrix} a_1 \\ \overline{a_2} \\ \overline{a_3} \\ a_4 \\ \overline{a_5} \\ a_6 \\ \overline{a_7} \\ a_8 \\ \overline{a_9} \end{bmatrix} \middle| a_i \in Q;\ 1 \leq i \leq 9 \right\}$$

be a super linear algebra of super column matrices over the field Q under the natural product $\times_n$.

Define $\eta : M \to P$ by

$$\eta \left( \begin{bmatrix} a_1 & a_2 & a_3 \\ \hline a_4 & a_5 & a_6 \\ a_7 & a_8 & a_9 \end{bmatrix} \right) = \begin{bmatrix} a_1 \\ a_2 \\ \overline{a_3} \\ a_4 \\ \overline{a_5} \\ a_6 \\ \overline{a_7} \\ a_8 \\ a_9 \end{bmatrix} ;$$

it is easily verified $\eta$ is a linear transformation of super linear algebras.

*Example 6.17:* Let

$$M = \left\{ \begin{bmatrix} a_1 & a_2 & a_3 \\ a_4 & a_5 & a_6 \\ \hline a_7 & a_8 & a_9 \\ a_{10} & a_{11} & a_{12} \\ a_{13} & a_{14} & a_{15} \\ a_{16} & a_{17} & a_{18} \end{bmatrix} \middle| \; a_i \in Q;\; 1 \leq i \leq 18 \right\}$$

be a super linear algebra of super matrices over the field Q under the natural product $\times_n$.

Let

$$P = \left\{ \begin{bmatrix} a_1 & a_3 & a_5 & a_6 & a_7 & a_{11} & a_{13} & a_{15} & a_{17} \\ a_2 & a_4 & a_8 & a_9 & a_{10} & a_{12} & a_{14} & a_{16} & a_{18} \end{bmatrix} \middle| \; a_i \in Q;\; 1 \leq i \leq 18 \right\}$$

be a super linear algebra of super matrices over the field Q under natural product $\times_n$.

Define $\eta : M \to P$ by

$$\eta \left( \begin{bmatrix} a_1 & a_2 & a_3 \\ a_4 & a_5 & a_6 \\ \hline a_7 & a_8 & a_9 \\ a_{10} & a_{11} & a_{12} \\ a_{13} & a_{14} & a_{15} \\ a_{16} & a_{17} & a_{18} \end{bmatrix} \right) = \begin{bmatrix} a_1 & a_3 & a_5 & a_6 & a_7 & a_{11} & a_{13} & a_{15} & a_{17} \\ a_2 & a_4 & a_8 & a_9 & a_{10} & a_{12} & a_{14} & a_{16} & a_{18} \end{bmatrix}$$

η is a linear transformation from M to P.

We now proceed onto define linear operator of super linear algebras of super matrices over the field F under natural product $\times_n$.

*Example 6.18:* Let

$$V = \left\{ \begin{bmatrix} a_1 & a_2 & a_3 & a_4 \\ a_5 & a_6 & a_7 & a_8 \\ a_9 & a_{10} & a_{11} & a_{12} \\ a_{13} & a_{14} & a_{15} & a_{16} \\ a_{17} & a_{18} & a_{19} & a_{20} \end{bmatrix} \middle| a_i \in Q; 1 \leq i \leq 20 \right\}$$

be a super linear algebra of super matrices under the natural product $\times_n$.

Consider $\eta : V \to V$ defined by

$$\eta \left( \begin{bmatrix} a_1 & a_2 & a_3 & a_4 \\ a_5 & a_6 & a_7 & a_8 \\ a_9 & a_{10} & a_{11} & a_{12} \\ a_{13} & a_{14} & a_{15} & a_{16} \\ a_{17} & a_{18} & a_{19} & a_{20} \end{bmatrix} \right) = \left( \begin{bmatrix} 0 & a_1 & a_2 & 0 \\ 0 & a_3 & a_4 & 0 \\ 0 & a_5 & a_6 & 0 \\ 0 & a_7 & a_8 & 0 \\ 0 & 0 & 0 & 0 \end{bmatrix} \right).$$

η is easily verified to be a linear operator on V. Consider μ : V → V defined by

$$\mu \left( \begin{bmatrix} a_1 & a_2 & a_3 & a_4 \\ a_5 & a_6 & a_7 & a_8 \\ a_9 & a_{10} & a_{11} & a_{12} \\ a_{13} & a_{14} & a_{15} & a_{16} \\ a_{17} & a_{18} & a_{19} & a_{20} \end{bmatrix} \right) = \left( \begin{bmatrix} a_1 & 0 & 0 & a_{11} \\ a_2 & 0 & 0 & a_{10} \\ a_3 & 0 & 0 & a_9 \\ a_4 & 0 & 0 & a_8 \\ a_5 & a_6 & a_7 & a_7 \end{bmatrix} \right).$$

µ is also a linear operator on V.

**Example 6.19:** Let

$$V = \left\{ \begin{bmatrix} a_1 & a_2 & a_3 \\ \hline a_4 & a_5 & a_6 \\ a_7 & a_8 & a_9 \\ a_{10} & a_{11} & a_{12} \\ \hline a_{13} & a_{14} & a_{15} \\ a_{16} & a_{17} & a_{18} \\ \hline a_{19} & a_{20} & a_{21} \end{bmatrix} \,\middle|\, a_i \in Q; 1 \le i \le 21 \right\}$$

be a super linear algebra of super column vectors defined over the field Q, under the natural product $\times_n$.

Define $\eta : V \to V$

$$\text{by } \eta \left( \begin{bmatrix} a_1 & a_2 & a_3 \\ \hline a_4 & a_5 & a_6 \\ a_7 & a_8 & a_9 \\ a_{10} & a_{11} & a_{12} \\ \hline a_{13} & a_{14} & a_{15} \\ a_{16} & a_{17} & a_{18} \\ \hline a_{19} & a_{20} & a_{21} \end{bmatrix} \right) = \begin{bmatrix} a_1 & a_2 & a_3 \\ 0 & 0 & 0 \\ \hline a_4 & a_5 & a_6 \\ 0 & 0 & 0 \\ a_7 & a_8 & a_9 \\ \hline 0 & 0 & 0 \\ a_{10} & a_{11} & a_{12} \end{bmatrix}.$$

$\eta$ is a linear operator on V.

We can also define linear function which is a matter of routine. However we give examples of them.

**Example 6.20:** let $V = \{(x_1\ x_2\ |\ x_3\ |\ x_4\ x_5\ x_6) \mid x_i \in Q; 1 \le i \le 6\}$ be a super linear algebra of super row matrices over the field Q. Define $f : V \to Q$ such that $f((x_1\ x_2\ |\ x_3\ |\ x_4\ x_5\ x_6)) = x_1 + x_2 + x_6$ is a linear functional on V.

*Example 6.21:* Let

$$M = \left\{ \begin{bmatrix} a_1 & a_2 & a_3 & a_4 \\ a_5 & a_6 & a_7 & a_8 \\ a_9 & a_{10} & a_{11} & a_{12} \\ a_{13} & a_{14} & a_{15} & a_{16} \\ a_{17} & a_{18} & a_{19} & a_{20} \\ a_{21} & a_{22} & a_{23} & a_{24} \\ a_{25} & a_{26} & a_{27} & a_{28} \end{bmatrix} \middle| a_i \in Q; 1 \leq i \leq 28 \right\}$$

be a super linear algebra of super matrices over the field Q under the natural product $\times_n$.

Define $f : M \to Q$ by

$$f \left( \begin{bmatrix} a_1 & a_2 & a_3 & a_4 \\ a_5 & a_6 & a_7 & a_8 \\ a_9 & a_{10} & a_{11} & a_{12} \\ a_{13} & a_{14} & a_{15} & a_{16} \\ a_{17} & a_{18} & a_{19} & a_{20} \\ a_{21} & a_{22} & a_{23} & a_{24} \\ a_{25} & a_{26} & a_{27} & a_{28} \end{bmatrix} \right) = a_1 + 3a_3 + 9a_{27},$$

f is a linear functional on M.

Interested reader can develop other related properties as all properties can be derived with appropriate modifications provided the situation is feacible. We can define super matrix coefficients polynomials

Let $F_R^S[x] = \left\{ \sum_{i=0}^{\infty} a_i x^i \middle| a_i \in (x_1 \mid x_2 \ x_3 \mid \ldots \mid x_{n-1} \ x_n) \in M = \right.$ {collection of all super row matrices of same type with entries

from Q or Z or R}. $F_R^S[x]$ is defined as the super row matrix coefficient polynomials or polynomials in the variable x with row super matrix coefficients.

We will illustrate this situation by an example.

*Example 6.22:* Let

$$F_R^S[x] = \left\{ \sum a_i x^i \ \middle| \ a_i \in (x_1 \mid x_2 \ x_3 \mid x_4 \ x_5 \ x_6 \mid x_7) \in M \right.$$

= {all super row matrices of the type $(x_1 \mid x_2 \ x_3 \mid x_4 \ x_5 \ x_6 \mid x_7)$ with $x_j \in R$, $1 \leq j \leq 7$}} be the polynomials with super row matrix coefficients.

*Example 6.23:* Let

$$F_R^S[x] = \left\{ \sum a_i x^i \ \middle| \ a_i \in (x_1 \mid x_2 \mid x_3 \mid x_4) \in P = \right.$$

{all super row matrices of the form $(x_1 \mid x_2 \mid x_3 \mid x_4)$ with $x_j \in Z$, $1 \leq j \leq 4$}} be the super row matrix coefficient polynomials.

We will illustrate how a polynomial looks like, its degree and operations on them with super row matrix coefficients.

Let $p(x) = (3 \mid 2 \mid -7 \mid 2 \ 0) + (7 \mid -4 \mid 0 \mid 3 \ 1)x^3 + (0 \mid 0 \mid 2 \mid -7 \ -9)x^5$ be a super row matrix coefficient polynomial. Clearly degree of $p(x)$ is 5.

$q(x) = (3 \mid 2 \ 0) + (-7 \mid -2 \ 9) \ x + (9 \mid 0 \ 0)x^2 + (0 \mid 3 \ 1)x^3 + (1 \mid 1 \ 1)x^4$ is a super row matrix coefficient polynomial in the variable x and degree of $q(x)$ is four.

Clearly $0(x) = (0 \mid 0 \ 0) + (0 \mid 0 \ 0)x + (0 \mid 0 \ 0)x^2 + \ldots + (0 \mid 0 \ 0)x^n$.

Now we can add two polynomial with super row matrices if and only if all the coefficients are from the same type of super row matrices.

Clearly we cannot add p(x) with q(x). However q(x) + 0(x) can be added and q(x) + 0(x) = q(x).

We will illustrate addition of two super matrix polynomials.

Let m(x) = (1 1 | 0 2 3| 7 5 0 1) + (0 1 | 2 0 1 | 0 0 1 1) x + (0 x | 1 0 0 | 8 0 0 5)$x^2$ + (0 0 | 0 0 1 | 20 1 2)$x^3$ and n(x) = (0 1 | 2 2 2 | 3 1 2 0) + (6 2 | 0 0 0 | 2 1 0 0)x + (0 1 | 1 0 1| 0 2 0 1)$x^2$ + (0 4 | 4 2 –1 | 0 7 2 1)$x^3$ + (1 2 | 0 1 4 | 3 0 1 4)$x^4$ be two super row matrix coefficient polynomials in the variable x.

m(x) + (n(x)) = (1 2 | 2 4 5 | 1 0 6 2 1) + (6 3 | 2 0 1| 2 1 1 1)x + (0 9 | 2 0 1 | 8 2 0 6)$x^2$ + (0 4 | 4 2 0 | 2 7 3 3)$x^3$ + (1 2 | 0 1 4 | 3 0 1 4)$x^4$. Thus we see addition of two super row matrix coefficient polynomials is again a polynomial with super row matrices coefficients. Infact the set of super row matrix coefficient polynomials under addition is a group. Further it is a commutative group under addition.

We will give examples of such groups.

***Example 6.24:*** Let $F_R^S = \left\{ \sum_{i=0}^{\infty} a_i x^i \;\middle|\; a_i = (x_1 \mid x_2 \mid x_3 \; x_4 \; x_5 \mid x_6 \; x_7) \right.$

$\in$ {to the collection of all super row matrices of same type with entries from R}, +} be an abelian group of infinite order. This has subgroups.

***Example 6.25:*** Let $F_R^S = \left\{ \sum_{i=0}^{\infty} a_i x^i \;\middle|\; a_i = (x_1 \; x_2 \; x_3 \; x_4 \mid x_5 \; x_6) \in N \right.$

= {all super row matrices of the same type as ($x_1 \; x_2 \; x_3 \; x_4 \mid x_5 \; x_6$) with entries from R}, +} be a group under addition.

Now we proceed onto give examples of super column matrix coefficient polynomials.

$$p(x) = \begin{bmatrix} 3/2 \\ 0 \\ 1 \\ 1 \\ 4/5 \end{bmatrix} + \begin{bmatrix} 0 \\ 1 \\ 2/0 \\ 1 \\ 4/0 \end{bmatrix} x + \begin{bmatrix} 7 \\ 0 \\ 8/5 \\ 0 \\ 1/8 \end{bmatrix} x^3 + \begin{bmatrix} 8 \\ 0 \\ 7/0 \\ 0 \\ 1/9 \end{bmatrix} x^4 + \begin{bmatrix} 1/0 \\ 0 \\ 0 \\ 0 \\ 9/2 \end{bmatrix} x^6$$

is a super column matrix coefficient polynomial of degree six.

$$\text{Consider } q(x) = \begin{bmatrix} 3/1 \\ 1 \\ 1/8 \end{bmatrix} + \begin{bmatrix} 0 \\ 2 \\ 3 \\ 4/0 \end{bmatrix} x + \begin{bmatrix} 9/2 \\ 3 \\ 1/-8 \end{bmatrix} x^2 + \begin{bmatrix} 1/2 \\ 3 \\ 4/0 \end{bmatrix} x^3 + \begin{bmatrix} 9/2 \\ 3 \\ 4/9 \end{bmatrix} x^8$$

is a column matrix coefficient polynomial of degree 8 in the variable x.

Now we show how addition of column matrix coefficient polynomials are carried out in case of same type of column matrices. For if these column matrices are different type certainly we cannot add any two column matrix coefficient polynomials.

Let p(x) = $\begin{bmatrix} 3 \\ 2 \\ 1 \\ 0 \\ 1 \\ 5 \end{bmatrix} + \begin{bmatrix} 0 \\ 1 \\ 2 \\ 3 \\ 4 \\ 7 \end{bmatrix}$ x + $\begin{bmatrix} 2 \\ 0 \\ 1 \\ 0 \\ 0 \\ -8 \end{bmatrix}$ x$^2$ + $\begin{bmatrix} 0 \\ 1 \\ 2 \\ 0 \\ 7 \\ 9 \end{bmatrix}$ x$^3$ + $\begin{bmatrix} 9 \\ 2 \\ 0 \\ 1 \\ 1 \\ 8 \end{bmatrix}$ x$^5$ and

q(x) = $\begin{bmatrix} 9 \\ 0 \\ 1 \\ 2 \\ 5 \\ 0 \end{bmatrix} + \begin{bmatrix} 8 \\ 7 \\ 0 \\ 0 \\ 7 \\ 6 \end{bmatrix}$ x + $\begin{bmatrix} 1 \\ 2 \\ 2 \\ 9 \\ 0 \\ 7 \end{bmatrix}$ x$^3$ + $\begin{bmatrix} 9 \\ 0 \\ 1 \\ 2 \\ 7 \\ 8 \end{bmatrix}$ x$^5$.

p(x) + q(x) = $\begin{bmatrix} 3 \\ 2 \\ 1 \\ 0 \\ 1 \\ 5 \end{bmatrix} + \begin{bmatrix} 9 \\ 0 \\ 1 \\ 2 \\ 5 \\ 0 \end{bmatrix} + \left( \begin{bmatrix} 0 \\ 1 \\ 2 \\ 3 \\ 4 \\ 7 \end{bmatrix} + \begin{bmatrix} 8 \\ 7 \\ 0 \\ 0 \\ 7 \\ 6 \end{bmatrix} \right)$ x

+ $\left( \begin{bmatrix} 2 \\ 0 \\ 1 \\ 0 \\ 0 \\ -8 \end{bmatrix} \right)$ x$^2$ + $\left( \begin{bmatrix} 0 \\ 1 \\ 2 \\ 0 \\ 7 \\ 9 \end{bmatrix} + \begin{bmatrix} 1 \\ 2 \\ 2 \\ 9 \\ 0 \\ 7 \end{bmatrix} \right)$ x$^3$ + $\left( \begin{bmatrix} 9 \\ 2 \\ 0 \\ 1 \\ 1 \\ 8 \end{bmatrix} + \begin{bmatrix} 9 \\ 0 \\ 1 \\ 2 \\ 7 \\ 8 \end{bmatrix} \right)$ x$^5$

$$= \begin{bmatrix} 12 \\ 2 \\ \dfrac{2}{2} \\ \dfrac{6}{5} \end{bmatrix} + \begin{bmatrix} 8 \\ 8 \\ \dfrac{2}{3} \\ \dfrac{12}{13} \end{bmatrix} x + \begin{bmatrix} 2 \\ 0 \\ \dfrac{1}{0} \\ \dfrac{0}{-8} \end{bmatrix} x^2 + \begin{bmatrix} 1 \\ 3 \\ \dfrac{4}{9} \\ \dfrac{7}{16} \end{bmatrix} x^3 + \begin{bmatrix} 18 \\ 2 \\ \dfrac{1}{3} \\ \dfrac{8}{16} \end{bmatrix} x^5.$$

This is the way addition of super column matrix coefficient polynomials are added. Thus addition is performed. Infact the collection of all super column matrix coefficient polynomials with same type of super column matrix coefficients is an abelian group under addition.

We shall illustrate this situation by some simple examples.

**Example 6.26:** Let $F_C^S[x] = \left\{ \sum_{i=0}^{\infty} \begin{bmatrix} a_1 \\ a_2 \\ \overline{a_3} \\ a_4 \\ \overline{a_5} \\ a_6 \\ a_7 \end{bmatrix} x^i \text{ with } \begin{bmatrix} a_1 \\ a_2 \\ \overline{a_3} \\ a_4 \\ \overline{a_5} \\ a_6 \\ a_7 \end{bmatrix} \in M \right.$

= {collection of all super column matrices of the same type with $a_i \in Z$, $1 \leq i \leq 7$}} be an abelian group of super column matrix coefficient polynomials in the variable x.

***Example 6.27:*** Let $F_C^S[x] = \left\{ \sum_{i=0}^{\infty} a_i x^i \;\middle|\; a_i = \begin{bmatrix} x_1 \\ x_2 \\ \hline x_3 \\ x_4 \\ \hline x_5 \\ x_6 \end{bmatrix} \in M = \left\{ \begin{bmatrix} x_1 \\ x_2 \\ x_3 \\ x_4 \\ x_5 \\ x_6 \end{bmatrix} \right. \right.$

with $x_j \in Z$, $1 \leq j \leq 6\}$, $+\}$ be an abelian group under addition.

***Example 6.28:*** Let

$$F_{3\times 5}^S = \left\{ \sum_{i=0}^{\infty} a_i x^i \;\middle|\; a_i = \begin{bmatrix} d_1 & d_2 & d_3 & d_4 & d_5 \\ \hline d_6 & d_7 & d_8 & d_9 & d_{10} \\ d_{11} & d_{12} & d_{13} & d_{14} & d_{15} \end{bmatrix} \in M = \right.$$

{all $3 \times 5$ matrices with entries from Z, $d_j \in Z$; $1 \leq j \leq 15$}} be an abelian group under addition.

***Example 6.29:*** Let $F_{4\times 4}^S = \left\{ \sum_{i=0}^{\infty} a_i x^i \;\middle|\; a_i = \begin{bmatrix} x_1 & x_2 & x_3 & x_4 \\ \hline x_5 & x_6 & x_7 & x_8 \\ x_9 & x_{10} & x_{11} & x_{12} \\ x_{13} & x_{14} & x_{15} & x_{16} \end{bmatrix} \right.$

with $x_j \in Q$; $1 \leq j \leq 16\}$ be an abelian group under addition. Clearly we see $F_{4\times 4}^S$ has subgroups.

$$P = \left\{ \sum_{i=0}^{\infty} a_i x^i \;\middle|\; a_i = \begin{bmatrix} p_1 & p_2 & p_3 & p_4 \\ \hline p_5 & p_6 & p_7 & p_8 \\ p_9 & p_{10} & p_{11} & p_{12} \\ p_{13} & p_{14} & p_{15} & p_{16} \end{bmatrix} \text{ with } p_j \in Z, \right.$$

$1 \leq j \leq 16\}$

is a subgroup of G under addition.

Now we proceed onto give the semigroup structure under the natural product $\times_n$.

Let $F_C^S = \left\{ \sum_{i=0}^{\infty} a_i x^i \;\middle|\; a_i = \begin{bmatrix} x_1 \\ \hline x_2 \\ \vdots \\ x_{20} \end{bmatrix} \text{ with } x_j \in Z, 1 \leq j \leq 20 \right\}$

be the collection of super column matrix coefficient polynomial. $F_C^S$ is a semigroup under the natural product $\times_n$. Infact $F_C^S$ is a commutative semigroup.

Suppose

$$p(x) = \begin{bmatrix} a_1 \\ \hline a_2 \\ \vdots \\ a_{20} \end{bmatrix} + \begin{bmatrix} b_1 \\ \hline b_2 \\ \vdots \\ b_{20} \end{bmatrix} x + \begin{bmatrix} c_1 \\ \hline c_2 \\ \vdots \\ c_{20} \end{bmatrix} x^3 \text{ and}$$

$$q(x) = \begin{bmatrix} d_1 \\ \hline d_2 \\ \vdots \\ d_{20} \end{bmatrix} + \begin{bmatrix} e_1 \\ \hline e_2 \\ \vdots \\ e_{20} \end{bmatrix} x^2 + \begin{bmatrix} m_1 \\ \hline m_2 \\ \vdots \\ m_{20} \end{bmatrix} x^4$$

are in $F_C^S$, then

$$p(x) \times_n q(x) = \begin{bmatrix} a_1d_1 \\ a_2d_2 \\ \vdots \\ a_{20}d_{20} \end{bmatrix} + \begin{bmatrix} b_1d_1 \\ b_2d_2 \\ \vdots \\ b_{20}d_{20} \end{bmatrix} x + \begin{bmatrix} a_1e_1 \\ a_2e_2 \\ \vdots \\ a_{20}e_{20} \end{bmatrix} x^2 + \begin{bmatrix} b_1e_1 \\ b_2e_2 \\ \vdots \\ b_{20}e_{20} \end{bmatrix} x^3$$

$$+ \begin{bmatrix} c_1d_1 \\ c_2d_2 \\ \vdots \\ c_{20}d_{20} \end{bmatrix} x^3 + \begin{bmatrix} a_1m_1 \\ a_2m_2 \\ \vdots \\ a_{20}m_{20} \end{bmatrix} x^4 + \begin{bmatrix} c_1e_1 \\ c_2e_2 \\ \vdots \\ c_{20}e_{20} \end{bmatrix} x^5 + \begin{bmatrix} c_1m_1 \\ c_2m_2 \\ \vdots \\ c_{20}m_{20} \end{bmatrix} x^7$$

$$+ \begin{bmatrix} b_1m_1 \\ b_2m_2 \\ \vdots \\ b_{20}m_{20} \end{bmatrix} x^5$$

$$= \begin{bmatrix} a_1d_1 \\ a_2d_2 \\ \vdots \\ a_{20}d_{20} \end{bmatrix} + \begin{bmatrix} b_1d_1 \\ b_2d_2 \\ \vdots \\ b_{20}d_{20} \end{bmatrix} x + \begin{bmatrix} a_1e_1 \\ a_2e_2 \\ \vdots \\ a_{20}e_{20} \end{bmatrix} x^2 + \begin{bmatrix} b_1e_1 + c_1d_1 \\ b_2e_2 + c_2d_2 \\ \vdots \\ b_me_m + c_md_m \end{bmatrix} x^3$$

$$+ \begin{bmatrix} a_1m_1 \\ a_2m_2 \\ \vdots \\ a_{20}m_{20} \end{bmatrix} x^4 + \begin{bmatrix} c_1e_1 + b_1m_1 \\ c_2e_2 + b_2m_2 \\ \vdots \\ c_{20}e_{20} + b_{20}m_{20} \end{bmatrix} x^5 + \begin{bmatrix} c_1m_1 \\ c_2m_2 \\ \vdots \\ c_{20}m_{20} \end{bmatrix} x^7 \text{ is in } F_C^S.$$

This way the natural product $\times_n$ is made on $F_C^S$.

We illustrate this situation by some examples.

*Example 6.30:* Let

$$F_{3\times 6}^S = \left\{\sum_{i=0}^{\infty} a_i x^i \;\middle|\; a_i = \left[\begin{array}{c|cc|ccc|c} x_1 & x_2 & x_3 & x_4 & x_5 & x_6 \\ \hline x_7 & x_8 & x_9 & x_{10} & x_{11} & x_{12} \\ x_{13} & x_{14} & x_{15} & x_{16} & x_{17} & x_{18} \end{array}\right]; x_j \in Z; 1 \leq j \leq 18\right\}$$

be a super $3 \times 6$ matrix coefficient polynomial semigroup under the natural product $\times_n$.

*Example 6.31:* Let

$$F_{3\times 3}^S = \left\{\sum_{i=0}^{\infty} a_i x^i \;\middle|\; a_i = \begin{bmatrix} x_1 & x_2 & x_3 \\ x_4 & x_5 & x_6 \\ \hline x_7 & x_8 & x_9 \end{bmatrix} \text{ where } x_j = Q; 1 \leq j \leq 9\right\}$$

be a suepr square matrix coefficient semigroup under natural product $\times_n$.

Let

$$p(x) = \begin{pmatrix} 3 & 2 & 0 \\ 1 & 0 & 1 \\ 0 & 2 & 3 \end{pmatrix} + \begin{pmatrix} 7 & 5 & 1 \\ 0 & 1 & 2 \\ 0 & 0 & 3 \end{pmatrix} x + \begin{pmatrix} 1 & 2 & 3 \\ 0 & 0 & 7 \\ 0 & 1 & 2 \end{pmatrix} x^2$$

$$+ \begin{pmatrix} 0 & 0 & 9 \\ 1 & 0 & 3 \\ 2 & 7 & 2 \end{pmatrix} x^4$$

and

$$q(x) = \begin{pmatrix} 4 & 0 & 2 \\ 1 & 5 & 6 \\ 7 & 0 & 2 \end{pmatrix} + \begin{pmatrix} 1 & 2 & 3 \\ 4 & 5 & 6 \\ 7 & 8 & 9 \end{pmatrix} x^2 + \begin{pmatrix} 0 & 3 & 1 \\ 2 & 1 & 0 \\ 3 & 4 & 5 \end{pmatrix} x^3$$

be in $F^S_{3\times 3}$.

To find

$$p(x) \times_n q(x); p(x) \times_n q(x) = \begin{pmatrix} 3 & 2 & | & 0 \\ 1 & 0 & | & 1 \\ \hline 0 & 2 & | & 3 \end{pmatrix} \times_n \begin{pmatrix} 4 & 0 & | & 2 \\ 1 & 5 & | & 6 \\ \hline 7 & 0 & | & 2 \end{pmatrix}$$

$$+ \begin{pmatrix} 3 & 2 & | & 0 \\ 1 & 0 & | & 1 \\ \hline 0 & 2 & | & 3 \end{pmatrix} \times_n \begin{pmatrix} 1 & 2 & | & 3 \\ 4 & 5 & | & 6 \\ \hline 7 & 8 & | & 9 \end{pmatrix} x^2 + \begin{pmatrix} 3 & 2 & | & 0 \\ 1 & 0 & | & 1 \\ \hline 0 & 2 & | & 3 \end{pmatrix} \times_n \begin{pmatrix} 0 & 3 & | & 1 \\ 2 & 1 & | & 0 \\ \hline 3 & 4 & | & 5 \end{pmatrix} x^3$$

$$+ \begin{pmatrix} 7 & 5 & | & 1 \\ 0 & 1 & | & 2 \\ \hline 0 & 0 & | & 3 \end{pmatrix} \times_n \begin{pmatrix} 4 & 0 & | & 2 \\ 1 & 5 & | & 6 \\ \hline 7 & 0 & | & 2 \end{pmatrix} x + \begin{pmatrix} 7 & 5 & | & 1 \\ 0 & 1 & | & 2 \\ \hline 0 & 0 & | & 3 \end{pmatrix} \times_n \begin{pmatrix} 1 & 2 & | & 3 \\ 4 & 5 & | & 6 \\ \hline 7 & 8 & | & 9 \end{pmatrix} x^3$$

$$+ \begin{pmatrix} 7 & 5 & | & 1 \\ 0 & 1 & | & 2 \\ \hline 0 & 0 & | & 3 \end{pmatrix} \times_n \begin{pmatrix} 0 & 3 & | & 1 \\ 2 & 1 & | & 0 \\ \hline 3 & 4 & | & 5 \end{pmatrix} x^4 + \begin{pmatrix} 1 & 2 & | & 3 \\ 0 & 0 & | & 7 \\ \hline 0 & 1 & | & 2 \end{pmatrix} \times_n \begin{pmatrix} 4 & 0 & | & 2 \\ 1 & 5 & | & 6 \\ \hline 7 & 0 & | & 2 \end{pmatrix} x^2$$

$$+ \begin{pmatrix} 1 & 2 & | & 3 \\ 0 & 0 & | & 7 \\ \hline 0 & 1 & | & 2 \end{pmatrix} \times_n \begin{pmatrix} 1 & 2 & | & 3 \\ 4 & 5 & | & 6 \\ \hline 7 & 8 & | & 9 \end{pmatrix} x^4 + \begin{pmatrix} 1 & 2 & | & 3 \\ 0 & 0 & | & 7 \\ \hline 0 & 1 & | & 2 \end{pmatrix} \times_n \begin{pmatrix} 0 & 3 & | & 1 \\ 2 & 1 & | & 0 \\ \hline 3 & 4 & | & 5 \end{pmatrix} x^5$$

$$+ \begin{pmatrix} 0 & 0 & | & 9 \\ 1 & 0 & | & 3 \\ \hline 2 & 7 & | & 2 \end{pmatrix} \times_n \begin{pmatrix} 4 & 0 & | & 2 \\ 1 & 5 & | & 6 \\ \hline 7 & 0 & | & 2 \end{pmatrix} x^4 + \begin{pmatrix} 0 & 0 & | & 9 \\ 1 & 0 & | & 3 \\ \hline 2 & 7 & | & 2 \end{pmatrix} \times_n \begin{pmatrix} 1 & 2 & | & 3 \\ 4 & 5 & | & 6 \\ \hline 7 & 8 & | & 9 \end{pmatrix} x^6$$

$$+ \begin{pmatrix} 0 & 0 & | & 9 \\ 1 & 0 & | & 3 \\ \hline 2 & 7 & | & 2 \end{pmatrix} \times_n \begin{pmatrix} 0 & 3 & | & 1 \\ 2 & 1 & | & 0 \\ \hline 3 & 4 & | & 5 \end{pmatrix} \times_n x^7$$

$$= \begin{pmatrix} 12 & 0 & 0 \\ 1 & 0 & 6 \\ \hline 0 & 0 & 6 \end{pmatrix} + \begin{pmatrix} 3 & 4 & 0 \\ 4 & 0 & 6 \\ \hline 0 & 16 & 27 \end{pmatrix} x^2 + \begin{pmatrix} 0 & 6 & 0 \\ 2 & 0 & 0 \\ \hline 0 & 8 & 15 \end{pmatrix} x^3$$

$$+ \begin{pmatrix} 28 & 0 & 2 \\ 0 & 5 & 12 \\ \hline 0 & 0 & 6 \end{pmatrix} x + \begin{pmatrix} 7 & 10 & 3 \\ 0 & 5 & 12 \\ \hline 0 & 0 & 27 \end{pmatrix} x^3 + \begin{pmatrix} 0 & 15 & 1 \\ 0 & 1 & 0 \\ \hline 0 & 0 & 15 \end{pmatrix} x^4$$

$$+ \begin{pmatrix} 1 & 4 & 9 \\ 0 & 0 & 42 \\ \hline 0 & 8 & 18 \end{pmatrix} x^4 + \begin{pmatrix} 0 & 6 & 3 \\ 0 & 0 & 0 \\ \hline 0 & 4 & 10 \end{pmatrix} x^5 + \begin{pmatrix} 0 & 0 & 18 \\ 1 & 0 & 18 \\ \hline 14 & 0 & 4 \end{pmatrix} x^4$$

$$+ \begin{pmatrix} 0 & 0 & 27 \\ 4 & 0 & 18 \\ \hline 14 & 56 & 18 \end{pmatrix} x^6 + \begin{pmatrix} 0 & 0 & 9 \\ 2 & 0 & 0 \\ \hline 6 & 28 & 10 \end{pmatrix} x^7 + \begin{pmatrix} 4 & 0 & 6 \\ 0 & 0 & 42 \\ \hline 0 & 0 & 4 \end{pmatrix} x^2$$

$$= \begin{pmatrix} 12 & 0 & 0 \\ 1 & 0 & 6 \\ \hline 0 & 0 & 6 \end{pmatrix} + \begin{pmatrix} 28 & 0 & 2 \\ 0 & 5 & 12 \\ \hline 0 & 0 & 6 \end{pmatrix} x + \begin{pmatrix} 7 & 4 & 6 \\ 4 & 0 & 48 \\ \hline 0 & 16 & 31 \end{pmatrix} x^2 +$$

$$\begin{pmatrix} 7 & 16 & 3 \\ 2 & 5 & 12 \\ \hline 0 & 8 & 42 \end{pmatrix} x^3 + \begin{pmatrix} 1 & 19 & 28 \\ 1 & 1 & 56 \\ \hline 14 & 8 & 37 \end{pmatrix} x^4 + \begin{pmatrix} 0 & 6 & 3 \\ 0 & 0 & 0 \\ \hline 0 & 4 & 10 \end{pmatrix} x^5$$

$$+ \begin{pmatrix} 0 & 0 & 27 \\ 4 & 0 & 18 \\ \hline 14 & 56 & 18 \end{pmatrix} x^6 + \begin{pmatrix} 0 & 0 & 9 \\ 2 & 0 & 0 \\ \hline 6 & 28 & 10 \end{pmatrix} x^7.$$

Thus we see $F^S_{3\times 3}[x]$ under natural product $\times_n$ is a semigroup. This semigroup has zero divisors and ideals. We can derive all related properties of this semigroup as a matter of routine.

Now we can also give these super matrix coefficient polynomials a ring structure. We just recall if $F^S_C[x]$ be a super column matrix polynomials, we know $F^S_C[x]$ under addition is an abelian group and under the natural product $\times_n$, $F^S_C[x]$ is a semigroup. Thus it is easily verified $(F^S_C[x], +, \times_n)$ is a commutative ring known as the super column matrix coefficient ring.

We will illustrate this situation by some simple examples.

*Example 6.32:* Let

$$F^S_C[x] = \left\{ \sum_{i=0}^{\infty} a_i x^i \;\middle|\; a_i = \begin{bmatrix} x_1 \\ x_2 \\ \hline x_3 \\ x_4 \\ \hline x_5 \\ x_6 \end{bmatrix} \text{ with } x_j \in Z, 1 \leq j \leq 6 \right\}$$

be the super column matrix coefficient polynomial ring under $+$ and $\times_n$.

**Example 6.33:** Let $(F_C^S[x], +, \times_n) = \left\{ \sum_{i=0}^{\infty} a_i x^i \;\middle|\; a_i = \begin{bmatrix} d_1 \\ \hline d_2 \\ \hline d_3 \\ \hline d_4 \\ \hline d_5 \\ \hline d_6 \\ \hline d_7 \\ \hline d_8 \end{bmatrix} \right.$ with

$d_j \in Q$; $1 \le j \le 8$, $+$, $\times_n\}$ be the super column matrix coefficient polynomial ring of infinite order. $F_C^S[x]$ has subrings which are not ideals, has zero divisors and idempotents only of a very special form which are only constant polynomials.

For instance $\alpha = \begin{bmatrix} 1 \\ \hline 0 \\ \hline 0 \\ \hline 1 \\ \hline 1 \\ \hline 1 \\ \hline 0 \\ \hline 1 \end{bmatrix} \in F_C^S[x]$ is such that $\alpha^2 = \alpha$.

Similarly $\alpha = \begin{bmatrix} 0 \\ \hline 1 \\ \hline 1 \\ \hline 1 \\ \hline 1 \\ \hline 1 \\ \hline 1 \\ \hline 0 \end{bmatrix} \in F_C^S[x]$ is such that $\alpha^2 = \alpha$.

Further we see $P = \left\{ \sum_{i=0}^{\infty} a_i x^i \mid a_i = \begin{bmatrix} x_1 \\ x_2 \\ x_3 \\ \hline x_4 \\ x_5 \\ \hline x_6 \\ \hline x_7 \\ x_8 \end{bmatrix} \right.$ with $x_j \in 3Z$,

$1 \leq j \leq 8, +, \times_n\} \subseteq F_C^S[x]$ is a subring of super column matrix coefficient polynomial ring. However P is not an ideal of $F_C^S[x]$.

However $F_C^S[x]$ has infinite number of zero divisors.

***Example 6.34:*** Let $F_C^S[x] = \left\{ \sum_{i=0}^{\infty} a_i x^i \mid a_i = \begin{bmatrix} d_1 \\ \hline d_2 \\ d_3 \end{bmatrix} \right.$ with $d_j \in Z_n$,

$\times_n, +, 1 \leq j \leq 3\}$ be a super column matrix coefficient polynomial ring.

Consider $p(x) = \begin{bmatrix} 0 \\ \hline a_1 \\ a_2 \end{bmatrix} + \begin{bmatrix} 0 \\ \hline b_1 \\ b_2 \end{bmatrix} x + \begin{bmatrix} 0 \\ \hline c_1 \\ c_2 \end{bmatrix} x^2,$

$q(x) = \begin{bmatrix} a_1 \\ \hline 0 \\ 0 \end{bmatrix} + \begin{bmatrix} b_1 \\ \hline 0 \\ 0 \end{bmatrix} x + \begin{bmatrix} d_1 \\ \hline 0 \\ 0 \end{bmatrix} x^3 + \begin{bmatrix} x_1 \\ \hline 0 \\ 0 \end{bmatrix} x^4 \in F_C^S[x].$

Clearly $p(x) \times_n q(x) = \begin{bmatrix} 0 \\ \hline 0 \\ \hline 0 \end{bmatrix}$.

Further if

$$P = \left\{ \sum_{i=0}^{\infty} a_i x^i \;\middle|\; a_i = \begin{bmatrix} x_1 \\ \hline 0 \\ \hline x_2 \end{bmatrix} \text{ with } x_1, x_2 \in Z, +, \times_n \right\} \subseteq F_C^S[x]$$

is a subring.

$$\text{Also } T = \left\{ \sum_{i=0}^{\infty} a_i x^i \;\middle|\; a_i = \begin{bmatrix} 0 \\ \hline y_1 \\ \hline 0 \end{bmatrix}, y_1 \in Z, +, \times_n \right\} \subseteq F_C^S[x]$$

is a subring.

We see $F_C^S[x] =$

$$P + T \text{ and } P \cap T = \begin{bmatrix} 0 \\ \hline 0 \\ \hline 0 \end{bmatrix}.$$

Further for every $\alpha \in P$ we have every $\beta \in T$ is such that

$$\alpha \times_n \beta = \begin{bmatrix} 0 \\ \hline 0 \\ \hline 0 \end{bmatrix}.$$

We see both P and T are also ideals of $F_C^S$.

Now we proceed onto define super row matrix coefficient polynomial ring. Consider $\{F_R^S[x], +, \times_n\}$ is a super row matrix coefficient polynomial ring which is commutative and of infinite order.

We will give examples of it.

***Example 6.35:*** Let $R = \{ F_R^S [x] = \sum_{i=0}^{\infty} a_i x^i ; a_i = (t_1\ t_2\ t_3\ |\ t_4\ |\ t_5\ t_6);$ $t_j \in Z, 1 \leq j \leq 6, +, \times_n\}$ be a super row matrix polynomial coefficient ring. R has zero divisors, units, ideals and subrings.

$$P = \left\{ \sum_{i=0}^{\infty} a_i x^i \,\middle|\, a_i = (0\ 0\ 0\ |\ d_1\ |\ d_2\ d_3)\ d_1, d_2, d_3 \in Z, +, \times_n \right\} \subseteq R$$

is a subring as well as an ideal of R.

$$T = \left\{ \sum_{i=0}^{\infty} a_i x^i \,\middle|\, a_i = (y_1\ y_2\ y_3\ |\ 0\ |\ 0\ \ 0)\ y_1, y_2, y_3 \in Z, +, \times_n \right\} \subseteq R$$

is a subring as well as an ideal of R.

We see every $a \in P$ and every $b \in T$ are such that $a \times_n b = (0\ 0\ 0\ |\ 0\ |\ 0\ 0)$.

Also $(1\ -1\ 1\ |\ -1\ |\ -1\ 1) = p$ is such that $p^2 = (1\ 1\ 1\ |\ 1\ |\ 1\ 1)$, only units of this form are in R.

***Example 6.36:*** Let $F_R^S [x] = \left\{ \sum_{i=0}^{\infty} a_i x^i \,\middle|\, a_i = (p_1\ p_2\ |\ p_3\ |\ p_4\ p_5\ |\ p_6 \right.$ $|\ p_7\ p_8\ |\ p_9)$ with $p_j \in Q; 1 \leq j \leq 9, +, \times_n\}$ be the super row matrix coefficient ring of polynomials. Clearly $F_R^S [x]$ has ideals, subrings which are not ideals and zero divisors.

For take $M = \left\{ \sum_{i=0}^{\infty} a_i x^i \,\middle|\, a_i = (d_4\ 0\ |\ 0\ |\ 0\ 0\ |\ d_5\ |\ d_1\ d_2\ |\ d_3) \right.$ with $d_j \in Q; 1 \leq j \leq 5, +, \times_n\} \subseteq F_R^S [x]$; M is an ideal. Consider

$$T = \left\{\sum_{i=0}^{\infty} a_i x^i \;\middle|\; a_i = (y_1\; y_2 \mid y_3 \mid y_4\; y_5 \mid y_6 \mid y_7\; y_8\; y_9) \text{ with } y_j \in Z; 1 \leq j \leq 9, +, \times_n \right\} \subseteq F_R^S[x]$$; T is a only a subring and not an ideal of $F_R^S[x]$.

It is easily verified $F_R^S[x]$ has zero divisors.

***Example 6.37:*** Let $F_R^S[x] = \left\{\sum_{i=0}^{\infty} a_i x^i \;\middle|\; a_i = (m_1 \mid m_2); m_1, m_2 \in R, +, \times_n \right\}$ be a super row matrix coefficient polynomial ring. Clearly $F_R^S[x]$ has units, zero divisors, idempotents all of them are only constant polynomials. For $\alpha = (1 \mid -1)$ is a unit as $\alpha^2 = (1 \mid 1)$ and $\beta = (0 \mid 1)$ and $b_1 = (1 \mid 0)$ are all idempotents. We also see $P = \left\{\sum_{i=0}^{\infty} a_i x^i \;\middle|\; a_i = (t \mid s), t, s, \in Z, +, \times_n \right\} \subseteq F_R^S[x]$ is a subring and not an ideal of $F_R^S[x]$.

Take $M = \left\{\sum_{i=0}^{\infty} a_i x^i \;\middle|\; a_i = (t \mid 0); t \in R, +, \times_n \right\} \subseteq F_R^S[x]$ is an ideal of $F_R^S[x]$. $N = \left\{\sum_{i=0}^{\infty} a_i x^i \;\middle|\; a_i = (0 \mid s)\; s \in R, +, \times_n \right\} \subseteq R_R^S[x]$ is also an ideal of $F_R^S[x]$. Every $\alpha$ in M and every $\beta$ in M are such that $\alpha \times_n \beta = (0 \mid 0)$.

Next we proceed onto describe the concept of super matrix coefficient polynomial rings.

Let $F_{m \times n}^S[x] = \left\{\sum_{i=0}^{\infty} a_i x^i \;\middle|\; a_i = (m_{ij}), m \times n \text{ super matrices of} \right.$ same type of every $a_i$ and $1 \leq i \leq m$ and $1 \leq j \leq n$ with $m \neq n$ and $m_{ij} \in Z$ (or Q or R), $+, \times_n\}$ be the super matrix coefficient polynomial ring.

We will illustrate this by some examples.

***Example 6.38:*** Let $\{F^S_{2\times 4}[x], +, \times_n\} = \left\{\sum_{i=0}^{\infty} a_i x^i \mid a_i = \begin{bmatrix} d_1 & d_2 & d_3 & d_4 \\ d_5 & d_6 & d_7 & d_8 \end{bmatrix} d_j \in z, 1 \leq j \leq 8, +, \times_n \right\}$ be the super row vector coefficient polynomial ring. $F^S_{2\times 4}$ has zero divisors, units, idempotents ideals and subrings.

$\alpha = \begin{bmatrix} 1 & 1 & -1 & 1 \\ -1 & 1 & -1 & -1 \end{bmatrix}$ is such that $\alpha^2 = \begin{bmatrix} 1 & 1 & 1 & 1 \\ 1 & 1 & 1 & 1 \end{bmatrix}$ is a unit.

Consider $\beta = \begin{bmatrix} 1 & 0 & 1 & 1 \\ 0 & 1 & 0 & 0 \end{bmatrix}$ in $F^S_{2\times 4}[x]$. We see $\beta^2 = \beta$ so $\beta$ is an idempotent in $F^S_{2\times 4}[x]$ only if it is a constant polynomial and the super row vector takes its entries only as 0 or 1.

Similarly an element $a \in F^S_{2\times 4}[x]$ is a unit only if a is a constant polynomial and all its entries are from the set $\{-1, 1\}$. Consider

$$P = \left\{\sum_{i=0}^{\infty} a_i x^i \mid a_i = \begin{bmatrix} t_1 & t_2 & t_3 & t_4 \\ t_5 & t_6 & t_7 & t_8 \end{bmatrix}; t_j \in 5Z, \right.$$
$$\left. 1 \leq j \leq 8, +, \times_n \right\}$$

is a subring which is also an ideal of $F^S_R[x]$.

**Example 6.39:** Let

$$\{F^S_{5\times 3}[x], +, \times_n\} = \left\{ \sum_{i=0}^{\infty} a_i x^i \;\middle|\; a_i = \begin{bmatrix} y_1 & y_2 & y_3 \\ y_4 & y_5 & y_6 \\ \hline y_7 & y_8 & y_9 \\ \hline y_{10} & y_{11} & y_{12} \\ y_{13} & y_{14} & y_{15} \end{bmatrix} \text{ with } y_j \in Q, \right.$$

$$\left. 1 \leq j \leq 15, \times_n, + \right\}$$

be the super column vector coefficient polynomial ring. $F^S_{5\times 3}[x]$ has subrings which are not ideals, ideals, units, zero divisors and idempotents.

Consider $M = \left\{ \sum_{i=0}^{\infty} a_i x^i \;\middle|\; a_i = \begin{bmatrix} m_1 & m_2 & m_3 \\ m_4 & m_5 & m_6 \\ \hline m_7 & m_8 & m_9 \\ \hline m_{10} & m_{11} & m_{12} \\ m_{13} & m_{14} & m_{15} \end{bmatrix} \text{ with } m_j \in \right.$

$\left. z, 1 \leq j \leq 15, +, \times_n \right\} \subseteq F^S_{5\times 3}[x]$

is only a subring and not an ideal of $F^S_{5\times 3}[x]$.

Take $p = \begin{bmatrix} 1 & 0 & 1 \\ 1 & 1 & 1 \\ \hline 0 & 0 & 1 \\ \hline 1 & 1 & 1 \\ 0 & 0 & 1 \end{bmatrix}$ in $F^S_{5\times 3}[x]$ is such that $p^2 = p$, that is an idempotent. All idempotents are only constant polynomials that is $5 \times 3$ super column vectors.

Consider $t = \begin{bmatrix} 1 & -1 & 1 \\ \hline 1 & 1 & -1 \\ \hline -1 & 1 & 1 \\ \hline 1 & -1 & -1 \\ \hline -1 & -1 & -1 \end{bmatrix}$ in $F_{5\times 3}^S[x]$.

We see $t^2 = \begin{bmatrix} 1 & 1 & 1 \\ \hline 1 & 1 & 1 \\ \hline 1 & 1 & 1 \\ \hline 1 & 1 & 1 \\ \hline 1 & 1 & 1 \end{bmatrix}$ the unit; that t is a unit.

Suppose $y \in F_{5\times 3}^S[x]$ is to be a unit then we see it should be a constant polynomial and the $5 \times 3$ matrix must take its entries from the set $\{1, -1\}$.

*Example 6.40:* Let

$$F_{5\times 4}^S[x] = \left\{ \sum_{i=0}^{\infty} a_i x^i \;\middle|\; a_i = \begin{bmatrix} m_1 & m_6 & m_{11} & m_{16} \\ m_2 & m_7 & m_{12} & m_{17} \\ m_3 & m_8 & m_{13} & m_{18} \\ m_4 & m_9 & m_{14} & m_{19} \\ m_5 & m_{10} & m_{15} & m_{20} \end{bmatrix} \right.$$

with $m_j \in Z$; $1 \leq j \leq 20$; $+, \times_n\}$

be the super matrix coefficient polynomial ring.

Take $M = \left\{ \sum_{i=0}^{\infty} a_i x^i \mid a_i = \begin{bmatrix} m_1 & m_6 & m_{11} & m_{16} \\ m_2 & m_7 & m_{12} & m_{17} \\ m_3 & m_8 & m_{13} & m_{18} \\ m_4 & m_9 & m_{14} & m_{19} \\ m_5 & m_{10} & m_{15} & m_{20} \end{bmatrix} \right.$

with $m_j \in 5Z$, $1 \leq j \leq 20$, $+, \times_n\} \subseteq F_{5\times 4}^S [x]$;
M is an ideal of $F_{5\times 4}^S [x]$.
Suppose

$P = \left\{ \sum_{i=0}^{\infty} a_i x^i \mid a_i = \begin{bmatrix} m_1 & 0 & 0 & m_9 \\ 0 & m_3 & m_4 & 0 \\ 0 & m_5 & m_6 & 0 \\ 0 & m_7 & m_8 & 0 \\ m_2 & 0 & 0 & m_{10} \end{bmatrix} \right.$

with $m_j \in 3Z$, $1 \leq j \leq 10$, $+, \times_n\} \subseteq F_{5\times 4}^S [x]$

be an ideal of $F_{5\times 4}^S [x]$.

Take $T = \left\{ \sum_{i=0}^{\infty} a_i x^i \mid a_i = \begin{bmatrix} 0 & m_4 & m_5 & 0 \\ m_1 & 0 & 0 & m_7 \\ m_2 & 0 & 0 & m_8 \\ m_3 & 0 & 0 & m_9 \\ 0 & m_{10} & m_6 & 0 \end{bmatrix} \; m_j \in 13Z; \right.$

$1 \leq j \leq 10, +, \times_n\} \subseteq F_{5\times 4}^S [x]$

be an ideal, we see for every $\alpha$ in P is such that for every $\beta$ in

$$T \alpha \times_n \beta = \begin{bmatrix} 0 & 0 & 0 & 0 \\ \hline 0 & 0 & 0 & 0 \\ \hline 0 & 0 & 0 & 0 \\ 0 & 0 & 0 & 0 \\ \hline 0 & 0 & 0 & 0 \end{bmatrix}.$$

However $P \cap T = \begin{bmatrix} 0 & 0 & 0 & 0 \\ \hline 0 & 0 & 0 & 0 \\ \hline 0 & 0 & 0 & 0 \\ 0 & 0 & 0 & 0 \\ \hline 0 & 0 & 0 & 0 \end{bmatrix}$ but $P \cup T \neq F_{3 \times 4}^S [x]$.

*Example 6.41:* Let

$$F_{5\times 7}^S [x] = \left\{ \sum_{i=0}^{\infty} a_i x^i \;\middle|\; a_i = \begin{bmatrix} m_1 & m_2 & m_3 & m_4 & m_5 & m_6 & m_7 \\ m_8 & m_9 & m_{10} & m_{11} & m_{12} & m_{13} & m_{14} \\ m_{15} & m_{16} & m_{17} & m_{18} & m_{19} & m_{20} & m_{21} \\ m_{22} & m_{23} & m_{24} & m_{25} & m_{26} & m_{27} & m_{28} \\ m_{29} & m_{30} & m_{31} & m_{32} & m_{33} & m_{34} & m_{35} \end{bmatrix} \text{ with } m_j \in Q; \right.$$

$$1 \leq j \leq 35, +, \times_n \}$$

be the super matrix coefficient polynomial ring $F_{5\times 7}^S [x]$ has zero divisors, units, idempotents, ideals and subrings which are not ideals.

$$m = \begin{bmatrix} 1 & 1 & 1 & 1 & -1 & 1 & -1 \\ \hline -1 & -1 & 1 & 1 & 1 & -1 & 1 \\ -1 & 1 & -1 & 1 & -1 & -1 & -1 \\ \hline 1 & -1 & 1 & -1 & 1 & 1 & -1 \\ 1 & 1 & -1 & 1 & 1 & -1 & 1 \end{bmatrix} \text{ in } F_{5\times 7}^S[x]$$

is such that

$$m^2 = \begin{bmatrix} 1 & 1 & 1 & 1 & 1 & 1 & 1 \\ 1 & 1 & 1 & 1 & 1 & 1 & 1 \\ 1 & 1 & 1 & 1 & 1 & 1 & 1 \\ 1 & 1 & 1 & 1 & 1 & 1 & 1 \\ 1 & 1 & 1 & 1 & 1 & 1 & 1 \end{bmatrix} \in F_{5\times 7}^S[x].$$

Further if $p = \begin{bmatrix} 1 & 0 & 0 & 1 & 0 & 1 & 1 \\ 0 & 1 & 1 & 0 & 1 & 0 & 0 \\ 0 & 1 & 1 & 0 & 0 & 1 & 1 \\ 1 & 0 & 1 & 0 & 1 & 0 & 1 \\ 1 & 1 & 1 & 1 & 0 & 1 & 0 \end{bmatrix} \in F_{5\times 7}^S[x],$

then $p^2 = p$ is an idempotent of $F_{5\times 7}^S[x]$.

Take

$$P = \left\{ \sum_{i=0}^{\infty} a_i x^i \ \Bigg| \ a_i = \begin{bmatrix} m_1 & 0 & 0 & m_6 & 0 & m_{11} & 0 \\ m_2 & 0 & 0 & m_7 & 0 & m_{12} & 0 \\ m_3 & 0 & 0 & m_8 & 0 & m_{13} & 0 \\ m_4 & 0 & 0 & m_9 & 0 & m_{14} & 0 \\ m_5 & 0 & 0 & m_{10} & 0 & m_{15} & 0 \end{bmatrix} \right. \text{ with}$$

$m_j \in Z; 1 \leq j \leq 15, +, \times_n\}$

to be a subring and not an ideal of $F_{5 \times 7}^S [x]$.

Take

$$S = \left\{ \sum_{i=0}^{\infty} a_i x^i \ \Bigg| \ a_i = \begin{bmatrix} 0 & m_1 & m_2 & 0 & m_{11} & 0 & m_{16} \\ 0 & m_3 & m_4 & 0 & m_{12} & 0 & m_{17} \\ 0 & m_5 & m_6 & 0 & m_{13} & 0 & m_{18} \\ 0 & m_7 & m_8 & 0 & m_{14} & 0 & m_{19} \\ 0 & m_9 & m_{10} & 0 & m_{15} & 0 & m_{20} \end{bmatrix} \right.$$

with $m_j \in Q; 1 \leq j \leq 20, +, \times_n\}$ is an ideal of $F_{5 \times 7}^S [x]$.

We see every polynomial $p(x) \in P$ is such that for every

$$q(x) \in S \text{ we have } p(x) \times_n q(x) = \begin{bmatrix} 0 & 0 & 0 & 0 & 0 & 0 & 0 \\ 0 & 0 & 0 & 0 & 0 & 0 & 0 \\ 0 & 0 & 0 & 0 & 0 & 0 & 0 \\ 0 & 0 & 0 & 0 & 0 & 0 & 0 \\ 0 & 0 & 0 & 0 & 0 & 0 & 0 \end{bmatrix}.$$

Now we describe the super square matrix coefficient polynomial ring.

Let $F^S_{m \times m}[x]$ be the collection of a super square matrix coefficient polynomial ring under $+$ and $\times_n$. We will illustrate this by some examples.

*Example 6.42:* Let

$$F^S_{4 \times 4}[x] = \left\{ \sum_{i=0}^{\infty} a_i x^i \mid a_i = \begin{bmatrix} m_1 & m_2 & m_3 & m_4 \\ \hline m_5 & m_6 & m_7 & m_8 \\ m_9 & m_{10} & m_{11} & m_{12} \\ \hline m_{13} & m_{14} & m_{15} & m_{16} \end{bmatrix} \right. ;$$

$$\left. m_j \in Z, 1 \le j \le 16, +, \times_n \right\}$$

be the super square matrix coefficient polynomial ring. $F^S_{4 \times 4}[x]$ has zero divisors units, idempotents, ideals and subrings.

Take $\alpha = \begin{bmatrix} 1 & -1 & 1 & 1 \\ \hline -1 & 1 & 1 & 0 \\ -1 & 0 & 1 & 0 \\ \hline 0 & 1 & 1 & 1 \end{bmatrix} \in F^S_{4 \times 4}[x]$, we see $\alpha^2 \ne \alpha$ and $\alpha$

is not a unit.

Take $\beta = \begin{bmatrix} -1 & 1 & 1 & 1 \\ \hline -1 & 1 & -1 & 1 \\ \hline 1 & 1 & 1 & -1 \\ 1 & 1 & -1 & 1 \end{bmatrix} \in F^S_{4 \times 4}[x];$

$\beta^2 = \begin{bmatrix} 1 & 1 & 1 & 1 \\ \hline 1 & 1 & 1 & 1 \\ 1 & 1 & 1 & 1 \\ \hline 1 & 1 & 1 & 1 \end{bmatrix} .$

Take $p = \begin{bmatrix} 1 & 0 & 1 & 1 \\ 0 & 1 & 0 & 1 \\ 1 & 0 & 1 & 0 \\ 1 & 1 & 0 & 0 \end{bmatrix}$ in $F_{4\times 4}^S[x]$; we see $p^2 = p$ is an idempotent.

Now $S = \left\{ \sum_{i=0}^{\infty} a_i x^i \;\middle|\; a_i = \begin{bmatrix} d_1 & 0 & d_2 & 0 \\ 0 & d_3 & 0 & d_4 \\ d_5 & 0 & d_6 & 0 \\ 0 & d_7 & 0 & d_8 \end{bmatrix}, d_j \in Z, \right.$

$1 \le j \le 8 \} \subseteq F_{4\times 4}^S[x]$ is an ideal of $F_{4\times 4}^S[x]$.

$M = \left\{ \sum_{i=0}^{\infty} a_i x^i \;\middle|\; a_i = \begin{bmatrix} 0 & d_1 & 0 & d_2 \\ d_8 & 0 & d_3 & 0 \\ 0 & d_4 & 0 & d_5 \\ d_6 & 0 & d_7 & 0 \end{bmatrix} \right.$ with $d_j \in Z$; $1 \le j \le$

$8, +, \times_n \} \subseteq F_{4\times 4}^S[x]$ is also an ideal of $F_{4\times 4}^S[x]$. We see every polynomial $p(x)$ is $S$ is such that for every polynomial $q(x)$ in $M$ we have $p(x) \times_n q(x) = \begin{bmatrix} 0 & 0 & 0 & 0 \\ 0 & 0 & 0 & 0 \\ 0 & 0 & 0 & 0 \\ 0 & 0 & 0 & 0 \end{bmatrix}$.

Now we have the following theorems, the proofs of which are left as an exercise to the reader.

**THEOREM 6.1:** *Let $F_C^S[x]$ (or $F_R^S[x]$ or $F_{m\times n}^S[x]$ ($m \ne n$) or $F_{n\times n}^S[x]$) be super matrix coefficient polynomial ring.*

1. *Every constant polynomial with entries from the set {1, –1} is a unit under $\times_n$.*
2. *Every constant polynomial with entries from the set {0, 1} is an idempotent under $\times_n$.*

**THEOREM 6.2:** *Every super matrix coefficient polynomial ring has ideals.*

**THEOREM 6.3:** *Every super matrix coefficient polynomial ring over Q or R has subrings which are not ideals.*

**THEOREM 6.4:** *Every super matrix coefficient polynomial ring has infinite number of zero divisors.*

Now we can define super vector space of polynomials over R or Q.

Suppose we take $V = F_C^S[x]$ to be an abelian group under addition. V is a vector space over the reals or rationals.

We will give examples of them.

*Example 6.43:* Let

$$V = \left\{ F_C^S[x] = \sum_{i=0}^{\infty} a_i x^i \text{ with } a_i = \begin{bmatrix} t_1 \\ t_2 \\ t_3 \\ t_4 \\ t_5 \end{bmatrix} ; t_j \in Q, 1 \leq j \leq 5, + \right\}$$

be a vector space over Q. V is called the super column matrix coefficient polynomial vector space or super polynomial vector space over Q.

*Example 6.44:* Let

$$V = \left\{ F_C^S[x] = \sum_{i=0}^{\infty} a_i x^i \text{ with } a_i = \begin{bmatrix} t_1 \\ \hline t_2 \\ t_3 \end{bmatrix}; t_j \in Q, 1 \leq j \leq 3, + \right\}$$

be an abelian group under '+'. V is a super column matrix coefficient polynomial vector space over Q.

$$P_1 = \left\{ \sum_{i=0}^{\infty} a_i x^i \,\Big|\, a_i = \begin{bmatrix} t_1 \\ \hline 0 \\ 0 \end{bmatrix}; t_1 \in Q, + \right\} \subseteq V$$

is a subspace of V over Q.

$$P_2 = \left\{ \sum_{i=0}^{\infty} a_i x^i \,\Big|\, a_i = \begin{bmatrix} 0 \\ \hline t_2 \\ 0 \end{bmatrix}; t_2 \in Q, + \right\} \subseteq V$$

is a subspace of V over Q.

$$P_3 = \left\{ \sum_{i=0}^{\infty} a_i x^i \,\Big|\, a_i = \begin{bmatrix} 0 \\ \hline 0 \\ t_3 \end{bmatrix}; t_3 \in Q, + \right\} \subseteq V$$

is a subspace of V over Q.

We see $V = P_1 + P_2 + P_3$ and

$$P_i \cap P_j = \begin{bmatrix} 0 \\ \hline 0 \\ 0 \end{bmatrix} \text{ if } i \neq j; 1 \leq i, j \leq 3.$$

Thus V is a direct sum of subspaces $P_1$, $P_2$ and $P_3$ over Q.

*Example 6.45:* Let

$$M = \left\{ \sum_{i=0}^{\infty} a_i x^i \;\middle|\; a_i = \begin{bmatrix} m_1 \\ \hline m_2 \\ m_3 \\ m_4 \\ \hline m_5 \\ \hline m_6 \\ m_7 \end{bmatrix} ; m_j \in Q; 1 \leq j \leq 7, + \right\}$$

be a super column matrix coefficient polynomial group under '+'.

M is a super column matrix coefficient vector space over Q.

Take

$$P_1 = \left\{ \sum_{i=0}^{\infty} a_i x^i \;\middle|\; a_i = \begin{bmatrix} m_1 \\ \hline m_2 \\ 0 \\ 0 \\ \hline 0 \\ \hline 0 \\ m_3 \end{bmatrix}, m_1, m_2, m_3 \in Q, + \right\} \subseteq M,$$

$$P_2 = \left\{ \sum_{i=0}^{\infty} a_i x^i \;\middle|\; a_i = \begin{bmatrix} m_1 \\ \hline 0 \\ m_2 \\ m_3 \\ \hline 0 \\ \hline 0 \\ m_4 \end{bmatrix} \text{ with } m_j \in Q, 1 \leq j \leq 4, + \right\} \subseteq M$$

and

$$P_3 = \left\{ \sum_{i=0}^{\infty} a_i x^i \,\middle|\, a_i = \begin{bmatrix} m_1 \\ \hline 0 \\ 0 \\ 0 \\ \hline m_2 \\ \hline m_3 \\ \hline m_4 \end{bmatrix} \text{ with } m_j \in Q,\ 1 \le j \le 4 \right\} \subseteq M.$$

Clearly $M \subseteq P_1 + P_2 + P_3$ but $P_i \cap P_j \ne \begin{bmatrix} 0 \\ \hline 0 \\ 0 \\ 0 \\ \hline 0 \\ \hline 0 \\ \hline 0 \end{bmatrix}$ if $i \ne j$,

$1 \le i, j \le 3$. $P_1$, $P_2$ and $P_3$ are subspaces of M over Q. Clearly M is the pseudo direct sum of vector subspaces of M.

*Example 6.46:* Let

$$V = \left\{ \sum_{i=0}^{\infty} a_i x^i \,\middle|\, a_i = (d_1 \mid d_2 \mid d_3\ d_4\ d_5\ d_6\ d_7 \mid d_8);\ d_j \in Q, \ 1 \le j \le 8 \right\}$$

be a super row matrix coefficient polynomial vector space over Q.

Consider

$$P_1 = \left\{ \sum_{i=0}^{\infty} a_i x^i \,\middle|\, a_i = (0 \mid 0 \mid d_1, d_2, 0, 0, 0 \mid 0),\ d_1, d_2 \in Q \right\} \subseteq V,$$

$$P_2 = \left\{\sum_{i=0}^{\infty} a_i x^i \;\middle|\; a_i = (d_1 \mid d_2 \mid 0, 0, 0, 0, 0 \mid 0),\, d_1, d_2 \in Q\right\} \subseteq V,$$

$$P_3 = \left\{\sum_{i=0}^{\infty} a_i x^i \;\middle|\; a_i = (0 \mid 0 \mid 0\ 0\ d_1\ d_2\ 0 \mid 0);\, d_1, d_2 \in Q\right\} \subseteq V$$

and

$$P_4 = \left\{\sum_{i=0}^{\infty} a_i x^i \;\middle|\; a_i = (0 \mid 0 \mid 0 \ldots 0\ d_1 \mid d_2) \text{ with } d_1, d_2 \in Q\right\} \subseteq V$$

be a super row matrix coefficient polynomials subvector space of V over Q.

Clearly $V = P_1 + P_2 + P_3 + P_4$ and

$$P_i \cap P_j = (\,0 \mid 0 \mid 0\ 0\ 0\ 0\ 0 \mid 0) \text{ if } i \neq j,\, 1 \leq i, j \leq 4.$$

Clearly V is a direct sum of subspaces.

*Example 6.47:* Let

$$V = \left\{\sum_{i=0}^{\infty} a_i x^i \;\middle|\; a_i = (t_1\ t_2 \mid t_3\ t_4 \mid t_5);\, t_j \in Q,\, 1 \leq j \leq 5\right\}$$

be the super row matrix coefficient polynomial vector space over Q.

Take

$$M_1 = \left\{\sum_{i=0}^{\infty} a_i x^i \;\middle|\; a_i = (0\ t_1 \mid t_2\ 0 \mid 0),\, t_1, t_2 \in Q\right\} \subseteq V,$$

$$M_2 = \left\{ \sum_{i=0}^{\infty} a_i x^i \mid a_i = (t_1, t_2 \mid 0\ 0 \mid 0),\ t_1, t_2 \in Q \right\} \subseteq V,$$

$$M_3 = \left\{ \sum_{i=0}^{\infty} a_i x^i \mid a_i = (0\ t_1 \mid 0\ t_2 \mid 0),\ t_1, t_2 \in Q \right\} \subseteq V$$

and

$$M_4 = \left\{ \sum_{i=0}^{\infty} a_i x^i \mid a_i = (0\ t_1 \mid 0\ 0 \mid t_2),\ t_1, t_2 \in Q \right\} \subseteq V$$

be super row matrix coefficient polynomial vector subspace of V over Q. We see $V \subseteq M_1 + M_2 + M_3 + M_4$, however
$M_i \cap M_j \neq (0\ 0 \mid 0\ 0 \mid 0)$ if $i \neq j$; $1 \leq i, j \leq 4$.

Thus V is only a pseudo direct sum of vector subspaces over Q.

*Example 6.48:* Let

$$V = \left\{ \sum_{i=0}^{\infty} a_i x^i \mid a_i = (m_1 \mid m_2\ m_3\ m_4 \mid m_5\ m_6 \mid m_7) \right.$$
$$\text{with } m_j \in R;\ 1 \leq j \leq 7 \}$$

be a super row matrix coefficient polynomial vector space over Q.

Consider

$$M_1 = \left\{ \sum_{i=0}^{\infty} a_i x^i \mid a_i = (0 \mid m_1\ 0\ m_2 \mid 0\ m_3 \mid 0) \right.$$
$$\text{with } m_1, m_2, m_3 \in R \} \subseteq V$$

and

$$M_2 = \left\{ \sum_{i=0}^{\infty} a_i x^i \mid a_i = (m_1 \mid 0\ m_2\ 0 \mid m_3\ 0 \mid m_4) \right.$$
$$\text{with } m_1, m_2, m_3, m_4 \in R \} \subseteq V$$

to be super row matrix coefficient polynomial vector subspaces of V over Q.

We see $V = M_1 + M_2$ and $M_1 \cap M_2 = (\,0\,|\,0\,0\,0\,|\,0\,0\,|\,0)$.

Infact we see every $q(x) \in M_1$ is orthogonal with every other $p(x) \in M_2$.

Thus $M_1$ is the orthogonal complement of $M_2$ and vice versa. Now we give examples of super m × n matrix coefficient polynomial vector spaces over Q or R.

*Example 6.49:* Let

$$V = \left\{ \sum_{i=0}^{\infty} a_i x^i \,\middle|\, a_i = \begin{bmatrix} d_1 & d_4 & d_7 & d_{10} & d_{13} & d_{16} \\ d_2 & d_5 & d_8 & d_{11} & d_{14} & d_{17} \\ d_3 & d_6 & d_9 & d_{12} & d_{15} & d_{18} \end{bmatrix} \text{ with } d_j \in Q, \right.$$

$$1 \leq j \leq 18\}$$

be the super 3 × 7 row vector coefficient polynomial vector space over Q.

*Example 6.50:* Let

$$W = \left\{ \sum_{i=0}^{\infty} a_i x^i \;\middle|\; a_i = \begin{bmatrix} m_1 & m_2 & m_3 & m_4 \\ m_5 & m_6 & m_7 & m_8 \\ m_9 & m_{10} & m_{11} & m_{12} \\ \hline m_{13} & m_{14} & m_{15} & m_{16} \\ m_{17} & m_{18} & m_{19} & m_{20} \\ m_{21} & m_{22} & m_{23} & m_{24} \\ \hline m_{25} & m_{26} & m_{27} & m_{28} \\ \hline m_{29} & m_{30} & m_{31} & m_{32} \\ \hline m_{33} & m_{34} & m_{35} & m_{36} \\ \hline m_{37} & m_{38} & m_{39} & m_{40} \end{bmatrix} \right.$$

with $m_j \in R$, $1 \le j \le 40\}$

be a super column vector coefficient polynomial vector space over Q.

*Example 6.51:* Let

$$V = \left\{ \sum_{i=0}^{\infty} a_i x^i \;\middle|\; a_i = \begin{bmatrix} m_1 & m_2 & m_3 & m_4 & m_5 \\ \hline m_6 & m_7 & m_8 & m_9 & m_{10} \\ m_{11} & m_{12} & m_{13} & m_{14} & m_{15} \\ m_{16} & m_{17} & m_{18} & m_{19} & m_{20} \\ m_{21} & m_{22} & m_{23} & m_{24} & m_{25} \\ m_{26} & m_{27} & m_{28} & m_{29} & m_{30} \\ m_{31} & m_{32} & m_{33} & m_{34} & m_{35} \end{bmatrix} \right.$$

with $m_j \in Q$, $1 \le j \le 35\}$

be a super $7 \times 5$ matrix coefficient polynomial vector space over Q.

Now we see for all these one can easily find a basis subspaces etc.

*Example 6.52:* Let

$$V = \left\{ \sum_{i=0}^{\infty} a_i x^i \;\middle|\; a_i = \right.$$

$$\begin{bmatrix} m_1 & m_8 & m_{15} & m_{22} & m_{29} & m_{36} & m_{43} & m_{50} \\ m_2 & m_9 & m_{16} & m_{23} & m_{30} & m_{37} & m_{44} & m_{51} \\ m_3 & m_{10} & m_{17} & m_{24} & m_{31} & m_{38} & m_{45} & m_{52} \\ m_4 & m_{11} & m_{18} & m_{25} & m_{32} & m_{39} & m_{46} & m_{53} \\ m_5 & m_{12} & m_{19} & m_{26} & m_{33} & m_{40} & m_{47} & m_{54} \\ m_6 & m_{13} & m_{20} & m_{27} & m_{34} & m_{41} & m_{48} & m_{55} \\ m_7 & m_{14} & m_{21} & m_{28} & m_{35} & m_{42} & m_{49} & m_{56} \end{bmatrix}$$ with

$m_j \in Q, 1 \leq j \leq 56\}$

be a super $7 \times 8$ matrix coefficient polynomial vector space over Q. Consider

$$M_1 = \left\{ \sum_{i=0}^{\infty} a_i x^i \;\middle|\; a_i = \begin{bmatrix} m_1 & 0 & 0 & 0 & 0 & 0 & 0 & 0 \\ m_2 & 0 & 0 & 0 & 0 & 0 & 0 & 0 \\ m_3 & 0 & 0 & 0 & 0 & 0 & 0 & 0 \\ m_4 & 0 & 0 & 0 & 0 & 0 & 0 & 0 \\ m_5 & 0 & 0 & 0 & 0 & 0 & 0 & 0 \\ m_6 & 0 & 0 & 0 & 0 & 0 & 0 & 0 \\ m_7 & 0 & 0 & 0 & 0 & 0 & 0 & 0 \end{bmatrix} \right.$$

with $m_j \in Q, 1 \leq j \leq 7\} \subseteq V$,

$$M_2 = \left\{ \sum_{i=0}^{\infty} a_i x^i \,\middle|\, a_i = \left[\begin{array}{c|cc|ccc|cc} 0 & m_1 & m_2 & 0 & 0 & 0 & 0 & 0 \\ \hline 0 & m_3 & m_4 & 0 & 0 & 0 & 0 & 0 \\ 0 & m_5 & m_6 & 0 & 0 & 0 & 0 & 0 \\ 0 & m_7 & m_8 & 0 & 0 & 0 & 0 & 0 \\ 0 & m_9 & m_{10} & 0 & 0 & 0 & 0 & 0 \\ 0 & m_{11} & m_{12} & 0 & 0 & 0 & 0 & 0 \\ \hline 0 & m_{13} & m_{14} & 0 & 0 & 0 & 0 & 0 \end{array}\right] \right.$$

with $m_j \in Q$, $1 \leq j \leq 14\} \subseteq V$,

$$M_3 = \left\{ \sum_{i=0}^{\infty} a_i x^i \,\middle|\, a_i = \left[\begin{array}{c|cc|ccc|cc} 0 & 0 & 0 & m_1 & m_2 & m_3 & 0 & 0 \\ \hline 0 & 0 & 0 & m_4 & m_5 & m_6 & 0 & 0 \\ 0 & 0 & 0 & m_7 & m_8 & m_9 & 0 & 0 \\ 0 & 0 & 0 & m_{10} & m_{11} & m_{12} & 0 & 0 \\ 0 & 0 & 0 & m_{13} & m_{14} & m_{15} & 0 & 0 \\ 0 & 0 & 0 & m_{16} & m_{17} & m_{18} & 0 & 0 \\ \hline 0 & 0 & 0 & m_{19} & m_{20} & m_{21} & 0 & 0 \end{array}\right] \right.$$

with $m_j \in Q$, $1 \leq j \leq 21\} \subseteq V$

and

$$M_4 = \left\{ \sum_{i=0}^{\infty} a_i x^i \,\middle|\, a_i = \left[\begin{array}{c|cc|ccc|cc} 0 & 0 & 0 & 0 & 0 & 0 & m_1 & m_2 \\ \hline 0 & 0 & 0 & 0 & 0 & 0 & m_3 & m_4 \\ 0 & 0 & 0 & 0 & 0 & 0 & m_5 & m_6 \\ 0 & 0 & 0 & 0 & 0 & 0 & m_7 & m_8 \\ 0 & 0 & 0 & 0 & 0 & 0 & m_9 & m_{10} \\ 0 & 0 & 0 & 0 & 0 & 0 & m_{11} & m_{12} \\ \hline 0 & 0 & 0 & 0 & 0 & 0 & m_{13} & m_{14} \end{array}\right] \right.$$

with $m_j \in Q$, $1 \leq j \leq 14\} \subseteq V$.

Clearly $M_1$, $M_2$, $M_3$ and $M_4$ are super matrix coefficient polynomial vector subspaces of V and $M_1 + M_2 + M_3 + M_4$ and

$$M_i \cap M_j = \begin{bmatrix} 0 & 0 & 0 & 0 & 0 & 0 & 0 & 0 \\ 0 & 0 & 0 & 0 & 0 & 0 & 0 & 0 \\ 0 & 0 & 0 & 0 & 0 & 0 & 0 & 0 \\ \hline 0 & 0 & 0 & 0 & 0 & 0 & 0 & 0 \\ 0 & 0 & 0 & 0 & 0 & 0 & 0 & 0 \\ 0 & 0 & 0 & 0 & 0 & 0 & 0 & 0 \\ 0 & 0 & 0 & 0 & 0 & 0 & 0 & 0 \end{bmatrix}, 1 \leq i, j \leq 4.$$

Thus V is the direct sum of subspaces.

*Example 6.53:* Let

$$P = \left\{ \sum_{i=0}^{\infty} a_i x^i \;\middle|\; a_i = \begin{bmatrix} m_1 & m_2 & m_3 & m_4 & m_5 & m_6 & m_7 \\ m_8 & m_9 & m_{10} & m_{11} & m_{12} & m_{13} & m_{14} \\ m_{15} & m_{16} & m_{17} & m_{18} & m_{19} & m_{20} & m_{21} \\ m_{22} & m_{23} & m_{24} & m_{25} & m_{26} & m_{27} & m_{28} \end{bmatrix} \right.$$

with $m_i \in R$, $1 \leq j \leq 28$ }

be a 4 × 7 super row vector matrix coefficient polynomial vector space over Q.

$$\text{Let } B_1 = \left\{ \sum_{i=0}^{\infty} a_i x^i \;\middle|\; a_i = \begin{bmatrix} m_1 & m_5 & 0 & 0 & 0 & 0 & 0 \\ m_2 & m_6 & 0 & 0 & 0 & 0 & 0 \\ m_3 & m_7 & 0 & 0 & 0 & 0 & 0 \\ m_4 & m_8 & 0 & 0 & 0 & 0 & 0 \end{bmatrix} \right.;$$

$m_i \in Q$, $1 \leq i \leq 8$} $\subseteq P$,

$$B_2 = \left\{ \sum_{i=0}^{\infty} a_i x^i \;\middle|\; a_i = \begin{bmatrix} m_1 & 0 & m_5 & 0 & 0 & 0 & 0 \\ m_2 & 0 & m_6 & 0 & 0 & 0 & 0 \\ m_3 & 0 & m_7 & 0 & 0 & 0 & 0 \\ m_4 & 0 & m_8 & 0 & 0 & 0 & 0 \end{bmatrix} \right. ;$$

$$m_i \in Q, 1 \leq i \leq 8 \} \subseteq P,$$

$$B_3 = \left\{ \sum_{i=0}^{\infty} a_i x^i \;\middle|\; a_i = \begin{bmatrix} m_1 & 0 & 0 & m_5 & 0 & 0 & 0 \\ m_2 & 0 & 0 & m_6 & 0 & 0 & 0 \\ m_3 & 0 & 0 & m_7 & 0 & 0 & 0 \\ m_4 & 0 & 0 & m_8 & 0 & 0 & 0 \end{bmatrix} \right. ;$$

$$m_i \in Q, 1 \leq i \leq 8 \} \subseteq P,$$

$$B_4 = \left\{ \sum_{i=0}^{\infty} a_i x^i \;\middle|\; a_i = \begin{bmatrix} m_1 & 0 & 0 & 0 & m_5 & 0 & 0 \\ m_2 & 0 & 0 & 0 & m_6 & 0 & 0 \\ m_3 & 0 & 0 & 0 & m_7 & 0 & 0 \\ m_4 & 0 & 0 & 0 & m_8 & 0 & 0 \end{bmatrix} \right. ;$$

$$m_i \in Q, 1 \leq i \leq 8 \} \subseteq P,$$

$$B_5 = \left\{ \sum_{i=0}^{\infty} a_i x^i \;\middle|\; a_i = \begin{bmatrix} m_1 & 0 & 0 & 0 & 0 & m_5 & 0 \\ m_2 & 0 & 0 & 0 & 0 & m_6 & 0 \\ m_3 & 0 & 0 & 0 & 0 & m_7 & 0 \\ m_4 & 0 & 0 & 0 & 0 & m_8 & 0 \end{bmatrix} \right. ;$$

$$m_i \in Q, 1 \leq i \leq 8 \} \subseteq P$$

and

$$B_6 = \left\{ \sum_{i=0}^{\infty} a_i x^i \,\middle|\, a_i = \begin{bmatrix} m_1 & 0 & 0 & 0 & 0 & 0 & m_5 \\ m_2 & 0 & 0 & 0 & 0 & 0 & m_6 \\ m_3 & 0 & 0 & 0 & 0 & 0 & m_7 \\ m_4 & 0 & 0 & 0 & 0 & 0 & m_8 \end{bmatrix}; \right.$$

$$\left. m_i \in Q, 1 \leq i \leq 8 \right\} \subseteq P$$

be super matrix coefficient polynomial subspaces of P over the field Q.

Clearly $B_i \cap B_j \neq 0$, if $i \neq j$, $1 \leq i, j \leq 6$.

Also $P \subseteq B_1 + B_2 + \ldots + B_6$ so P is a pseudo direct sum of vector subspaces $B_1, B_2, \ldots, B_6$.

We can define orthogonal subspaces and orthogonal complements also.

*Example 6.54:* Let

$$V = \left\{ \sum_{i=0}^{\infty} a_i x^i \,\middle|\, a_i = \begin{bmatrix} m_1 \\ \hline m_2 \\ \hline m_3 \\ m_4 \\ \hline m_5 \\ \hline m_6 \\ m_7 \\ \hline m_8 \end{bmatrix} \text{ with } m_i \in Q, 1 \leq i \leq 8 \right\}$$

be a super column matrix coefficient polynomial vector space over the field Q.

Consider

$$M_1 = \left\{ \sum_{i=0}^{\infty} a_i x^i \ \middle| \ a_i = \begin{bmatrix} m_1 \\ \hline m_2 \\ \hline 0 \\ \hline 0 \\ \hline 0 \\ \hline 0 \\ \hline 0 \\ \hline 0 \end{bmatrix} \text{ with } m_1, m_2 \in Q \right\} \subseteq V$$

is a vector subspace of V over Q.

$$M_2 = \left\{ \sum_{i=0}^{\infty} a_i x^i \ \middle| \ a_i = \begin{bmatrix} 0 \\ \hline 0 \\ \hline m_1 \\ \hline m_2 \\ \hline m_3 \\ \hline m_4 \\ \hline m_5 \\ \hline m_6 \end{bmatrix} \text{ with } m_i \in Q; 1 \leq i \leq 6 \right\} \subseteq V$$

is a vector subspace of V over Q. Clearly $M_1^{\perp} = M_2$ and vice versa.

Consider

$$M_3 = \left\{ \sum_{i=0}^{\infty} a_i x^i \;\middle|\; a_i = \begin{bmatrix} 0 \\ \hline 0 \\ \hline m_1 \\ m_2 \\ m_3 \\ \hline 0 \\ \hline 0 \\ \hline 0 \end{bmatrix}, m_1, m_2, m_3 \in Q \right\} \subseteq V,$$

$M_3$ is also a vector subspace of V over Q and for every $x \in M_1$ and for every $y \in M_3$ we see $x \times_n y = (0)$ however $M_1^\perp \neq M_3$ for $M_1 + M_3 \neq V$ however $M_1 + M_2 = V$ and $M_1^\perp = M_2$ and $M_2^\perp = M_1$.

Thus we see we can have subspaces in V orthogonal to $M_1$ but they need not be the orthogonal complement of M, in V over Q.

Now we proceed onto define semivector space of super matrix coefficient polynomials defined over the semifield $Z^+ \cup \{0\}$ or $Q^+ \cup \{0\}$ or $R^+ \cup \{0\}$.

We just describe them in the following.

$$\text{Let } P = \left\{ \sum_{i=0}^{\infty} a_i x^i \;\middle|\; a_i = \begin{bmatrix} d_1 \\ \hline d_2 \\ d_3 \\ d_4 \\ d_5 \end{bmatrix} \text{ with } d_j \in Z^+ \cup \{0\}; 1 \leq j \leq 5 \right\}$$

be a semigroup under addition known as the super column matrix coefficient polynomial semigroup. P is a semivector space over the semifield $S = Z^+ \cup \{0\}$.

$$W = \left\{ \sum_{i=0}^{\infty} a_i x^i \;\middle|\; a_i = (m_1 \mid m_2 \mid m_3 \; m_4) \text{ with } m_j \in Q^+ \cup \{0\}, \right.$$

$1 \leq j \leq 4\}$ is a super column matrix coefficient polynomial semivector space over the semifield $S = Q^+ \cup \{0\}$ (or $Z^+ \cup \{0\}$).

$$T = \left\{ \sum a_i x^i \;\middle|\; a_i = \begin{bmatrix} d_1 & d_2 & d_3 \\ \hline d_4 & d_5 & d_6 \\ d_7 & d_8 & d_9 \\ \hline d_{10} & d_{11} & d_{12} \\ d_{13} & d_{14} & d_{15} \\ d_{16} & d_{17} & d_{18} \\ \hline d_{19} & d_{20} & d_{21} \\ \hline d_{22} & d_{23} & d_{24} \\ \hline d_{25} & d_{26} & d_{27} \end{bmatrix} \text{ with } d_j \in R^+ \cup \{0\}, \right.$$

$1 \leq j \leq 27\}$ is a super column vector coefficient polynomial semivector space over the semifield $S = Z^+ \cup \{0\}$ (or $Q^+ \cup \{0\}$ or $R^+ \cup \{0\}$).

$$M = \left\{ \sum_{i=0}^{\infty} a_i x^i \;\middle|\; a_i = \begin{bmatrix} t_1 & t_2 & t_3 & t_4 & t_5 & t_6 & t_7 \\ t_8 & t_9 & t_{10} & t_{11} & t_{12} & t_{13} & t_{14} \\ t_{15} & t_{16} & t_{17} & t_{18} & t_{19} & t_{20} & t_{21} \end{bmatrix} \; t_i \in \right.$$

$Q^+ \cup \{0\}, 1 \leq i \leq 21\}$ is the super row vector coefficient polynomial semivector space over the semifield $S = Q^+ \cup \{0\}$ or $Z^+ \cup \{0\}$. However M is not a semivector space over $R^+ \cup \{0\}$.

$$B = \left\{ \sum_{i=0}^{\infty} a_i x^i \;\middle|\; a_i = \begin{bmatrix} y_1 & y_2 & y_3 & y_4 & y_5 & y_6 & y_7 & y_8 \\ y_9 & y_{10} & y_{11} & y_{12} & y_{13} & y_{14} & y_{15} & y_{16} \\ y_{17} & y_{18} & y_{19} & y_{20} & y_{21} & y_{22} & y_{23} & y_{24} \\ y_{25} & y_{26} & y_{27} & y_{28} & y_{29} & y_{30} & y_{31} & y_{32} \end{bmatrix} \right.$$

$$y_i \in Z^+ \cup \{0\}, 1 \leq i \leq 32\}$$

is a super 4 × 8 matrix coefficient polynomial semivector space over the semifield $Z^+ \cup \{0\}$.

$$C = \left\{ \sum_{i=0}^{\infty} a_i x^i \;\middle|\; a_i = \begin{bmatrix} m_1 & m_2 & m_3 & m_4 & m_5 & m_6 \\ m_7 & m_8 & m_9 & m_{10} & m_{11} & m_{12} \\ m_{13} & m_{14} & m_{15} & m_{16} & m_{17} & m_{18} \\ m_{19} & m_{20} & m_{21} & m_{22} & m_{23} & m_{24} \\ m_{25} & m_{26} & m_{27} & m_{28} & m_{29} & m_{30} \\ m_{31} & m_{32} & m_{33} & m_{34} & m_{35} & m_{36} \end{bmatrix} \right.$$

$$\left. m_j \in Z^+ \cup \{0\},\; 1 \leq i \leq 36 \right\}$$

is a super square matrix coefficient polynomial semivector space over the semifield $Z^+ \cup \{0\}$.

The authors by examples show how subsemivector space direct sum etc looks like.

*Example 6.55:* Let

$$M = \left\{ \sum_{i=0}^{\infty} a_i x^i \;\middle|\; a_i = \right.$$

$$\begin{bmatrix} m_1 & m_2 & m_3 & m_4 & m_5 & m_6 & m_7 & m_8 \\ m_9 & m_{10} & m_{11} & m_{12} & m_{13} & m_{14} & m_{15} & m_{16} \\ m_{17} & m_{18} & m_{19} & m_{20} & m_{21} & m_{22} & m_{23} & m_{24} \end{bmatrix}$$

$$\left. m_j \in Q^+ \cup \{0\},\; 1 \leq i \leq 24 \right\}$$

be a super row vector polynomial coefficient semivector space defined over the semifield $Z^+ \cup \{0\}$. M is a semilinear algebra under the natural product.

Consider

$$p(x) = \begin{bmatrix} 8 & 0 & 1 & 7 & 0 & 3 & 8 & 1 \\ 0 & 0 & 2 & 0 & 8 & 1 & 0 & 1 \\ 1 & 5 & 0 & 1 & 0 & 0 & 0 & 1 \end{bmatrix} +$$

$$\begin{bmatrix} 6 & 7 & 0 & 1 & 6 & 0 & 9 & 1 \\ 0 & 9 & 1 & 2 & 0 & 1 & 6 & 2 \\ 2 & 0 & 2 & 8 & 1 & 1 & 7 & 3 \end{bmatrix} x +$$

$$\begin{bmatrix} 1 & 0 & 1 & 8 & 1 & 2 & 0 & 7 \\ 1 & 2 & 4 & 0 & 3 & 2 & 0 & 0 \\ 2 & 0 & 3 & 1 & 1 & 0 & 0 & 2 \end{bmatrix} x^2$$

and

$$q(x) = \begin{bmatrix} 0 & 1 & 2 & 1 & 0 & 5 & 7 & 2 \\ 2 & 0 & 1 & 3 & 7 & 0 & 5 & 1 \\ 1 & 4 & 8 & 1 & 0 & 0 & 4 & 0 \end{bmatrix} +$$

$$\begin{bmatrix} 2 & 8 & 1 & 3 & 4 & 3 & 3 & 1 \\ 1 & 0 & 5 & 0 & 1 & 2 & 2 & 2 \\ 0 & 1 & 4 & 6 & 0 & 0 & 1 & 3 \end{bmatrix} x^2 +$$

$$\begin{bmatrix} 0 & 0 & 1 & 3 & 0 & 0 & 6 & 2 \\ 7 & 2 & 0 & 0 & 1 & 0 & 0 & 1 \\ 0 & 0 & 0 & 0 & 0 & 1 & 1 & 0 \end{bmatrix} x^3 \text{ be in M. } p(x).q(x) \in M.$$

$$T = \left\{ \sum_{i=0}^{\infty} a_i x^i \;\middle|\; a_i = \begin{bmatrix} 0 & m_1 & m_2 & 0 & 0 & 0 & m_7 & m_{10} \\ 0 & m_3 & m_4 & 0 & 0 & 0 & m_8 & m_{11} \\ 0 & m_5 & m_6 & 0 & 0 & 0 & m_9 & m_{12} \end{bmatrix} \right.$$
$$\left. m_i \in Q^+ \cup \{0\},\ 1 \leq i \leq 12 \right\} \subseteq M$$

and

$$P = \left\{ \sum_{i=0}^{\infty} a_i x^i \;\middle|\; a_i = \begin{bmatrix} m_1 & 0 & 0 & m_4 & m_5 & m_6 & 0 & 0 \\ m_2 & 0 & 0 & m_7 & m_8 & m_9 & 0 & 0 \\ m_3 & 0 & 0 & m_{10} & m_{11} & m_{12} & 0 & 0 \end{bmatrix} \right.$$
$$\left. m_i \in Q^+ \cup \{0\}, 1 \le i \le 12 \right\} \subseteq M,$$

T and P are super row vector polynomial coefficient semivector subspace of M over the semifield $Q^+ \cup \{0\}$;

we see $M = T + P$ with $T \cap P = \begin{bmatrix} 0 & 0 & 0 & 0 & 0 & 0 & 0 & 0 \\ 0 & 0 & 0 & 0 & 0 & 0 & 0 & 0 \\ 0 & 0 & 0 & 0 & 0 & 0 & 0 & 0 \end{bmatrix}$.

Further for every $x \in T$ we have a $y \in P$;

$$\text{with } x \times_n y = \begin{bmatrix} 0 & 0 & 0 & 0 & 0 & 0 & 0 & 0 \\ 0 & 0 & 0 & 0 & 0 & 0 & 0 & 0 \\ 0 & 0 & 0 & 0 & 0 & 0 & 0 & 0 \end{bmatrix}.$$

*Example 6.56:* Let

$$W = \left\{ \sum_{i=0}^{\infty} a_i x^i \;\middle|\; a_i = \begin{bmatrix} d_1 & d_2 & d_3 & d_4 \\ d_5 & d_6 & d_7 & d_8 \\ d_9 & d_{10} & d_{11} & d_{12} \\ d_{13} & d_{14} & d_{15} & d_{16} \\ d_{17} & d_{18} & d_{19} & d_{20} \\ d_{21} & d_{22} & d_{23} & d_{24} \\ d_{25} & d_{26} & d_{27} & d_{28} \\ d_{29} & d_{30} & d_{31} & d_{32} \end{bmatrix} \; d_j \in R^+ \cup \{0\}, \right.$$
$$\left. 1 \le i \le 32 \right\}$$

be a super column vector polynomial coefficient semivector space (linear algebra) over the semifield $S = Z^+ \cup \{0\}$.

Now we can define like wise semivector spaces of super matrices and study those structures.

**Chapter Seven**

# APPLICATIONS OF THESE ALGEBRAIC STRUCTURES WITH NATURAL PRODUCT

We define natural product on matrices and that results in the compatability in column matrices and m × n (m ≠ n) matrices.

Several algebraic structures using them are developed.

The notion of super matrix coefficient polynomials are defined and described. These new structures will certainly find nice and appropriate applications in due course of time. Of course it is a very difficult thing to define products on super matrices in the way defined in [8, 19]; however because of this new 'natural product' we can define product on super matrices provided they are of the same type.

One has to find the uses of these new structures in eigen value problems, mathematical models, coding theory and finite element analysis methods.

Further these natural product on matrices works like the usual product on the real line and the matrix product of row matrices. If the concept of matrices is a an array of number than certainty the natural product seems to be appropriate so in due course of time researchers will find nice applications of them.

**Chapter Eight**

# SUGGESTED PROBLEMS

In this chapter we suggest over 100 problems. Some of them can be taken up as research problems. These problems however makes the reader to understand these new notions introduced in this book.

1. Find some interesting properties enjoyed by polynomials with matrix coefficients.

2. For the row matrix coefficient polynomial semigroup

   $$S[x] = \left\{ \sum_{i=0}^{\infty} a_i x^i \,\middle|\, a_i = (x_1, x_2, x_3, x_4) \text{ with } x_j \in R,\ 1 \leq j \leq 4 \right\}$$

   (i) Find zero divisors in $S[x]$.
   (ii) Can $S[x]$ have ideals?
   (iii) Can $S[x]$ have subrings which are not ideals?
   (iv) Can $S[x]$ have idempotents?

3. Let $p(x) = (3, 2, 1, 5) + (-2, 0, 1, 3)x + (7, 8, 4, -6)x^2$ be a row matrix coefficient polynomial in the variable x.
   (i) Find roots of $p(x)$.
   (ii) If $\alpha$ and $\beta$ are roots of $p(x)$ find a row matrix coefficient polynomial whose roots are $\alpha^2 + \beta^2$ and $\alpha^2 \beta^2$.

(iii) If the row matrix coefficients are from Z will α and β be in Z × Z × Z × Z?

4. Give some nice properties enjoyed by the semigroup of square matrix coefficient polynomials.

5. Solve the equation $\begin{bmatrix} 1 & 1 \\ 1 & 1 \end{bmatrix} x^3 - \begin{bmatrix} 2 & 7 & 8 \\ 1 & 6 & 4 \end{bmatrix} = 0$.

6. Suppose $p(x) = \begin{bmatrix} 9 \\ 2 \\ 1 \\ 3 \\ 7 \end{bmatrix} - \begin{bmatrix} 8 \\ 9 \\ 2 \\ 3 \\ -1 \end{bmatrix} x + \begin{bmatrix} 3 \\ 7 \\ 0 \\ 8 \\ 1 \end{bmatrix} x^2$. Solve for x.

7. Let $p(x) = (1, 2, 3)x^3 - (2, 4, 5)x^2 + (1, 0, 2)x - (3, 8, 1)$.
Solve for x.
Does $q(x) = (1, 6, 9)x + (2, 1, 3)$ divide $p(x)$?

8. Find the properties enjoyed by the group of square matrix coefficient polynomials in the variable x.

9. Let $p(x) = \begin{pmatrix} 2 & 1 & 0 \\ 1 & 5 & 6 \end{pmatrix} + \begin{pmatrix} 7 & 8 & 0 \\ 1 & 1 & 1 \end{pmatrix} x^2 + \begin{pmatrix} 3 & 1 & 2 \\ 0 & 1 & 5 \end{pmatrix} x$.

Is p(x) solvable as a quadratic equation?

10. Let $p(x) = \begin{bmatrix} 7 \\ 8 \\ 9 \\ 3 \end{bmatrix} + \begin{bmatrix} 1 \\ 2 \\ 3 \\ 4 \end{bmatrix} x^2 + \begin{bmatrix} 0 \\ 2 \\ 5 \\ 7 \end{bmatrix} x^3$. Find the derivatives,

the coefficients are from Z.
Can p(x) be integrated with respect to x? Justify.

11. Suppose $p(x) = \sum_{i=0}^{9} a_i x^i$ where $a_i \in \{V_{3\times 6} = \{\text{all } 3 \times 6$ matrices with coefficients from Z$\}$, $0 \le i \le 9\}$; prove $p(x)$ cannot be integrated and the resultant coefficients will not be in $V_{3\times 6}$.

12. Prove $V_R = \{(a_1, a_2, \ldots, a_{12}) \mid a_i \in R\}$ has zero divisors under product.

13. Prove $V_{3\times 3} = \left\{ \begin{pmatrix} a_1 & a_2 & a_3 \\ a_4 & a_5 & a_6 \\ a_7 & a_8 & a_9 \end{pmatrix} \middle| a_i \in Z, 1 \le i \le 9 \right\}$ is a semigroup under multiplication.
   (i) Is $V_{3\times 3}$ a commutative semigroup?
   (ii) Find ideals in $V_{3\times 3}$.
   (iii) Can $V_{3\times 3}$ have subsemigroups which are not ideals?

14. Can $V_R[x] = \left\{ \sum_{i=0}^{\infty} a_i x^i \middle| a_i = (x_1, x_2, \ldots, x_8) \text{ with } x_j \in Z, 1 \le j \le 8 \right\}$, a semigroup under product have zero divisors?
   (i) Find ideals of $V_R[x]$?
   (ii) Find subsemigroups which are not ideals in $V_R[x]$.
   (iii) Can $V_R[x]$ be a S-semigroup?

15. Let $V_{7\times 7}[x] = \left\{ \sum_{i=0}^{\infty} a_i x^i \middle| a_i\text{'s are } 7 \times 7 \text{ matrices with entries from } R \right\}$ be a semigroup under the matrix multiplication.
   (i) Is $V_{7\times 7}[x]$ a commutative semigroup?
   (ii) Find right ideals in $V_{7\times 7}[x]$ which are not left ideals and vice versa.
   (iii) Find two sided ideals of $V_{7\times 7}[x]$.

16. Distinguish between $R[x]$ and $V_{3\times 3}[x]$.

17. Distinguish between Q [x] and

$$V_R [x] = \left\{ \sum_{i=0}^{\infty} a_i x^i \,\middle|\, a_i = (x_1, x_2, x_3, x_4); x_j \in Q; 1 \le j \le 4 \right\}.$$

18. What are the benefits of natural product in matrices?

19. Prove $V_{n \times n} [x]$ under natural product is a commutative semigroup.

20. Prove $V_{3 \times 7} [x]$ is a commutative semigroup under natural product $\times_n$.

21. Prove natural product $\times_n$ and the usual product of row matrices on $V_R [x]$ are identical.

22. Show $V_C [x]$ under natural product is a semigroup with zero divisors.

23. Obtain some nice properties enjoyed by

$$V_C = \left\{ \begin{bmatrix} a_1 \\ a_2 \\ \vdots \\ a_m \end{bmatrix} \,\middle|\, a_i \in Q; 1 \le i \le m \right\} \text{ under the natural}$$

product, $\times_n$.

24. Show $V_{5 \times 2} = \{\text{all } 5 \times 2 \text{ matrices with entries from } Q\}$ under natural product $\times_n$ is a semigroup.

   (i) Find zero divisors in $V_{5 \times 2}$.

   (ii) Show all elements of $V_{5 \times 2}$ are not invertible in general.

   (iii) Show $V_{5 \times 2}$ has subsemigroups which are not ideals.

   (iv) Find ideals of $V_{5 \times 2}$.

(v) Can $V_{5\times 2}$ have idempotents justify?

25. Show $(V_{3\times 3}, \times_n)$ and $(V_{3\times 3}, \times)$ are distinct as semigroups.

   (i) Can they be isomorphic?
   (ii) Find any other stricking difference between them.

26. Can the set of $5 \times 5$ diagonal matrices with entries from Q under the natural product and the usual product be same?

27. Prove $(V_{2\times 2}, +, \times_n)$ is a commutative ring.

28. Prove $(V_{3\times 3}, +, \times)$ is a non commutative ring with

$$\begin{pmatrix} 1 & 0 & 0 \\ 0 & 1 & 0 \\ 0 & 0 & 1 \end{pmatrix} \text{ as unit.}$$

29. Find the differences between a ring of matrices under natural product and usual matrix product.

30. Prove $(V_{2\times 7}, +, \times_n)$ is a commutative ring with identity.

31. Let $S = (V_{5\times 2}, +, \times_n)$ be a ring.
   (i) Find subrings of S.

   (ii) Is S a Smarandache ring?

   (iii) Can S have S-subrings?

   (iv) Can S have subrings which are not S-ring?

   (v) Find ideals in S.

   (vi) Find subrings in S which are not S-ideals.

   (vii) Find zero divisors in S.

32. Find some special properties enjoyed by $(V_C, +, \times_n)$.

33. Distinguish between $(V_R, +, \times_n)$ and $(V_R, +, \times)$.

34. Find the difference between the rings $(V_{3\times3}, +, \times)$ and $(V_{3\times3}, +, \times_n)$.

35. Let $M = (V_R^+, +, \times_n)$ be a semiring where
    $V_R^+ = \{(x_1, x_2, \ldots, x_n) \mid x_i \in R^+ \cup \{0\}, 1 \leq i \leq n\}$.

    (i) Is M a semifield?
    (ii) Is M a S-semiring?
    (iii) Find subsemiring in M.
    (iv) Show every subsemiring need not a be S-subsemiring.
    (v) Find zero divisors in M.
    (vi) Can M have idempotents?

36. Let $P = \{V_C = \left\{\begin{bmatrix} a_1 \\ a_2 \\ a_3 \\ a_4 \\ a_5 \end{bmatrix}\right\} \mid a_i \in Q^+ \cup \{0\}; 1 \leq i \leq 5\}$ be a semiring under $+$ and $\times_n$.

    (i) Find ideals of P.
    (ii) Is P a S-semiring?
    (iii) Can P have S-subsemiring?
    (iv) Find S-ideals if any in P.
    (v) Find zero divisors in P.

37. Obtain some special properties enjoyed by the semiring of $7 \times 1$ column matrices with entries from $R^+ \cup \{0\}$.

38. Mention some of the special features enjoyed by the semiring of $5 \times 8$ matrices with $+$ and $\times_n$; the entries are from $R^+ \cup \{0\}$.

39. Is $P = \left\{ \begin{bmatrix} x_1 \\ x_2 \\ x_3 \\ x_4 \\ x_5 \end{bmatrix} \middle| x_i \in R^+; 1 \leq i \leq 5 \right\} \cup \left\{ \begin{bmatrix} 0 \\ 0 \\ 0 \\ 0 \\ 0 \end{bmatrix} \right\}$ a semifield under $+$ and $\times_n$?

40. Can $M = \left\{ \begin{bmatrix} a & b & c & d & e \\ a_1 & b_1 & c_1 & d_1 & e_1 \end{bmatrix} \middle| a_1, b_1, c_1, d_1, e_1, a, b, c, d, e \in Q^+ \right\} \cup \left\{ \begin{bmatrix} 0 & 0 & 0 & 0 & 0 \\ 0 & 0 & 0 & 0 & 0 \end{bmatrix} \right\}, +, \times_n\}$ be a semifield?

41. Find for $S = \left\{ \begin{bmatrix} a_1 & a_2 \\ a_3 & a_4 \\ a_5 & a_6 \\ a_7 & a_8 \end{bmatrix} \middle| a_i \in Q^+ \cup \{0\}, 1 \leq i \leq 8, +, \times_n \right\}$

   the semiring.

   (i) Ideals which are not S-ideals.
   (ii) Subsemirings which are not S-semirings.
   (iii) Zero divisors.
   (iv) Subsemirings which are not ideals.
   (v) Is S a S-semiring?

42. Let $P = \left\{ \begin{bmatrix} a_1 & a_2 & a_3 & a_4 \\ a_5 & a_6 & a_7 & a_8 \\ a_9 & a_{10} & a_{11} & a_{12} \\ a_{13} & a_{14} & a_{15} & a_{16} \end{bmatrix} \middle| a_i \in Q^+ \cup \{0\}, \right.$

$1 \leq i \leq 16\}$ be a semiring under + and natural product $\times_n$.

   (i) Is P a S-ring?
   (ii) Can P have zero divisors?
   (iii) Show P is commutative.

   (iv) If $\times_n$ replaced by usual matrix product will P be a semiring? Justify your claim.
   (v) Find S-ideals in P.
   (vi) Find subsemirings which are not S-subsemirings.

43. Let $V = \{(x_1, x_2, \ldots, x_n) \mid x_i \in R, 1 \leq i \leq n\}$ be a vector space over F. Find Hom (V, V).

44. Let $P = \left\{ \begin{bmatrix} x_1 \\ x_2 \\ \vdots \\ x_{10} \end{bmatrix} \middle| x_i \in Q; 1 \leq i \leq 10 \right\}$ be a linear algebra

over Q.
   (i) Find dimension of P over Q.
   (ii) Find a basis of P over Q.
   (iii) Write P as a direct sum of subspaces.
   (iv) Write P as a pseudo direct sum of subspaces.
   (v) Find a linear operator on P which is invertible.

45. Give an example of a natural Smarandache special field.

46. What is the difference between a natural Smarandache special field and the field?

47. Obtain the special properties enjoyed by S-special strong column matrix linear algebra.

48. Obtain the special and distinct features of S-special strong $3 \times 3$ matrix linear algebra.

49. Find differences between Smarandache vector spaces and Smarandache special strong vector spaces.

50. Let $P = \left\{ \begin{bmatrix} a_1 & a_2 \\ a_3 & a_4 \end{bmatrix} \middle| a_i \in R,\ 1 \le i \le 4 \right\}$ be the $3 \times 3$ square matrix of natural special Smarndache field.

    (i) What are the special properties enjoyed by P?
    (ii) Can P have zero divisors?

51. Obtain some interesting properties about S-special strong column matrix vector spaces constructed over R.

52. Find some applications of S-special strong $m \times n$ ($m \ne n$) linear algebras constructed over Q.

53. Define some nice types of inner products on vector spaces using the natural product $\times_n$.

54. Can linear functionals be defined on S-special super $n \times m$ ($m \ne n$) matrix vector spaces?

55. Let $V = \left\{ \begin{bmatrix} a_1 & a_2 & \ldots & a_{12} \\ a_{13} & a_{14} & \ldots & a_{24} \\ a_{25} & a_{26} & \ldots & a_{36} \end{bmatrix} \middle| a_i \in Q,\ 1 \le i \le 36 \right\}$ be a

    S-special strong vector space over the S-field.

$$F_{3\times 12} = \left\{ \begin{bmatrix} x_1 & x_2 & \cdots & x_{12} \\ x_{13} & x_{14} & \cdots & x_{24} \\ x_{25} & x_{26} & \cdots & x_{36} \end{bmatrix} \middle| \; x_i \in Q, \; 1 \leq i \leq 36 \right\}.$$

(i) Find a basis for V.
(ii) What is the dimension of V over $F_{3\times 12}$?
(iii) Write V as a direct sum of subspaces.
(iv) Write V as a pseudo direct sum of subspaces.

56. Let V be a S-special strong vector space of $n \times n$ matrices over the S-field $F_C$ of $n \times n$ matrices with elements from the field Q.
    (i) Find a basis for V.
    (ii) Write V as a direct sum of subspaces.
    (iii) Write V as a pseudo direct sum of subspaces.
    (iv) Find a linear operator on V.
    (v) Does every subspace W of V have $W^\perp$?
    (vi) Write V as $W + W^\perp$;

57. Obtain some interesting properties about orthogonal subspaces.

58. Find some interesting properties related with S-special row matrix linear algebras.

59. Study the special properties enjoyed by S-special strong $m \times n$ matrix linear algebras ($m \neq n$).

60. Let $S = \left\{ \begin{bmatrix} a_1 \\ a_2 \\ a_3 \\ a_4 \\ a_5 \end{bmatrix} \middle| \; a_i \in Z^+ \cup \{0\}, \; 1 \leq i \leq 5 \right\}$ be a semivector space of column matrices over the semifield $F = Z^+ \cup \{0\}$.

(i) Find basis for S.
(ii) Write S as a direct sum of subsemivector spaces.
(iii) Write S as a pseudo direct sum of subsemivector spaces.
(iv) Can S be a semilinear algebra?

61. Obtain some interesting properties enjoyed by semivector space of column matrices V over the semifield $S = Q^+ \cup \{0\}$.

62. Enumerate the special properties enjoyed by the semivector space of $m \times n$ matrices ($m \neq n$) over the semifield $F = Q^+ \cup \{0\}$.

63. Bring out the differences between the semivector space of column matrices over $Q^+ \cup \{0\}$ and vector space of column matrices over Q.

64. Find some special properties enjoyed by semivector space of $m \times n$ matrices over the semifield $Z^+ \cup \{0\} = S$.

65. Let $V = \left\{ \begin{bmatrix} a_1 & a_2 & a_3 & a_4 & a_5 \\ a_6 & a_7 & \ldots & \ldots & a_{10} \\ a_{11} & a_{12} & \ldots & \ldots & a_{15} \\ a_{16} & a_{17} & \ldots & \ldots & a_{20} \end{bmatrix} \middle| a_i \in Z^+ \cup \{0\}, 1 \le i \le 20 \right\}$ be a semivector space over the semifield $S = Z^+ \cup \{0\}$.

(i) Can V be made into a semilinear algebra over S?

(ii) Find a basis for V.

(iii) Find semivector subspaces of V so that V can be written as a direct sum of semivector subspaces.

(iv) Write V as a pseudo direct sum of semivector subspaces.

(v) Write $V = W \oplus W^\perp$, $W^\perp$ the orthogonal complement of W.

66. Let $F_C^S = \left\{ \left[ \begin{array}{c} a_1 \\ a_2 \\ \hline a_3 \\ a_4 \\ \hline a_5 \\ \hline a_6 \end{array} \right] \middle| a_i \in Q, 1 \leq i \leq 6, \times_n \right\}$ be a semigroup.

(i) Find ideals in $F_C^S$.
(ii) Can $F_C^S$ have subsemigroups which are not ideals?
(iii) Prove $F_C^S$ has zero divisors.
(iv) Find units in $F_C^S$.
(v) Is $F_C^S$ a S-semigroup?

67. Let $F_R^S = \{(a_1\ a_2\ |\ a_3\ a_4\ |\ a_5) \mid a_i \in Z, 1 \leq i \leq 5, \times_n\}$ be a semigroup.

(i) Find subsemigroups which are not ideals in $F_R^S$
(ii) Find zero divisors in $F_R^S$.
(iii) Can $F_R^S$ have units?
(iv) Is $x = (1, -1\ |\ 1\ -1\ |-1)$ a unit in $F_R^S$.

68. Let $F_{5\times 3}^S = \left\{ \left[ \begin{array}{ccc} a_1 & a_2 & a_3 \\ \hline a_4 & a_5 & a_6 \\ a_7 & a_8 & a_9 \\ a_{10} & a_{11} & a_{12} \\ \hline a_{13} & a_{14} & a_{15} \end{array} \right] \middle| a_i \in Q, 1 \leq i \leq 15, \times_n \right\}$ be a

semigroup of super matrices.

(i) Find units in $F^S_{5\times 3}$.

(ii) Is $F^S_{5\times 3}$ a S-semigroup?

(iii) Can $F^S_{5\times 3}$ have S-subsemigroups?

(iv) Can $F^S_{5\times 3}$ have S-ideals?

(v) Does $F^S_{5\times 3}$ have S-zero divisors/

(vi) Can $F^S_{5\times 3}$ have S-idempotents?

(vii) Show $F^S_{5\times 3}$ have only finite number of idempotents.

69. Let $F^S_{4\times 4} = \left\{ \begin{bmatrix} a_1 & a_2 & a_3 & a_4 \\ a_5 & a_6 & a_7 & a_8 \\ a_9 & a_{10} & a_{11} & a_{12} \\ a_{13} & a_{14} & a_{15} & a_{16} \end{bmatrix} \,\middle|\, a_i \in Q,\ 1 \leq i \leq 16 \right\}$ be a

semigroup of super square matrices.

(i) Find zero divisors in $F^S_{4\times 4}$.

(ii) Can $F^S_{4\times 4}$ have S-zero divisors?

(iii) Can $F^S_{4\times 4}$ have S-idempotents?

(iv) Find the main complement of $\begin{bmatrix} 3 & 0 & 1 & 2 \\ 1 & 0 & 0 & 5 \\ 0 & 3 & 15 & 7 \\ 0 & 1 & 0 & 0 \end{bmatrix}$.

(v) Is $F^S_{4\times 4}$ a S-semigroup?

70. Let $F^S_{3\times 7} = \left\{ \begin{bmatrix} a_1 & a_2 & a_3 & a_4 & a_5 & a_6 & a_7 \\ a_8 & a_9 & a_{10} & a_{11} & a_{12} & a_{13} & a_{14} \\ a_{15} & a_{16} & a_{17} & a_{18} & a_{19} & a_{20} & a_{21} \end{bmatrix} \,\middle|\, a_i \in Z, \right.$

$1 \leq i \leq 21,\ \times_n \}$ be a super row vector semigroup.

- (i) Find ideal of $F^S_{3\times 7}$.
- (ii) Is $F^S_{3\times 7}$ a S-semigroup?
- (iii) Can $F^S_{3\times 7}$ have S-ideals?
- (iv) Show $F^S_{3\times 7}$ can have only finite number of idempotents.
- (v) Show $F^S_{3\times 7}$ has no units.

71. Obtain some interesting properties about $(F^S_C, \times_n)$.

72. Find some applications of the semigroup $(F^S_{m\times m}; m \neq n, \times_n)$.

73. Find the difference between $(F^S_{n\times n}, \times)$ and $(F^S_{n\times n}, \times_n)$.

74. Find some special and distinct features enjoyed by $F^S_{9\times 8}$.

75. Prove $\{F^S_{n\times m}, n \neq m, +, \times_n\}$ is a commutative ring of infinite order.

76. Let $F^S_{2\times 5} = \left\{ \begin{bmatrix} a_1 & a_2 & a_3 & a_4 & a_5 \\ a_6 & a_7 & a_8 & a_9 & a_{10} \end{bmatrix} \middle| a_i \in Q, 1 \leq i \leq 10, +, \times_n \right\}$ be a ring of super row vectors.

- (i) Find ideals in $F^S_{2\times 5}$.
- (ii) Is $F^S_{2\times 5}$ a S-ring?
- (iii) Prove $F^S_{2\times 5}$ has ideals.
- (iv) Can $F^S_{2\times 5}$ have S-ideals?
- (v) Can $F^S_{2\times 5}$ have S-zero divisors and S-idempotents?

77. Let $F^S_{8\times 3} = \left\{ \begin{bmatrix} \begin{array}{|ccc|} \hline a_1 & a_2 & a_3 \\ \hline a_4 & a_5 & a_6 \\ \hline a_7 & a_8 & a_9 \\ \hline a_{10} & a_{11} & a_{12} \\ \hline a_{13} & a_{14} & a_{15} \\ \hline a_{16} & a_{17} & a_{18} \\ \hline a_{19} & a_{20} & a_{21} \\ \hline a_{22} & a_{23} & a_{24} \\ \hline \end{array} \end{bmatrix} \middle| a_i \in Q, 1 \leq i \leq 24, +, \times_n \right\}$ be a super column vector ring.

 (i) Prove $F^S_{8\times 3}$ has zero divisors.
 (ii) Prove $F^S_{8\times 3}$ has units.
 (iii) Can $F^S_{8\times 3}$ have S-units?
 (iv) Is $F^S_{8\times 3}$ a S-ring?
 (v) Prove $F^S_{8\times 3}$ has idempotents?

78. Let $F^S_C = \left\{ \begin{bmatrix} a_1 \\ a_2 \\ a_3 \\ a_4 \\ a_5 \\ a_6 \\ a_7 \\ a_8 \\ a_9 \\ a_{10} \end{bmatrix} \middle| a_i \in R, 1 \leq i \leq 10, +, \times_n \right\}$ be a super column matrix ring.

 (i) Find the number of idempotents in $F^S_C$.
 (ii) Show $F^S_C$ has infinite number of zero divisors but only finite number of idempotents.

(iii) Can $F_C^S$ have S-idempotents?

(iv) Can $F_C^S$ have S-units?

79. Let $F_{3\times 3}^S = \left\{ \begin{bmatrix} a_1 & a_2 & a_3 \\ a_4 & a_5 & a_6 \\ a_7 & a_8 & a_9 \end{bmatrix} \middle| a_i \in Q, 1 \le i \le 8, +, \times_n \right\}$ be a square super matrix ring.

    (i) Prove $F_{3\times 3}^S$ is a commutative ring.

    (ii) Find ideals in $F_{3\times 3}^S$.

(iii) Is $F_{3\times 3}^S$ a S-ring?

(iv) Show $F_{3\times 3}^S$ has only finite number of idempotents.

(v) Can $F_{3\times 3}^S$ have S-ideals?

80. Enumerate some special features enjoyed by super matrix rings $F_C^S$ (or $F_R^S$ or $F_{n\times n}^S$ or $F_{n\times m}^S$; $m \ne n$).

81. Find some applications of the rings mentioned in problem (80).

82. Prove $R = \{(a_1 \mid a_2 \mid \ldots \mid a_n) \mid a_i \in Q^+ \cup \{0\}, 1 \le i \le n\}$ is a semigroup under $+$.

    (i) Can R have ideals?

    (ii) Can R have S-zero divisors?

83. Let $P = \{(a_1 \mid a_2 \mid a_3 \, a_4 \, a_5 \mid a_6 \, a_7) \mid a_i \in Z^+ \cup \{0\}, 1 \le i \le 7\}$ be a semigroup under $\times_n$. Find the special properties enjoyed by these semigroups.

84. Let $T = \left\{ \begin{bmatrix} a_1 \\ a_2 \\ a_3 \\ a_4 \\ a_5 \end{bmatrix} \middle| a_i \in Z^+ \cup \{0\}, 1 \le i \le 5, \times_n \right\}$ be a

semigroup under $\times_n$.
   (i) Prove T is a commutative semigroup.
   (ii) Can T have S-zero divisors?
   (iii) Show T can have only finite number of idempotents.
   (iv) Show T can have no units.
   (v) Can T have S-ideals?

85. Let $W = \left\{ \begin{bmatrix} a_1 & a_2 & a_3 & a_4 & a_5 \\ a_6 & a_7 & a_8 & a_9 & a_{10} \\ a_{11} & a_{12} & a_{13} & a_{14} & a_{15} \end{bmatrix} \middle| a_i \in Q^+ \cup \{0\}, 1 \le i \right.$

$\le 15, \times_n\}$ be a semigroup.
   (i) Find the number of idempotents in W.
   (ii) Is W a S-semigroup?
   (iii) Find units in W.
   (iv) Show all elements are not units in W.

86. Let $M = \left\{ \begin{bmatrix} a_1 & a_2 & a_3 & a_4 & a_5 & a_6 \\ a_7 & a_8 & a_9 & a_{10} & a_{11} & a_{12} \\ a_{13} & a_{14} & ... & ... & ... & a_{18} \\ a_{19} & a_{20} & ... & ... & ... & a_{24} \\ a_{25} & a_{26} & ... & ... & ... & a_{30} \end{bmatrix} \middle| a_i \in Q^+ \cup \{0\}, \right.$

$1 \le i \le 30, \times_n\}$ be a semigroup.
   (i) Is M a S-semigroup?

(ii) Does M contain S-zero divisors?

(iii) Prove M has only finite number of idempotents.

(iv) Can M have S-idempotents?

(v) Can M have S-units?

87. Let $M = F_{2\times 3}^S = \left\{ \begin{bmatrix} a_1 & a_2 & a_3 \\ a_4 & a_5 & a_6 \end{bmatrix} \mid a_i \in Z^+ \cup \{0\}, 1 \le i \le 6, +, \times_n \right\}$ be a semiring.

   (i) Prove M is not a semifield.

   (ii) Find subsemirings in M.

88. Obtain some interesting properties enjoyed by column super matrix semirings with entries from $Q^+ \cup \{0\}$.

89. Distinguish between a super square matrix ring and a super square matrix semiring.

90. Let $P = \left\{ \begin{bmatrix} a_1 & a_6 \\ a_2 & a_7 \\ a_3 & a_8 \\ a_4 & a_9 \\ a_5 & a_{10} \end{bmatrix} \mid a_i \in R^+ \cup \{0\}, 1 \le i \le 10, +, \times_n \right\}$ be a semiring.

   (i) Find subsemirings of P.

   (ii) Is P a S-semiring?

   (iii) Can P have S-ideals?

   (iv) Can P have S-idempotents?

   (v) Can P have S-zero divisors?

91. Let T = $\left\{ \begin{bmatrix} a_1 & a_2 & a_3 & a_4 & a_5 & a_6 & a_7 \\ a_8 & a_9 & \cdots & \cdots & \cdots & \cdots & a_{14} \\ a_{15} & a_{16} & \cdots & \cdots & \cdots & \cdots & a_{21} \\ a_{22} & a_{23} & \cdots & \cdots & \cdots & \cdots & a_{28} \\ a_{29} & a_{30} & \cdots & \cdots & \cdots & \cdots & a_{35} \\ a_{36} & a_{37} & \cdots & \cdots & \cdots & \cdots & a_{42} \\ a_{43} & a_{44} & \cdots & \cdots & \cdots & \cdots & a_{49} \end{bmatrix} \middle| a_i \in Z^+ \cup \{0\}, 1 \le i \le 49, +, \times_n \right\}$ be a semiring.

   (i) Show T has only finite number of idempotents in it.

   (ii) Find zero divisors of T.

   (iii) Find idempotents of T.

   (iv) Can T have S-zero divisors?

   (v) Is T a S-semiring?

92. Find some interesting properties of super matrix semirings.

93. Let M = $\left\{ \begin{bmatrix} a_1 & a_2 & a_3 & a_4 \\ a_5 & a_6 & a_7 & a_8 \\ a_9 & a_{10} & a_{11} & a_{12} \\ a_{13} & a_{14} & a_{15} & a_{16} \end{bmatrix} \cup \begin{bmatrix} 0 & 0 & 0 & 0 \\ 0 & 0 & 0 & 0 \\ 0 & 0 & 0 & 0 \\ 0 & 0 & 0 & 0 \end{bmatrix} \right.$ where

   $a_i \in Q^+, 1 \le i \le 16, +, \times_n \}$.

   (i) Is M a semifield?

   (ii) Is M a S-semiring?

94. Can $P = \left\{ \begin{bmatrix} a_1 \\ a_2 \\ a_3 \\ \hline a_4 \\ a_5 \\ \hline a_6 \end{bmatrix} \cup \begin{bmatrix} 0 \\ 0 \\ 0 \\ \hline 0 \\ 0 \\ \hline 0 \end{bmatrix} \text{ where } a_i \in Q^+, 1 \le i \le 6, +, \times_n \right\}$ be a semifield?

95. Is $T = \{(a_1 \mid a_2 \ a_3 \mid a_4 \ a_5 \ a_6) \cup (0 \mid 0 \ 0 \mid 0 \ 0 \ 0) \mid a_i \in Q^+, 1 \le i \le 6, +, \times_n\}$ a semifield?

    Can T have subsemifields?

96. Let $W = \left\{ \begin{bmatrix} a_1 & a_2 & a_3 & a_4 \\ a_5 & a_6 & a_7 & a_8 \\ \hline a_9 & . & . & a_{12} \\ a_{13} & . & . & a_{16} \\ a_{17} & . & . & a_{20} \\ \hline a_{21} & . & . & a_{24} \\ a_{25} & . & . & a_{28} \\ a_{29} & . & . & a_{32} \\ a_{33} & . & . & a_{36} \end{bmatrix} \cup \begin{bmatrix} 0 & 0 & 0 & 0 \\ 0 & 0 & 0 & 0 \\ \hline 0 & 0 & 0 & 0 \\ 0 & 0 & 0 & 0 \\ 0 & 0 & 0 & 0 \\ \hline 0 & 0 & 0 & 0 \\ 0 & 0 & 0 & 0 \\ 0 & 0 & 0 & 0 \\ 0 & 0 & 0 & 0 \end{bmatrix} \text{ where } \right.$

    $a_i \in R^+, 1 \le i \le 36, +, \times_n\}$ be a semifield of super column vectors. Find subsemifield of W. Can W have subsemirings?

97. Find applications of matrices with natural product.

98. Prove super matrices of same type under natural product is a semigroup with zero divisors.

99. Let $X = \left\{ \begin{bmatrix} a_1 & a_2 & a_3 & a_4 \\ a_5 & a_6 & a_7 & a_8 \\ a_9 & a_{10} & a_{11} & a_{12} \\ a_{13} & a_{14} & a_{15} & a_{16} \end{bmatrix} \right.$ where $a_i \in R^+ \cup \{0\}$, $1 \le i \le 16\}$ be a super square matrix linear algebra over the semifield $S = R^+ \cup \{0\}$.

(i) Find a basis of X over S.

(ii) Is X finite dimensional?

(iii) Write X as a direct sum of semivector subspaces.

(iv) Write X as a pseudo direct sum of semivector subspaces.

(v) Let $V \subseteq X$ be a subspace find $V^\perp$ so that $X = V + V^\perp$ a complement?

100. Let $V = \left\{ \begin{bmatrix} a_1 & a_2 & a_3 \\ a_4 & a_5 & a_6 \\ a_7 & a_8 & a_9 \end{bmatrix} \right.$ $a_i \in Q^+ \cup \{0\}$, $1 \le i \le 9\}$ be a super matrix semivector space over the semifield $S = Z^+ \cup \{0\}$.

(i) Can $W = \left\{ \begin{bmatrix} a_1 & a_2 & a_3 \\ a_4 & a_5 & a_6 \\ a_7 & a_8 & a_9 \end{bmatrix} \right.$ $a_i \in Z^+ \cup \{0\}$, $1 \le i \le 9\} \subseteq V$ have a orthogonal complement space? Justify your claim.

101. Let $V = \left\{ \begin{bmatrix} a_1 & a_2 & a_3 & a_4 & a_5 & a_6 \\ a_7 & a_8 & a_9 & a_{10} & a_{11} & a_{12} \\ a_{13} & a_{14} & a_{15} & a_{16} & a_{17} & a_{18} \\ a_{19} & a_{20} & a_{21} & a_{22} & a_{23} & a_{24} \end{bmatrix} \middle| a_i \in Q^+ \cup \{0\}, \right.$

$1 \le i \le 24 \}$ be a super matrix semilinear algebra over $S = Z^+ \cup \{0\}$.

(i) Is V finite dimensional?

(ii) Find subspaces of V so that V can be written as a direct sum of super matrix semivector subspaces.

(iii) Write V as a pseudo direct sum of super matrix semilinear algebra.

(iv) Let $M = \left\{ \begin{bmatrix} a_1 & 0 & 0 & a_6 & a_7 & a_8 \\ a_2 & 0 & 0 & a_9 & a_{10} & a_{11} \\ a_3 & 0 & 0 & a_{12} & a_{13} & a_{14} \\ 0 & a_4 & a_5 & 0 & 0 & 0 \end{bmatrix} \middle| a_i \in Q^+ \cup \right.$

$\{0\}, 1 \le i \le 14\} \subseteq V$ be a super matrix semilinear subalgebra of V.
a) How many complements exists for M?

b) Write down the main complement of M.

102. Let $V = \left\{ \begin{bmatrix} a_1 & a_2 & a_3 & a_4 & a_5 & a_6 \\ a_7 & a_8 & a_9 & a_{10} & a_{11} & a_{12} \\ a_{13} & a_{14} & a_{15} & a_{16} & a_{17} & a_{18} \end{bmatrix} \middle| a_i \in Q^+ \cup \{0\}, \right.$

$1 \le i \le 18\}$ be a super row vector semilinear algebra over the semifield $S = Q^+ \cup \{0\}$.

(i) Find the dimension of V over S.
(ii) Can V have more than one basis over S?

(iii) Can V have linearly independent elements whose number (cardinality) is greater than that of cardinality of the basis of V over S?

103. Let M = $\left\{ \begin{bmatrix} a_1 & a_2 & a_3 & a_4 \\ \hline a_5 & a_6 & a_7 & a_8 \\ a_9 & a_{10} & a_{11} & a_{12} \\ \hline a_{13} & a_{14} & a_{15} & a_{16} \\ a_{17} & a_{18} & a_{19} & a_{20} \\ a_{21} & a_{22} & a_{23} & a_{24} \\ \hline a_{25} & a_{26} & a_{27} & a_{28} \\ a_{29} & a_{30} & a_{31} & a_{32} \\ a_{33} & a_{34} & a_{35} & a_{36} \\ a_{37} & a_{38} & a_{39} & a_{40} \end{bmatrix} \right.$ where $a_i \in Z^+ \cup \{0\}$, $1 \leq i \leq 40\}$ be a super column vector semilinear algebra over the semifield $S = Z^+ \cup \{0\}$.

(i) Find a basis for M over S.

(ii) Write M as a pseudo direct sum of super column vector semilinear subalgebra over S.

(iii) Write M as a direct sum of super column vector semilinear subalgebras over S.

(iv) Write $M = W + W^\perp$ where $W^\perp$ is the orthogonal complement of W.

104. Let $P = \left\{ \sum_{i=0}^{\infty} a_i x^i \mid a_i = (m_1 \, m_2 \mid m_3 \, m_4 \, m_5 \mid m_6 \, m_7 \mid m_8) \right.$ where $m_i \in Q^+ \cup \{0\}$, $1 \leq i \leq 8\}$ be a super row matrix semilinear algebra over the semifield $S = Z^+ \cup \{0\}$.

(i) Write P as a direct sum of semilinear subalgebras.

(ii) Let $M = \left\{ \sum_{i=0}^{\infty} a_i x^i \,\middle|\, a_i = (0\ 0\ |\ m_1\ m_2\ m_3\ |\ 0\ 0\ |\ 0) \right.$ $m_1$, $m_2$, $m_3 \in Q^+ \cup \{0\}\} \subseteq P$ be a super row matrix semi linear subalgebra over S. Find $M^\perp$ so that $P = M + M^\perp$.

(iii) Let $T = \left\{ \sum_{i=0}^{\infty} a_i x^i \,\middle|\, a_i = (d_1\ d_2\ |\ 0\ 0\ 0\ |\ d_3\ d_4\ |\ 0), d_j \in \right.$ $Z^+ \cup \{0\}, 1 \leq j \leq 4\} \subseteq P$ be a semi sublinear algebra of P over S. Can we find a orthogonal complement of T of P over S so that $T + T^\perp = P$?

105. Let $S = \left\{ \sum_{i=0}^{\infty} a_i x^i \,\middle|\, a_i = \begin{bmatrix} d_1 \\ d_2 \\ d_3 \\ d_4 \\ \hline d_5 \\ d_6 \\ d_7 \\ \hline d_8 \\ d_9 \\ \hline d_{10} \\ d_{11} \end{bmatrix}, d_j \in Z^+ \cup \{0\}, 1 \leq j \leq 11 \right\}$ be

a super column matrix semilinear algebra over the semifield $F = Z^+ \cup \{0\}$.

(i) Find a basis for S over F.

(ii) Write S as a direct sum of linear sbalgebras.

(iii) Write S as $P + P^\perp$ so that for every $x \in P$ we have every $y \in P^\perp$ with $x \times_n y = (0)$.

106. Let $T = \left\{ \sum_{i=0}^{\infty} a_i x^i \,\middle|\, a_i = \begin{bmatrix} m_1 & m_2 & m_3 \\ \hline m_4 & m_5 & m_6 \\ \hline m_7 & m_8 & m_9 \\ \hline m_{10} & m_{11} & m_{12} \\ \hline m_{13} & m_{14} & m_{15} \\ \hline m_{16} & m_{17} & m_{18} \\ \hline m_{19} & m_{20} & m_{21} \\ \hline m_{22} & m_{23} & m_{24} \\ \hline m_{25} & m_{26} & m_{27} \\ \hline m_{28} & m_{29} & m_{30} \\ \hline m_{31} & m_{32} & m_{33} \end{bmatrix} \; m_i \in Z^+ \cup \{0\}, \right.$

$1 \leq i \leq 33\}$ be a super column vector semilinear algebra over $S = Z^+ \cup \{0\}$.

(i) Find a basis for T over S.

(ii) Can T have more than one basis?

(iii) Write T as direct sum.

(iv) Write T as pseudo direct sum.

107. Let $M = \left\{ \sum_{i=0}^{\infty} d_i x^i \,\middle|\, d_i = \begin{bmatrix} a_1 & a_2 & a_3 & a_4 \\ a_5 & a_6 & a_7 & a_8 \\ \hline a_9 & a_{10} & a_{11} & a_{12} \\ a_{13} & a_{14} & a_{15} & a_{16} \\ \hline a_{17} & a_{18} & a_{19} & a_{20} \\ \hline a_{21} & a_{22} & a_{23} & a_{24} \end{bmatrix} \; a_i \in Z^+ \cup \right.$

$\{0\}$, $1 \leq i \leq 24\}$ be a super matrix semilinear algebra over the semifield $S = Z^+ \cup \{0\}$.

(i) Find a basis of M over S.

(ii) Prove M has subspaces which are complements to each other in m.

108. Obtain some unique properties enjoyed by super matrix coefficient polynomial rings.

109. Find some special features of super matrix coefficient polynomial semivector spaces over the semifield $S = Z^+ \cup \{0\}$.

110. Describe any special feature enjoyed by super matrix coefficient polynomial semivector space over the semifield $S = Z^+ \cup \{0\}$.

111. Give some applications of super matrix coefficient polynomial semilinear algebras defined over a semifield.

# FURTHER READING

# INDEX

**D**

# ABOUT THE AUTHORS

**Dr.W.B.Vasantha Kandasamy** is an Associate Professor in the Department of Mathematics, Indian Institute of Technology Madras, Chennai. In the past decade she has guided 13 Ph.D. scholars in the different fields of non-associative algebras, algebraic coding theory, transportation theory, fuzzy groups, and applications of fuzzy theory of the problems faced in chemical industries and cement industries. She has to her credit 646 research papers. She has guided over 68 M.Sc. and M.Tech. projects. She has worked in collaboration projects with the Indian Space Research Organization and with the Tamil Nadu State AIDS Control Society. She is presently working on a research project funded by the Board of Research in Nuclear Sciences, Government of India. This is her $63^{rd}$ book.

On India's 60th Independence Day, Dr.Vasantha was conferred the Kalpana Chawla Award for Courage and Daring Enterprise by the State Government of Tamil Nadu in recognition of her sustained fight for social justice in the Indian Institute of Technology (IIT) Madras and for her contribution to mathematics. The award, instituted in the memory of Indian-American astronaut Kalpana Chawla who died aboard Space Shuttle Columbia, carried a cash prize of five lakh rupees (the highest prize-money for any Indian award) and a gold medal.
She can be contacted at vasanthakandasamy@gmail.com
Web Site: http://mat.iitm.ac.in/home/wbv/public_html/
or http://www.vasantha.in

---

**Dr. Florentin Smarandache** is a Professor of Mathematics at the University of New Mexico in USA. He published over 75 books and 200 articles and notes in mathematics, physics, philosophy, psychology, rebus, literature. In mathematics his research is in number theory, non-Euclidean geometry, synthetic geometry, algebraic structures, statistics, neutrosophic logic and set (generalizations of fuzzy logic and set respectively), neutrosophic probability (generalization of classical and imprecise probability). Also, small contributions to nuclear and particle physics, information fusion, neutrosophy (a generalization of dialectics), law of sensations and stimuli, etc. He got the 2010 Telesio-Galilei Academy of Science Gold Medal, Adjunct Professor (equivalent to Doctor Honoris Causa) of Beijing Jiaotong University in 2011, and 2011 Romanian Academy Award for Technical Science (the highest in the country). Dr. W. B. Vasantha Kandasamy and Dr. Florentin Smarandache got the 2011 New Mexico Book Award for Algebraic Structures. He can be contacted at smarand@unm.edu